\documentclass[11pt,leqno]{amsart}
\usepackage[letterpaper,textwidth=6.1in,textheight=8.8in,centering]{geometry}
\usepackage{amsmath,amssymb,mathtools,mathrsfs,bm,enumitem,microtype}
\usepackage{needspace}
\usepackage{etoolbox}
\usepackage[T1]{fontenc}
\usepackage{lmodern}
\usepackage{xcolor}
\usepackage[colorlinks=true,allcolors=blue,linktoc=page]{hyperref}
\usepackage[nameinlink,noabbrev]{cleveref}
\hypersetup{
  pdftitle={Differential Monodromy, q-Opers, and Universal (G,q)-Wronskians},
  pdfauthor={Anton M. Zeitlin}
}
\numberwithin{equation}{section}
\makeatletter
\renewcommand\paragraph{\@startsection{paragraph}{4}{\z@}%
  {.8\baselineskip \@plus .2\baselineskip}%
  {-.5em}%
  {\normalfont\normalsize\bfseries}}
\makeatother
\theoremstyle{plain}
\newtheorem{theorem}{Theorem}[section]
\newtheorem{proposition}[theorem]{Proposition}
\newtheorem{lemma}[theorem]{Lemma}
\newtheorem{corollary}[theorem]{Corollary}

\theoremstyle{definition}
\newtheorem{definition}[theorem]{Definition}
\AtBeginEnvironment{definition}{\pushQED{\qed}}
\AtEndEnvironment{definition}{\popQED}
\theoremstyle{plain}
\newtheorem{hypothesis}{Hypothesis}
\theoremstyle{definition}
\newtheorem{example}[theorem]{Example}

\theoremstyle{remark}
\newtheorem{remark}[theorem]{Remark}
\newcommand{\g}{\mathfrak g}
\newcommand{\gh}{\mathfrak h}
\newcommand{\n}{\mathfrak n}

\newcommand{\qaff}{\mathsf q}
\title[Differential Monodromy and $q$-Opers]{Differential Monodromy, $q$-Opers, and Universal $(G,q)$-Wronskians}
\author[A. M. Zeitlin]{Anton M. Zeitlin}
\address{School of Mathematics,\newline
  Georgia Institute of Technology,\newline
  686 Cherry Street,\newline
  Atlanta, GA 30332-0160, USA}

\email{\href{mailto:zeitlin@gatech.edu}{zeitlin@gatech.edu}}
\date{}
\begin{document}
\begin{abstract}
We give a geometric formulation of the $q$DE/ODE/IM correspondence for families of differential operators associated with simple Lie algebras of ADE type.  Under analytic hypotheses, shifted monodromy defines a $Z$-twisted $(G,q)$-oper.  Its underlying $q$-connection is represented by the constant Cartan element $Z$ in the Frobenius frame at zero and has Steinberg form in the frame of solutions normalized by prescribed asymptotics in a sector at infinity.  We identify the central connection between these frames with a universal meromorphic $(G,q)$-Wronskian.  Its generalized minors give $Q$-functions satisfying the ADE $QQ$-relations.  When $q$ is not a root of unity, we parametrize shifted monodromy by quantum-affine $q$-characters evaluated on shifted $Q$-ratios, using the same character formulas that describe normalized transfer-matrix eigenvalues.  We verify the analytic hypotheses for generic untwisted Feigin--Frenkel--Hernandez (FFH) opers, thereby relating these affine opers to $Z$-twisted $(G,q)$-opers.  We conclude by reviewing a construction from a forthcoming companion paper on the non-simply-laced correspondence.  Fractional rotation of twisted vacuum FFH opers gives rise to a new class of difference opers, the $(G,\eta,\sigma)$-opers, whose connections combine the spectral shift $u\mapsto\eta u$ with a Dynkin diagram automorphism $\sigma$.  The minors of their universal Wronskians satisfy the $QQ$-relations of the Langlands-dual quantum affine types.
\end{abstract}
\maketitle

\tableofcontents

\section{Introduction}\label{sec:introduction}

\subsection{From spectral determinants to a group-valued correspondence}\label{sec:intro-spectral-problems}

A basic problem in quantum integrable systems is to determine the joint spectrum of a commuting family of Hamiltonians.  In many models these Hamiltonians are generated by transfer matrices depending on a spectral parameter.  The $TQ$ relation provides a functional approach to this spectral problem \cite{Baxter1972}: on a common eigenvector of the transfer matrices and a commuting $Q$-operator, it relates the transfer eigenvalue to shifted values of the corresponding scalar $Q$-function.  With the required analyticity and asymptotic conditions, this functional equation can be used to determine the transfer eigenvalues and hence the Hamiltonian spectrum.  In particular, regularity of the transfer eigenvalue at zeros of $Q$ gives equations for the positions of these zeros, called the \emph{Bethe equations}.  In conformal field theory, Bazhanov, Lukyanov, and Zamolodchikov constructed such commuting transfer and Baxter $Q$-operators for the quantum Korteweg--de Vries (KdV) hierarchy \cite{BLZ1996,BLZ1997}.  These operators act in representations of the Virasoro algebra; the asymptotic expansions of the transfer operators recover the commuting quantum KdV Hamiltonians.

Their two $Q$-operators also satisfy a bilinear relation involving spectral shifts, called the quantum Wronskian identity \cite{BLZ1997}.  This is the rank-one $QQ$-relation.  In higher rank, the fundamental $QQ$-relations couple pairs of $Q$-functions indexed by the vertices of a Dynkin diagram: the difference of two products of shifted functions at one vertex is related to a product of functions at neighboring vertices.  Bethe equations follow by evaluating suitable shifts of the relations at zeros and eliminating the second $Q$-function, when the required denominators are nonzero.  Frenkel and Hernandez generalized the Baxter relations to arbitrary untwisted quantum affine algebras \cite{FH2015} and proved representation-theoretic analogues of the higher-rank $QQ$-relations \cite{FrenkelHernandez2018}.

In 1998, Dorey and Tateo observed that the spectral determinants of anharmonic oscillators, whose zeros are the Schr\"odinger eigenvalues, satisfy functional equations associated with integrable quantum field theories \cite{DoreyTateo1999a}.  Bazhanov, Lukyanov, and Zamolodchikov identified these determinants with vacuum eigenvalues of the $Q$-operators and extended the correspondence to general Virasoro highest weights \cite{BLZ2001}.  Dorey and Tateo related the corresponding transfer eigenvalues to Stokes data, which describe the changes between solution bases normalized in different angular sectors at infinity \cite{DoreyTateo1999b}.  The early determinant identities were also related to short exact sequences of quantum-affine representations \cite{Suzuki1999}.  These results concern second-order Schr\"odinger equations and the rank-one Lie algebra $\mathfrak{sl}_2$.

In the rank-one case, the spectral determinants are obtained from connection coefficients of the Schr\"odinger equation.  A solution selected by its decay in a sector at infinity is expanded in a basis normalized at zero; the resulting coefficients are the $Q$-functions.  A \emph{Symanzik rotation} relates solutions in different sectors by rotating the differential variable and rescaling the spectral parameter so that the equation retains its form.  Linear relations among the resulting solutions give the Stokes coefficients and the $TQ$ equation.  Thus, in these examples, the spectral functions of a quantum integrable model can be studied by means of a finite-order linear differential equation.  This is the ordinary differential equation/integrable model (ODE/IM) correspondence; we refer to \cite{DDTReview2007} for its development.  The analytic background includes the theory of sectorial solutions and Stokes multipliers, together with the analysis of spectral determinants and exact connection formulas \cite{Sibuya1975,Voros1983}.

The passage to higher rank led to $N$th-order scalar equations for $\mathfrak{sl}_N$ \cite{Suzuki2000,DDT2000} and to scalar or pseudodifferential equations for the other classical Lie algebras \cite{DDMST2007}.  Feigin and Frenkel proposed a common conjectural framework for these spectral problems in terms of affine opers for the Langlands dual affine algebra \cite{FeiginFrenkel2011}.  These are connections with additional flag data.  Their scalar realizations recover the earlier rank-one and type-$A$ equations; we explain this passage in Section~\ref{sec:intro-FFH-opers}.

The formulation of the classical-type equations as first-order systems in fundamental representations was developed in \cite{Sun2012}.  Masoero, Raimondo, and Valeri studied the differential equations associated with a class of these affine connections, first for simply-laced Lie algebras \cite{MRV2016}.  They represented the same connection in each fundamental representation of the underlying finite-dimensional Lie algebra, giving a first-order matrix differential equation for each vertex of the Dynkin diagram.  In each equation, they constructed a solution distinguished by its asymptotic behavior at infinity and expanded it in a basis normalized at zero.  The resulting connection coefficients are the $Q$-functions attached to that vertex.

They related these families of $Q$-functions by proving tensor identities among the distinguished solutions.  Under intertwining maps, these identities relate the exterior product of two solutions related by Symanzik rotations at one vertex to the tensor product of distinguished solutions at neighboring vertices.  These tensor identities form the $\Psi$-system.  Its scalar type-$A$ version was proved by Wronskian identities in \cite{DDMST2007}, where extensions to the other classical types were proposed.  Since the rotations also shift the spectral parameter, expansion in the normalized bases at zero yields solutions of the fundamental $QQ$-relations described above.  The non-simply-laced cases were treated using Langlands-dual affine connections \cite{MRV2017}.

The bilinear form of these relations suggests looking for a matrix whose generalized minors are the $Q$-functions.  The Pl\"ucker identity supplies the relation between minors, and the covariance of the differential equation supplies the spectral shift \cite{FominZelevinsky1999,KZ2021}.  For $SL_N$, generalized minors are ordinary matrix minors; for other simple groups, they are matrix coefficients associated with highest-weight vectors and their Weyl-group translates in fundamental representations.  We compare bases of solutions normalized at zero and infinity and show that their transition matrices in the fundamental representations come from a single group element.

We carry out this construction for a connected, simply connected, simple complex group $G$ of type $A$, $D$, or $E$, with Lie algebra $\g$.  We consider a spectral family of $\g$-valued differential operators
\[
 \mathcal L_u=\partial_z+\mathscr A(z;u)
\]
parametrized by the spectral parameter $u\in U\subset\mathbb C^\times$.

Near zero, the leading powers are determined by a fixed element $\ell$ of a Cartan subalgebra $\gh\subset\g$.  We impose the normalization in a specified local $G$-valued gauge $T_0(z;u)$, defined on a punctured sector at zero.  For a chosen vector $v$ of weight $\mu$ in a fundamental representation $\rho_V:G\to GL(V)$, the corresponding solution $y$ of the linear problem $\partial_z y+\rho_V(\mathscr A(z;u))y=0$ is normalized by
\[
 \rho_V(T_0(z;u))y(z;u)
 =z^{-\mu(\ell)}(v+o(1))
 \qquad(z\to0).
\]
Here $\rho_V$ also denotes the induced Lie algebra representation.  We require the gauged solutions to have prescribed convergent expansions in possibly nonintegral powers of $z$ and powers of $\log z$, called \emph{Frobenius expansions}.  Within each allowed expansion class, the leading vector fixes the solution uniquely.  Powers and logarithms use fixed branches on the universal cover of the punctured plane.

The frame formed by these solutions, normalized at zero, is denoted by $\mathcal Y_{\mathrm F,V}$; the subscript $\mathrm F$ refers to this Frobenius normalization.

Near infinity, the solutions have exponential factors whose growth or decay depends on the direction of approach.  We prescribe their asymptotics on angular sectors.  We write $\mathcal L_u^T:=T\mathcal L_uT^{-1}$ for the differential operator after a change of frame $T$.  The irregular behavior considered here is governed by a \emph{cyclic element}: after a prescribed change of frame $T_\infty$, the equation has the form
\[
 \mathcal L_u^{T_\infty}=\partial_z+p(z;u)\Lambda+R(z;u),
 \qquad
 \Lambda=\sum_{i\in I}f_i+v_\theta,
 \qquad R(z;u)=o\bigl(p(z;u)\bigr)
\]
as $z\to\infty$ on the chosen sectors.  Here $I$ indexes the simple roots, $f_i$ are nonzero negative simple-root vectors, $v_\theta$ is a nonzero vector for the highest root $\theta$, and $p$ is a scalar function.  Thus $p\Lambda$ is the leading coefficient in this frame.  The name \emph{cyclic} is illustrated by $SL_N$: in the defining representation, $\sum_i f_i$ sends each coordinate line $\mathbb C e_j$ to $\mathbb C e_{j+1}$ for $j<N$, while the highest-root term closes the cycle by sending $\mathbb C e_N$ to $\mathbb C e_1$.  The nonzero coefficients along this cycle need not all equal one.

In every ADE type, the cyclic structure of $\Lambda$ is governed by a Coxeter element, a product of the simple reflections, each occurring once.  Its order is the Coxeter number $h$.  The centralizer of $\Lambda$ in $\g$ is a Cartan subalgebra.  There is a conjugation preserving this subalgebra that acts on it as a Coxeter element and sends $\Lambda$ to $e^{-2\pi i/h}\Lambda$ \cite{Kostant1959}.  Consequently, the spectrum of $\Lambda$ in every representation is invariant under multiplication by $e^{2\pi i/h}$.

In a representation $V$, an eigenvalue $\nu$ of $\rho_V(\Lambda)$ gives the exponential factor $e^{-\nu S(z;u)}$, where $\partial_zS=p$.  Our hypotheses impose $p(z;u)\sim z^{1/h-1}$, so $S(z;u)\sim hz^{1/h}$.  A full turn in $z$ multiplies this leading phase by $e^{2\pi i/h}$ and hence permutes the leading exponential factors by the same cyclic symmetry.  The term \emph{cyclic type at infinity} refers to the leading coefficient $p\Lambda$ and its Coxeter symmetry.

In the gauge $T_\infty$, each normalized solution, divided by its exponential factor $e^{-\nu S(z;u)}$, tends to a prescribed eigenvector of $\Lambda$.  We assume that these solutions form a frame on a sector of the universal cover whose angular opening in the variable $z^{1/h}$ exceeds $\pi$.  On this sector, the prescribed asymptotics determine the frame uniquely.  This is the sectorial frame $\mathcal Y_{\mathrm{as},V}$, normalized at infinity.  The normalization uses a fixed phase $S$, including its additive constant and all lower-order terms that do not tend to zero.  Such terms, including possible logarithms, affect the normalized solutions and how they change under analytic continuation.  Classical background on sectorial asymptotic solutions is given in \cite{Wasow1965,Sibuya1975}.

To compare the local frames, we continue both to a common connected domain in $z$ avoiding the singularities.  Preservation of the invariant tensors defining $G$ then lets us assemble the fundamental frames at each end into one $G$-valued frame.  At zero, we require the direct sum of the normalized fundamental frames to preserve a finite family of such tensors; at infinity, preservation follows from the sectorial normalization.  The complete analytic assumptions are stated in Section~\ref{sec:analytic-input}.  We denote the resulting unique $G$-valued frames by $\mathcal Y_{\mathrm F}$ at zero and $\mathcal Y_{\mathrm{as}}$ at infinity.  Their matrices in $V$ are $\mathcal Y_{\mathrm F,V}$ and $\mathcal Y_{\mathrm{as},V}$.  These frames respect tensor products, duals, and intertwiners.  Their meromorphic transition is
\[
 \mathcal Y_{\mathrm{as}}(z;u)=\mathcal Y_{\mathrm F}(z;u)C(u),
 \qquad C:U\longrightarrow G.
\]
Since the two frames solve the same differential equation, $C$ is independent of $z$.  We call it the \emph{central connection matrix}.  Its fundamental highest-weight columns give the connection coefficients of the distinguished solutions discussed above.

To relate $C$ to a difference equation, we use covariance of the differential family under the map
\[
 \mathsf S_1(z,u)=(e^{2\pi i}z,qu),\qquad q\in\mathbb C^\times.
\]
Here $e^{2\pi i}z$ denotes the point on the universal cover obtained by increasing the argument of $z$ by $2\pi$.  The pullback $\mathsf S_1^*$ on solutions is called \emph{shifted monodromy}.  For $n\in\mathbb Z$, the iterate $\mathsf S_n=\mathsf S_1^n$ has pullback $\mathsf S_n^*$.  The inverse of $\mathsf S_1^*$ gives transport from $u$ to $qu$.  Before computing it in the normalized frames, we describe the bundle map and the two flag conditions that enter the construction.

\subsection{Differential monodromy and \texorpdfstring{$(G,q)$}{(G,q)}-opers}\label{sec:intro-qopers}

We assume $qU=U$.  A \emph{$q$-connection} on the trivial principal $G$-bundle $P=U\times G$ is a meromorphic family of $G$-equivariant isomorphisms
\[
 \mathcal A_u:P_u\longrightarrow P_{qu}.
\]
Here $G$ acts on each fiber on the right.  Equivariance means $\mathcal A_u(ph)=\mathcal A_u(p)h$ for $h\in G$, so the family is determined by a meromorphic function $A_q:U\to G$ through
\[
 \mathcal A_u(u,g)=(qu,A_q(u)g).
\]
In a representation $\rho_V:G\to GL(V)$, the associated vector bundle is $U\times V$, and the induced maps are $(u,v)\mapsto(qu,\rho_V(A_q(u))v)$.  A section $s:U\to V$ preserved by these maps satisfies the \emph{$q$-difference equation}
\[
 s(qu)=\rho_V(A_q(u))s(u).
\]
Thus a $q$-connection compares fibers at spectral values related by multiplication by $q$.  A meromorphic change of trivialization $s(u)\mapsto\rho_V(\eta(u))s(u)$, with $\eta:U\to G$, changes its matrix by
\[
 A_q(u)\longmapsto\eta(qu)A_q(u)\eta(u)^{-1}.
\]
This is a \emph{$q$-gauge transformation}.  A $q$-connection equivalent in this way to a constant $Z\in G$ is called \emph{$Z$-twisted}.

To define a $(G,q)$-oper, we add a Borel reduction to this bundle.  For a Borel subgroup $B\subset G$, a \emph{meromorphic $B$-reduction} is a principal $B$-subbundle $Q\subset P$ on a dense open subset of $U$.  In the trivialization $P_u=\{u\}\times G$, its fiber $Q_u$ is one right $B$-orbit, hence has the form $\{u\}\times gB$ for some $g\in G$.  The resulting points $F(u)=gB$ are required to define a meromorphic map $F:U\to G/B$, where $G/B$ is the flag variety.  Replacing $g$ by $gb$, $b\in B$, leaves this point unchanged.  For $G=SL_N$ and $B$ upper triangular, $F(u)$ is a complete flag of subspaces of $\mathbb C^N$; the coset $gB$ represents the spans of the first $a$ columns of $g$, $1\le a<N$.

Intrinsically, $Q$ determines a meromorphic section $s_Q:U\to P/B$ of the flag bundle.  In the chosen trivialization, $s_Q(u)=(u,F(u))$.  We reserve $Q$ for the principal subbundle, $F$ for its representing map, and $F(u)$ for the individual flag; the map $F$ depends on the trivialization.

The $q$-connection sends the reduction fiber $Q_u=\{u\}\times gB$ to $\{qu\}\times A_q(u)gB$.  Its induced action on the corresponding flag is therefore
\[
 F(u)=gB\longmapsto A_q(u)F(u):=A_q(u)gB.
\]
This gives two flags in the fiber at $qu$: the chosen flag $F(qu)$ and the image of the flag at $u$.

We use opposite Borel subgroups $B_\pm$, with common Cartan torus $H$, $\operatorname{Lie}(H)=\gh$, and unipotent radicals $N_\pm$.  A \emph{$(G,q)$-oper} is a $q$-connection together with a meromorphic $B_+$-reduction $Q_+\subset P$, represented in the chosen trivialization by a map $F_+:U\to G/B_+$, for which the pair of flags
\[
 \bigl(F_+(qu),A_q(u)F_+(u)\bigr)
\]
has relative position equal to a fixed Coxeter element $c$ for generic $u$ \cite{FKSZ2020}.  Explicitly, if locally $F_+(u)=g_+(u)B_+$, this condition is
\[
 g_+(qu)^{-1}A_q(u)g_+(u)\in B_+\dot cB_+.
\]
For the ordered word $c=s_{i_1}\cdots s_{i_r}$, we use its pinned lift
\[
 \dot c=\widehat s_{i_1}^{-1}\cdots\widehat s_{i_r}^{-1}\in N_G(H).
\]
The pinning fixes a simple-root homomorphism $SL_2\to G$ at each vertex $i$. The element $\widehat s_i^{-1}$ is the image of $\left(\begin{smallmatrix}0&-1\\1&0\end{smallmatrix}\right)$ under this homomorphism. The double-coset condition is independent of the representatives $g_+$ and specifies the incidence relations between the two flags.  For $SL_2$, a flag is a line, and the condition says that the two lines are distinct.

A \emph{Miura $(G,q)$-oper} is a $(G,q)$-oper $(P,\mathcal A,Q_+)$ together with a meromorphic $B_-$-reduction $Q_-\subset P$ preserved by $\mathcal A$. This additional reduction is called the \emph{Miura structure}. Preservation means that, for generic $u$,
\[
 \mathcal A_u\bigl((Q_-)_u\bigr)=(Q_-)_{qu}.
\]
For its representing map $F_-:U\to G/B_-$, this is equivalent to the flag equation
\[
 A_q(u)F_-(u)=F_-(qu).
\]
In this trivialization the triple $(A_q,F_+,F_-)$ therefore records the Miura $(G,q)$-oper. The Coxeter relative-position condition concerns $F_+$, while the preservation equation concerns $F_-$.

We now apply these definitions to the normalized differential frames.  In a representation $V$, every differential solution at parameter $u$ has unique expressions
\[
 y(z;u)=\mathcal Y_{\mathrm{as},V}(z;u)v
       =\mathcal Y_{\mathrm F,V}(z;u)w,
 \qquad w=\rho_V(C(u))v.
\]
The vectors $v,w\in V$ are independent of $z$; they are the coefficients of the solution in the sectorial and Frobenius bases, respectively.  For a family of solutions they may depend on $u$.  Each normalized frame thus identifies the solution space at $u$ with $V$, and $C(u)$ relates these two identifications.

Shifted monodromy $\mathsf S_1^*$ carries solutions at $qu$ to solutions at $u$.  Compatibility with the prescribed Frobenius expansions gives
\[
 \mathsf S_1^*\mathcal Y_{\mathrm F}
 =\mathcal Y_{\mathrm F} Z^{-1},
 \qquad Z:=e^{2\pi i\ell}.
\]
In the Frobenius basis, $\mathsf S_1^*$ sends a coefficient vector $w$ at $qu$ to $\rho_V(Z)^{-1}w$ at $u$.  We use its inverse as the $q$-connection from $u$ to $qu$, with constant matrix $Z$.  Its matrix in the sectorial basis is
\begin{equation}\label{eq:intro-central-q-connection}
 \mathscr M(u)
 :=(\mathsf S_1^*\mathcal Y_{\mathrm{as}})^{-1}\mathcal Y_{\mathrm{as}}
 =C(qu)^{-1}ZC(u).
\end{equation}
Both matrices in the defining ratio solve the differential equation at $u$, so the ratio is independent of $z$.  Thus $\mathsf S_1^*\mathcal Y_{\mathrm{as}}=\mathcal Y_{\mathrm{as}}\mathscr M(u)^{-1}$.

On a family of solutions, shifted monodromy acts on the sectorial coefficient function by $v(u)\mapsto\rho_V(\mathscr M(u)^{-1})v(qu)$.  For a section preserved by this $q$-connection, its sectorial and Frobenius coefficient functions satisfy
\[
 v(qu)=\rho_V(\mathscr M(u))v(u),\qquad
 w(qu)=\rho_V(Z)w(u),\qquad
 w(u)=\rho_V(C(u))v(u).
\]

\Needspace{13\baselineskip}
We use the sectorial frame to define the $B_+$-reduction and the Frobenius frame to define the $B_-$-reduction preserved by the constant twist $Z$.  In their respective trivializations, these reductions are $U\times B_+$ and $U\times B_-$.  We denote the simple roots by $\alpha_i$ and the corresponding fundamental weights by $\omega_i$.  In the fundamental representation $\rho_i:G\to GL(V_i)$, with $V_i=V(\omega_i)$, the corresponding coefficient lines are $\mathbb Cv_{\omega_i}$ and $\mathbb Cv_{w_0\omega_i}$, where $w_0$ is the longest Weyl element.  Their common stabilizers, taken over all fundamental representations, are $B_+$ and $B_-$, respectively.  The solution frames carry these chosen lines to lines of differential solutions, generated by
\[
 \begin{aligned}
 \Phi_0^{(i)}(z;u)
 &:=\mathcal Y_{\mathrm{as},V_i}(z;u)v_{\omega_i},\\
 \chi_{w_0\omega_i}(z;u)
 &:=\mathcal Y_{\mathrm F,V_i}(z;u)v_{w_0\omega_i}.
 \end{aligned}
\]
The first is the distinguished sectorial solution in $V_i$, and the second is the chosen lowest-weight Frobenius column.  Their weight labels refer to the coefficient vectors.  At fixed $u$, the compatible families of solution lines $(\mathbb C\Phi_0^{(i)})_i$ and $(\mathbb C\chi_{w_0\omega_i})_i$ realize the two Borel reductions just chosen.

The Frobenius solution above has sectorial coefficient
\[
 \mathcal Y_{\mathrm{as},V_i}(z;u)^{-1}
       \chi_{w_0\omega_i}(z;u)
 =\rho_i(C(u))^{-1}v_{w_0\omega_i}.
\]
Since $gB_-$ represents the family $(\rho_i(g)\mathbb Cv_{w_0\omega_i})_i$, the two reductions are represented in the sectorial trivialization by
\[
 F_+(u)=B_+\in G/B_+,
 \qquad F_-(u)=C(u)^{-1}B_-\in G/B_-.
\]
Here $B_\pm$ as points of $G/B_\pm$ denote the identity cosets $1B_\pm$.  The stabilizers of the flags $F_+(u)$ and $F_-(u)$ are $B_+$ and $C(u)^{-1}B_-C(u)$, respectively.  The Frobenius reduction is preserved because its representing flag map satisfies, using $Z\in H\subset B_-$,
\[
 \begin{aligned}
 \mathscr M(u)F_-(u)
 &=C(qu)^{-1}ZC(u)\bigl(C(u)^{-1}B_-\bigr)\\
 &=C(qu)^{-1}B_-=F_-(qu).
 \end{aligned}
\]
This supplies the preserved reduction required for a Miura structure.  The remaining assertion is that the sectorial flag map $F_+$ satisfies the $q$-oper condition.  We prove this by showing that the normalized sectorial frame puts $\mathscr M$ in \emph{Steinberg form}.  The sector at infinity determines a bipartite Coxeter word $c$. The normalization of the leading eigenvectors identifies the Coxeter factor of the transport from $u$ to $qu$ as $h_{\mathrm{as}}(u)\dot c$. Here $\dot c$ is the pinned lift of this word and $h_{\mathrm{as}}(u)\in H$ is determined by the increment of the scalar phase under rotation.
\Needspace{7\baselineskip}
\begin{theorem}\label{thm:intro-main}
For a simply-laced spectral family $\mathcal L_u$ satisfying Hypotheses~\ref{hyp:frob}--\ref{hyp:ccm}, the normalized Frobenius and sectorial frames define a meromorphic $G$-valued central connection matrix $C$.  In the normalized sectorial frame, $\mathscr M$ of \eqref{eq:intro-central-q-connection} satisfies
\[
 \mathscr M(u)\in U_c\,h_{\mathrm{as}}(u)\dot c,
 \qquad U_c=N_+\cap\dot cN_-\dot c^{-1}.
\]
The triple $(\mathscr M,F_+,F_-)$ represents a $Z$-twisted Miura $(G,q)$-oper, with sectorial oper flag $F_+$ and preserved Frobenius flag $F_-$.
\end{theorem}
The unipotent factor has one coordinate for each simple reflection in $c$ \cite{Steinberg1965}.  This form implies the required relative position of $F_+(qu)$ and $\mathscr M(u)F_+(u)$.

The Steinberg form follows by comparing the sectorial frame with $\mathsf S_1^*\mathcal Y_{\mathrm{as}}$ on an overlapping sector.  Asymptotic limits in two directions constrain the transition to the subgroup $U_c$.  This argument precedes the differential bilinear identities: the matrix $\mathscr M$ determines relations among the iterates of the distinguished solutions, whose expansion in the Frobenius basis gives the scalar functional relations.

In type $A_{N-1}$, the same $q$-connection also admits a cyclic realization for every $N\ge 2$. In the defining representation of $SL_N$, the pullbacks $\Phi_j^{(1)}:=\mathsf S_j^*\Phi_0^{(1)}$ of the distinguished sectorial solution all solve the differential equation at $u$. Under the same hypotheses the first $N$ of these solutions form a basis for generic $u$, giving the cyclic frame
\[
 \mathcal Y_{\mathrm{cyc}}(z;u)
   =\bigl(\Phi_0^{(1)}(z;u),\ldots,\Phi_{N-1}^{(1)}(z;u)\bigr).
\]
The spans of its first $a$ columns, for $1\le a<N$, form the cyclic oper flag. Its first line is the sectorial line $\mathbb C\Phi_0^{(1)}$, although its higher-dimensional subspaces can differ from those of the sectorial flag.

The pullback $\mathsf S_1^*$ takes each column to the next. Expressing $\Phi_N^{(1)}$ in this basis gives the last column of the companion matrix of the pullback. Its inverse represents the forward $q$-connection in the cyclic frame. Rescaling the last column gives a frame of determinant one. In this frame the forward matrix has Steinberg form for the linear Coxeter element $s_{N-1}\cdots s_1$.

\subsection{The universal \texorpdfstring{$(G,q)$}{(G,q)}-Wronskian and the full \texorpdfstring{$QQ$}{QQ}-system}\label{sec:intro-Wronskian-QQ}

The Steinberg form of $\mathscr M$ gives difference equations for the central connection itself.  The subgroup $U_c\subset N_+$ fixes each fundamental highest-weight vector $v_{\omega_i}$, so
\[
 \mathscr M(u)s_{\mathrm{as}}(u)^{-1}v_{\omega_i}=v_{\omega_i},
 \qquad s_{\mathrm{as}}(u):=h_{\mathrm{as}}(u)\dot c.
\]
Together with $C(qu)\mathscr M(u)=ZC(u)$, this gives the one-step $(G,q)$-Wronskian equations
\begin{equation}\label{eq:intro-one-step-Wronskian}
 Z^{-1}C(qu)v_{\omega_i}
 =C(u)s_{\mathrm{as}}(u)^{-1}v_{\omega_i},
 \qquad i\in I.
\end{equation}
A meromorphic section $\mathscr W:U\to G$ satisfying the same equations with $C$ replaced by $\mathscr W$ is called a generalized $(G,q)$-Wronskian for the prescribed twist and Coxeter data \cite{KZ2021}.  Thus the normalized central connection is such a $(G,q)$-Wronskian.

Its scalar connection coefficients are generalized minors.  For extremal weights $\lambda,\mu$ in $V_i$, the minor $\Delta_{\lambda,\mu}(g)$ is the coefficient of the chosen $\lambda$-weight vector in $\rho_i(g)v_\mu$, extracted by the corresponding dual covector.  For $w$ in the Weyl group,
\begin{equation}\label{eq:intro-Q-as-defined-coefficient}
 \mathcal Q_w^{(i)}(u)
 =\Delta_{w\omega_i,\omega_i}(C(u)).
\end{equation}
This is the coefficient of the Frobenius solution with leading weight $w\omega_i$ in the distinguished sectorial solution in $V_i$.

To obtain the scalar relation at node $i$, we apply the Pl\"ucker identity for $i$ and use the $(G,q)$-Wronskian equation to replace a column by its value at $qu$.  The bipartite order determines which neighboring minors acquire the same shift.  We use a bipartition $I=I_0\sqcup I_1$ of the ADE Dynkin diagram, with every edge joining vertices of different colors, and write $\epsilon_i\in\{0,1\}$ for the color of $i$.  The chosen Coxeter element is
\[
 c=c_1c_0,\qquad c_a:=\prod_{i\in I_a}s_i.
\]
Reflections of the same color commute.  Generalized minors satisfy the Pl\"ucker identity
\[
 \begin{aligned}
 &\Delta_{\omega_i,\omega_i}(g)
       \Delta_{s_i\omega_i,s_i\omega_i}(g)
 -\Delta_{\omega_i,s_i\omega_i}(g)
       \Delta_{s_i\omega_i,\omega_i}(g)\\
 &\hspace{25mm}=\prod_{j\sim i}\Delta_{\omega_j,\omega_j}(g),
 \qquad g\in G,
 \end{aligned}
\]
where $j\sim i$ denotes adjacency in the Dynkin diagram \cite{FominZelevinsky1999,KZ2021}.

We normalize the connection coefficients by $Q_w^{(i)}(u)=\exp\!\bigl[\pi i(\epsilon_i-\tfrac12)(w\omega_i)(\ell)\bigr]\mathcal Q_w^{(i)}(u)$.  With $Q_+^{(i)}:=Q_1^{(i)}$, $Q_-^{(i)}:=Q_{s_i}^{(i)}$, and $\mathfrak z_i:=e^{\pi i\alpha_i(\ell)}$, this substitution gives
\begin{equation}\label{eq:intro-QQ-system}
 \mathfrak z_i Q_+^{(i)}(u)Q_-^{(i)}(qu)
 -\mathfrak z_i^{-1}Q_+^{(i)}(qu)Q_-^{(i)}(u)
 =\varphi_i(u)\prod_{j\sim i}Q_+^{(j)}\!\left(q^{\epsilon_i}u\right).
\end{equation}
The fundamental weights $\omega_j:H\to\mathbb C^\times$ give the coefficient
\[
 \varphi_i(u)=\omega_i(h_{\mathrm{as}}(u))^{-1}
               \prod_{j\sim i}\omega_j(h_{\mathrm{as}}(u))^{\epsilon_i}.
\]
  The two products on the left are the two terms of a $2\times2$ determinant of connection coefficients at $u$ and $qu$; the product on the right is supplied by the generalized Pl\"ucker identity.  In rank one, the product is empty and the relation is the quantum Wronskian identity.  The Weyl translates of the same Pl\"ucker identity give the relations for pairs $w,ws_i$ for which multiplication by $s_i$ increases the Weyl-group length by one.  Together, these relations form the \emph{full $QQ$-system} \cite[Definition~3.18]{KZ2021}; equation~\eqref{eq:intro-QQ-system} is the case $w=1$.  Evaluating the fundamental relation at a zero of $Q_+^{(i)}$ and at its $q^{-1}$-shift, under the required regularity and nonvanishing conditions, and eliminating $Q_-^{(i)}$ gives the Bethe equation for that zero.

The full $QQ$-system and Bethe equations involve scalar coefficients of the fundamental columns $\mathscr W(u)v_{\omega_i}$.  Specifying these columns and requiring the one-step equations~\eqref{eq:intro-one-step-Wronskian}, with $C$ replaced by $\mathscr W$, need not determine the full $G$-valued section.  Indeed, right multiplication $\mathscr W\mapsto\mathscr W n$ by a meromorphic $n:U\to N_+$ leaves all these columns unchanged, since $n$ fixes every $v_{\omega_i}$.  Restricting to those $n$ for which $\mathscr W n$ still satisfies the one-step equations gives the fixed $(G,q)$-Wronskian orbit of Definition~\ref{def:Wronskian-right-orbit}.  The remaining problem is to recover the other columns from these data and further iterates of the differential solutions under shifted monodromy.

The recursive construction in Section~5 imposes further equations comparing the columns at $u$ and $q^mu$, as specified in~\eqref{eq:W-41-u}.  On the right-hand side, the product of $m$ inverse Coxeter factors has Weyl-group element $c^{-m}$.  For the chosen bipartite Coxeter word in types $D,E$, and also in type $A_{N-1}$ with even $N=h$, one has
\[
 c^{-d}=w_0,\qquad d=h/2.
\]
Thus the equation at $m=d$ relates the shifted highest-weight column to the lowest-weight column in each fundamental representation.  If two representatives in the fixed orbit satisfy all these shifted equations, their right-unipotent ratio must fix the lowest-weight vectors as well as the highest-weight vectors.  It therefore lies in $N_+\cap N_-=\{1\}$, giving uniqueness.  The \emph{recursive construction of the universal $(G,q)$-Wronskian} produces this representative in the fixed orbit by $d$ successive Bruhat factorizations.  These shifted equations characterize the \emph{universal} meromorphic $(G,q)$-Wronskian in the fixed orbit.  They also force the extremal minors $\Delta_{\omega_i,w_0\omega_i}(\mathscr W)$ and $\Delta_{w_0\omega_i,\omega_i}(\mathscr W)$ to be nonzero meromorphic functions for every $i$, placing $\mathscr W$ in the two longest Bruhat cells.

In type $A_{N-1}$, the cyclic frame of Section~1.2 gives an independent reconstruction for every $N$.  With $e_1$ the first standard basis vector, $\mathbf c(u)=C(u)e_1$ is the Frobenius coefficient column of the distinguished sectorial solution in the defining representation of $SL_N$.  Its iterates under shifted monodromy have Frobenius coefficient columns
\[
 \mathbf w_j(u)=Z^{-j}\mathbf c(q^ju),\qquad 0\le j<N.
\]
The analytic hypotheses imply that these form a basis over the field of meromorphic functions on $U$.  The matrix $P^Z=(\mathbf w_0,\ldots,\mathbf w_{N-1})$ gives the coefficients of this cyclic frame in the Frobenius basis.  We write $P$ for the coefficient matrix of the same frame in the sectorial basis.  Since $C$ takes sectorial coefficients to Frobenius coefficients,
\[
 P^Z=CP,\qquad C=P^ZP^{-1}.
\]
Applying the same last-column determinant normalization to $P^Z$ and $P$ leaves this quotient unchanged.  This reconstruction applies for every $N$, including odd $N=h$, since it uses the cyclic basis rather than the identity $c^{h/2}=w_0$.

The recursive and cyclic constructions give the following identification.
\Needspace{7\baselineskip}
\begin{theorem}\label{thm:intro-reconstruction}
Under Hypotheses~\ref{hyp:frob}--\ref{hyp:ccm}, if the Coxeter number $h$ is even, the central connection $C$ is the universal meromorphic $(G,q)$-Wronskian.  In each fundamental representation $V_i$, the coefficients of the normalized distinguished sectorial solution in the Frobenius frame form the column $C(u)v_{\omega_i}$. The $(G,q)$-Wronskian equations determine $C(u)$ uniquely from the values of these columns at $u,qu,\ldots,q^{h/2}u$.

For $G=SL_N$ with arbitrary $N$, $C$ is the unique isomorphism of $q$-connections identifying the sectorial and Frobenius realizations of the cyclic frame generated by the distinguished solution under shifted monodromy.

In every ADE type, the normalized fundamental highest-weight minors of $C$ satisfy the full $QQ$-system.
\end{theorem}

\subsection{\texorpdfstring{$T$}{T}-functions and \texorpdfstring{$TQ$}{TQ}-relations}\label{sec:intro-T-functions}

In this subsection we explain how the $Q$-functions of the central connection $C$ determine the coefficients of shifted monodromy.  By Theorem~\ref{thm:intro-main}, removing the prescribed Cartan factor gives an element of the Steinberg section,
\[
 S(u):=h_{\mathrm{as}}(u)^{-1}\mathscr M(u)\in U_c\dot c.
\]
Its principal minors in the fundamental representations $V_i$ are the \emph{fundamental $T$-functions},
\begin{equation}\label{eq:intro-T-coordinate}
 T_i(u):=\Delta_{\omega_i,\omega_i}\bigl(S(u)\bigr).
\end{equation}
These functions will give coordinates for $S(u)$ and hence recover $\mathscr M(u)=h_{\mathrm{as}}(u)S(u)$.

To express the $T$-functions through the $Q$-functions, we write the $q$-connection in \emph{Miura form}. This means choosing a trivialization in which the oper flag is $B_+$ and the preserved Frobenius flag is $B_-$ (the Borels are interchanged in \cite[Corollary~3.4]{FKSZ2020}). In the sectorial frame these flags are $B_+$ and $C(u)^{-1}B_-$, respectively. They are opposite for generic $u$ under our analytic hypotheses, so $C$ has a meromorphic Gauss factorization $C(u)=v_-(u)n(u)$ with $v_-(u)\in B_-$ and $n(u)\in N_+$. The change of coordinates $g\mapsto n(u)g$ fixes the oper flag and sends the Frobenius flag to $B_-$, since $n(u)C(u)^{-1}=v_-(u)^{-1}$. In these coordinates the map $P_u\to P_{qu}$ is left multiplication by
\[
 A_-(u)=n(qu)\mathscr M(u)n(u)^{-1}=v_-(qu)^{-1}Zv_-(u)\in B_- .
\]

For $b\in B_-$, its Cartan component $[b]_0$ is its image in $B_-/N_-=H$. The factor $n(u)$ fixes $v_{\omega_j}$, and the unipotent part of $v_-(u)$ adds only terms of lower weight. Thus the principal minor $\Delta_j=\Delta_{\omega_j,\omega_j}$ satisfies $\Delta_j(C(u))=\omega_j([v_-(u)]_0)$. Projecting the Miura form to $H$ therefore gives
\[
 \omega_j([A_-(u)]_0)
 =\omega_j(Z)\frac{\Delta_j(C(u))}{\Delta_j(C(qu))}.
\]

Forgetting the preserved reduction gives the underlying $q$-oper. In coordinates this is the \emph{$q$-Miura transformation}. After the prescribed Cartan factor has been removed, it expresses the principal minors of the Steinberg representative as Laurent polynomials in the Cartan coordinates of the connection in Miura form at successive spectral shifts. Frenkel and Reshetikhin conjectured that these polynomials are the fundamental $q$-characters~\cite{FR1998}. The $q$-character of a finite-dimensional representation of the quantum affine algebra records the simultaneous generalized eigenvalues of its commuting Cartan currents. Each eigenvalue is encoded by a Laurent monomial whose variables carry a Dynkin node and a spectral parameter. Its coefficient is the dimension of the corresponding generalized eigenspace. The fundamental $q$-character at node $i$ is the $q$-character of the fundamental representation at that node, with a specified spectral parameter.

Francone and Leclerc first proved this identification for all simply laced types by passing through the theory of cluster algebras and the generalized Baxter relations~\cite[Theorem~30]{FL2025}. In \cite{ZeitlinFRWronskians} we give a direct proof by identifying the variations of the preserved flag with the quantum screening operators. Incorporating a Dynkin diagram automorphism into the Miura transformation extends the argument to types $B$, $C$, $F_4$, and $G_2$.

The same character polynomials enter the Baxter relations for quantum integrable models.  In the construction of \cite{FH2015}, a finite-dimensional quantum-affine representation determines a transfer matrix, a family of operators depending on a spectral parameter.  These transfer matrices commute with one another and with the Baxter $Q$-operators.  On a common eigenvector, the latter give scalar $Q$-functions, and the transfer-matrix eigenvalue is obtained by evaluating the corresponding $q$-character on ratios of these functions at shifted parameters, with the scalar factors prescribed by the model.

The character formula uses the pinned lift $\dot c$ at each shift. Both standard flags are fixed by $H$, so a change of torus coordinates at $q^ku$ and $q^{k+1}u$ can remove the Cartan factor from $h_{\mathrm{as}}(q^ku)\dot c$. Starting at $u$, the torus elements $t_k=t_k(u)$, $k\in\mathbb Z$, are determined by
\[
 t_0=1,\qquad t_{k+1}=c(t_k)h_{\mathrm{as}}(q^ku)^{-1},
 \qquad c(t)=\dot c t\dot c^{-1}.
\]
Indeed, $t_{k+1}h_{\mathrm{as}}(q^ku)\dot c\,t_k^{-1}=\dot c$. In these coordinates the map from $q^ku$ to $q^{k+1}u$ is left multiplication by $A_k^{\mathrm M}=t_{k+1}A_-(q^ku)t_k^{-1}$. Evaluating the fundamental weight $\omega_j$ on the Cartan component and using $\Delta_j(C)=\omega_j([v_-]_0)$ gives the Miura variable at node $j$,
\begin{equation}\label{eq:intro-Miura-Q-ratio}
 \mathsf q_{j,k}(u):=\omega_j([A_k^{\mathrm M}]_0)
 =\omega_j(t_{k+1}Zt_k^{-1})\,
   \frac{\Delta_j(C(q^ku))}{\Delta_j(C(q^{k+1}u))}
 =\omega_j(t_{k+1}Zt_k^{-1})\,
   \frac{Q_+^{(j)}(q^ku)}{Q_+^{(j)}(q^{k+1}u)} .
\end{equation}
The normalization of $Q_+^{(j)}$ multiplies the principal minor by the same constant at both arguments. For each $j,k$ the resulting $\mathsf q_{j,k}$ is meromorphic on $U$. If $h_{\mathrm{as}}=1$, then $t_k=1$ for every $k$ and $\mathsf q_{j,k}(u)=\mathsf q_{j,0}(q^ku)$. In general these can differ, since starting the sequence at $q^ku$ imposes $t_0(q^ku)=1$ there.

For the fixed Coxeter word $c=s_{i_1}\cdots s_{i_r}$, where $r=\operatorname{rank}G$, the Miura form factors in these variables as
\[
 A_k^{\mathrm M}
 =P_{i_1}(\mathsf q_{i_1,k}(u))\cdots P_{i_r}(\mathsf q_{i_r,k}(u)),
 \qquad P_j(x)=\exp(x^{-1}f_j)\,\check\alpha_j(x).
\]
Here $f_j$ is the negative simple-root generator and $\check\alpha_j:\mathbb C^\times\to H$ is the corresponding coroot. At the original point $u$, the recurrence gives $A_-(u)=h_{\mathrm{as}}(u)A_0^{\mathrm M}$.

For the character formula, $q$ has infinite order and the quantum affine deformation parameter $\hbar$ satisfies $\hbar^{-2}=q$. The color $\epsilon_j\in\{0,1\}$ of vertex $j$ in the fixed bipartition determines the spectral label $b_{j,k}=\hbar^{2k+\epsilon_j}$. The $q$-character $\chi_i$ of the fundamental representation at node $i$ has highest monomial $Y_{i,b_{i,0}}$. The evaluation $\operatorname{ev}_C$ substitutes $\mathsf q_{j,k}(u)$ for $Y_{j,b_{j,k}}$ in this character. Since $\chi_i$ is a finite Laurent polynomial, its evaluation uses only finitely many spectral shifts.

\Needspace{10\baselineskip}
To recover $S(u)$ from these character evaluations, we use the root coordinates for the same Coxeter word. The corresponding roots are
\[
 \beta_a=s_{i_1}\cdots s_{i_{a-1}}(\alpha_{i_a}).
\]
For the root subgroup $U_{\beta_a}$, we choose the one-parameter map $X_a:\mathbb C\to U_{\beta_a}$ so that, with the fixed Coxeter lift and weight vectors,
\[
 \Delta_{\omega_{i_a},\omega_{i_a}}\bigl(X_a(x)\dot c\bigr)=x.
\]
The following theorem combines these root coordinates with the character formula to reconstruct shifted monodromy from the $Q$-functions.

\Needspace{14\baselineskip}
\begin{theorem}\label{thm:intro-TQ-relations}
Under Hypotheses~\ref{hyp:frob}--\ref{hyp:ccm}, the normalized element $S$ and the element $\mathscr M$ representing the original $q$-connection have the following parametrization.
\begin{enumerate}[label=(\alph*)]
\item The principal minors $(\Delta_{\omega_i,\omega_i})_{i\in I}$ give an algebraic isomorphism $U_c\dot c\simeq\mathbb A^r$.  The $T$-functions therefore determine $S$ uniquely, with
\begin{equation}\label{eq:intro-Steinberg-T-parameterization}
 S(u)=X_1\bigl(T_{i_1}(u)\bigr)\cdots
       X_r\bigl(T_{i_r}(u)\bigr)\dot c.
\end{equation}
The original connection is recovered by $\mathscr M(u)=h_{\mathrm{as}}(u)S(u)$.
\item If $q$ is not a root of unity, these coefficients are given by
\begin{equation}\label{eq:intro-TQ-parameterization}
 T_i(u)=\operatorname{ev}_C(\chi_i)(u).
\end{equation}
Consequently, each $T_i$ is a Laurent polynomial in the shifted $Q$-ratios~\eqref{eq:intro-Miura-Q-ratio}, with the prescribed twist and Cartan factors.
\end{enumerate}
\end{theorem}

For $G=SL_2$, the fundamental character is $Y_{1,1}+Y_{1,\hbar^2}^{-1}$, so its two monomials use the shifts $k=0,1$.  The pinned lift is $\dot c=\left(\begin{smallmatrix}0&-1\\1&0\end{smallmatrix}\right)$. We write $Z=\operatorname{diag}(\zeta,\zeta^{-1})$, $\phi(u):=\varphi_1(u)$, and $Q(u):=Q_+^{(1)}(u)$, where the superscript labels the unique Dynkin node.  Then
\[
 h_{\mathrm{as}}(u)=\operatorname{diag}(\phi(u)^{-1},\phi(u)),
 \qquad
 S(u)=\begin{pmatrix}T_1(u)&-1\\1&0\end{pmatrix},
\]
with
\[
 \begin{aligned}
 T_1(u)
 &=\mathsf q_{1,0}(u)+\mathsf q_{1,1}(u)^{-1}\\
 &=\zeta\phi(u)\frac{Q(u)}{Q(qu)}
   +\zeta^{-1}\frac{\phi(u)^2}{\phi(qu)}\frac{Q(q^2u)}{Q(qu)}.
 \end{aligned}
\]
Thus the coefficient for the step $u\mapsto qu$ is expressed through $Q(u)$, $Q(qu)$, and $Q(q^2u)$.

When the normalized shifted $Q$-ratios on a common eigenvector of a quantum model give the same assignment to the $q$-character variables as $\operatorname{ev}_C$, Theorem~\ref{thm:intro-TQ-relations} identifies the $T$-functions with the corresponding fundamental transfer-matrix eigenvalues.  This agreement includes the spectral parameters and the scalar factors from the model and our Cartan normalization.

\subsection{Feigin--Frenkel--Hernandez opers}\label{sec:intro-FFH-opers}

Theorem~\ref{thm:intro-main} constructs a $(G,q)$-oper over the spectral domain from a differential family in $z$.  Our principal application starts with affine differential opers.  We recall the differential notion here, on the $z$-curve: for a finite-dimensional group $G$, a differential $G$-oper consists of a principal $G$-bundle, a differential connection, and a Borel reduction.  Away from the singularities, in a local trivialization adapted to the reduction, the connection must have a nonzero component in every negative simple-root direction and no other negative-root component.  This condition describes how the flag changes under differential parallel transport.  The formalism goes back to Drinfeld and Sokolov \cite{DrinfeldSokolov1985}; its intrinsic definition is due to Beilinson and Drinfeld \cite{BeilinsonDrinfeld2005}. For the role of opers in quantum integrable systems, see also \cite{FrenkelGaudin2005}. For the Gaudin model, the critical-level construction of Bethe vectors in \cite{FeiginFrenkelReshetikhin1994} led to the description of Bethe solutions by differential Miura opers \cite{Frenkel2004}.

The affine extension proposed in \cite{FeiginFrenkel2011} was intended to encode the joint spectrum of quantum KdV Hamiltonians through affine Gaudin models.  Frenkel and Hernandez established the corresponding $QQ$ relations in quantum-affine representation theory and developed their connection with this description by opers \cite{FrenkelHernandez2018}.  We consider the resulting differential objects, called Feigin--Frenkel--Hernandez (FFH) opers.  They are affine Lie algebra-valued connections with prescribed singularities at zero and infinity.  For higher states, they also have finitely many additional regular singular points \cite{MRhigher2020,MR2024}.  These points are required to be \emph{apparent}: continuation of every local solution once around such a point returns the same solution.  The solutions may have poles there, but they have no branching or logarithmic terms.  Thus the local monodromy is the identity.  In rank one, the use of additional poles with trivial local monodromy to describe higher states already appears in the work of Bazhanov, Lukyanov, and Zamolodchikov \cite{BLZhigher2003}.

To apply Theorem~\ref{thm:intro-main}, we express the affine equation as a family of differential equations in $z$ parametrized by $u$.  The FFH equation has two realizations.  In the \emph{current realization}, the affine loop coordinate $\lambda_{\mathrm{aff}}$ is an independent variable, and the differential part is $z\partial_z+k\lambda_{\mathrm{aff}}\partial_{\lambda_{\mathrm{aff}}}$, where $0<k<1$ is a parameter of the oper.  The substitution $\lambda_{\mathrm{aff}}=z^ku$ gives the \emph{loop realization}, in which this differential part becomes $z\partial_z$ at fixed $u$.  In a finite-dimensional representation of $G$, the result is an ordinary first-order system with spectral parameter $u$.  Current differential monodromy at fixed $\lambda_{\mathrm{aff}}$ becomes shifted monodromy $\mathsf S_1^*$ in loop coordinates, with $q=e^{-2\pi ik}$.  In the loop realization, zero is regular singular and infinity has the cyclic irregular behavior described in Section~\ref{sec:intro-spectral-problems}.

The loop system also recovers the scalar equations discussed in Section~\ref{sec:intro-spectral-problems}.  In rank one, eliminating a component gives a second-order equation, which takes Schr\"odinger form after a scalar gauge transformation.  In type $A_{N-1}$, elimination in the defining $N$-dimensional representation gives a scalar equation of order $N$ \cite{FeiginFrenkel2011}.

We call an FFH oper \emph{generic} when its residue at zero is semisimple and the Frobenius coefficient recursions there are invertible in every fundamental representation; the precise conditions are given in Definition~\ref{def:generic-FFH-data}.  Section~\ref{sec:FFH-frames} verifies Hypotheses~\ref{hyp:frob}--\ref{hyp:ccm} for these untwisted ADE opers, using the analytic results of \cite{MRV2016,MRhigher2020,MR2024,CGM2025}.  The general construction therefore gives the following specialization.

\Needspace{7\baselineskip}
\begin{corollary}\label{cor:intro-FFH-specialization}
For $0<k<1$ and an untwisted simply-laced FFH oper generic in the preceding sense, Theorems~\ref{thm:intro-main} and~\ref{thm:intro-reconstruction} apply with $q=e^{-2\pi ik}$. The $q$-connection associated with its current differential monodromy has Steinberg form in the normalized sectorial frame:
\[
 \mathscr M(u)\in U_c\,h_{\mathrm{as}}(u)\dot c.
\]
\end{corollary}

For FFH opers, a separate calculation using only the normalization at infinity shows that, in type $A_{N-1}$, every block of $N$ consecutive iterates of the distinguished solution in the defining representation forms a cyclic frame whose entries are entire in the spectral parameter $u$; see Section~\ref{sec:FFH-typeA-basis}.  Its determinant is nonzero at every complex spectral value, and the calculation also applies when the Frobenius exponents are resonant.

\subsection{Relation to previous work}

Here we relate our construction to earlier work on differential equations, functional relations, and opers.  We begin with the massive ODE/IM correspondence, where connection coefficients and bases of consecutive solutions play a central role.  For modified affine Toda equations, this correspondence uses a flat linear system in $z$ and $\bar z$. For the modified sinh-Gordon equation, Lukyanov and Zamolodchikov identify the sine- and sinh-Gordon $Q$-functions with connection coefficients obtained by expanding the distinguished solution at infinity in the basis at zero. A combined angular rotation and spectral shift generates further solutions. In the basis of two consecutive solutions, the Stokes coefficients are bilinear expressions in the shifted $Q$-functions \cite[Sections~2--4]{LukyanovZamolodchikov2010}. A related rank-one construction is reviewed in \cite[Sections~4.4--4.6]{DDNT2020}, where matrices of companion form relate adjacent Stokes frames and determinants of pairs of subdominant solutions give the $T$-functions.

Consecutive rotated solutions also enter higher-rank linear problems.  The modified Tzitz\'eica--Bullough--Dodd equation, of twisted affine type $A_2^{(2)}$, has an $\mathfrak{sl}_3$-valued linear problem and a third-order scalar reduction. Three consecutive rotated solutions form a basis, and Wronskian identities give functional relations for the connection coefficients \cite[Sections~2 and~4]{DFNT2013}. For untwisted type $A$, the conformal limit recovers the massless scalar equations, and a Stokes coefficient is expressed as a sum of products of shifted $Q$-ratios \cite[Sections~3.2 and~4.2--4.3]{AdamopoulouDunning2014}. A Lie-theoretic formulation uses tensor identities among rotated subdominant solutions in fundamental representations to obtain $QQ$- and Bethe relations in the conformal limit. In this formulation, the general massive correspondence is stated as a conjecture \cite[Sections~3--5]{ItoLocke2015}.

For the affine opers considered here, we use the analytic existence results of \cite{MRV2016,MR2024,CGM2025} and the finite-dimensional identities between $\Lambda$-eigenvectors of \cite{MRV2016,MR2024}.  In this setting, full central connection matrices in the fundamental representations are defined in \cite{MR2024}, where the $QQ$-relations are derived from selected connection coefficients.  The further study of the full matrices and the general $TQ$-relations is left outside the scope of that work \cite{MR2024}.  In the related rank-one affine Gaudin setting, \cite{MasoeroMukhinRaimondo2025} constructs $Q$-functions from connection coefficients of differential opers depending on a spectral parameter and proves their $QQ$- and $TQ$-relations.

For the ADE differential families studied here, we construct the normalized $G$-valued frames and identify their central connection by the recursive and cyclic constructions of Theorem~\ref{thm:intro-reconstruction}.  In type $A$, the determinant formulas of Section~\ref{sec:typeA-differential-realization} use bases of consecutive solutions, with the iterates defined here by shifted monodromy.  The $q$-connection with matrix $\mathscr M$, together with the sectorial and Frobenius flags, defines a $Z$-twisted Miura $(G,q)$-oper.  The normalized differential solutions determine distinguished frames in which the $q$-connection matrices have Steinberg form.  Remark~\ref{rem:FFH-normalization-comparison} distinguishes the normalization factors from the Stokes transition needed to compare the full sectorial frames.

The group-theoretic normalization uses the centralizer torus of $\Lambda$ and its Coxeter action, described in \cite{Kostant1959,Steinberg1959}.  Group-valued sectorial solutions and Stokes factors are studied in \cite{Boalch2002}.  For the $tt^*$-Toda equations of topological--antitopological fusion, \cite{GuestHo2019} places canonical partial monodromy in a Steinberg cross-section.  Connections between $tt^*$ and ODE/IM also arise through the supersymmetric thermodynamic Bethe ansatz \cite{CecottiFendleyIntriligatorVafa1992} and its Airy realization \cite{Fendley1999}.  The conformal limit from $tt^*$-Toda connections to holomorphic differential systems is discussed in \cite{Guest2021}.

On the $q$-difference side, the Drinfeld--Sokolov reduction and its Miura transformation were developed in \cite{FrenkelReshetikhinSTS1998,SS1998}.  Earlier work associated additive difference operators, called discrete Miura opers, with solutions of rational spin-chain Bethe equations \cite{MukhinVarchenko2005}.  The intrinsic geometric definition of $(G,q)$-opers was introduced for $SL_N$ in \cite{KSZ2021} and for general simple simply connected $G$ in \cite{FKSZ2020}.  The $q$DE/IM correspondence formulated in \cite{FKSZ2020} relates $(G,q)$-opers to the $QQ$-system and Bethe equations.  That work also anticipated the broader $q$DE/ODE/IM correspondence, formulated more explicitly in \cite{KZ2021}, which proposed relating $(G,q)$-oper connections to differential monodromy.  The construction in this paper realizes this proposal in ADE under the analytic hypotheses.  The preprint \cite{ESV2020} discusses extended $Q$-systems and fused flags in simply-laced types within this $(G,q)$-oper framework.

Generalized $(G,q)$-Wronskians were introduced in \cite{KZ2021}.  For a general simple simply connected group $G$, \cite{KZ2021} derived the full $QQ$-system from identities between generalized minors of $(G,q)$-Wronskians and pointed out their relation to cluster mutations on double Bruhat cells.  This construction was subsequently used in \cite{GHL2024} to define quantum $(G,c)$-Wronskians over a Grothendieck ring. For rational input satisfying the hypotheses of our recursive construction, we recover the universal normalization of \cite{KZ2021}.  In type $A$, the cyclic construction gives the column matrices of \cite{KSZ2021,KZ2021}, after matching the twist, Coxeter lift, and diagonal normalization.

On the representation-theoretic side, the $q$-characters introduced in \cite{FR1998} encode the finite-dimensional representations entering the transfer matrices.  The infinite-dimensional representations of quantum-affine Borel algebras constructed in \cite{HernandezJimbo2012} provide the representation-theoretic setting for the Baxter operators in \cite{FH2015}.  The use of ODE/IM relations to study these Borel representations already appears in \cite{Sun2012}, which proposes matching differential connection coefficients with normalized $q$-characters and proves a determinant identity for the latter in type $A$.  The substitution formulas of Frenkel and Hernandez \cite{FH2015} express transfer eigenvalues through shifted Baxter functions.  Applying them to the differential $Q$-functions requires the normalization matching stated above.

The description of Bethe equations by opers also extends beyond the present differential construction.  Wronskian relations for differential Miura opers describe the Bethe equations of the inhomogeneous Gaudin model \cite{BrinsonSageZeitlin2024}, while toroidal $q$-opers encode Bethe equations for XXZ-type models with quantum toroidal symmetry \cite{KoroteevZeitlinToroidal2023}.  The broader relations between opers and quantum integrable systems include quantum/classical dualities for opers and their deformations, studied in \cite{KoroteevZeitlinZoo2024}.

The XXZ Bethe equations also link $q$-opers to the quantum equivariant $K$-theory of quiver varieties \cite{KSZ2021,KoroteevPushkarSmirnovZeitlin2021}, a relation reviewed in \cite[Sections~3 and~4]{ZeitlinBethe2024}.  This link has applications to 3d mirror symmetry, where twisted $q$-opers and $QQ$-systems enter the proof of self-duality for Hilbert schemes of points on $\mathbb C^2$ \cite{KoroteevZeitlinMirror2023}.

\subsection{Organization of the paper}

The introduction sets out most of the paper's principal notation and provides a map of its main results and the connections between them.

Section~2 fixes the Lie-theoretic conventions, $\Lambda$-eigenvectors, and tensor data.  Section~\ref{sec:analytic-input} states the five analytic hypotheses, derives their frame and tensor identities, proves Steinberg form in the sectorial frame, and verifies the hypotheses for FFH opers, including the normalization and sectorial-basis comparison for their transition matrices.  Section~\ref{sec:central-functional} identifies the two flag reductions of the Miura $(G,q)$-oper and derives the generalized-minor, $QQ$, and Bethe relations.

Section~5 constructs universal $(G,q)$-Wronskians by finitely many shifted equations when the Coxeter number is even.  Section~\ref{sec:consequences} explains the fundamental $q$-characters and evaluates them on ratios of the differential $Q$-functions. The resulting principal minors determine the Steinberg representative and hence the matrix of shifted monodromy.  Section~\ref{sec:ADE-central-connection} first develops the algebraic type-$A$ construction from successive pullbacks, including its Steinberg form and comparison with the recursive construction. It then applies these results to differential solutions, derives the scalar Baxter equation, and combines the recursive and cyclic constructions of the central connection in one ADE theorem with the $q$-character parametrization of shifted monodromy.  Section~\ref{sec:non-simply-laced-sketch} summarizes the non-simply-laced correspondence, emphasizing affine Langlands duality, fractional monodromy, and universal Wronskians with a diagram automorphism.  For even Coxeter number, the \hyperref[sec:appendix-distinguished-solutions]{Appendix} derives the shifted Wronskian equations from the distinguished solutions at infinity. These equations identify the central connection with the universal Wronskian and give the sectorial Steinberg form.

\section*{Acknowledgments}

The author is partially supported by NSF grant DMS-2526435 (formerly DMS-2203823) and by the Simons Foundation through its Travel Support for Mathematicians program.

\paragraph{Use of AI tools.}
The author developed the mathematical ideas and original proofs of this paper before using ChatGPT on the manuscript. During revision, the author used ChatGPT for language editing, revisions to the exposition, and help with organizing and checking proofs. The principal mathematical changes at that stage were a direct proof of Theorem~\ref{prop:analytic-sectorial-Steinberg} and refinements to the five analytic hypotheses. The \hyperref[sec:appendix-distinguished-solutions]{Appendix} retains the argument with distinguished solutions from the earlier proof and uses the resulting Wronskian equations to deduce Steinberg form.

The author also used ChatGPT to search the literature, look for simpler proofs of the background results in Section~2, and assist with \LaTeX{} typesetting. The author takes full responsibility for the mathematical content and final text.

\section{Lie-theoretic setup and tensor compatibility}

This section is purely expository.  We first fix the conventions for Bruhat factorization and generalized minors.  We then normalize the eigenvectors used in the sectorial frame of Section~\ref{sec:intro-spectral-problems} and record the resulting Coxeter action.  The section closes with a finite tensor criterion for assembling matrices in the fundamental representations into a single element of $G$.  Proofs are included where they help clarify the results, their normalizations, or their role in the subsequent constructions.

\subsection{Lie-theoretic data}\label{sec:lie-theoretic-data}
Throughout this section, $G$ is a connected simply-connected simple complex algebraic group of type $ADE$.  As in Section~\ref{sec:intro-qopers}, $B_\pm$ are opposite Borel subgroups, $H=B_+\cap B_-$ is their common maximal torus, and $N_\pm$ are their unipotent radicals.  Their Lie algebras are
\[
 \g=\operatorname{Lie}(G),\qquad
 \gh=\operatorname{Lie}(H),\qquad
 \n_\pm=\operatorname{Lie}(N_\pm),\qquad
 \mathfrak b_\pm=\operatorname{Lie}(B_\pm)=\gh\oplus\n_\pm.
\]  For $G=SL_N$, we take $B_+$ and $B_-$ to be the upper and lower triangular subgroups, respectively.  Then $H$ is the diagonal torus, and $N_\pm$ consist of the corresponding triangular matrices with diagonal entries equal to $1$.

All root, weight, and Weyl-group labels refer to this choice of $H$ and $B_+$, which fixes the based root datum.  We use $I=\{1,\dots,r\}$ to index the simple roots $\alpha_i$, the coroots $\check\alpha_i$, and the fundamental weights $\omega_i$.  We denote the set of positive roots by $\Delta_+$, the set of negative roots by $\Delta_-:=-\Delta_+$, and the root space of weight $\alpha$ by $\g_\alpha$.  In the real Cartan subspace
\[
 \gh_{\mathbb R}:=\operatorname{span}_{\mathbb R}\{\check\alpha_i:i\in I\},
\]
the open dominant chamber consists of the elements $x$ with $\alpha_i(x)>0$ for every $i$; its negative is the open antidominant chamber.  Since $G$ is simply-laced, its Cartan matrix is given by
\begin{equation}
 a_{ij}=2\delta_{ij}-B_{ij},
\end{equation}
where $B_{ij}\in\{0,1\}$ is the adjacency matrix.  The quotient
\[
 W:=N_G(H)/H
\]
is the Weyl group.  We denote by $s_i\in W$ the simple reflection corresponding to $\alpha_i$, by $w_0$ the longest Weyl element, and by $i\mapsto i^\star$ the opposition involution, so that $w_0\omega_i=-\omega_{i^\star}$.  We write $\ell_W:W\to\mathbb Z_{\ge0}$ for the length function and $h$ for the Coxeter number.

The unipotent factor in the Bruhat cell of $w\in W$ will be denoted by $U_w$.  For any representative $\dot w\in N_G(H)$, it is
\begin{equation}\label{eq:weighted-Steinberg-cells}
 U_w:=N_+\cap\dot wN_-\dot w^{-1}.
\end{equation}
This subgroup is independent of the representative, since $H$ normalizes $N_-$.

A \emph{pinning} of $G$ consists of a Borel subgroup $B_+$, a maximal torus $H\subset B_+$, and a nonzero vector $e_i\in\g_{\alpha_i}$ for each simple root $\alpha_i$.  The vector $e_i$ fixes the parametrization $t\mapsto\exp(te_i)$ of the corresponding root subgroup.  The opposite Chevalley generators $f_i\in\g_{-\alpha_i}$ are determined by $[e_i,f_j]=\delta_{ij}\check\alpha_i$.  We denote the highest root by $\theta$, its coroot by $\theta^\vee$, and the half-sum of the positive coroots by $\rho^\vee$.

For each root $\alpha$, $U_\alpha$ denotes the root subgroup with Lie algebra
$\g_\alpha$.  We use the following Bruhat and root-subgroup facts over every
field extension $F/\mathbb C$; see
\cite[Theorem~21.68, Remark~21.69, Proposition~21.77, Lemma~21.78 and
Theorem~21.80]{Milne2017}.
\begin{proposition}\label{prop:Bruhat-root-factorization}
\begin{enumerate}[label=(\alph*),topsep=3pt,itemsep=4pt]
\item For either choice of sign, the Bruhat decompositions are
\[
 G(F)=\bigsqcup_{w\in W}B_\pm(F)\dot wB_\pm(F)
     =\bigsqcup_{w\in W}B_\mp(F)\dot wN_\pm(F).
\]
If $\ell_W(vw)=\ell_W(v)+\ell_W(w)$, then
\[
 (B_+(F)\dot vB_+(F))(B_+(F)\dot wB_+(F))
 =B_+(F)\dot{vw}B_+(F).
\]
\item For a reduced word $w=s_{j_1}\cdots s_{j_m}$, the roots
\[
 \beta_a=s_{j_1}\cdots s_{j_{a-1}}\alpha_{j_a}
 \quad(1\le a\le m)
\]
are distinct and give precisely
$\{\alpha\in\Delta_+:w^{-1}\alpha\in\Delta_-\}$.  Multiplication induces
an isomorphism
\[
 U_{\beta_1}\times\cdots\times U_{\beta_m}\xrightarrow{\sim}U_w.
\]
In particular,
\[
 \operatorname{Lie}(U_w)=
 \bigoplus_{\substack{\alpha\in\Delta_+\\w^{-1}\alpha\in\Delta_-}}\g_\alpha,
 \qquad \dim U_w=\ell_W(w).
\]
For $V_w:=N_+\cap\dot wN_+\dot w^{-1}$, multiplication in either order
gives an isomorphism $U_w\times V_w\simeq N_+$ and
$V_w\times U_w\simeq N_+$.  Each factor is the ordered product of the
root subgroups it contains; any fixed order of those roots gives an
isomorphism of varieties.  Multiplication also gives
\[
 U_w\times B_+\xrightarrow{\sim}B_+\dot wB_+,
 \qquad (u,b)\longmapsto u\dot wb.
\]
All these isomorphisms commute with extension of $F$.
\item For every $w\in W$, multiplication gives the unique factorization
\[
 \dot wB_-\dot w^{-1}
 =H\,(N_+\cap\dot wN_-\dot w^{-1})
       (N_-\cap\dot wN_-\dot w^{-1}).
\]
The last two factors may be taken in the opposite order as well.
\end{enumerate}
\end{proposition}

We use the fundamental representations to describe the flags of Section~\ref{sec:intro-qopers} by their highest-weight lines and to detect the Bruhat cells by matrix coefficients.

For a dominant weight $\nu$, $V(\nu)$ denotes the irreducible representation with highest weight $\nu$.  If $V$ is a finite-dimensional representation, we denote its action by $\rho_V:G\to\operatorname{GL}(V)$.  For each $i$, the fundamental representation is $V_i=V(\omega_i)$, with action $\rho_i:=\rho_{V_i}$ and a chosen nonzero highest-weight vector $v_{\omega_i}\in V_i$.  Since $G$ is simply connected, the character lattice is $X^*(H)=\bigoplus_i\mathbb Z\omega_i$.
The following highest-weight realization of the flag variety is standard; see \cite[Sections~21(i) and~22(a)]{Milne2017}.
\begin{proposition}\label{prop:highest-weight-stabilizers}
The subgroup $P_i:=\operatorname{Stab}_G(\mathbb Cv_{\omega_i})$ is the maximal parabolic associated with the simple roots other than $\alpha_i$, and $\bigcap_iP_i=B_+$.  The map
\begin{equation}\label{eq:Borel-fundamental-embedding}
 G/B_+\hookrightarrow\prod_{i\in I}\mathbb P(V_i),\qquad
 gB_+\longmapsto\bigl([gv_{\omega_i}]\bigr)_{i\in I},
\end{equation}
is a closed embedding.  Replacing highest-weight lines by lowest-weight lines gives the corresponding closed embedding of $G/B_-$.
\end{proposition}

The weights in $W\omega_i$ are the \emph{extremal weights} of $V_i$.  Their weight spaces are one-dimensional.  To fix the signs of the generalized minors introduced in Section~\ref{sec:intro-Wronskian-QQ}, we choose extremal vectors and covectors using the inverse-Tits representatives determined by the pinning
\[
 \widehat s_i:=\exp(e_i)\exp(-f_i)\exp(e_i),
 \qquad
 \widehat w:=\widehat s_{j_1}\cdots\widehat s_{j_m}
 \quad\text{if }w=s_{j_1}\cdots s_{j_m}\text{ is reduced}.
\]
By the braid relations, $\widehat w$ does not depend on the reduced expression; see \cite[Section~1.4, equations~(1.8)--(1.9)]{FominZelevinsky1999}.  The covector $v_{\omega_i}^*$ takes the value $1$ on $v_{\omega_i}$ and vanishes on the other weight spaces.  The corresponding extremal vectors and covectors are
\[
 v_{w\omega_i}:=\rho_i(\widehat w)v_{\omega_i},\qquad
 v_{w\omega_i}^*:=v_{\omega_i}^*\rho_i(\widehat w^{-1}).
\]
The corresponding \emph{generalized minors} are the matrix coefficients
\[
 \Delta_{w\omega_i,\,v\omega_i}(g)
 :=v_{w\omega_i}^*\rho_i(g)v_{v\omega_i},
 \qquad w,v\in W.
\]
The chosen lifts of elements of the Weyl parabolic subgroup generated by $s_j$, $j\ne i$, fix $v_{\omega_i}$, since $e_jv_{\omega_i}=f_jv_{\omega_i}=0$ for these nodes.  Thus the indices of a minor specify weights rather than particular Weyl elements.  We use the displayed normalization throughout; in particular, $v_{s_i\omega_i}=-f_iv_{\omega_i}$ and $v_{w\omega_i}^*(v_{w\omega_i})=1$.  For $SL_N$, the fundamental representations are exterior powers of the defining representation, and these coefficients are ordinary matrix minors, with signs fixed by the extremal vectors.  For the general construction, see \cite[Definition~1.4]{FominZelevinsky1999}.

The following properties of generalized minors are standard; see \cite[Sections~1.4 and~2.3]{FominZelevinsky1999}.  We state them over a field extension $F$ of $\mathbb C$, with the groups and representations extended to $F$.
\begin{proposition}\label{prop:minor-covariance}
For every $i\in I$, $x\in G(F)$, $w,v\in W$, and $h_1,h_2\in H(F)$,
\begin{align*}
 \Delta_{w\omega_i,v\omega_i}(x)
 &=\Delta_{\omega_i,\omega_i}(\widehat w^{-1}x\widehat v),\\
 \Delta_{w\omega_i,v\omega_i}(h_1xh_2)
 &=(w\omega_i)(h_1)(v\omega_i)(h_2)
       \Delta_{w\omega_i,v\omega_i}(x).
\end{align*}
For $n_\pm\in N_\pm(F)$,
\[
 \Delta_{w\omega_i,\omega_i}(xn_+)
 =\Delta_{w\omega_i,\omega_i}(x),\qquad
 \Delta_{\omega_i,v\omega_i}(n_-x)
 =\Delta_{\omega_i,v\omega_i}(x).
\]
The quotient homomorphism $\pi_H:B_\pm(F)\to H(F)$ has kernel $N_\pm(F)$.  For $b\in B_\pm(F)$,
\[
 \Delta_{w\omega_i,w\omega_i}(b)=(w\omega_i)(\pi_H(b)).
\]
In particular, each diagonal extremal minor restricts to a multiplicative character on either Borel subgroup.
\end{proposition}

The Gauss-factorization criterion below is \cite[Corollary~2.5]{FominZelevinsky1999}; the two longest-cell criteria follow by translation by $\widehat w_0$.
\begin{proposition}\label{prop:Gauss-minor-criteria}
Multiplication is an isomorphism of varieties from $N_-\times H\times N_+$ onto the open subset $B_-N_+\subset G$, called the \emph{Gauss cell}.  For $x\in G(F)$,
\begin{align*}
 x\in B_-(F)N_+(F)
 &\Longleftrightarrow
 \Delta_{\omega_i,\omega_i}(x)\ne0\quad(i\in I),\\
 x\in B_+(F)\dot w_0B_+(F)
 &\Longleftrightarrow
 \Delta_{w_0\omega_i,\omega_i}(x)\ne0\quad(i\in I),\\
 x\in B_-(F)\dot w_0B_-(F)
 &\Longleftrightarrow
 \Delta_{\omega_i,w_0\omega_i}(x)\ne0\quad(i\in I),
\end{align*}
where $\dot w_0$ is any representative of $w_0$.  On the Gauss cell the factorization $x=n_-hn_+$ is unique, and
\[
 \Delta_{\omega_i,\omega_i}(x)=\omega_i(h).
\]
\end{proposition}
\begin{proof}
For $w=w_0$, Proposition~\ref{prop:Bruhat-root-factorization}(b) writes $B_+\dot w_0B_+$ as $N_+\dot w_0B_+$.  Left multiplication by $\dot w_0^{-1}$ gives the Gauss cell $N_-B_+=N_-HN_+$, with its unique factorization.  In its factorization, Proposition~\ref{prop:minor-covariance} gives $\Delta_{\omega_i,\omega_i}(x)=\omega_i(h)\ne0$.  Conversely, Proposition~\ref{prop:Bruhat-root-factorization}(a) gives $x=b_-\widehat w n_+$.  If $w\omega_i\ne\omega_i$, the difference $\omega_i-w\omega_i$ is a nonzero sum of simple roots with nonnegative coefficients.  Since $B_-$ can only lower weights, the corresponding principal minor is zero.  Nonvanishing of all principal minors forces $w$ to fix every fundamental weight, so $w=1$.  This proves the first equivalence.  Applying it to $\widehat w_0^{-1}x$ and to $x\widehat w_0$, and using Proposition~\ref{prop:minor-covariance}, gives the remaining two equivalences.
\end{proof}

The finite Weyl-group facts needed below are recorded together; see
\cite[Chapters~1 and~3]{Humphreys1990}.
\begin{proposition}\label{prop:finite-Weyl-data}
The symmetric Cartan matrix $2I-B$ is positive definite.  The Weyl group
acts simply transitively on the open Weyl chambers.  All Coxeter elements
are conjugate and have order $h$, and
\[
 \ell_W(w_0)=|\Delta_+|=rh/2.
\]
\end{proposition}

\subsubsection{The cyclic element and integer-indexed eigenvectors}
The sectorial normalization at infinity described in Section~\ref{sec:intro-spectral-problems} uses a regular semisimple element constructed from the principal nilpotent and a highest-root vector.  We use the principal nilpotent $\mathring f=\sum_i f_i\in\n_-$ and a chosen highest-root vector $v_\theta\in\g_\theta$.  For this choice of $\mathring f\in\n_-$, the principal grading element is $\mathbf h=-\rho^\vee$.  It satisfies
\begin{equation}\label{eq:principal-grading-relations}
 [\mathbf h,\mathring f]=\mathring f,
 \qquad
 [\mathbf h,v_\theta]=-(h-1)v_\theta.
\end{equation}
The cyclic element and its rotation factor are
\begin{equation}\label{eq:Kostant-cyclic-element}
 \Lambda:=\mathring f+v_\theta,
 \qquad
 \gamma_h:=e^{2\pi i/h}.
\end{equation}
The element $\Lambda$ is called a \emph{cyclic element}.  For $SL_N$, the matrix units are denoted by $E_{ab}$, and our choices are $f_i=E_{i+1,i}$ and $v_\theta=aE_{1N}$, $a\ne0$.  In the coordinate basis $e_1,\ldots,e_N$, we then have
\[
 \Lambda e_j=e_{j+1}\quad(1\le j<N),\qquad
 \Lambda e_N=ae_1,\qquad \Lambda^N=a\,\mathrm{id}.
\]
Thus $\Lambda$ carries the coordinate lines around one cycle; for $a=1$, it cyclically permutes the coordinate vectors themselves.  A Lie-algebra element is \emph{regular semisimple} if it is conjugate into a Cartan subalgebra and no root vanishes on the resulting Cartan element.  Its centralizer in $G$ is then a maximal torus.  Throughout, $\gamma_h^{a\mathbf h}$ means $\exp(2\pi ia\mathbf h/h)$.

\Needspace{12\baselineskip}
The following is Kostant's cyclic-element theorem in our choice of positive roots; see \cite[Lemma~6.3, Corollary~6.4, Section~6.7 and Corollary~8.6]{Kostant1959}.
\begin{theorem}\label{thm:Kostant-cyclic-Cartan}
For every nonzero $v_\theta\in\g_\theta$, the element $\Lambda=\mathring f+v_\theta$ is regular semisimple.  The centralizer and principal rotation are
\[
 H_\Lambda:=Z_G(\Lambda),\qquad Q:=\gamma_h^{-\mathbf h}.
\]
$H_\Lambda$ is a maximal torus normalized by $Q$, and the Weyl-group element represented by $Q$ is a Coxeter element.  Moreover, $\Lambda\in\operatorname{Lie}(H_\Lambda)$ satisfies
\[
 \operatorname{Ad}_Q\Lambda=\gamma_h^{-1}\Lambda.
\]
\end{theorem}

Conjugating $H$ to $H_\Lambda$ carries the usual weight decomposition to one that diagonalizes $\Lambda$ in every representation.  The principal rotation $Q$ permutes the eigenspaces as follows \cite[Lemma~3.1]{MRV2016}.  If $\Lambda\xi=\nu\xi$ in a representation $V$, the grading relations give
\begin{equation}\label{eq:rotated-cyclic-vector}
 \Lambda\rho_V(Q)^n\xi=\gamma_h^n\nu\,\rho_V(Q)^n\xi,
 \qquad n\in\mathbb Z.
\end{equation}
Indeed, $\operatorname{Ad}_{Q^n}\Lambda=\gamma_h^{-n}\Lambda$.

\subsubsection{The positive adjacency eigenvector}
The coefficients of a conjugate $\Lambda_H\in\gh$ of $\Lambda$ will come from the positive eigenvector of the Dynkin adjacency matrix $B$.  We recall the form of Perron--Frobenius used to obtain this vector.

For a square matrix $M$, its \emph{spectral radius} is
\[
 \rho(M):=\max\{ |\lambda|:\lambda\in\operatorname{Spec}(M)\}.
\]
A real symmetric matrix with nonnegative entries has an associated undirected graph: the vertices are its indices, and distinct vertices $i,j$ are joined when $M_{ij}>0$.  We use the following symmetric form of the Perron--Frobenius theorem; see \cite[Chapter~8]{HornJohnson1985}.
\begin{theorem}\label{thm:symmetric-Perron-Frobenius}
For a real symmetric matrix $M$ with nonnegative entries and connected associated graph, $\rho(M)$ is a simple eigenvalue of $M$.  It has an eigenvector with strictly positive coordinates, unique up to multiplication by a positive scalar.  If $M\ne0$, then $\rho(M)>0$.
\end{theorem}

The eigenvalue $\rho(M)$ is called the \emph{Perron eigenvalue}, and a corresponding positive eigenvector is called a \emph{Perron eigenvector}.  The adjacency matrix $B$ has connected graph because the Dynkin diagram is connected.  We normalize its Perron eigenvector $p=(p_i)_{i\in I}$ by $p_1=1$.  When the rank is greater than one, $B\ne0$ and its Perron eigenvalue is positive.  In type $A_1$, $B=[0]$ and we take $p_1=1$, with eigenvalue $0$.  We denote by $\kappa_h:=2\cos(\pi/h)$ the number that the next lemma identifies with the Perron eigenvalue of $B$.  The lemma also uses $p$ to construct a Cartan representative of $\Lambda$.

\subsubsection{The Cartan representative and its positive chamber}
The bipartition of the diagram introduced in Section~\ref{sec:intro-Wronskian-QQ} relates the eigenvalue equation for $B$ to the action of a Coxeter element.  We use this relation to specify the Cartan representative and its chamber.  Our bipartition is
\begin{equation}\label{eq:fixed-bipartition}
 I=I_0\sqcup I_1,
 \qquad
 \epsilon_i:=
 \begin{cases}
 0,&i\in I_0,\\
 1,&i\in I_1,
 \end{cases}
 \qquad
 \epsilon_j=1-\epsilon_i\quad\text{whenever }i\sim j,
\end{equation}
where $1\in I_0$.
The bipartite Coxeter factors and the chosen square root of $\gamma_h$ are
\[
 c_0:=\prod_{i\in I_0}s_i,\qquad
 c_1:=\prod_{i\in I_1}s_i,\qquad
 \zeta_h:=e^{\pi i/h}=\gamma_h^{1/2}.
\]
The reflections within either color commute.  We take real parts with respect to $\gh_{\mathbb R}$.  The bipartite expression for the Cartan element conjugate to $\Lambda$ in terms of the Perron eigenvector appears in \cite[equations~(3.5)--(3.7), Remark~3.2 and Lemma~3.3]{MRV2016}.  We use it to choose a conjugating element $a_0$ and a positive chamber in a real Cartan subspace of $\operatorname{Lie}(H_\Lambda)$.
\begin{lemma}\label{lem:bipartite-Coxeter-identification}
The Perron eigenvalue of $B$ is $\kappa_h$.  For its normalized Perron eigenvector $p$, the associated Cartan element is
\begin{equation}\label{eq:cyclic-Cartan-PF-vector}
 \Lambda_H
 :=\sum_{i\in I_0}p_i\check\alpha_i
   +\zeta_h\sum_{i\in I_1}p_i\check\alpha_i.
\end{equation}
The element $\Lambda_H$ is regular, $c_0c_1(\Lambda_H)=\gamma_h^{-1}\Lambda_H$, and $\operatorname{Re}(e^{-\pi i/(2h)}\Lambda_H)$ lies in the interior of the dominant chamber.  One can choose the scalar multiple of $v_\theta$ and an element $a_0\in G$ such that
\[
 \operatorname{Ad}_{a_0^{-1}}(\mathring f+v_\theta)=\Lambda_H,
 \qquad
 a_0^{-1}\gamma_h^{-\mathbf h}a_0\in N_G(H),
\]
with the latter element having Weyl image $c_0c_1$.
\end{lemma}
\begin{proof}
We first relate the spectra of $B$ and $c_0c_1$.  This determines the Perron eigenvalue and the line spanned by $\Lambda_H$.

In rank greater than one, the adjacency matrix in the simple-coroot basis ordered by $I_0,I_1$ has the form
\[
 B=\begin{pmatrix}0&X\\X^{\mathsf T}&0\end{pmatrix}.
\]
The formula $s_i(x)=x-\alpha_i(x)\check\alpha_i$ gives
\begin{equation}\label{eq:bipartite-reflection-blocks}
 [c_0]=\begin{pmatrix}-\mathbf1_{I_0}&X\\0&\mathbf1_{I_1}\end{pmatrix},
 \qquad
 [c_1]=\begin{pmatrix}\mathbf1_{I_0}&0\\X^{\mathsf T}&-\mathbf1_{I_1}\end{pmatrix}.
\end{equation}
Thus
\[
 [c_0c_1]=
 \begin{pmatrix}XX^{\mathsf T}-\mathbf1_{I_0}&-X\\
 X^{\mathsf T}&-\mathbf1_{I_1}\end{pmatrix}.
\]
These matrices give the following correspondence for every $t\in\mathbb C^\times$:
\begin{equation}\label{eq:bipartite-spectral-correspondence}
 B\binom{x}{y}=(t+t^{-1})\binom{x}{y}
 \quad\Longleftrightarrow\quad
 (c_0c_1)\binom{x}{ty}=t^{-2}\binom{x}{ty}.
\end{equation}
Indeed, both eigenvalue equations are equivalent to
$Xy=(t+t^{-1})x$ and $X^{\mathsf T}x=(t+t^{-1})y$.  Thus scaling the second block by $t$ identifies the two eigenspaces.

We first use an arbitrary nonzero highest-root vector $e_\theta\in\g_\theta$ and the cyclic element $\Lambda_0=\mathring f+e_\theta$.  After identifying the relevant Coxeter eigenspace, we will rescale $e_\theta$ to obtain the prescribed Cartan representative.  Theorem~\ref{thm:Kostant-cyclic-Cartan} applies to $\Lambda_0$: its centralizer is a maximal torus, $Q=\gamma_h^{-\mathbf h}$ induces a Coxeter element, and $\operatorname{Ad}_Q\Lambda_0=\gamma_h^{-1}\Lambda_0$.  Thus $\gamma_h^{-1}$ is an eigenvalue of a Coxeter element.  By Proposition~\ref{prop:finite-Weyl-data}, the Coxeter elements are conjugate and have order $h$.  In particular, $\gamma_h^{-1}$ is an eigenvalue of $c_0c_1$.  Taking $t=\zeta_h$ in \eqref{eq:bipartite-spectral-correspondence} shows that $\kappa_h=2\cos(\pi/h)$ is an eigenvalue of $B$, so $\rho(B)\ge\kappa_h$.

Conversely, with $r_B:=\rho(B)$, Perron--Frobenius gives $r_B>0$, and positive definiteness of the finite-type Cartan matrix $2I-B$ gives $r_B<2$.  Thus $r_B=2\cos\vartheta$ for $0<\vartheta<\pi/2$.  Applying \eqref{eq:bipartite-spectral-correspondence} to a Perron eigenvector shows that $e^{-2i\vartheta}$ is an eigenvalue of $c_0c_1$.  Since $(c_0c_1)^h=1$, we have $\vartheta=k\pi/h$ for an integer $k\ge1$.  Hence
\[
 r_B=2\cos(k\pi/h)\le2\cos(\pi/h)=\kappa_h.
\]
It follows that $\rho(B)=\kappa_h$.  Its eigenspace is one-dimensional by Perron--Frobenius, so \eqref{eq:bipartite-spectral-correspondence} also proves that the $\gamma_h^{-1}$-eigenspace of $c_0c_1$ is one-dimensional.

For the decomposition $p=(p^{(0)},p^{(1)})$ according to the two colors, the spectral correspondence~\eqref{eq:bipartite-spectral-correspondence} at $t=\zeta_h$ gives
\[
 [c_0c_1]\binom{p^{(0)}}{\zeta_hp^{(1)}}
 =\gamma_h^{-1}\binom{p^{(0)}}{\zeta_hp^{(1)}}.
\]
This is $c_0c_1(\Lambda_H)=\gamma_h^{-1}\Lambda_H$.  In type $A_1$, $B=0$, $h=2$, $p_1=1$, and $c_0c_1=s_1=-1$ on the one-dimensional Cartan, so the eigenvalue assertions, simplicity, and this identity hold directly.

Moreover,
\[
 \operatorname{Re}(e^{-\pi i/(2h)}\Lambda_H)
 =\cos\!\left(\frac{\pi}{2h}\right)\sum_i p_i\check\alpha_i,
\]
and every simple root takes on this element the positive value
\[
 (2-\kappa_h)\cos\!\left(\frac{\pi}{2h}\right)p_i.
\]
Every positive root is a nonzero nonnegative integral combination of simple roots.  Thus no root vanishes on $\Lambda_H$, proving regularity as well as the assertion about the chamber.

It remains to realize this Cartan element by a cyclic element with the prescribed principal nilpotent.  There is an element $b\in G$ carrying $H$ to $Z_G(\Lambda_0)$.  By replacing $b$ on the right by a representative of a Weyl element, we may arrange that $b^{-1}Qb$ has Weyl image $c_0c_1$.  Then $\operatorname{Ad}_{b^{-1}}\Lambda_0$ lies in the one-dimensional $\gamma_h^{-1}$-eigenspace spanned by $\Lambda_H$.  Consequently
\[
 \operatorname{Ad}_{b^{-1}}\Lambda_0=t\Lambda_H
 \quad\text{for some }t\in\mathbb C^\times.
\]
Scaling $\Lambda_0$ by $t^{-1}$ gives the desired Cartan representative.  A torus conjugation restores the coefficient of the principal nilpotent.  The definitions $z_*=t^{-1}$ and $g_*:=\exp((\log z_*)\rho^\vee)$ use an arbitrary logarithm.  Since $\alpha_i(\rho^\vee)=1$ and $\theta(\rho^\vee)=h-1$, we have
\[
 z_*\operatorname{Ad}_{g_*}\Lambda_0
 =\mathring f+z_*^h e_\theta.
\]
With $v_\theta=z_*^h e_\theta$ and $a_0=g_*b$, the preceding two equalities give $\operatorname{Ad}_{a_0^{-1}}\Lambda=\Lambda_H$.  Also $g_*$ commutes with $Q$, so $a_0^{-1}Qa_0=b^{-1}Qb$ has the required Weyl image.  Changing the logarithm multiplies $g_*$ by a central element and changes neither $v_\theta$ nor the adjoint identities.  This completes the normalization.
\end{proof}

We henceforth use this choice of $v_\theta$.  Conjugation by $a_0$ carries $\gh_{\mathbb R}$ and its dominant chamber to a real Cartan subspace of $\operatorname{Lie}(H_\Lambda)$ and its chosen positive chamber.  The corresponding Borel subgroup is
\[
 B_\Lambda:=a_0B_+a_0^{-1},
\]
which contains $H_\Lambda$.  In $V_i$, the line $\rho_i(a_0)\mathbb Cv_{\omega_i}$ is the highest-weight line for $B_\Lambda$ and is a $\Lambda$-eigenline.

\subsubsection{Fundamental Pl\"ucker morphisms}
To compare the wedge of two consecutive eigenvectors in $V_i$ with the tensor product of eigenvectors at the neighboring nodes, we use an equivariant map between these representations.  The common highest weight is
\begin{equation}\label{eq:eta-i}
 \eta_i:=2\omega_i-\alpha_i
 =\sum_{j\in I}B_{ij}\omega_j.
\end{equation}
The following highest-weight decomposition is
\cite[Lemmas~1.1--1.2 and equation~(1.4)]{MRV2016}.
\begin{proposition}\label{prop:fundamental-Plucker-morphisms}
The irreducible module $V(\eta_i)$ occurs with multiplicity one in both
$\bigwedge^2V_i$ and $\bigotimes_jV_j^{\otimes B_{ij}}$; every other
irreducible summand has strictly smaller highest weight.  There is a
unique $G$-morphism
\begin{equation}\label{eq:mi-finite}
 m_i:\bigwedge^2V_i\longrightarrow
       \bigotimes_jV_j^{\otimes B_{ij}}
\end{equation}
that vanishes on every other irreducible summand and satisfies
\begin{equation}\label{eq:mi-highest-normalization}
 m_i(f_iv_{\omega_i}\wedge v_{\omega_i})
   =\bigotimes_j(v_{\omega_j})^{\otimes B_{ij}}.
\end{equation}
\end{proposition}
\begin{proof}
For the exterior square, every weight $\mu\ne\omega_i$ of $V_i$ satisfies
\begin{equation}\label{eq:fundamental-weight-lowering}
 \omega_i-\mu=\sum_j n_j\alpha_j,
 \qquad n_j\geq0,\qquad n_i\geq1.
\end{equation}
Indeed, $V_i$ is spanned by words in the simple lowering operators applied to $v_{\omega_i}$.  In any nonzero word of positive length, the rightmost operator must be $f_i$, since $f_jv_{\omega_i}=0$ for $j\ne i$.  The highest-weight space of $V_i$ is one-dimensional, so at least one factor of a nonzero weight wedge has weight below $\omega_i$.  Equation~\eqref{eq:fundamental-weight-lowering} therefore bounds its weight by $2\omega_i-\alpha_i=\eta_i$.  The weight space at this bound is the line spanned by $f_iv_{\omega_i}\wedge v_{\omega_i}$, and every $e_j$ kills this vector.  This proves the assertion for $\bigwedge^2V_i$.  In the neighboring tensor product, the unique top weight is $\sum_jB_{ij}\omega_j=\eta_i$, and its weight space is the line spanned by the tensor of highest-weight vectors.  Decomposing the two modules into irreducible summands now gives both assertions.
Multiplicity one and the prescribed highest-weight vectors determine the normalized morphism $m_i$ uniquely.
\end{proof}

We now apply $m_i$ to the distinguished eigenlines and their images under $Q$.
For each $i$, we use $\nu_i:=\gamma_h^{\epsilon_i/2}p_i$ and a nonzero vector $\xi_i$ in the highest-weight line of $B_\Lambda$ in $V_i$.  The iterates under $Q$ and their wedge and tensor products are, for $n\in\mathbb Z$,
\begin{equation}\label{eq:integral-cyclic-vectors}
 \xi_{i,n}:=\rho_i(Q)^n\xi_i,
 \qquad
 \xi_{i,n}^{\wedge}:=\xi_{i,n}\wedge\xi_{i,n+1},
 \qquad
 \xi_{i,n}^{\otimes}:=\bigotimes_{j\in I}\xi_{j,n+\epsilon_i}^{\otimes B_{ij}}.
\end{equation}
All vectors $\xi_{i,n}$ belong to the fixed representation $V_i$; the second subscript records the number of applications of the principal rotation $Q$.
To state spectral maximality for $\mathring f+(-1)^{\epsilon_i}v_\theta$ with a positive real eigenvalue, we also use $\psi_i:=\gamma_h^{\epsilon_i\mathbf h/2}\xi_i$.
The eigenpairs and the wedge--tensor identity below were established in \cite[Proposition~3.4 and Remark~3.5]{MRV2016}; see also \cite[Theorem~2.5]{MR2024}.  We use the displayed order of the wedge factors.

\Needspace{6\baselineskip}
\begin{theorem}\label{thm:Kostant-PF-Plucker}
For each $i$, $\nu_i$ is a simple eigenvalue of $\Lambda$ on $V_i$, with
\begin{equation}\label{eq:index-zero-eigenpair}
 \Lambda\xi_i=\nu_i\xi_i.
\end{equation}
The vector $\psi_i$ is an eigenvector of $\mathring f+(-1)^{\epsilon_i}v_\theta$ with eigenvalue $p_i$, whose real part is strictly greater than that of every other eigenvalue.  The eigenvalues $\nu_i$ satisfy
\begin{equation}\label{eq:oriented-cyclic-eigenvalues}
 \sum_{j\in I}B_{ij}\nu_j
 =\gamma_h^{-\epsilon_i}(1+\gamma_h)\nu_i.
\end{equation}
The vectors $\xi_{i,n}^{\wedge}$ and $\xi_{i,n}^{\otimes}$ lie in the respective multiplicity-one summands isomorphic to $V(\eta_i)$ in $\bigwedge^2V_i$ and $\bigotimes_jV_j^{\otimes B_{ij}}$.  They have the same $\Lambda$-eigenvalue
\begin{equation}\label{eq:integral-edge-eigenvalue}
 \gamma_h^n(1+\gamma_h)\nu_i
 =\gamma_h^{n+\epsilon_i}\sum_jB_{ij}\nu_j.
\end{equation}
The vectors $\xi_i$ may be rescaled simultaneously so that
\begin{equation}\label{eq:normalized-cyclic-Plucker-vector}
 m_i(\xi_{i,n}^{\wedge})=\xi_{i,n}^{\otimes}
 \qquad(i\in I,\ n\in\mathbb Z).
\end{equation}
\end{theorem}
\begin{proof}
For $a_0$ as in Lemma~\ref{lem:bipartite-Coxeter-identification}, we first use $\xi_i=\rho_i(a_0)v_{\omega_i}$; any other choice on the same line differs by a nonzero scalar.  Since the fundamental weights are dual to the simple coroots,
\[
 \omega_i(\Lambda_H)=\zeta_h^{\epsilon_i}p_i.
\]
Conjugation by $a_0$ gives \eqref{eq:index-zero-eigenpair}, and the highest-weight lines of $B_\Lambda$ are precisely $\mathbb C\xi_i$.  Equation~\eqref{eq:rotated-cyclic-vector} then gives the eigenvalue $\gamma_h^n\nu_i$ for $\xi_{i,n}$.

We verify the asserted spectral maximality, which also proves simplicity.  The grading identity
\[
 \Lambda\gamma_h^{-\epsilon_i\mathbf h/2}
 =\gamma_h^{\epsilon_i/2}\gamma_h^{-\epsilon_i\mathbf h/2}
   \bigl(\mathring f+(-1)^{\epsilon_i}v_\theta\bigr)
\]
shows that $\mathring f+(-1)^{\epsilon_i}v_\theta$ is conjugate to $\zeta_h^{-\epsilon_i}\Lambda$.  It therefore suffices to compare the real parts of the weights evaluated on $\zeta_h^{-\epsilon_i}\Lambda_H$.

The normalized Perron eigenvector satisfies
\begin{equation}\label{eq:PF-cyclic-eigenvalues}
 \sum_{j\in I}B_{ij}p_j=\kappa_hp_i,
 \qquad p_1=1.
\end{equation}
Evaluating simple roots on \eqref{eq:cyclic-Cartan-PF-vector} and using \eqref{eq:PF-cyclic-eigenvalues} gives
\[
 \alpha_j(\Lambda_H)=
 \begin{cases}
 (1-\zeta_h^2)p_j,&j\in I_0,\\
 (\zeta_h-\zeta_h^{-1})p_j,&j\in I_1.
 \end{cases}
\]
For a fixed $i$, every $\operatorname{Re}(\zeta_h^{-\epsilon_i}\alpha_j(\Lambda_H))$ is nonnegative, and the value is strictly positive when $j=i$.  Applying \eqref{eq:fundamental-weight-lowering} to any weight $\mu\ne\omega_i$, we obtain
\[
 \operatorname{Re}\!\left(
 \zeta_h^{-\epsilon_i}(\omega_i-\mu)(\Lambda_H)\right)>0.
\]
Thus the highest-weight eigenvalue of $\zeta_h^{-\epsilon_i}\Lambda_H$ has strictly greatest real part.  Its eigenspace is the one-dimensional highest-weight space.  The distinguished eigenvector of $\mathring f+(-1)^{\epsilon_i}v_\theta$ is $\psi_i=\gamma_h^{\epsilon_i\mathbf h/2}\xi_i$, with eigenvalue $p_i>0$.  At node $1$ this value is $1$.  This proves the promised spectral normalization.

The eigenvalue on a wedge or tensor of eigenvectors is the sum of the eigenvalues of its factors.  The adjacency equation~\eqref{eq:PF-cyclic-eigenvalues} gives the equality of these sums:
\[
 \sum_j B_{ij}\nu_j
 =\zeta_h^{1-\epsilon_i}\kappa_hp_i
 =\gamma_h^{-\epsilon_i}(1+\gamma_h)\nu_i.
\]
Together with \eqref{eq:rotated-cyclic-vector}, it gives the common eigenvalue \eqref{eq:integral-edge-eigenvalue}.

We next prove that the two vectors lie in the distinguished copies of $V(\eta_i)$ and are related by $m_i$.  For this comparison, we use
\[
 Q=\gamma_h^{-\mathbf h},\qquad S=a_0^{-1}Qa_0.
\]
By Lemma~\ref{lem:bipartite-Coxeter-identification}, $S$ has Weyl image $c_0c_1$.  It suffices first to consider $n=0$ and apply $a_0^{-1}$ in all the factors.  The two vectors then become
\[
 v_{\omega_i}\wedge Sv_{\omega_i},
 \qquad
 \bigotimes_{j\sim i}S^{\epsilon_i}v_{\omega_j}.
\]
If $i\in I_0$, then $c_1\omega_i=\omega_i$ and $c_0\omega_i=s_i\omega_i$.  Hence $Sv_{\omega_i}$ is a nonzero multiple of $f_iv_{\omega_i}$, and the displayed tensor is the tensor of highest-weight vectors.  The two vectors therefore belong to their respective $V(\eta_i)$ summands, and \eqref{eq:mi-highest-normalization} relates them by a nonzero scalar.

For $i\in I_1$, we remove the $c_0$ factor from the Weyl image of $S$ by applying the same lift $D^{-1}$ to the wedge and to each tensor factor, where $D=\prod_{j\in I_0}\widehat s_j$ is a lift of $c_0$.  Each root subgroup with index in $I_0$ fixes $v_{\omega_i}$, so $D^{-1}v_{\omega_i}=v_{\omega_i}$.  The Weyl element of $D^{-1}S$ is $c_1$.  Consequently
\[
 D^{-1}Sv_{\omega_i}\in\mathbb C^\times f_iv_{\omega_i},
 \qquad
 D^{-1}Sv_{\omega_j}\in\mathbb C^\times v_{\omega_j}
 \quad(j\sim i).
\]
Here $c_1\omega_i=s_i\omega_i$, whereas $c_1\omega_j=\omega_j$ for every neighbor $j\in I_0$; extremal weight spaces are one-dimensional.  Applying the same $D^{-1}$ to the wedge and to every tensor factor reduces this case to the highest-weight identity as well.  Equivariance of $m_i$ now gives, for both colors,
\[
 m_i(\xi_i\wedge Q\xi_i)
 =c_i\bigotimes_{j\sim i}Q^{\epsilon_i}\xi_j,
 \qquad c_i\in\mathbb C^\times.
\]
The coefficient is nonzero because the highest-weight wedge is nonzero and $m_i$ is an isomorphism on its $V(\eta_i)$ summand.  Equivariance gives the same coefficient after applying $Q^n$ in every factor:
\[
 m_i(\xi_{i,n}^{\wedge})=c_i\xi_{i,n}^{\otimes}
 \qquad(n\in\mathbb Z).
\]
This proves the asserted summand membership for every $n$.

Finally, rescaling $\xi_i$ by $a_i\in\mathbb C^\times$ multiplies $c_i$ by $a_i^2\prod_j a_j^{-B_{ij}}$.  For logarithms $\ell_i$ of $c_i^{-1}$, the normalization is determined by the system
\[
 (2I-B)x=\ell.
\]
The Cartan matrix is nonsingular, so this system has a solution.  Taking $a_i=e^{x_i}$ makes every coefficient equal to $1$.  Equivalently, the monomial map with exponent matrix $2I-B$ is a surjective homomorphism of complex tori with finite kernel.  We fix one such normalization; any other differs by an element of this kernel.

In type $A_1$, the neighboring tensor product is $\mathbb C$, and the rescaling amounts to taking a square root of $c_1^{-1}$.
\end{proof}

\subsubsection{The normalized diagonalizer and the pinned Coxeter lift}
The normalization of the vectors $\xi_i$ determines how the principal rotation $Q$ acts in the eigenbasis of $\Lambda$.  Under the identification of $H_\Lambda$ with $H$ by $a_0$, the Weyl image of $Q$ is the inverse of the bipartite Coxeter element
\begin{equation}\label{eq:fixed-datum-bipartite-Coxeter}
 c:=c_1c_0,
 \qquad c^{-1}=c_0c_1.
\end{equation}
As in Section~\ref{sec:intro-Wronskian-QQ}, we order the vertices $i_1,\ldots,i_r$ with $I_1$ before $I_0$.  The pinning fixes the product
\begin{equation}\label{eq:coxeter}
 c=s_{i_1}\cdots s_{i_r},\qquad
 \dot c:=\widehat s_{i_1}^{-1}\cdots\widehat s_{i_r}^{-1}.
\end{equation}
The factors within each color commute. In the convention for extremal vectors above, $\dot c^{-1}=\widehat{c^{-1}}$. The next lemma identifies $\dot c$ in the normalized eigenbasis.
\begin{lemma}\label{lem:Kostant-full-transport}
For $a_0$ as in Lemma~\ref{lem:bipartite-Coxeter-identification} and the normalized vectors $\xi_i$ of Theorem~\ref{thm:Kostant-PF-Plucker}, there is a unique $t\in H$ such that $a_\Lambda=a_0t$ satisfies
\begin{equation}\label{eq:Kostant-diagonalizer}
 \rho_i(a_\Lambda)v_{\omega_i}=\xi_i,
 \qquad i\in I.
\end{equation}
This element satisfies $\operatorname{Ad}_{a_\Lambda^{-1}}\Lambda=\Lambda_H$ and carries $(H,B_+)$ to $(H_\Lambda,B_\Lambda)$.  In these coordinates the inverse rotation is the pinned lift of \eqref{eq:coxeter},
\begin{equation}\label{eq:Kostant-Coxeter-lift}
 a_\Lambda^{-1}Q^{-1}a_\Lambda
 =a_\Lambda^{-1}\gamma_h^{\mathbf h}a_\Lambda=\dot c,
\end{equation}
and its powers give the rotated vectors
\begin{equation}\label{eq:Kostant-Coxeter-orbit-vector}
 \rho_i(a_\Lambda\dot c^{-n})v_{\omega_i}
 =\xi_{i,n},
 \qquad n\in\mathbb Z.
\end{equation}
\end{lemma}
\begin{proof}
Since $\rho_i(a_0)v_{\omega_i}$ and $\xi_i$ span the same highest-weight line for $B_\Lambda$, we have $\rho_i(a_0)v_{\omega_i}=b_i\xi_i$ for $b_i\in\mathbb C^\times$.  Because $G$ is simply connected, the fundamental weights form a basis of $X^*(H)$, so there is a unique $t\in H$ with $\omega_i(t)=b_i^{-1}$ for all $i$.  Then $a_\Lambda=a_0t$ satisfies \eqref{eq:Kostant-diagonalizer}.  The torus $H$ centralizes $\Lambda_H$ and normalizes $B_+$, so right multiplication by $t$ preserves the Cartan and Borel identifications.

To identify the rotation, we first compute the action of the pinned product on the same tensors that normalize the $\xi_i$.  For $b=0,1$, write $D_b=\prod_{j\in I_b}\widehat s_j$, so that $\dot c^{-1}=D_0D_1$.  For $i\in I_0$, we have $\dot c^{-1}v_{\omega_i}=-f_iv_{\omega_i}$.  For $i\in I_1$, the element $D_0^{-1}$ fixes $v_{\omega_i}$ and sends $\dot c^{-1}v_{\omega_i}$ to $-f_iv_{\omega_i}$.  It also sends $\dot c^{-1}v_{\omega_j}$ to $v_{\omega_j}$ for every neighbor $j\in I_0$.  Equivariance and the highest-weight normalization \eqref{eq:mi-highest-normalization} therefore give, for either color,
\[
 m_i(v_{\omega_i}\wedge\dot c^{-1}v_{\omega_i})
   =\bigotimes_{j\sim i}\dot c^{-\epsilon_i}v_{\omega_j}.
\]

By Lemma~\ref{lem:bipartite-Coxeter-identification}, $a_\Lambda^{-1}Q^{-1}a_\Lambda$ has Weyl image $c$, so it equals $d\dot c$ for some $d\in H$.  Applying $a_\Lambda^{-1}$ to \eqref{eq:normalized-cyclic-Plucker-vector} at $n=0$ gives the preceding tensor identity with $d\dot c$ in place of $\dot c$.  Since $(d\dot c)^{-1}v_{\omega_i}=\omega_i(d)^{-1}\dot c^{-1}v_{\omega_i}$, comparison gives
\[
 \omega_i(d)=\prod_{j\sim i}\omega_j(d)^{\epsilon_i}.
\]
For $i\in I_0$ the right-hand side is $1$.  Every neighbor of a vertex in $I_1$ belongs to $I_0$, so $\omega_i(d)=1$ for this color as well.  The fundamental weights separate points of $H$, hence $d=1$.  This proves \eqref{eq:Kostant-Coxeter-lift}.  Its $n$th power gives $a_\Lambda\dot c^{-n}=Q^na_\Lambda$, and \eqref{eq:Kostant-diagonalizer} gives \eqref{eq:Kostant-Coxeter-orbit-vector}.
\end{proof}

For this Coxeter element, \eqref{eq:weighted-Steinberg-cells} gives
\begin{equation}\label{eq:fixed-unipotent-Uc}
 U_c=N_+\cap\dot cN_-\dot c^{-1},\qquad \dim U_c=r.
\end{equation}
The principal rotation is represented by $\dot c^{-1}$.  This is the direction of pullback from $qu$ to $u$ in Section~\ref{sec:intro-qopers}, so transport from $u$ to $qu$ uses $\dot c$.

For even $h$, the recursive construction of the universal $(G,q)$-Wronskian described in Section~\ref{sec:intro-Wronskian-QQ} uses the step at which the Coxeter action takes the highest-weight lines to the opposite extremal lines.  The next corollary records the classical half-Coxeter identity for the bipartite ordering; see \cite{Steinberg1959} and \cite[Section~3.17 and Section~3.19, Exercise~2]{Humphreys1990}.
\begin{corollary}\label{cor:fixed-datum-half-Coxeter}
For even $h$, the bipartite Coxeter element $c=c_1c_0$ satisfies
\begin{equation}\label{eq:fixed-datum-half-Coxeter}
 c^{h/2}=w_0.
\end{equation}
Consequently $\dot c^{h/2}$ sends every fundamental highest-weight line to its $w_0$-extremal line.
\end{corollary}
\begin{proof}
The chamber comparison uses
\[
 t_\Lambda:=\operatorname{Re}\!\left(e^{-\pi i/(2h)}\Lambda_H\right)
 =\cos\!\left(\frac{\pi}{2h}\right)\sum_{i\in I}p_i\check\alpha_i.
\]
By Lemma~\ref{lem:bipartite-Coxeter-identification}, $t_\Lambda$ lies in the interior of the dominant chamber.  Since
$c\Lambda_H=\gamma_h\Lambda_H$ and $h$ is even, we have
$c^{h/2}\Lambda_H=-\Lambda_H$.  The Weyl action is real on the coroot span, and taking the indicated real parts therefore gives
$c^{h/2}t_\Lambda=-t_\Lambda$.  Thus the half-power takes the dominant chamber to the antidominant chamber.  There is a unique Weyl element taking the dominant chamber to the antidominant chamber, namely $w_0$.  This proves the equality.  The lift $\dot c^{h/2}$ has Weyl image $w_0$, so it sends the line of weight $\omega_i$ to the line of weight $w_0\omega_i$.
\end{proof}

The morphism $m_i$ also gives the generalized-minor identity used to derive the full $QQ$-system in Section~\ref{sec:intro-Wronskian-QQ}.  The calculation uses the pinned representatives $\widehat w$ that fix the normalization of the minors.  The following is the fundamental identity of \cite[Theorem~1.17, equation~(1.26)]{FominZelevinsky1999}.  The proof below checks its coefficient in our extremal-vector normalization.
\begin{proposition}\label{prop:fundamental-Plucker-minors}
For $w,v\in W$ satisfying
\(
 \ell_W(ws_i)=\ell_W(w)+1
\)
and
\(
 \ell_W(vs_i)=\ell_W(v)+1,
\)
the fixed extremal-vector and covector normalizations give
\begin{align}\label{eq:fundamental-Plucker-minors}
 &\Delta_{w\omega_i,v\omega_i}(g)
  \Delta_{ws_i\omega_i,vs_i\omega_i}(g)
 -\Delta_{ws_i\omega_i,v\omega_i}(g)
  \Delta_{w\omega_i,vs_i\omega_i}(g)\notag\\
 &\hspace{38mm}=
 \prod_{j\ne i}
 \Delta_{w\omega_j,v\omega_j}(g)^{-a_{ji}}
 =\prod_{j\sim i}\Delta_{w\omega_j,v\omega_j}(g).
\end{align}
\end{proposition}
\begin{proof}
With the inverse-Tits representatives, we have $v_{s_i\omega_i}=-f_iv_{\omega_i}$.  Thus, in the prescribed order of wedge factors,
\[
 m_i(v_{\omega_i}\wedge v_{s_i\omega_i})
 =m_i(f_iv_{\omega_i}\wedge v_{\omega_i})
 =\bigotimes_{j\sim i}v_{\omega_j}.
\]
If $\ell_W(ws_i)=\ell_W(w)+1$, we have $\widehat{ws_i}=\widehat w\widehat s_i$.  We apply $\widehat w$ to both sides of the preceding highest-weight identity.  By $G$-equivariance, we obtain the corresponding identity between extremal vectors:
\begin{equation}\label{eq:extremal-Plucker-normalization}
 m_i(v_{w\omega_i}\wedge v_{ws_i\omega_i})
 =\bigotimes_{j\sim i}v_{w\omega_j}.
\end{equation}
The same equality holds with $v$ in place of $w$.

The wedge of the two normalized extremal covectors is the pullback under $m_i$ of the tensor product of the neighboring extremal covectors.  At $w=1$, the wedge covector and the pullback of the tensor-product covector are linear forms on $\bigwedge^2V_i$ of weight $-\eta_i$; both take value $1$ on $v_{\omega_i}\wedge v_{s_i\omega_i}$.  The weight space is one-dimensional.  Applying $\widehat w$ gives the equality at $w$.

Applying $g$ to $v_{v\omega_i}\wedge v_{vs_i\omega_i}$ and pairing with this covector gives the $2\times2$ determinant on the left-hand side of \eqref{eq:fundamental-Plucker-minors}.  Equivariance of $m_i$ and \eqref{eq:extremal-Plucker-normalization} express it as the product of the neighboring generalized minors.  Since $-a_{ji}=B_{ij}$ in simply-laced type, the product involves precisely the adjacent vertices.
\end{proof}

\Needspace{10\baselineskip}
\subsection{Invariant tensors and \texorpdfstring{$G$}{G}-valued frames}\label{sec:finite-tensor-data}
The analytic construction described in Section~\ref{sec:intro-spectral-problems} gives a solution frame in each fundamental representation.  To combine these frames into a single $G$-valued frame, we test membership in the image of $G$ in the faithful direct sum of the fundamental representations:
\[
 E=\bigoplus_{i\in I}V_i,
 \qquad \rho:=\rho_{E}:G\longrightarrow\operatorname{GL}(E).
\]

\begin{lemma}\label{lem:faithful-fundamental-sum}
The representation $\rho$ is faithful.
\end{lemma}
\begin{proof}
Its differential $d\rho:\g\to\mathfrak{gl}(E)$ is nonzero, since $f_iv_{\omega_i}\ne0$ in each fundamental representation.  Its kernel is an ideal of the simple Lie algebra $\g$, so $\ker(d\rho)=0$.  Thus the algebraic subgroup $K_0:=\ker\rho$ has zero-dimensional Lie algebra and is finite.  For $k\in K_0$, the map $G\to K_0$, $g\mapsto gkg^{-1}$, is constant because $G$ is connected.  Hence $K_0\subset Z(G)\subset H$.

If $k\in K_0$, its action on the fundamental highest-weight lines gives $\omega_i(k)=1$ for every $i$.  Since $G$ is simply connected, the fundamental weights form a basis of the character lattice $X^*(H)$, and characters separate points of $H$.  Therefore $k=1$; see \cite[22.8 and Proposition~23.59]{Milne2017} for these highest-weight facts.
\end{proof}

For nonnegative integers $p,q$, the notation
\[
 T^{p,q}(E):=E^{\otimes p}\otimes(E^*)^{\otimes q},
 \qquad
 T^{p,q}(A):=A^{\otimes p}\otimes((A^{-1})^*)^{\otimes q}
\]
denotes the mixed tensor space and the action induced by $A\in\operatorname{GL}(E)$.  Exponent zero means the trivial one-dimensional factor.

An \emph{algebraic representation} of an algebraic group on a finite-dimensional vector space $W$ is a group homomorphism to $\operatorname{GL}(W)$ that is a morphism of algebraic varieties.

To characterize $\rho(G)$ by invariant tensors, we first recall Chevalley's description of a closed algebraic subgroup as the stabilizer of a line.  The two representation-theoretic facts below turn this description into a finite tensor criterion for reductive groups; see \cite[Theorem~4.27]{Milne2017} and \cite[Proposition~3.1(a)--(c)]{Deligne1982}.
\begin{theorem}\label{thm:Chevalley-tensor-realization}
\begin{enumerate}[label=(\alph*),topsep=3pt,itemsep=3pt]
\item For every closed algebraic subgroup $G'\subset\operatorname{GL}(E)$, there is a finite-dimensional algebraic representation $P_0$ of $\operatorname{GL}(E)$ and a line $D\subset P_0$ such that
\[
 G'=\{g\in\operatorname{GL}(E):gD=D\}.
\]
This is an equality of closed algebraic subgroups.
\item Every finite-dimensional algebraic representation of $\operatorname{GL}(E)$ embeds equivariantly in a finite direct sum of mixed tensor spaces $T^{p,q}(E)$.
\item In a finite-dimensional algebraic representation of a reductive complex algebraic group, every invariant subspace has an invariant complement.
\end{enumerate}
\end{theorem}

The following theorem makes explicit the finite family of tensors needed for the invariant-tensor criterion of \cite[Proposition~3.1(c)]{Deligne1982}.

\begin{theorem}\label{thm:finite-invariant-tensors}
For a closed reductive complex algebraic subgroup $G'\subset\operatorname{GL}(E)$, there is a finite family of $G'$-invariant tensors $t_\kappa\in T^{p_\kappa,q_\kappa}(E)$, $1\le\kappa\le m$, whose common stabilizer is $G'$.  More precisely, for every field extension $F/\mathbb C$,
\[
 G'(F)=\{A\in\operatorname{GL}(E\otimes F):
 T^{p_\kappa,q_\kappa}(A)(t_\kappa\otimes1)=t_\kappa\otimes1\text{ for all }\kappa\}.
\]
\end{theorem}
\begin{proof}
By Theorem~\ref{thm:Chevalley-tensor-realization}(a)--(b), $G'$ is the stabilizer of a line $D$ in a finite direct sum $P=\bigoplus_{a=1}^{s}T^{p_a,q_a}(E)$.  Part~(c) gives a $G'$-invariant complement $Q$ to $D$ in $P$.  The projection $e\in\operatorname{End}(P)$ onto $D$ along $Q$ is therefore fixed by $G'$ under the induced conjugation action.  Conversely, an element of $\operatorname{GL}(E)$ fixing $e$ preserves its image $D$ and therefore belongs to $G'$.

The natural identification $\operatorname{End}(P)\simeq P\otimes P^*$ and permutation of tensor factors give the $\operatorname{GL}(E)$-equivariant decomposition
\[
 \operatorname{End}(P)\simeq
 \bigoplus_{a,b=1}^{s}T^{p_a+q_b,\,q_a+p_b}(E).
\]
The components of $e$ in these finitely many summands are $G'$-invariant tensors $t_\kappa$.  Since every summand is preserved by $\operatorname{GL}(E)$, fixing all $t_\kappa$ is equivalent to fixing $e$.

For a field extension $F/\mathbb C$, the decomposition $P\otimes F=(D\otimes F)\oplus(Q\otimes F)$ is preserved by $G'(F)$.  Thus $G'(F)$ fixes $e\otimes1$ and its components $t_\kappa\otimes1$.  Conversely, any $A\in\operatorname{GL}(E\otimes F)$ fixing these components fixes $e\otimes1$ and therefore preserves its image $D\otimes F$.  The algebraic subgroup identity in Theorem~\ref{thm:Chevalley-tensor-realization}(a), after extension to $F$, then gives $A\in G'(F)$.
\end{proof}

For $G'=\rho(G)$, we choose once and for all a family $t_\kappa\in T^{p_\kappa,q_\kappa}(E)$ supplied by Theorem~\ref{thm:finite-invariant-tensors}.  Its members will be called the \emph{defining invariant tensors} of $G$ in $E$.  They give a finite test for membership in $\rho(G)$, and a matrix satisfying this test fixes every $\rho(G)$-invariant tensor.  The criterion is
\begin{equation}\label{eq:finite-defining-tensors}
 \rho(G)=\{A\in\operatorname{GL}(E):T^{p_\kappa,q_\kappa}(A)t_\kappa=t_\kappa\text{ for all }\kappa\}.
\end{equation}
The tensors themselves are fixed, rather than merely their spans.

\begin{proposition}\label{prop:finite-tensor-reconstruction}
For a field extension $F/\mathbb C$ and matrices $Y_i\in\operatorname{GL}(V_i\otimes F)$, $i\in I$, their direct sum is $Y_E:=\bigoplus_iY_i$.  There is a unique $g\in G(F)$ with $Y_i=\rho_i(g)$ for every $i$ if and only if
\[
 T^{p_\kappa,q_\kappa}(Y_E)(t_\kappa\otimes1)=t_\kappa\otimes1
 \qquad\text{for every }\kappa.
\]
\end{proposition}
\begin{proof}
Theorem~\ref{thm:finite-invariant-tensors} identifies the matrices satisfying these identities with $\rho(G)(F)$.  Faithfulness makes $\rho:G\to\rho(G)$ an isomorphism of algebraic groups, so the unique preimage $g$ has the stated fundamental matrices.  Conversely, every $\rho(g)$ fixes the defining tensors.
\end{proof}

For every other finite-dimensional algebraic representation $V$, the induced matrix is $Y_V:=\rho_V(g)$.  These matrices respect direct sums, tensor products and duals, and satisfy $Y_Wf=fY_V$ for every $G$-equivariant map $f:V\to W$.  In particular, their exterior-power and neighboring tensor-product actions commute with the maps $m_i$.

The following analytic form applies the criterion to meromorphic solution frames.
\begin{corollary}\label{cor:finite-tensor-analytic}
If the fundamental matrices $Y_i$ are meromorphic on a complex manifold $\Omega$, their determinants are not identically zero on any connected component, and $Y_E=\bigoplus_iY_i$ fixes every defining tensor, there is a unique meromorphic $G$-valued family $g$ with $Y_i=\rho_i(g)$ for every $i$.  Wherever $Y_E$ is holomorphic and invertible, $g$ is holomorphic.
\end{corollary}
\begin{proof}
Proposition~\ref{prop:finite-tensor-reconstruction} applies over the field of meromorphic functions on each connected component of $\Omega$.  On an open subset where $Y_E$ is holomorphic and invertible, its values lie in $\rho(G)$.  The inverse of the closed algebraic embedding $\rho:G\to\rho(G)$ is regular, so composition gives holomorphic $g$ there.
\end{proof}

\Needspace{14\baselineskip}
\section{Abstract analytic hypotheses and their verification}\label{sec:analytic-input}
In this section we establish the analytic input for the construction of $(G,q)$-opers from differential equations.  We first formulate five hypotheses on the normalized solutions at zero and infinity.  These allow us to construct the central connection and the $q$-connection associated with shifted monodromy, and to prove that the latter has Steinberg form in the normalized sectorial frame.  We then verify the hypotheses for untwisted simply-laced FFH opers under explicit genericity conditions on their parameters, derive determinant formulas for consecutive solutions in type $A$, and obtain the tensor identities used below.

\subsection{The differential family}\label{sec:general-hypotheses}\label{sec:general-analytic-setting}
Here we specify the differential family and the conditions under which its normalized solution frames at zero and infinity can be compared.  The spectral shift will then act on the solution spaces described by these frames.

For a fixed $q\in\mathbb C^\times$, the spectral parameter $u$ ranges over a connected domain $U\subset\mathbb C^\times$ invariant under $u\mapsto q^{\pm1}u$.  The spectral difference field is
\[
 K:=\mathcal M(U),\qquad \sigma(f)(u):=f(qu).
\]
Here $\mathcal M(U)$ is the field of meromorphic functions on $U$.  A \emph{difference field} is a field equipped with a specified automorphism; ours is $(K,\sigma)$.  Elements of $K$ are independent of the differential variable $z$.

The differential variable lies on the universal cover
\[
 \widetilde{\mathbb C^\times}
 =\{(r,\theta):r>0,\ \theta\in\mathbb R\},
 \qquad \Pi(r,\theta)=re^{i\theta}.
\]
We write $z$ for a point of this cover, with $|z|=r$, $\arg z=\theta$, and $\log z=\log r+i\theta$.  The function $\log z$ is a global holomorphic coordinate, and $z^\alpha=e^{\alpha\log z}$.  In differential and algebraic expressions, $z$ also denotes the projected coordinate $\Pi(z)$; thus $z\partial_z=\partial_{\log z}$.  Coefficients and solutions are functions on the cover and need not descend to $\mathbb C^\times$.

On a connected open subset $\mathscr X$ of this cover, our family of connections on the trivial principal $G$-bundle has the following expression in a fixed trivialization:
\begin{equation}\label{eq:general-connection}
 \nabla_u=d_z+\mathscr A(z;u)\,dz,\qquad
 \mathcal L_u=\partial_z+\mathscr A(z;u),\qquad
 \mathscr A\in\g\otimes\mathcal O(\mathscr X\times U).
\end{equation}
The notation $\mathcal O$ means holomorphic functions.  The set $\mathscr X$, fixed for all $u\in U$ and all representations $V$, is the \emph{common differential domain}.  It excludes the finite poles of the connection.  Continued solution frames are compared on this same set, in the fixed original trivialization.  We require $\mathscr X$ to be invariant under $\theta\mapsto\theta+2\pi n$ and to contain the chosen sectors of the forms
\[
 0<|z|<r_0,\quad a<\arg z<b,
 \qquad\text{and}\qquad
 |z|>R_0,\quad a'<\arg z<b'.
\]
Limits at zero and infinity mean $|z|\to0$ and $|z|\to\infty$ within these sectors, with $u$ fixed.  The coefficient $\mathscr A$ is holomorphic on $\mathscr X\times U$, although the connection may be singular at its boundary and at the excluded finite points.

In a representation $V$, a \emph{horizontal section}, or solution, is a function $y$ satisfying
\[
 \nabla_u y=0,
 \qquad\text{equivalently}\qquad
 \partial_zy+\rho_V(\mathscr A)y=0.
\]
A solution frame $Y_V$ is an invertible matrix with solution columns.  A holomorphic differential gauge $T$, with values in $G$, acts by
\[
 Y_V^T=\rho_V(T)Y_V,\qquad
 \mathscr A^T=T\mathscr A T^{-1}-(\partial_zT)T^{-1}.
\]
The transformed differential operator is $\mathcal L_u^T=T\mathcal L_uT^{-1}=\partial_z+\mathscr A^T$.
Holomorphic dependence on $z$ and $u$ means joint holomorphy on the stated product domain.  An ordinary point of the differential equation is one where the connection coefficient is holomorphic; every point of $\mathscr X$ has this property.  All vector and matrix norms below are fixed norms on the indicated finite-dimensional spaces.

\paragraph{Shifted monodromy in a local gauge.}
The maps $\mathsf S_n$ of Section~\ref{sec:intro-spectral-problems} have the explicit form
\begin{equation}\label{eq:sym-action}
 \mathsf S_n(z,u)=(e^{2\pi in}z,q^nu).
\end{equation}
The first entry denotes the deck translation
\[
 e^{2\pi in}z=(r,\theta+2\pi n),\qquad
 \Pi(e^{2\pi in}z)=\Pi(z).
\]
Pullback gives
\begin{equation}\label{eq:forward-Symanzik-map}
 \mathsf S_n^*Y=Y\circ\mathsf S_n,\qquad
 \mathsf S_n^*(fY)=f(q^nu)\mathsf S_n^*Y\quad(f\in K).
\end{equation}
The second identity expresses \emph{semilinearity} with respect to $\sigma^n$.  The two gauge factors in a pullback are evaluated at different points.  The expression in a local gauge $T$ therefore involves the comparison factor $\Gamma_{T,n}$ below.
\begin{equation}\label{eq:local-gauge-covariant-pullback}
 \begin{aligned}
  \Gamma_{T,n}(z;u)&:=T(z;u)T(e^{2\pi in}z;q^nu)^{-1},\\
  \rho_V(T)\,\mathsf S_n^*\bigl(\rho_V(T)^{-1}Y_V^T\bigr)
   &=\rho_V(\Gamma_{T,n})\,\mathsf S_n^*Y_V^T
     =\rho_V(T)\,\mathsf S_n^*Y_V.
 \end{aligned}
\end{equation}
Thus pullback is first taken in the original trivialization and then expressed in the gauge $T$, wherever both gauge factors are defined.  For a constant gauge $T$, the factor $\Gamma_{T,n}$ is the identity.

\begingroup
\flushbottom
\par\vspace{0pt plus \baselineskip}
\Needspace*{18\baselineskip}
\endgroup
\subsection{Analytic hypotheses for the differential family}\label{sec:five-analytic-hypotheses}
Here we state the five conditions used in the analytic construction.  We first specify the local frames at zero and infinity, then impose compatibility of their normalizations with shifted pullbacks.  The last condition allows us to compare the frames on a common domain.

\par\addvspace{\baselineskip}
\begingroup
\centering
\small
\renewcommand{\arraystretch}{1.25}
\begin{tabular}{@{}c@{\quad}p{.29\linewidth}@{\qquad}p{.60\linewidth}@{}}
\hline
 & \textbf{Hypothesis} & \textbf{Role in the construction} \\
\hline
\ref{hyp:frob} & Frobenius normalization & The prescribed expansions select a $G$-valued solution frame at zero. \\[3pt]
\ref{hyp:subdominant} & Sectorial normalization & The asymptotics select a $G$-valued solution frame on a sector at infinity. \\[3pt]
\ref{hyp:symanzik} & Covariance at zero & A full turn with $u\mapsto qu$ preserves the equation and determines the Cartan twist in the Frobenius frame. \\[3pt]
\ref{hyp:common-orbit} & Phase and gauge increments at infinity & Comparison on overlapping sectors determines the Cartan factor and the Coxeter lift in the sectorial frame. \\[3pt]
\ref{hyp:ccm} & Common continuation & The two frames extend to the same domain, where their comparison gives the central connection matrix. \\
\hline
\end{tabular}
\par
\endgroup
\par\addvspace{\baselineskip}

The normalization at zero is specified by convergent Frobenius expansions in a chosen gauge.  The gauge is defined on a punctured sector; a meromorphic extension through zero is not assumed.  Prescribing the leading vector fixes the scalar normalization of each column.

\begin{definition}\label{def:Frobenius-classes}
A prescribed class of Frobenius expansions specifies the allowed exponents of $z$, the allowed powers of $\log z$, the restrictions on the coefficients, and convergence on the chosen sector.  Nonintegral exponents are allowed, and all branches are fixed.  The class is preserved by multiplication by a nonzero scalar.

On a sector $\Sigma$ at zero and a nonempty open set $U_0\subset U$, \emph{Frobenius normalization data} consist of a gauge $T_0\in\mathcal O(\Sigma\times U_0,G)$, an element $\ell\in\gh$, and, for each fundamental representation $V=V_i$, prescribed classes of Frobenius expansions indexed by an ordered weight basis $(v_{\mu,a})$.  At an extremal weight, $v_{\mu,1}=v_\mu$ with the normalization fixed in Section~2.  A column $\chi_{\mu,a}$ has the prescribed normalization if its gauged column belongs to the prescribed class for $\mathcal L_u^{T_0}$ with label $(\mu,a)$ and
\begin{equation}\label{eq:frob-general}
 \chi_{\mu,a}(z;u)
   =\rho_V(T_0(z;u))^{-1}z^{-\mu(\ell)}
       \bigl(v_{\mu,a}+r_{\mu,a}(z;u)\bigr),
 \qquad r_{\mu,a}=o(1)\quad(z\to0).
\end{equation}
The matrix of these columns is denoted by $\mathcal Y_{\mathrm F,V}$:
\[
 \mathcal Y_{\mathrm F,V}(z;u)v_{\mu,a}=\chi_{\mu,a}(z;u).
\qedhere
\]
\end{definition}

\begin{hypothesis}\label{hyp:frob}
There are Frobenius normalization data as in Definition~\ref{def:Frobenius-classes} on a sector $\Sigma$ at zero and a nonempty open set $U_0\subset U$.  For every fundamental representation, fixed $u\in U_0$ and label $(\mu,a)$, there is a unique solution of $\mathcal L_u$ with the prescribed normalization.  The ordered solution columns form a holomorphic frame on $\Sigma\times U_0$.  Their direct sum $\mathcal Y_{\mathrm F,E}:=\bigoplus_i\mathcal Y_{\mathrm F,V_i}$ satisfies
\[
 T^{p_\kappa,q_\kappa}(\mathcal Y_{\mathrm F,E})t_\kappa=t_\kappa
 \qquad\text{for every defining tensor }t_\kappa.
\]
\end{hypothesis}

The following criterion verifies the tensor-preservation clause of the Frobenius hypothesis.  We will use it for FFH opers in Section~\ref{sec:FFH-hyp-frob}.

\Needspace{12\baselineskip}
\begin{proposition}\label{prop:Frobenius-tensor-prescription}\label{rem:Frobenius-tensor-prescription}
We assume that the prescribed Frobenius columns of Definition~\ref{def:Frobenius-classes} form holomorphic frames on $\Sigma\times U_0$.  For each $u\in U_0$, let $\mathcal R_u$ be a unital $\mathbb C$-algebra of scalar holomorphic germs at zero on $\Sigma$ with convergent Frobenius expansions and finite limits at zero.  The ideal of germs tending to zero is denoted by $\mathfrak m_u$.  We require that $\Psi(f_1,\ldots,f_d)\in\mathcal R_u$ whenever $f_1,\ldots,f_d\in\mathfrak m_u$ and $\Psi$ is a holomorphic germ at $0\in\mathbb C^d$.

We also assume that for each fundamental representation $V_i$ and each prescribed basis vector $v\in V_i$ of weight $\mu$, the column labelled by $v$ is the unique horizontal germ $y$ at this $u$ satisfying
\[
 z^{\mu(\ell)}\rho_i(T_0)y-v\in\mathfrak m_u\otimes_{\mathbb C} V_i .
\]
The same algebra $\mathcal R_u$ is used for all columns and all $i$.  Then the fundamental Frobenius frames are induced by a single holomorphic $G$-valued frame on $\Sigma\times U_0$.  In particular, they satisfy the tensor-preservation clause of Hypothesis~\ref{hyp:frob}.
\end{proposition}
\begin{proof}
We work first at a fixed $u$.  By Lemma~\ref{lem:faithful-fundamental-sum}, the fundamental representations identify $G$ with a closed subgroup of $L:=\prod_i GL(V_i)$.  We use this identification also for $\g\subset\operatorname{Lie}(L)$.  On the chosen branch of $\log z$, the tuple of prefactors $F_0=(F_{0,i})_i$ is given by
\[
 F_{0,i}=\rho_i(T_0)\mathcal Y_{\mathrm F,V_i}\,z^{\rho_i(\ell)}.
\]
The normalization of the columns puts every entry of $F_{0,i}-I_{V_i}$ in $\mathfrak m_u$.  Differentiating the gauged solution frames gives the equation
\[
 z\partial_z F_0=-(z\mathscr A^{T_0})F_0+F_0\ell .
\]

To compare $F_0$ with a $G$-valued prefactor, we construct a holomorphic projection near the identity of $L$.  By Theorem~\ref{thm:Chevalley-tensor-realization}(c), $\g$ has an $\operatorname{Ad}(G)$-invariant complement $\mathfrak s$ in $\operatorname{Lie}(L)$.  The differential of $(g,X)\mapsto g\exp X$ at $(1,0)$ is the isomorphism $(\xi,X)\mapsto\xi+X$.  The holomorphic inverse function theorem therefore gives a unique local factorization $p=g\exp X$ with $g\in G$ near $1$ and $X\in\mathfrak s$ near $0$.  The projection $\pi(p):=g$ satisfies $\pi(apb)=a\pi(p)b$ locally for $a,b\in G$ near $1$, since
\[
 a(g\exp X)b=(agb)\exp(\operatorname{Ad}_{b^{-1}}X).
\]
Differentiation in $a$ and $b$ gives
\[
 d\pi_p(\xi p+p\eta)=\xi\pi(p)+\pi(p)\eta,
 \qquad \xi,\eta\in\g .
\]

Since $F_0\to1$, the germ $\widetilde F_0:=\pi(F_0)$ is defined.  The differential identity with $\xi=-z\mathscr A^{T_0}$ and $\eta=\ell$ shows that $\widetilde F_0$ satisfies the same prefactor equation.  Each entry of $\widetilde F_0-I$ is a holomorphic function of the entries of $F_0-I$ with zero constant term, so it belongs to $\mathfrak m_u$.  For a basis vector $v$ of weight $\mu$ in $V_i$, the column
\[
 \widetilde y_v=\rho_i(T_0)^{-1}\widetilde F_{0,i}\,z^{-\rho_i(\ell)}v
\]
is thus horizontal and satisfies $z^{\mu(\ell)}\rho_i(T_0)\widetilde y_v-v=(\widetilde F_{0,i}-I_{V_i})v\in\mathfrak m_u\otimes V_i$.  Uniqueness identifies it with the prescribed column.  This holds on each fundamental basis, hence $F_0=\widetilde F_0\in G$ as germs.

The tuple of original frames is $T_0^{-1}F_0z^{-\ell}$, so it fixes every defining tensor near zero.  At each $u$ these tensor identities extend by holomorphy to the connected sector $\Sigma$.  By \eqref{eq:finite-defining-tensors}, the jointly holomorphic tuple consequently takes values in the embedded copy of $G$ on $\Sigma\times U_0$, as asserted.
\end{proof}

At infinity we use the diagonalizer $a_\Lambda$ of Lemma~\ref{lem:Kostant-full-transport}, for which $a_\Lambda^{-1}Q^{-1}a_\Lambda=\dot c$.  In each fundamental representation $V=V_i$, we choose an ordered $H$-weight basis $(e_b)$ containing the normalized vector $v_{\omega_i}$, with a basis in every weight space, including those of higher dimension.  The vectors $\xi_b=\rho_V(a_\Lambda)e_b$ form the ordered eigenbasis of $\rho_V(\Lambda)$ used below, with eigenvalues $\nu_b$.  We write
\[
 D:=\operatorname{Ad}_{a_\Lambda^{-1}}\Lambda\in\gh,
 \qquad \operatorname{Ad}_{\dot c}D=\gamma_hD.
\]
  The prescribed exponential and power factors of the gauged frame are included in $e^{-S\rho_V(D)}$, with the powers of $z$ coming from logarithmic terms in the scalar phase $S$.  The normalization below is a limit for the remaining prefactor.

\begin{definition}\label{def:leading-sectorial-normalization}
A leading sectorial normalization consists of a holomorphic $G$-gauge $T_\infty$, a fixed holomorphic phase $S$, and numbers
\[
 -\pi/h<a<b<0.
\]
For every relatively compact open set $U_0\subset U$, the data are defined on a tail of
\[
 \Sigma_\infty:=\{z:a-\pi<\arg z/h<b\}
\]
times $U_0$, with compatible choices on overlaps.  A solution frame $Y_V$ has this normalization if its $GL(V)$-valued \emph{normalized prefactor} tends to the identity:
\begin{equation}\label{eq:general-sectorial-columns}
 \begin{aligned}
 F_V(z;u)&:=\rho_V(a_\Lambda^{-1}T_\infty(z;u))Y_V(z;u)
                e^{S(z;u)\rho_V(D)},\\
 F_V(z;u)&\longrightarrow I_V
 \qquad(|z|\to\infty,\ z\in\Sigma_\infty).
 \end{aligned}
\end{equation}
The convergence is uniform on closed angular subsectors and compact spectral sets: for every closed interval $J\subset(a-\pi,b)$ and compact set $L\subset U$,
\[
 \lim_{r\to\infty}
 \sup_{\substack{|z|\ge r,\ \arg z/h\in J\\u\in L}}
 \left\|F_V(z;u)-I_V\right\|=0.
\]
The tail radius may depend on $J,L$.  Branches and the additive constant of $S$ are fixed throughout.
\end{definition}

\begin{hypothesis}\label{hyp:subdominant}
There are leading sectorial normalization data as above for which the following asymptotic statements hold as $|z|\to\infty$ in $\Sigma_\infty$, uniformly on closed angular subsectors and compact spectral sets:
\begin{equation}\label{eq:general-irregular-normal-form}
 \mathcal L_u^{T_\infty}=\partial_z+p(z;u)\Lambda+R(z;u),
 \qquad p=\partial_zS,\quad \|R\|=o(|p|),
\end{equation}
and
\begin{equation}\label{eq:abstract-principal-asymptotic-type}
 p(z;u)z^{1-1/h}\longrightarrow1,\qquad
 S(z;u)=hz^{1/h}+o(z^{1/h}).
\end{equation}
In every fundamental representation there is a jointly holomorphic solution frame on the sectorial tails specified in Definition~\ref{def:leading-sectorial-normalization}, with normalization~\eqref{eq:general-sectorial-columns} and compatible choices on overlaps.
\end{hypothesis}

\begin{remark}\label{rem:leading-normalization-scope}
The coefficient estimates and the existence of the normalized solution frame are separate requirements of Hypothesis~\ref{hyp:subdominant}.  The normalization is the leading limit~\eqref{eq:general-sectorial-columns}, without a complete asymptotic expansion of its prefactor.
\end{remark}

The angular width of $\Sigma_\infty$ makes the leading limit strong enough to prove uniqueness and tensor preservation.

\begin{lemma}\label{lem:leading-sectorial-tensors}
The frames of Hypothesis~\ref{hyp:subdominant} are unique.  Their direct sum in $E$ preserves the defining tensors of $G$.
\end{lemma}
\begin{proof}
Throughout the proof, $u$ is fixed and $D_V=\rho_V(D)$.  For two frames $Y_V$ and $\widehat Y_V$ with the prescribed normalization, the change-of-basis matrix is
\[
 C_V:=Y_V^{-1}\widehat Y_V,\qquad \widehat Y_V=Y_VC_V.
\]
Both frames solve $\partial_zY=-\rho_V(\mathscr A)Y$.  Consequently,
\[
 \partial_zC_V
 =Y_V^{-1}\rho_V(\mathscr A)\widehat Y_V
    -Y_V^{-1}\rho_V(\mathscr A)\widehat Y_V=0.
\]
Thus $C_V$ may depend on $u$, but not on $z$.  The normalized prefactors $F_V$ and $\widehat F_V$ of $Y_V$ and $\widehat Y_V$, defined by~\eqref{eq:general-sectorial-columns}, both tend to $I_V$.  Therefore
\[
 e^{-SD_V}C_Ve^{SD_V}=F_V^{-1}\widehat F_V\longrightarrow I_V.
\]
In the fixed eigenbasis of $D_V$, this says
\[
 e^{S(z;u)(\nu_j-\nu_i)}(C_V)_{ij}\longrightarrow\delta_{ij}.
\]
If $\nu_i=\nu_j$, then $(C_V)_{ij}=\delta_{ij}$.  If $\nu_i\ne\nu_j$, the interval $(a-\pi,b)$ has length greater than $\pi$, so it contains an angle $\theta$ for which
$\operatorname{Re}(e^{i\theta}(\nu_j-\nu_i))>0$.
Along the ray $z=re^{ih\theta}$ on the cover,
\[
 S(z;u)=hr^{1/h}e^{i\theta}+o(r^{1/h}),\qquad r\to\infty.
\]
The exponential multiplier then grows in modulus, whereas the displayed matrix entry must tend to zero.  Since $(C_V)_{ij}$ is independent of $z$, it is zero.  Hence $C_V=I_V$, proving uniqueness, also in the presence of repeated eigenvalues.

For tensor preservation, we apply the same exponential growth argument to the constant coordinates of an invariant tensor in the induced solution frame.  The direct-sum frame on $E$ is $Y_E=\bigoplus_iY_{V_i}$.  For a defining tensor $t_\kappa$, we abbreviate its mixed tensor space and the induced action by
\[
 W_\kappa:=E^{\otimes p_\kappa}\otimes(E^*)^{\otimes q_\kappa},\qquad
 \tau_\kappa:=T^{p_\kappa,q_\kappa}:GL(E)\longrightarrow GL(W_\kappa).
\]
This action is defined for every invertible matrix; no membership of $Y_E$ in $G$ is assumed.  The induced coefficient, constant leading operator, and solution frame are
\[
 \mathscr A_\kappa:=d\tau_\kappa(\rho_E(\mathscr A)),\qquad
 D_\kappa:=d\tau_\kappa(\rho_E(D)),\qquad
 Y_\kappa:=\tau_\kappa(Y_E),
\]
so $\partial_zY_\kappa=-\mathscr A_\kappa Y_\kappa$.  Since $t_\kappa$ is $G$-invariant and $\mathscr A,D\in\g$,
\[
 \mathscr A_\kappa t_\kappa=0,\qquad D_\kappa t_\kappa=0.
\]
Thus the vector $c_\kappa:=Y_\kappa^{-1}t_\kappa$ satisfies
\[
 \partial_zc_\kappa=Y_\kappa^{-1}\mathscr A_\kappa t_\kappa=0.
\]
The normalized prefactor on $E$ is $F_E:=\bigoplus_iF_{V_i}$.  Its induced tensor action satisfies
\[
 \tau_\kappa(F_E)
 =\tau_\kappa\!\left(\rho_E(a_\Lambda^{-1}T_\infty)\right)
       Y_\kappa e^{SD_\kappa}\longrightarrow I_{W_\kappa}.
\]
The first factor fixes $t_\kappa$, because $a_\Lambda^{-1}T_\infty$ takes values in $G$.  Hence
\[
 \begin{aligned}
 \tau_\kappa(F_E)^{-1}
 &=e^{-SD_\kappa}Y_\kappa^{-1}
       \tau_\kappa\!\left(\rho_E(a_\Lambda^{-1}T_\infty)\right)^{-1},\\
 \tau_\kappa(F_E)^{-1}t_\kappa
 &=e^{-SD_\kappa}Y_\kappa^{-1}t_\kappa
  =e^{-SD_\kappa}c_\kappa.
 \end{aligned}
\]
Since $\tau_\kappa(F_E)\to I_{W_\kappa}$, its inverse also tends to $I_{W_\kappa}$, and therefore
\[
 e^{-SD_\kappa}c_\kappa\longrightarrow t_\kappa.
\]
Since $D\in\gh$, the induced operator $D_\kappa$ is semisimple.  If $P_\lambda$ denotes its projection onto the $\lambda$-eigenspace, then
\[
 e^{-\lambda S(z;u)}P_\lambda c_\kappa
       \longrightarrow P_\lambda t_\kappa.
\]
Since $D_\kappa t_\kappa=0$, the tensor $t_\kappa$ lies in the zero eigenspace.  The limit for $\lambda=0$ therefore gives
\[
 P_0c_\kappa=P_0t_\kappa=t_\kappa.
\]  For $\lambda\ne0$, one has $P_\lambda t_\kappa=0$, and the interval $(a-\pi,b)$ contains an angle $\theta$ with $\operatorname{Re}(e^{i\theta}\lambda)<0$.  Along $z=re^{ih\theta}$, the modulus of $e^{-\lambda S(z;u)}$ tends to infinity.  The last limit thus forces $P_\lambda c_\kappa=0$.  Summing the eigenspace components gives
\[
 c_\kappa=t_\kappa,\qquad\text{hence}\qquad
 T^{p_\kappa,q_\kappa}(Y_E)t_\kappa=t_\kappa.
\]
This holds for every defining tensor and every $u$ in the parameter domain.
\end{proof}

The two local frames use different ordered bases for their columns.  We distinguish their coordinate spaces so that we can express the change of coordinates after continuation.  To identify the distinguished column at infinity, we recall the eigenpair from Theorem~\ref{thm:Kostant-PF-Plucker}:
\[
 \xi_i=\rho_i(a_\Lambda)v_{\omega_i},\qquad
 \rho_i(\Lambda)\xi_i=\nu_i\xi_i,
\]
where $\nu_i$ is a simple eigenvalue.  This eigenpair specifies the leading term of the distinguished column in the sectorial frame.

\begin{definition}\label{def:normalized-local-frames}\label{def:reference-columns}
For a fundamental representation $V=V_i$, the spaces $V^{\mathrm F}$ and $V^{\mathrm{as}}$ are labelled copies of $V$, with ordered bases $(v_{\mu,a})$ and $(e_b)$.  The \emph{normalized Frobenius frame} of Hypothesis~\ref{hyp:frob} is regarded as a map $\mathcal Y_{\mathrm F,V}(z;u):V^{\mathrm F}\xrightarrow{\sim}V$.  The \emph{normalized sectorial frame} of Hypothesis~\ref{hyp:subdominant} is
\[
 \mathcal Y_{\mathrm{as},V}(z;u):V^{\mathrm{as}}\xrightarrow{\sim}V,
 \qquad\mathcal Y_{\mathrm{as},V}(z;u)e_b=y_b(z;u).
\]
Both are expressed in the original trivialization.  Their expressions in the local gauges are
\begin{equation}\label{eq:local-gauges-returned-frames}
 \mathcal Y_{\mathrm F,V}^{T_0}:=\rho_V(T_0)\mathcal Y_{\mathrm F,V},\qquad
 \mathcal Y_{\mathrm{as},V}^{T_\infty}:=\rho_V(T_\infty)\mathcal Y_{\mathrm{as},V}.
\end{equation}
The \emph{reference column} with label $b$ is $\mathcal Y_{\mathrm{as},V}e_b$.

The reference column associated with the highest-weight vector $v_{\omega_i}$ is the distinguished solution
\begin{equation}\label{eq:general-subdominant}
 \begin{aligned}
 \Phi_0^{(i)}&:=\mathcal Y_{\mathrm{as},V_i}v_{\omega_i},\\
 \rho_i(T_\infty)\Phi_0^{(i)}
 &=e^{-\nu_iS(z;u)}\bigl(\xi_i+o(1)\bigr).
 \end{aligned}
\end{equation}
The limit is taken on closed angular subsectors of $\Sigma_\infty$.
\end{definition}

Hypothesis~\ref{hyp:frob} and Lemma~\ref{lem:leading-sectorial-tensors} give tensor preservation for the two local frames.  Proposition~\ref{prop:finite-tensor-reconstruction} therefore reconstructs unique $G$-valued frames, denoted by $\mathcal Y_{\mathrm F}$ and $\mathcal Y_{\mathrm{as}}$, on their respective local domains.
The normalized sectorial frame has the $G$-valued prefactor
\begin{equation}\label{eq:group-sectorial-prefactor}
 F(z;u):=a_\Lambda^{-1}T_\infty(z;u)\mathcal Y_{\mathrm{as}}(z;u)
             e^{S(z;u)D}.
\end{equation}
For $Y_V=\mathcal Y_{\mathrm{as},V}$ in each fundamental representation, the prefactor of \eqref{eq:general-sectorial-columns} is $F_V=\rho_V(F)$.  Consequently $F\to1$ with the uniformity prescribed in Definition~\ref{def:leading-sectorial-normalization}.

Having specified both local frames, we first impose compatibility with shifted monodromy at zero.  For a sector $\Sigma_n\subset\mathscr X$ at zero and a nonempty open set $U_n\subset U$, the source and target products are
\[
 D_n:=\Sigma_n\times U_n,\qquad
 \mathsf S_n(D_n)=(e^{2\pi in}\Sigma_n)\times(q^nU_n).
\]
When the Frobenius data are defined on both products and $u\in U_n$, the gauge-corrected pullback~\eqref{eq:local-gauge-covariant-pullback} takes a gauged column at $q^nu$ on $e^{2\pi in}\Sigma_n$ to a column at $u$ on $\Sigma_n$.

\begin{hypothesis}\label{hyp:symanzik}
For every $n\in\mathbb Z$, the following hold.
\begin{enumerate}[label=(\alph*),topsep=3pt,itemsep=2pt]
\item The operator coefficient satisfies
\[
 \mathscr A(e^{2\pi in}z;q^nu)=\mathscr A(z;u).
\]
\item There is a product $D_n$ as above on which Hypothesis~\ref{hyp:frob} holds, as well as on $\mathsf S_n(D_n)$, with the chosen gauge $T_0$, the same $\ell$ and indexed weight bases, and the branches induced by the fixed $\log z$ on the cover.  For every fundamental representation, label $(\mu,a)$ and fixed $u\in U_n$, the gauge-corrected pullback maps the prescribed class of Frobenius expansions for $\mathcal L_{q^nu}^{T_0}$ into the class with the same label for $\mathcal L_u^{T_0}$.
\item For every compact set $L\subset U_n$,
\[
 \lim_{\substack{z\to0\\z\in\Sigma_n}}
 \sup_{u\in L}
 \left\|\rho_{E}(\Gamma_{T_0,n}(z;u))
       -I_{E}\right\|=0.
\]
\end{enumerate}
\end{hypothesis}

By covariance of the equation, the shifted pullback of a distinguished column again solves the equation at $u$.  We denote these solutions by
\begin{equation}\label{eq:abstract-integral-orbit}
 \Phi_n^{(i)}:=\mathsf S_n^*\Phi_0^{(i)},\qquad n\in\mathbb Z.
\end{equation}
Initially, $\Phi_n^{(i)}$ is defined on the inverse image under $\mathsf S_n$ of the normalization domain of $\Phi_0^{(i)}$.  After the continuation in Hypothesis~\ref{hyp:ccm} below, the same formula defines these solutions on $\mathscr X$, meromorphically in $u$.
To compare their asymptotics with the original sectorial frame, we first compare the phase and gauge for one shift on the common sector.  This is the overlap of the sector in Hypothesis~\ref{hyp:subdominant} with its preimage under $\mathsf S_1$, with angular range
\[
 a-\pi<\arg z/h<b-2\pi/h.
\]

\begin{hypothesis}\label{hyp:common-orbit}
There is a holomorphic function $\delta:U\to\mathbb C$ such that, as $|z|\to\infty$ within this overlap,
\begin{equation}\label{eq:analytic-phase-rotation}
 S(e^{2\pi i}z;qu)-\gamma_hS(z;u)\longrightarrow\delta(u),
\end{equation}
and
\begin{equation}\label{eq:analytic-gauge-rotation}
 \mathcal G(z;u):=a_\Lambda^{-1}T_\infty(z;u)
 T_\infty(e^{2\pi i}z;qu)^{-1}a_\Lambda
 \longrightarrow\dot c^{-1}.
\end{equation}
The limits are uniform on closed subsectors and compact spectral sets, with the fixed branches and phase constants of Hypothesis~\ref{hyp:subdominant}.
\end{hypothesis}

The factor $\mathcal G=a_\Lambda^{-1}\Gamma_{T_\infty,1}a_\Lambda$ is the gauge comparison from \eqref{eq:local-gauge-covariant-pullback}, expressed in the fixed eigenbasis.  Having imposed the phase and gauge limits at infinity, we require the two local normalizations to continue to a common domain.
\begin{hypothesis}\label{hyp:ccm}
The fundamental local frames in the original trivialization admit single-valued analytic continuation to the same domain $\mathscr X$, meromorphic in $u\in U$. On each $D_n$ and $\mathsf S_n(D_n)$, the continued Frobenius frame restricts to the normalized Frobenius frame prescribed in Hypothesis~\ref{hyp:symanzik}(b).
\end{hypothesis}

\Needspace{8\baselineskip}
\subsection{Normalized frames and differential monodromy}\label{sec:general-analytic-consequences}
In this subsection we construct the central connection and the $q$-connection associated with shifted monodromy from the normalized frames.  We then prove Steinberg form in the sectorial frame and derive the column and tensor identities used below.

\paragraph{The central connection.}
The continued frames give two bases of the same solution space.  Their comparison can be computed at any ordinary point.  For a simply connected open set $D\subset\mathscr X$ and a fixed spectral value $u$, the space $\mathscr S_V(D,u)$ consists of the holomorphic maps $y:D\to V$ satisfying $\partial_zy+\rho_V(\mathscr A(z;u))y=0$.  A solution is uniquely determined by its value at any $z_0\in D$, and every value in $V$ occurs.  Thus evaluation $y\mapsto y(z_0)$ identifies this horizontal solution space with $V$:
\begin{equation}\label{eq:sectorial-evaluation-isomorphism}
 \operatorname{ev}_{z_0}:\mathscr S_V(D,u)\xrightarrow{\sim}V,\qquad z_0\in D.
\end{equation}
Intertwiners commute with the differential equation.  Moreover, evaluation respects direct sums, tensor products, and duals, so the same is true of the horizontal solution spaces and their continuation.

The local tensor identities persist under the continuation in Hypothesis~\ref{hyp:ccm}.  Corollary~\ref{cor:finite-tensor-analytic} therefore continues $\mathcal Y_{\mathrm F}$ and $\mathcal Y_{\mathrm{as}}$ as $G$-valued frames, meromorphic in $u$ and holomorphic wherever their matrices on $E$ are holomorphic and invertible.  In every finite-dimensional algebraic representation $V$, we henceforth use the induced frames
\[
 \mathcal Y_{\mathrm F,V}:=\rho_V(\mathcal Y_{\mathrm F}),\qquad
 \mathcal Y_{\mathrm{as},V}:=\rho_V(\mathcal Y_{\mathrm{as}}).
\]
The transition between the $G$-valued frames is
\[
 C(u):=\mathcal Y_{\mathrm F}(z;u)^{-1}\mathcal Y_{\mathrm{as}}(z;u).
\]
We call this transition matrix the \emph{central connection matrix}, or simply the \emph{central connection}.  Both frames solve the same equation in the original trivialization, and therefore
\[
 \partial_zC=\mathcal Y_{\mathrm F}^{-1}\mathscr A\mathcal Y_{\mathrm{as}}
             -\mathcal Y_{\mathrm F}^{-1}\mathscr A\mathcal Y_{\mathrm{as}}=0.
\]
Thus $C=C(u)\in G(K)$, with
\begin{equation}\label{eq:group-frame-central-ratio}
 \mathcal Y_{\mathrm{as}}(z;u)=\mathcal Y_{\mathrm F}(z;u)C(u),\qquad
 \rho_V(C(u)):V^{\mathrm{as}}\longrightarrow V^{\mathrm F}.
\end{equation}
Where the two local gauge expressions are both defined, we can also write this transition matrix as
\begin{equation}\label{eq:central-local-gauge-transition}
 \rho_V(C(u))=(\mathcal Y_{\mathrm F,V}^{T_0})^{-1}
       \rho_V(T_0T_\infty^{-1})\mathcal Y_{\mathrm{as},V}^{T_\infty}.
\end{equation}
The factor $T_0T_\infty^{-1}$ is needed here because the two frames were normalized in different local gauges.

\begin{remark}\label{rem:analytic-comparison-scope}
Both $G$-valued frames trivialize the same solution system.  If we change these trivializations by $\mathcal Y_{\mathrm F}\mapsto\mathcal Y_{\mathrm F}g_{\mathfrak f}$ and $\mathcal Y_{\mathrm{as}}\mapsto\mathcal Y_{\mathrm{as}}g_{\mathfrak a}$, the central connection matrix transforms by $C\mapsto g_{\mathfrak f}^{-1}Cg_{\mathfrak a}$.  Proposition~\ref{prop:finite-tensor-reconstruction} identifies changes of the fundamental frames preserving the defining tensors with elements of $G$.
\end{remark}

\paragraph{Frobenius covariance.}
\begin{proposition}\label{prop:Frobenius-Symanzik-covariance}
Under Hypotheses~\ref{hyp:frob}, \ref{hyp:symanzik}, and~\ref{hyp:ccm}, in each fundamental representation $V=V_i$ and on the products $D_n$ of Hypothesis~\ref{hyp:symanzik},
\begin{equation}\label{eq:Frobenius-Symanzik-weight}
 \chi_{\mu,a}(e^{2\pi in}z;q^nu)
 =e^{-2\pi in\mu(\ell)}\chi_{\mu,a}(z;u),\qquad n\in\mathbb Z.
\end{equation}
Equivalently,
\begin{equation}\label{eq:Frobenius-Symanzik-frame}
 \mathcal Y_{\mathrm F,V}(e^{2\pi in}z;q^nu)
 =\mathcal Y_{\mathrm F,V}(z;u)\rho_V(e^{-2\pi i\ell})^n.
\end{equation}
Both identities extend by analytic continuation wherever their columns are defined.
\end{proposition}
\begin{proof}
On $D_n$, for the gauged column $\chi_{\mu,a}^{T_0}=\rho_V(T_0)\chi_{\mu,a}$, covariance and \eqref{eq:local-gauge-covariant-pullback} give
\[
 \rho_V(T_0)\mathsf S_n^*\chi_{\mu,a}
 =\rho_V(\Gamma_{T_0,n})\chi_{\mu,a}^{T_0}(e^{2\pi in}z;q^nu).
\]
Hypothesis~\ref{hyp:symanzik} implies that the right-hand side lies in the same prescribed class of Frobenius expansions.  Its leading coefficient is $e^{-2\pi in\mu(\ell)}v_{\mu,a}$, because $\Gamma_{T_0,n}\to1$.  After dividing by this scalar, we may apply the uniqueness assertion of Hypothesis~\ref{hyp:frob}.  Multiplying the resulting equality by $\rho_V(T_0)^{-1}$ proves the identity for $\chi_{\mu,a}$.  We obtain the frame identity by applying this argument to each basis vector, including the basis vectors in weight spaces of higher multiplicity.  Faithfulness gives the corresponding identity for the local $G$-valued frame and hence for every induced representation.  The compatibility in Hypothesis~\ref{hyp:ccm} identifies these local frames with restrictions of the continued Frobenius frame, so the identities extend by analytic continuation.
\end{proof}

In the Frobenius frame, shifted monodromy has matrix $Z^{-1}$, where
\[
 Z=e^{2\pi i\ell}.
\]

\begin{remark}\label{rem:Frobenius-eigenvector-vs-general-solution}
For a general solution $y(z;u)=\mathcal Y_{\mathrm F,V}(z;u)c(u)$, formula~\eqref{eq:Frobenius-Symanzik-frame} gives the coefficient transformation
\[
 \mathsf S_n^*y(z;u)=\mathcal Y_{\mathrm F,V}(z;u)\rho_V(Z)^{-n} c(q^nu).
\]
Thus $\mathsf S_n^*$ acts on the entire coefficient column and shifts its spectral argument.  The scalar phase in~\eqref{eq:Frobenius-Symanzik-weight} describes its action on each normalized Frobenius column; multiplying that column by a function of $u$ also shifts that function.
\end{remark}

\paragraph{Shifted monodromy and the $q$-connection.}
By covariance, $\mathsf S_1^*\mathcal Y_{\mathrm{as}}$ and $\mathcal Y_{\mathrm{as}}$ are solution frames for the same equation at $u$, so their ratio is independent of $z$.  The inverse of their comparison represents the $q$-connection from $u$ to $qu$ and has matrix
\begin{equation}\label{eq:infinity-frame-transition}
 \mathscr M(u):=(\mathsf S_1^*\mathcal Y_{\mathrm{as}})(z;u)^{-1}
                    \mathcal Y_{\mathrm{as}}(z;u).
\end{equation}
Equivalently,
\[
 \mathsf S_1^*\mathcal Y_{\mathrm{as}}
 =\mathcal Y_{\mathrm{as}}\mathscr M(u)^{-1}.
\]
Thus shifted monodromy acts on sectorial coefficient functions by $a(u)\mapsto\rho_V(\mathscr M(u)^{-1})a(qu)$.  The comparison with the Frobenius frame also uses Hypothesis~\ref{hyp:frob}.  By Proposition~\ref{prop:Frobenius-Symanzik-covariance}, the $q$-connection has matrix $Z$ in the Frobenius frame,
\[
 (\mathsf S_1^*\mathcal Y_{\mathrm F}(z;u))^{-1}
       \mathcal Y_{\mathrm F}(z;u)=Z.
\]
The matrices $Z$ and $\mathscr M$ therefore represent the same $q$-connection in Frobenius and sectorial coordinates, respectively.  Their inverses, $Z^{-1}$ and $\mathscr M^{-1}$, describe the pullback from $qu$ to $u$.

The iterated $q$-connection from $u$ to $q^nu$ has matrix
\begin{equation}\label{eq:infinity-frame-iterate}
 \mathscr M^{[0]}(u)=1,\qquad
 \mathscr M^{[n]}(u)
 =\mathscr M(q^{n-1}u)\cdots\mathscr M(u),
 \quad n>0.
\end{equation}
For negative indices, the convention is $\mathscr M^{[-n]}(u)=\mathscr M^{[n]}(q^{-n}u)^{-1}$, $n>0$.

In the pulled-back frame, $v_{\omega_i}$ selects the distinguished column $\Phi_n^{(i)}$.  The formula below expresses that solution in the fixed sectorial frame at $u$, where its coordinate vector is $\rho_i(\mathscr M^{[n]}(u)^{-1})v_{\omega_i}$.

\begin{proposition}\label{prop:infinity-frame-transition}
Under Hypotheses~\ref{hyp:frob}--\ref{hyp:symanzik} and~\ref{hyp:ccm}, the central connection and the $q$-connection satisfy
\begin{equation}\label{eq:central-infinity-transport}
 C(qu)\mathscr M(u)=Z C(u),
\end{equation}
and, for every $n\in\mathbb Z$,
\begin{align}
 \Phi_n^{(i)}(z;u)
 &=\rho_i(\mathcal Y_{\mathrm{as}}(z;u))
       \rho_i(\mathscr M^{[n]}(u)^{-1})v_{\omega_i},
       \label{eq:actual-infinity-frame-column}\\
 C(q^nu)\mathscr M^{[n]}(u)&=Z^nC(u).
       \label{eq:central-infinity-iterate}
\end{align}
\end{proposition}
\begin{proof}
Pulling back \eqref{eq:group-frame-central-ratio} and using \eqref{eq:Frobenius-Symanzik-frame} gives
\[
 \mathcal Y_{\mathrm{as}}(e^{2\pi i}z;qu)
 =\mathcal Y_{\mathrm F}(z;u)Z^{-1} C(qu).
\]
Substitution in the defining equation for $\mathscr M$ yields
\[
 \mathscr M(u)=C(qu)^{-1}Z C(u),
\]
which proves \eqref{eq:central-infinity-transport}.  Iterating the $q$-connection gives \eqref{eq:central-infinity-iterate}.  Its inverse form is
\[
 C(u)\mathscr M^{[n]}(u)^{-1}=Z^{-n}C(q^nu).
\]
Together with Frobenius covariance, this gives the coefficient formula for $\Phi_n^{(i)}$.  The negative-index convention gives the same identities for $n<0$.
\end{proof}

\begin{remark}\label{rem:general-hypotheses-scope}
These comparisons were made with both frames in the same differential gauge.  They also give the transformation rule under a common change of gauge.  A comparison of coefficient coordinates at fixed $u$ is a \emph{change of frame}, while \emph{differential parallel transport} means continuation along paths in $z$.

If we apply a common differential gauge $Q$, the normalizing gauges become $T_0Q^{-1}$ and $T_\infty Q^{-1}$.  The action of $\mathsf S_n^*$ on $Y^Q=\rho_V(Q)Y$ is then $\rho_V(Q)\mathsf S_n^*(\rho_V(Q)^{-1}Y^Q)$, as in \eqref{eq:local-gauge-covariant-pullback}.  In either local normalization gauge its expression remains unchanged, because
\[
 (TQ^{-1})\,\Gamma_{Q,n}\,
     (T\circ \mathsf S_n\,(Q\circ \mathsf S_n)^{-1})^{-1}=\Gamma_{T,n}.
\]
It follows that $C$ and the pullback in the sectorial frame remain unchanged:
\[
 (\mathcal Y_{\mathrm F}^Q)^{-1}\mathcal Y_{\mathrm{as}}^Q=C,
 \qquad
 (\mathcal Y_{\mathrm{as}}^Q)^{-1}
       Q\,\mathsf S_1^*\bigl(Q^{-1}\mathcal Y_{\mathrm{as}}^Q\bigr)=\mathscr M^{-1}.
\]
Thus a common differential gauge, with the compensated pullback above, leaves $C$ and $\mathscr M$ unchanged.  Right multiplication of either solution frame changes its basis of solutions.  The resulting transformation of $C$ is given in Remark~\ref{rem:analytic-comparison-scope}.
\end{remark}

\paragraph{Steinberg form from sectorial normalization.}
The constant factor $\dot c$ in sectorial transport comes from the principal rotation.  The phase increment in Hypothesis~\ref{hyp:common-orbit} supplies the accompanying Cartan factor.
\begin{definition}\label{def:asymptotic-Coxeter-factors}
For the function $\delta$ of Hypothesis~\ref{hyp:common-orbit} and the fixed Cartan element $D$, the holomorphic maps $h_{\mathrm{as}}:U\to H$ and $s_{\mathrm{as}}:U\to N_G(H)$ are
\begin{equation}\label{eq:asymptotic-Coxeter-factors}
 h_{\mathrm{as}}(u):=e^{\delta(u)D},\qquad
 s_{\mathrm{as}}(u):=h_{\mathrm{as}}(u)\dot c.
\end{equation}
\end{definition}
The Weyl image of $s_{\mathrm{as}}(u)$ is $c$. The sectorial normalization also determines the position of $\mathscr M(u)$ relative to this lift.

\begin{theorem}\label{prop:analytic-sectorial-Steinberg}\label{prop:sectorial-Steinberg}
Under Hypotheses~\ref{hyp:subdominant}--\ref{hyp:ccm}, the $q$-connection matrix in the normalized sectorial frame satisfies
\begin{equation}\label{eq:reference-sectorial-Steinberg}
 \mathscr M(u)\in U_c(K)s_{\mathrm{as}}(u)
 =U_c(K)h_{\mathrm{as}}(u)\dot c.
\end{equation}
This holds in every ADE type.
\end{theorem}

It suffices to show that the factor $n(u):=\mathscr M(u)s_{\mathrm{as}}(u)^{-1}$ belongs to $U_c=N_+\cap\dot cN_-\dot c^{-1}$. Comparing the original and rotated frames will give directions in which conjugation by a Cartan element contracts $n(u)$ to the identity. The following lemma identifies the unipotent subgroup determined by each such direction.

\begin{lemma}\label{lem:sectorial-direction-contraction}
Suppose $H_\theta=\operatorname{Re}(e^{i\theta}D)$ is regular in the real Cartan subalgebra, and let $N_\theta^+$ be the unipotent subgroup generated by the root subgroups for $\alpha(H_\theta)>0$. For any complex function $T(r)=re^{i\theta}+o(r)$ as $r\to\infty$, an element $g\in G$ belongs to $N_\theta^+$ if and only if
\[
 e^{-T(r)D}g e^{T(r)D}\longrightarrow 1 .
\]
\end{lemma}
\begin{proof}
In ordered root coordinates conjugation by $e^{-T(r)D}$ multiplies the coordinate of the root $\alpha$ by $e^{-T(r)\alpha(D)}$. Its modulus is $\exp(-r\alpha(H_\theta)+o(r))$, which tends to zero on every positive root for $H_\theta$. Hence every element of $N_\theta^+$ contracts to the identity.

Conversely, if the conjugates of $g$ tend to $1$, they eventually belong to the open Gauss cell $N_\theta^- H N_\theta^+$, where $N_\theta^-$ is generated by the opposite root subgroups. The cell is invariant under conjugation by $H$, so $g$ belongs to it as well. In its unique factorization $g=v_-h_0v_+$, continuity of the Gauss factors implies $h_0=1$ and $e^{-T(r)D}v_-e^{T(r)D}\to1$. For each root of $N_\theta^-$ the modulus of the multiplier tends to infinity. Convergence to the identity forces every root coordinate of $v_-$ to vanish. Thus $g=v_+\in N_\theta^+$.
\end{proof}

\begin{proof}[Proof of Theorem~\ref{prop:analytic-sectorial-Steinberg}]
\emph{Comparison of the sectorial frames.}
We work at a regular spectral value $u$ for which both frames are defined, and write $\omega=\gamma_h=e^{2\pi i/h}$ and $\Delta=2\pi/h$. A prime denotes evaluation at $(e^{2\pi i}z,qu)$. If $t=S(z;u)$, the two sectorial frames expressed in the same gauge $a_\Lambda^{-1}T_\infty(z;u)$ are
\[
 F e^{-tD}\qquad\text{and}\qquad \mathcal G F'e^{-t'D},
\]
where $F$ is the normalized prefactor in~\eqref{eq:group-sectorial-prefactor} and $\mathcal G$ is the gauge comparison in~\eqref{eq:analytic-gauge-rotation}. On the overlap $a-\pi<\arg z/h<b-\Delta$, Hypotheses~\ref{hyp:subdominant} and~\ref{hyp:common-orbit} give
\[
 F,F'\longrightarrow 1,\qquad
 \mathcal G\longrightarrow\dot c^{-1},\qquad
 t'-\omega t\longrightarrow\delta(u)
\]
as $|z|\to\infty$ in closed subsectors. By definition~\eqref{eq:infinity-frame-transition}, the inverse of the rotated frame times the original frame is
\[
 \mathscr M(u)=e^{t'D}(\mathcal G F')^{-1}F e^{-tD}.
\]

\emph{Contraction on the overlap.}
Since $\operatorname{Ad}_{\dot c}D=\omega D$, we have $e^{-tD}\dot c^{-1}e^{\omega tD}=\dot c^{-1}$. Substituting the frame comparison and $s_{\mathrm{as}}(u)=e^{\delta(u)D}\dot c$ in $n(u)=\mathscr M(u)s_{\mathrm{as}}(u)^{-1}$ gives
\[
 e^{-\omega tD}n(u)e^{\omega tD}
 =e^{(t'-\omega t)D}(\mathcal G F')^{-1}F\dot c^{-1}e^{-\delta(u)D}
 \longrightarrow 1 .
\]
Along a ray $z=(r/h)^he^{ih\vartheta}$ in the overlap, the leading asymptotic~\eqref{eq:abstract-principal-asymptotic-type} gives $t=re^{i\vartheta}+o(r)$. The exponent $\omega t$ therefore has direction $\theta=\vartheta+\Delta$ and leading term $re^{i\theta}$. Lemma~\ref{lem:sectorial-direction-contraction} now gives $n(u)\in N_\theta^+$ for every regular direction in
\[
 a-\pi+\Delta<\theta<b .
\]

\emph{Identification of the two unipotent subgroups.}
For any $\theta\in(a,b)$, both $\theta$ and $\theta'=\theta-\pi+\Delta$ lie in this interval because $\Delta\le\pi$. The simple-root values of $D=\Lambda_H$ in~\eqref{eq:cyclic-Cartan-PF-vector} give
\[
 \alpha_i(H_\theta)=
 \begin{cases}
 |\alpha_i(D)|\sin(\theta+\pi/h),&i\in I_0,\\
 -|\alpha_i(D)|\sin\theta,&i\in I_1.
 \end{cases}
\]
Both expressions are positive since $-\pi/h<a<\theta<b<0$. Thus $H_\theta$ is strictly dominant and $N_\theta^+=N_+$. Since $cD=\omega D$ and the Weyl action preserves the real Cartan subalgebra, we have
\[
 H_{\theta'}=\operatorname{Re}(e^{i(\theta-\pi+\Delta)}D)=-cH_\theta .
\]
The positive roots for $H_{\theta'}$ are therefore $c\Delta_-$, giving $N_{\theta'}^+=\dot cN_-\dot c^{-1}$. Contraction in these two directions puts $n(u)$ in $U_c$. In rank one $\Delta=\pi$ and $c=w_0$, so the directions coincide and both give $N_+$.

Finally $n(u)$ is meromorphic, since $\mathscr M(u)$ is meromorphic and $s_{\mathrm{as}}(u)$ is holomorphic. The subgroup $U_c$ is closed, so membership at every regular spectral value gives $n\in U_c(K)$. This proves~\eqref{eq:reference-sectorial-Steinberg}.
\end{proof}

The unipotent factor in the Steinberg form drops out of the fundamental highest-weight columns.  When $h$ is even, it also drops out of the iterated columns in the range below.

\begin{lemma}\label{lem:analytic-finite-columns}\label{prop:finite-infinity-transport}
One has $\mathscr M^{-1} v_{\omega_i}=s_{\mathrm{as}}^{-1}v_{\omega_i}$ in every fundamental representation.  If $h$ is even and $d=h/2$, then
\begin{equation}\label{eq:finite-infinity-transport-products}
 \mathscr M^{[m]}(u)^{-1}v_{\omega_i}
 =s_{\mathrm{as}}(u)^{-1}\cdots
   s_{\mathrm{as}}(q^{m-1}u)^{-1}v_{\omega_i},
 \qquad 0\le m\le d.
\end{equation}
The product is the identity for $m=0$.
\end{lemma}
\begin{proof}
The factorization is $\mathscr M^{-1}=s_{\mathrm{as}}^{-1}u_+$ with $u_+\in U_c$.  Since $U_c\subset N_+$ fixes all fundamental highest-weight vectors, the first identity follows.  In the iterated formula, each unipotent factor must fix the vector obtained by applying the inverse Coxeter factors to its right.  The reduced word for $w_0$ gives this property.  Indeed, Corollary~\ref{cor:fixed-datum-half-Coxeter} gives $c^d=w_0$.  The word obtained by repeating the reduced Coxeter word $d$ times has length $dr=|\Delta_+|=\ell_W(w_0)$ and is reduced.  Its inversion roots therefore give
\[
 c^j I(c^{-1})\subset\Delta_+\quad(0\le j<d),
 \qquad I(c^{-1})=\{\alpha\in\Delta_+:c^{-1}\alpha\in-\Delta_+\}.
\]
The roots in $I(c^{-1})$ are precisely the roots of $U_c$.  Consequently $\dot c^ju_+\dot c^{-j}\in N_+$, so $u_+$ fixes $\dot c^{-j}v_{\omega_i}$ for $j<d$.  Every length-$j$ product of the $s_{\mathrm{as}}^{-1}$ lies in $N_G(H)$ and sends $v_{\omega_i}$ to a nonzero vector of weight $c^{-j}\omega_i$.  Induction in
$\mathscr M^{[j+1]}(u)^{-1}=\mathscr M(u)^{-1}\mathscr M^{[j]}(qu)^{-1}$ now gives~\eqref{eq:finite-infinity-transport-products}.
\end{proof}

The formula for $\Phi_n^{(i)}$ in the sectorial frame turns these matrix identities into identities among solutions.

\begin{corollary}\label{cor:analytic-rotated-columns}\label{lem:actual-cyclic-column-comparison}\label{lem:infinity-normalized-columns}
For $n=0,1$, and also for $0\le n\le h/2$ when $h$ is even, there is a unique $a_{i,n}\in K^\times$ such that
\[
 \Phi_n^{(i)}
 =a_{i,n}\mathcal Y_{\mathrm{as},V_i}\dot c^{-n}v_{\omega_i}
 =a_{i,n}\mathcal Y_{\mathrm{as},V_i}
            \rho_i(a_\Lambda^{-1})\xi_{i,n}.
\]
Equivalently, these are the column identities
\begin{equation}\label{eq:infinity-normalized-columns}
 \mathscr M^{[n]}(u)^{-1}v_{\omega_i}=a_{i,n}(u)\dot c^{-n}v_{\omega_i}.
\end{equation}
Here $a_{i,0}=1$ and $a_{i,1}=\omega_i(h_{\mathrm{as}})^{-1}$.
\end{corollary}
\begin{proof}
The ordered products in Lemma~\ref{lem:analytic-finite-columns} normalize $H$ and have Weyl component $c^{-n}$.  Their actions on $v_{\omega_i}$ are therefore nonzero multiples of $\dot c^{-n}v_{\omega_i}$.  Equation~\eqref{eq:actual-infinity-frame-column} gives the result; the equality with $\xi_{i,n}$ is Lemma~\ref{lem:Kostant-full-transport}.
\end{proof}

\paragraph{The differential tensor identity.}
The column identities for $n=0,1$ transport the Pl\"ucker relation of Section~2 to the distinguished solutions.  Their normalization factors give the scalar
\begin{equation}\label{eq:asymptotic-tensor-multiplier}
 \varphi_i(u):=\frac{a_{i,1}(u)}{\prod_{j\sim i}a_{j,1}(u)^{\epsilon_i}}.
\end{equation}
\begin{proposition}\label{prop:general-Phi-from-Kostant}\label{prop:derived-Phi-system}
Hypotheses~\ref{hyp:subdominant}, \ref{hyp:symanzik}, \ref{hyp:common-orbit}, and~\ref{hyp:ccm} imply
\begin{equation}\label{eq:hyp-fundamental-Phi-system}
 m_i(\Phi_0^{(i)}\wedge\Phi_1^{(i)})
 =\varphi_i(u)\bigotimes_{j\sim i}\Phi_{\epsilon_i}^{(j)}.
\end{equation}
For every $n\in\mathbb Z$,
\begin{equation}\label{eq:hyp-integral-Phi-system}
 m_i(\Phi_n^{(i)}\wedge\Phi_{n+1}^{(i)})
 =\varphi_i(q^nu)\bigotimes_{j\sim i}\Phi_{n+\epsilon_i}^{(j)}.
\end{equation}
\end{proposition}
\begin{proof}
At a regular point $(z,u)$, the group element $\mathcal Y_{\mathrm{as}}(z;u)a_\Lambda^{-1}$ defines $L_V:=\rho_V(\mathcal Y_{\mathrm{as}}(z;u)a_\Lambda^{-1})$ in every representation $V$.  These maps commute with $G$-equivariant maps and respect tensor products.  For $V=V_i$, $L_{V_i}$ sends the chosen eigenbasis to the values of the corresponding normalized reference solutions.  Corollary~\ref{lem:actual-cyclic-column-comparison} gives
\[
 \Phi_0^{(i)}=L_{V_i}\xi_i,\qquad
 \Phi_1^{(i)}=a_{i,1}L_{V_i}\xi_{i,1}.
\]
With $W_i=\bigotimes_{j\sim i}V_j$, application of $L_{W_i}$ to the Pl\"ucker identity \eqref{eq:normalized-cyclic-Plucker-vector} for the leading eigenvectors at $n=0$ gives
\begin{align*}
 m_i(\Phi_0^{(i)}\wedge\Phi_1^{(i)})
 &=a_{i,1}L_{W_i}m_i(\xi_i\wedge\xi_{i,1})\\
 &=a_{i,1}L_{W_i}\bigotimes_{j\sim i}\xi_{j,\epsilon_i}\\
 &=\frac{a_{i,1}}{\prod_{j\sim i}a_{j,1}^{\epsilon_i}}
      \bigotimes_{j\sim i}\Phi_{\epsilon_i}^{(j)}.
\end{align*}
Since the differential point was arbitrary, this proves the first identity at every regular point.  We obtain the second by applying $\mathsf S_n^*$: semilinearity shifts the scalar coefficient, and $m_i$ is constant.  This also applies for negative $n$.
\end{proof}

\paragraph{Independence of the iterated columns.}
The column comparison of Corollary~\ref{lem:actual-cyclic-column-comparison} transfers independence of the leading eigenvectors to the iterated solutions.  For $1\le n\le h/2$, these eigenvectors have distinct eigenvalues, giving the following statement in types $D,E$.

\begin{corollary}\label{cor:finite-block-independence}
For $G$ of type $D$ or $E$ and a differential family satisfying Hypotheses~\ref{hyp:subdominant}, \ref{hyp:symanzik}, \ref{hyp:common-orbit}, and~\ref{hyp:ccm}, the solutions at every fundamental node $i$
\[
 \Phi_1^{(i)},\ldots,\Phi_{h/2}^{(i)}
\]
are linearly independent over $K$.
\end{corollary}
\begin{proof}
Corollary~\ref{lem:actual-cyclic-column-comparison} expresses these solutions as nonzero $K$-multiples of the images of $\xi_{i,1},\ldots,\xi_{i,h/2}$ under the common invertible matrix $\mathcal Y_{\mathrm{as},V_i}\rho_i(a_\Lambda^{-1})$.  These eigenvectors have pairwise distinct eigenvalues $\gamma_h\nu_i,\ldots,\gamma_h^{h/2}\nu_i$, since $\nu_i\ne0$.  They are therefore linearly independent over $\mathbb C$, and hence over $K$.  The invertible matrix and the nonzero scalar factors preserve their independence.
\end{proof}

\subsection{Verification for FFH opers}\label{sec:FFH-frames}\label{sec:FFH-verification}
In this subsection we verify the five analytic hypotheses for untwisted simply-laced Feigin--Frenkel--Hernandez (FFH) opers.  We use the analytic basis theorems of \cite{MR2024}, show that the normalized Frobenius frames are induced by a single $G$-valued frame, and compute the changes of phase and gauge at infinity under shifted monodromy.  We begin with the current realization and pass to an ordinary differential equation with a spectral parameter.

The parameters satisfy
\[
 0<k<1,
 \qquad q=e^{-2\pi ik}.
\]
As in Section~2, $\mathring f\in\n_-$ is principal nilpotent, $v_\theta\in\g_\theta$ is a highest-root vector, and $\theta^\vee\in\gh$ is the highest-root coroot.

\subsubsection{The FFH connection and its two realizations}
Here we present the current and loop realizations of the FFH connection, using the explicit representative in \cite[Theorem~6.18]{MR2024}.  Throughout this subsection, the FFH connection denotes this representative of the oper.

The parameters are $\ell\in\gh$, a finite index set $J$, distinct points $w_j\in\mathbb C^\times$, and elements $X(j)\in\n_+$ for $j\in J$.  The common differential domain is
\[
 \widetilde X_J
 :=\widetilde{\mathbb C^\times}\setminus\bigcup_{j\in J}\Pi^{-1}(w_j).
\]
All powers of $z$ are taken on this cover.

\begin{definition}\label{def:FFH-oper}
The \emph{current realization} of an untwisted simply-laced FFH connection is the differential operator on functions of $z$ and the affine loop coordinate $\lambda_{\mathrm{aff}}$ given by
\begin{equation}\label{eq:FFH-current-connection}
 z\mathcal L_{\mathrm{cur}}
 =z\partial_z+k\lambda_{\mathrm{aff}}\partial_{\lambda_{\mathrm{aff}}}
  +\mathring f+\ell+(z+\lambda_{\mathrm{aff}})v_\theta
  +\sum_{j\in J}\frac{z}{z-w_j}
       \left(-\theta^\vee+X(j)\right).
\end{equation}
The parameters at the finite poles are constrained so that the local monodromy about every $w_j$ is the identity, for every spectral value and in every finite-dimensional representation.
\end{definition}
The case $J=\varnothing$ gives the ground-state family, while $J\ne\varnothing$ gives a higher-state family.  The additional poles are \emph{apparent singularities}: they are regular singular points with identity local monodromy, so continuation around each pole returns every local solution to itself.  Solutions may have poles there, but have no branching or logarithmic terms.  A meromorphic change of frame removes such a singularity.

The vector field $z\partial_z+k\lambda_{\mathrm{aff}}\partial_{\lambda_{\mathrm{aff}}}$ annihilates $u=z^{-k}\lambda_{\mathrm{aff}}$.  In the coordinates $(z,u)$ this vector field is therefore $z\partial_z$ at fixed $u$.

\begin{definition}\label{def:FFH-loop-realization}
After the substitution $\lambda_{\mathrm{aff}}=z^ku$, the \emph{loop realization} of \eqref{eq:FFH-current-connection} is the family
\begin{equation}\label{eq:ffh}
 \mathcal L(z;u)
 =\partial_z+\frac1z
 \left(
 \mathring f+\ell+(z+z^ku)v_\theta
 +\sum_{j\in J}\frac{z}{z-w_j}
       \left(-\theta^\vee+X(j)\right)
 \right).
\end{equation}
\end{definition}
For a finite-dimensional $G$-module $V$, applying $\rho_V$ to the Lie-algebra coefficients in \eqref{eq:ffh} gives the operator $\mathcal L_V(z;u)$.  We call $\mathcal L_V(z;u)y=0$, with differentiation at fixed $u$, the \emph{FFH equation in the loop realization}.  The two local gauges $T_0$ and $T_\infty$ act on this equation, and the two frames are compared in its original gauge.

To compare monodromy in the two realizations, we use this change of coordinates on families of solutions that depend holomorphically on the spectral parameter.
\begin{definition}\label{def:FFH-solution-modules}
The coefficient rings of the current and loop solution modules are, respectively,
\begin{equation}\label{eq:Oprime}
 \mathcal O'
 =\{F(z^{-k}\lambda_{\mathrm{aff}}):F\in\mathcal O(\mathbb C)\},\qquad
 \mathcal O(\mathbb C_u).
\end{equation}
The operator $\mathcal L_{\mathrm{cur},V}$ denotes the current realization \eqref{eq:FFH-current-connection} in $V$.  The notation $\mathcal O$ denotes jointly holomorphic functions.  The global solution modules, with entire spectral dependence, are
\[
 \begin{aligned}
 \mathcal A_V^{\mathrm{cur}}
 &:=\{\psi\in\mathcal O(\widetilde X_J\times\mathbb C_{\lambda_{\mathrm{aff}}},V):
                 \mathcal L_{\mathrm{cur},V}\psi=0\},\\
 \mathcal A_V^{\mathrm{loop}}
 &:=\{\widetilde\psi\in\mathcal O(\widetilde X_J\times\mathbb C_u,V):
                 \mathcal L_V(z;u)\widetilde\psi=0\}.
 \end{aligned}
\]
An entire function $F$ of the spectral parameter acts on these modules by
\[
 \begin{aligned}
 (F\psi)(z;\lambda_{\mathrm{aff}})
   &=F(z^{-k}\lambda_{\mathrm{aff}})\psi(z;\lambda_{\mathrm{aff}}),\\
 (F\widetilde\psi)(z;u)
   &=F(u)\widetilde\psi(z;u).
 \end{aligned}
\qedhere
\]
\end{definition}
The current vector field annihilates $u$, so multiplication by $F(u)$ preserves current solutions.  To compare these modules, we evaluate the affine loop coordinate at $z^ku$.  On holomorphic functions this gives the map
\begin{equation}\label{eq:loop-current-evaluation}
 (\mathscr E_V\psi)(z;u)
 :=\psi(z;z^ku).
\end{equation}
The identification of solutions is given in \cite[Lemma~3.4]{MR2024}.  In the present coordinates it takes the following form.

\Needspace{12\baselineskip}
\begin{proposition}\label{lem:loop-current-comparison}
The map $\mathscr E_V$ restricts to an isomorphism
$\mathcal A_V^{\mathrm{cur}}\xrightarrow{\sim}\mathcal A_V^{\mathrm{loop}}$,
which is linear over the ring of entire functions of $u$.  For every such function $F$,
\[
 \mathscr E_V(F\psi)=F\,\mathscr E_V\psi.
\]
Its inverse is
\[
 (\mathscr E_V^{-1}\widetilde\psi)(z;\lambda_{\mathrm{aff}})
 :=\widetilde\psi(z;z^{-k}\lambda_{\mathrm{aff}}).
\]
Evaluation commutes with constant intertwiners, tensor products, and analytic continuation around the apparent singularities.
\end{proposition}
\begin{proof}
The chain rule gives
\[
 z\partial_z(\mathscr E_V\psi)
 =\mathscr E_V\bigl((z\partial_z
       +k\lambda_{\mathrm{aff}}\partial_{\lambda_{\mathrm{aff}}})\psi\bigr).
\]
Substitution in the remaining coefficients yields
$z\mathcal L_V\mathscr E_V=\mathscr E_V\,z\mathcal L_{\mathrm{cur},V}$.
The inverse substitution preserves holomorphy and proves bijectivity on solutions.  Under the substitution $\lambda_{\mathrm{aff}}=z^ku$, the current coefficient $F(z^{-k}\lambda_{\mathrm{aff}})$ becomes $F(u)$, which gives the asserted linearity.  The other compatibilities follow directly from substitution.
\end{proof}

Apparentness is imposed only at the additional poles $w_j$.  Under loop evaluation, continuation around zero at fixed $\lambda_{\mathrm{aff}}$ becomes the shifted monodromy $\mathsf S_1^*$.
\begin{corollary}\label{cor:FFH-current-monodromy}
The positive continuation $\mathcal M_{\mathrm{cur}}$ around zero at fixed $\lambda_{\mathrm{aff}}$ corresponds under loop evaluation to shifted monodromy through
\[
 \mathscr E_V\mathcal M_{\mathrm{cur}}\mathscr E_V^{-1}=\mathsf S_1^*.
\]
\end{corollary}
\begin{proof}
Since $(e^{2\pi i}z)^k(qu)=z^ku$, a positive turn at fixed $\lambda_{\mathrm{aff}}$ sends $(z,u)$ to $(e^{2\pi i}z,qu)=\mathsf S_1(z,u)$.  Pullback by this map gives the stated conjugacy.
\end{proof}

\Needspace{15\baselineskip}
\subsubsection{The generic FFH locus}
Here we describe the parameter locus on which the Frobenius normalization is well defined.

A connected analytic family is an \emph{FFH parameter family} if its ADE type, $k$, and labelled set $J$ are fixed, and $(\ell,w_j,X(j))$ varies in a connected solution locus of the trivial-monodromy equations.  The normalization at zero uses the following semisimplicity and nonresonance conditions.
\begin{definition}\label{def:generic-FFH-data}
Within a fixed FFH parameter family, an FFH oper is \emph{generic} if the following conditions hold for every fundamental representation $V_i$.  The residue
\[
 R_{0,i}:=\rho_i(\mathring f+\ell)
\]
is semisimple, and, for every residue eigenvalue $\omega(\ell)$ and every $(m,n)\in\mathbb Z_{\ge0}^2\setminus\{(0,0)\}$, the operator
\[
 R_{0,i}-\omega(\ell)+m+kn
\]
is invertible on $V_i$.  The \emph{generic locus} is the subset of the fixed parameter family on which these conditions hold.
\end{definition}
These are precisely the conditions of \cite[Definition~3.13]{MR2024} in the untwisted simply-laced case.  The second condition excludes the eigenvalue differences that would make a generalized Frobenius coefficient equation singular at an exponent $m+kn$.

For fixed $k$, these conditions define an open subset of the FFH parameter family.  To check openness of the semisimplicity condition, first note that $\mathring f+\ell$ is regular for every $\ell$.  Indeed,
\[
 t\operatorname{Ad}_{t^{\rho^\vee}}(\mathring f+\ell)
 =\mathring f+t\ell\longrightarrow\mathring f
 \qquad(t\longrightarrow0).
\]
The dimension of the centralizer is at least $r$ and is upper semicontinuous.  Since the centralizer of the principal nilpotent $\mathring f$ has dimension $r$, the residue $\mathring f+\ell$ is regular for every $\ell$.  It is semisimple precisely when $\ell$ is regular, an open condition.  The excluded resonance hyperplanes are locally finite, because $m+kn$ tends to infinity with $m+n$.  Nonemptiness and density in a given higher-state parameter family are separate questions.

The conditions above give the Frobenius normalization at zero.  The following theorem combines this normalization with the estimates at infinity, covariance, and continuation.  Its proof occupies the next five parts of this subsection.

\begin{theorem}\label{thm:FFH-verification}
Every untwisted simply-laced FFH oper on the generic locus of Definition~\ref{def:generic-FFH-data} satisfies Hypotheses~\ref{hyp:frob}--\ref{hyp:ccm} on a connected q-stable spectral domain $U$.  The iterates of shifted monodromy use $q=e^{-2\pi ik}$.
\end{theorem}

\subsubsection{Hypothesis \ref{hyp:frob}: Frobenius normalization}\label{sec:FFH-hyp-frob}

Here we express the local Frobenius bases of the FFH equation associated with \eqref{eq:ffh} in the normalization required by Hypothesis~\ref{hyp:frob}.  These bases are supplied by \cite[Proposition~5.10 and Corollary~5.11]{MRhigher2020} and \cite[Proposition~3.14]{MR2024}.  We use the same constant gauge $T_0=g_0\in N_-$ in every fundamental representation.  On the generic locus, $g_0$ carries the semisimple residue $\mathring f+\ell$ at zero to its Cartan projection $\ell$, by the standard conjugacy of maximal tori in a Borel subgroup \cite[16.29(d) and Theorem~16.33]{Milne2017}.

In the Frobenius expansion in this gauge, a current monomial $z^m\lambda_{\mathrm{aff}}^n$ becomes $z^{m+kn}u^n$.\nopagebreak[4]
\begin{lemma}\label{lem:FFH-Frobenius-basis}
On the generic locus of Definition~\ref{def:generic-FFH-data}, for every fundamental representation $V=V_i$ and every vector $v_{\omega,a}$ of weight $\omega$ in the fixed ordered basis, there is a unique solution $\chi_{\omega,a}$ of the FFH equation in $V$ with the expansion
\begin{equation}\label{eq:declared-FFH-Frobenius-class}
 \rho_V(g_0)\chi_{\omega,a}(z;u)
 =z^{-\omega(\ell)}\sum_{m,n\ge0}c_{\omega,a;m,n}z^{m+kn}u^n,
 \qquad c_{\omega,a;0,0}=v_{\omega,a}.
\end{equation}
Before loop evaluation, the series multiplying $z^{-\omega(\ell)}$ is a convergent power series in $(z,\lambda_{\mathrm{aff}})$ near $(0,0)$.  The solutions are entire in $u$ and form a local frame $\mathcal Y_{\mathrm F,V}$ for every $u\in\mathbb C$.
\end{lemma}

\paragraph{Frobenius frames and defining tensors.}
To assemble the fundamental Frobenius bases into a single $G$-valued frame, we use Proposition~\ref{prop:Frobenius-tensor-prescription} with $T_0=g_0$.  The convergent prefactors in Lemma~\ref{lem:FFH-Frobenius-basis} have entries in the common algebra
\[
 \mathcal R_u:=\{f(z,z^ku):f\in\mathbb C\{x,y\}\},
 \qquad
 \mathfrak m_u:=\{f(z,z^ku):f\in\mathbb C\{x,y\},\ f(0,0)=0\}
\]
at each fixed $u$.  Here $\mathbb C\{x,y\}$ denotes the algebra of power series convergent near $(0,0)$.  Since $k>0$, the limit of $f(z,z^ku)$ is $f(0,0)$.  For entries of $\mathfrak m_u$, composition with a holomorphic germ at the origin gives another convergent germ in $(x,y)$ before substitution, so $\mathcal R_u$ has the closure required by the proposition.  The corrections to the leading vectors lie in $\mathfrak m_u\otimes V_i$.

The criterion also requires uniqueness at fixed $u$ within this algebra.  Multiplying the gauged operator $\mathcal L^{g_0}(z;u)$ by $z$ gives
\[
 z\partial_z+\ell+B(z,z^ku),\qquad
 B\in\g\otimes\mathbb C\{x,y\},\quad B(0,0)=0.
\]
For two horizontal germs $y_1,y_2$ with the same prescribed leading vector of weight $\omega$, the difference of their prefactors $w=z^{\omega(\ell)}\rho_i(g_0)(y_1-y_2)$ belongs to $\mathfrak m_u\otimes V_i$ and satisfies
\[
 \bigl(z\partial_z+\rho_i(\ell)-\omega(\ell)
                  +\rho_i(B(z,z^ku))\bigr)w=0.
\]
After collecting terms of equal exponent, the expansion has exponents in the locally finite set $\{m+kn:m,n\ge0\}$.  If $w\ne0$, its first nonzero coefficient $w_s$ has $s>0$.  The terms from $B$ raise the exponent, so comparison at $s$ gives
\[
 \bigl(s+\rho_i(\ell)-\omega(\ell)\bigr)w_s=0.
\]
The operator is invertible by Definition~\ref{def:generic-FFH-data}, which contradicts $w_s\ne0$.  This also applies at $u=0$, where $\mathcal R_0=\mathbb C\{z\}$.

Proposition~\ref{prop:Frobenius-tensor-prescription} now gives a single holomorphic $G$-valued Frobenius frame.  Together with the normalized local bases of Lemma~\ref{lem:FFH-Frobenius-basis}, this verifies Hypothesis~\ref{hyp:frob}.

\subsubsection{Hypothesis \ref{hyp:subdominant}: Sectorial normalization}\label{sec:FFH-hyp-sectorial}

Here we construct the normalized sectorial frame at infinity by reducing the coefficient of the FFH equation associated with \eqref{eq:ffh} to a scalar multiple of $\Lambda$ plus a remainder integrable along rays.  The scalar coefficient will determine the exponential phase.  Throughout the estimates, the ADE type, $k$, and labelled set $J$ are fixed, while $(\ell,w_j,X(j))$ varies in the FFH parameter family defined above.

The reduction in the next proposition is obtained in \cite[proof of Proposition~3.20]{MR2024}.  We state it for the untwisted simply-laced loop realization \eqref{eq:ffh}, with the convention $\mathcal L^T=T\mathcal L T^{-1}$, making the common $G$-valued gauge and the uniformity on compact parameter sets explicit.

\begin{proposition}\label{prop:FFH-gauge-reduction}
For the FFH loop operator \eqref{eq:ffh}, there is a holomorphic $G$-valued gauge $T_\infty$ and a holomorphic $\g$-valued function $R$ such that
\begin{equation}\label{eq:uniform-tensorial-FFH-normal-form}
 \mathcal L^{T_\infty}
 =\partial_z+p(z;u)\Lambda+R(z;u).
\end{equation}
Here the scalar coefficient is
\begin{equation}\label{eq:FFH-scalar-coefficient}
 \begin{aligned}
 L_k&:=\left\lfloor\frac1{h(1-k)}\right\rfloor,\\
 p(z;u)&:=\sum_{a=0}^{L_k}
       \binom{1/h}{a}u^a z^{1/h-1-a(1-k)}.
 \end{aligned}
\end{equation}
In the faithful representation $E$, there is a decay exponent $\tau>0$ such that
\begin{equation}\label{eq:FFH-reduction-estimates}
 \|\rho_E(R(z;u))\|\le C|z|^{-1-\tau},\qquad
 \left\|\rho_E\!\left(T_\infty(z;u)z^{\rho^\vee/h}\right)-I_E\right\|
 \le C|z|^{-\tau}.
\end{equation}
For every finite angular interval $[\alpha,\beta]$, with $\alpha<\beta$, compact subset of the spectral plane, and compact subset of the fixed FFH parameter family, the gauge and remainder are holomorphic on a neighborhood of the corresponding product with
\[
 \{z:\ |z|>r,\quad \alpha\le\arg z\le\beta\},
\]
for sufficiently large $r$, and the estimates hold there with one constant $C$.  The choices agree on overlapping domains.  In particular, uniformly on these sets,
\begin{equation}\label{eq:FFH-reduction-limits}
 \lim_{|z|\to\infty}z^{1-1/h}p(z;u)=1,\qquad
 \lim_{|z|\to\infty}
 \rho_E\!\left(T_\infty(z;u)z^{\rho^\vee/h}\right)=I_E.
\end{equation}
\end{proposition}

The gauge used in the cited proof is a product of a torus element and a nilpotent exponential in $G$.  It therefore acts in every finite-dimensional representation $V$, where the reduced coefficient is $p\rho_V(\Lambda)+\rho_V(R)$.

The exponential phase is obtained by integrating $p$.  We denote the coefficients and exponents in \eqref{eq:FFH-scalar-coefficient} by
\[
 b_a:=\binom{1/h}{a},\qquad \beta_a:=\frac1h-a(1-k)\ge0,
 \qquad 0\le a\le L_k.
\]
With zero additive constant and the fixed branch of $\log z$, the primitive is
\begin{equation}\label{eq:FFH-exact-phase-primitive}
 \mathcal P_k(z;u):=\sum_{a=0}^{L_k}b_au^aJ_{\beta_a}(z),
 \qquad
 J_\beta(z)=
 \begin{cases}z^\beta/\beta,&\beta>0,\\ \log z,&\beta=0.\end{cases}
\end{equation}
Thus $\partial_z\mathcal P_k=p$.  The leading term is $hz^{1/h}$, and a term with $\beta_a=0$ contributes $b_au^a\log z$.

To construct the solution columns from this phase, we use the fixed ordered eigenbasis of $\rho_i(\Lambda)$, including the bases of repeated eigenspaces.  The distinguished simple eigenvalues and their eigenlines were specified in Section~2.  For an eigenvalue $\nu$, the leading term of the phase orders the exponential rates by $\operatorname{Re}(e^{i\arg(z)/h}\nu)$.  The basis theorem uses angular intervals on which this order is fixed for distinct eigenvalues.
\begin{definition}\label{def:spectrally-separating-interval}
For a finite-dimensional $G$-module $V$, $\Sigma_V$ denotes the spectrum of $\rho_V(\Lambda)$.  A Stokes direction is a value of $\varphi=\arg(z)/h$ for which two distinct eigenvalues have equal real parts after multiplication by $e^{i\varphi}$.  A closed interval $\mathcal I=[a,b]$ with nonempty interior \emph{separates the eigenvalues of $\rho_V(\Lambda)$} if
\begin{equation}\label{eq:explicit-spectral-separation}
 \operatorname{Re}(e^{i\varphi}\nu)
 =\operatorname{Re}(e^{i\varphi}\nu')
 \quad\Longleftrightarrow\quad\nu=\nu'
 \qquad(\varphi\in\mathcal I,\ \nu,\nu'\in\Sigma_V).
\end{equation}
It separates the spectra of $\rho_i(\Lambda)$ if this holds in every $V_i$.
\end{definition}
Only finitely many directions modulo $\pi$ are excluded by the fundamental spectra, so a common separating interval exists.  It is independent of the varying parameters in this fixed FFH family.

The reduced-system result used below is \cite[Theorem~3.26]{CGM2025}; its holomorphic parameter statement follows from \cite[Theorems~2.21, 2.27 and Corollaries~2.22, 2.28]{CGM2025}.  Part~(a) gives the full frame for Hypothesis~\ref{hyp:subdominant}.  Part~(b) gives a wider sector for a distinguished column, which will be used to compare consecutive solutions in type~$A$.
\begin{proposition}\label{prop:sectorial-reduced-system-input}
For a constant semisimple endomorphism $D$ of a finite-dimensional complex vector space $V$, the eigenbasis $(\xi_j)$ and eigenvalues $(\nu_j)$ are fixed.  For $t$ in a complex parameter domain, the scalar coefficient has the form
\[
 p(z;t)=z^{\beta_0}+\sum_{m=1}^{L}c_m(t)z^{\beta_m},
 \qquad \beta_0>\beta_1>\cdots>\beta_L\ge-1,
 \qquad \alpha:=\beta_0+1>0,
\]
where the exponents are fixed and the coefficients are holomorphic.  A term with exponent $-1$ is allowed to have zero coefficient.  The primitive $\mathcal P$ is obtained by integrating each power, using $\log z$ at exponent $-1$ and zero additive constant.  The remainder $\mathcal R(z;t)$ is jointly holomorphic near infinity on the universal cover of $\mathbb C^\times$ and satisfies
\[
 \|\mathcal R(z;t)\|=O(|z|^{-1-\tau}),\qquad \tau>0,
\]
uniformly on every closed finite angular sector and compact parameter set.  For the equation
\[
 \partial_zY=(-pD+\mathcal R)Y,
\]
the following hold for every fixed interval $a<b$.
\begin{enumerate}[label=(\alph*),topsep=3pt,itemsep=2pt]
\item If
\[
 \operatorname{Re}(e^{i\varphi}(\nu_j-\nu_k))\ne0
 \quad\text{whenever }\nu_j\ne\nu_k,\quad \varphi\in[a,b],
\]
there is a unique solution basis $(y_j)$ with
\[
 y_j(z;t)=(\xi_j+O(|z|^{-\tau}))e^{-\nu_j\mathcal P(z;t)}
\]
on closed subsectors of $a/\alpha<\arg z<(b+\pi)/\alpha$.
\item If $\nu_*$ is simple, $\xi_*\ne0$, $D\xi_*=\nu_*\xi_*$, and
\[
 \operatorname{Re}(e^{i\varphi}(\nu_*-\nu))>0
 \quad\bigl(\nu\in\operatorname{Spec}(D)\setminus\{\nu_*\},\ \varphi\in[a,b]\bigr),
\]
there is a unique solution with
\[
 y_*(z;t)=(\xi_*+O(|z|^{-\tau}))e^{-\nu_*\mathcal P(z;t)}
\]
on closed subsectors of $(a-\pi)/\alpha<\arg z<(b+\pi)/\alpha$.
\end{enumerate}
The solutions are jointly holomorphic near infinity, and the estimates are locally uniform in $t$.  Repeated eigenvalues are permitted in~(a), with the prescribed basis in each eigenspace.
\end{proposition}

Applying part~(a) to the reduced FFH operator gives the following frame in the original gauge, as in \cite[Proposition~3.20(1)]{MR2024} and \cite[Theorems~3.26(1) and~3.27(1)]{CGM2025}.  We use the phase $\mathcal P_k$ and the positive decay exponent $\tau$ fixed above.  For each bounded set of spectral parameters and each compact subset of the fixed FFH parameter family, the asymptotic estimates below hold beyond a common radius.

\begin{lemma}\label{lem:uniform-tensorial-FFH-reduction}
For an interval $\mathcal I=[a,b]$, $a<b$, separating the spectra of $\rho_i(\Lambda)$, the FFH equation associated with \eqref{eq:ffh} in each $V=V_i$ has a unique solution frame $Y_V=(y_b)_b$, indexed by the fixed $\Lambda$-eigenvectors $\xi_b$ with eigenvalues $\nu_b$, such that
\begin{equation}\label{eq:subdom}
 y_b(z;u)
 =\rho_V\!\left(z^{\rho^\vee/h}\right)
   \bigl(\xi_b+O(|z|^{-\tau})\bigr)
   e^{-\nu_b\mathcal P_k(z;u)}
\end{equation}
as $|z|\to\infty$ on closed subsectors of
\[
 ha<\arg z<h(b+\pi).
\]
The frame is jointly holomorphic in $z,u$ near infinity and entire in $u$ at every fixed ordinary point of the differential equation reached by continuation.  The estimates are locally uniform on the product of the spectral plane and the fixed FFH parameter family.  Repeated eigenvalues are allowed, with the fixed basis retained in each eigenspace.
\end{lemma}
\begin{proof}
Under the change of unknown $Y_V^{\mathrm{red}}=\rho_V(T_\infty)Y_V$, Proposition~\ref{prop:FFH-gauge-reduction} gives
\[
 \partial_zY_V^{\mathrm{red}}
 =\bigl(-p\rho_V(\Lambda)-\rho_V(R)\bigr)Y_V^{\mathrm{red}},
 \qquad \rho_V(R)=O(|z|^{-1-\tau}).
\]
The scalar coefficient~\eqref{eq:FFH-scalar-coefficient} has leading exponent $1/h-1$ and primitive $\mathcal P_k$.  Proposition~\ref{prop:sectorial-reduced-system-input}(a), with $D=\rho_V(\Lambda)$, $\mathcal R=-\rho_V(R)$ and $\alpha=1/h$, gives the required reduced basis.  Returning to the original gauge uses
\[
 Y_V=\rho_V(T_\infty)^{-1}Y_V^{\mathrm{red}},\qquad
 z^{-\rho^\vee/h}T_\infty^{-1}=1+O(|z|^{-\tau}),
\]
and proves~\eqref{eq:subdom}.  The parameter assertion in the proposition gives joint holomorphy and locally uniform estimates.  Uniqueness identifies the solutions on overlapping parameter domains.  Continuation by the original differential equation gives entire spectral dependence at every fixed ordinary point of the differential equation.
\end{proof}

For Hypothesis~\ref{hyp:subdominant} we use the full frame of Lemma~\ref{lem:uniform-tensorial-FFH-reduction}.  A common separating interval inside $(-\pi/h,0)$, translated by $-\pi$ in the basis theorem, gives the angular range required by the hypothesis.  There are numbers
\[
 -\pi/h<a'<a<b<b'<0
\]
such that $[a',b']$ separates the spectra of $\rho_i(\Lambda)$.  We fix
$\mathcal I_{\mathrm{ref}}=[a-\pi,b-\pi]$ and use the sectorial frame obtained from Lemma~\ref{lem:uniform-tensorial-FFH-reduction} on the larger interval $[a'-\pi,b'-\pi]$.  Its normalization holds on an open sector containing
\begin{equation}\label{eq:FFH-reference-wide-sector}
 \Sigma_{\mathrm{ref}}^{\mathrm{cl}}
 =\{z:\ h(a-\pi)\le\arg z\le hb,\quad |z|\text{ sufficiently large}\}.
\end{equation}
The angles depend only on the fixed spectra of $\rho_i(\Lambda)$.  The tail radius and the estimates are locally uniform on compact subsets of the product of the spectral plane and the fixed FFH parameter family.

\begin{corollary}\label{cor:FFH-sectorial-hypothesis}
The FFH equation satisfies Hypothesis~\ref{hyp:subdominant} with
\[
 T_\infty\text{ as in Proposition~\ref{prop:FFH-gauge-reduction}},\qquad S=\mathcal P_k,
\]
and with the normalization sector \eqref{eq:FFH-reference-wide-sector}.
\end{corollary}
\begin{proof}
Proposition~\ref{prop:FFH-gauge-reduction} gives the limits \eqref{eq:FFH-reduction-limits}, and the phase formula gives
\[
 \mathcal P_k(z;u)=hz^{1/h}+o(z^{1/h}).
\]
For $Y_i:=Y_{V_i}$, Lemma~\ref{lem:uniform-tensorial-FFH-reduction} gives
\[
 \rho_i(a_\Lambda^{-1}T_\infty)Y_i
       e^{\mathcal P_k\rho_i(\Lambda_H)}\longrightarrow I_{V_i}.
\]
The gauge reduction also gives $R/p=O(|z|^{-\tau-1/h})\to0$.  Taking the direct sum gives the leading normalization in $E$.  The phase includes all nondecaying terms, including any logarithm, and there is no separate residual power.  These statements verify Hypothesis~\ref{hyp:subdominant}.
\end{proof}

For the type-$A$ determinant calculation in Section~\ref{sec:FFH-typeA-basis}, we also need an asymptotic estimate for one column on a wider sector.  Part~(b) of Proposition~\ref{prop:sectorial-reduced-system-input} gives this estimate when the corresponding exponential is subdominant.  This is the subdominant-solution theorem of \cite[Theorem~3.15(2) and Proposition~3.20(2)]{MR2024}, in the form of \cite[Theorem~3.26(2)]{CGM2025}.

\begin{lemma}\label{lem:FFH-subdominant-solution}
For a finite-dimensional $G$-module $V$, an eigenvector $\xi_*\ne0$ of $\rho_V(\Lambda)$ with simple eigenvalue $\nu_*$, and real numbers $a<b$ satisfying
\begin{equation}\label{eq:FFH-subdominant-interval}
 \operatorname{Re}\!\left(e^{i\varphi}(\nu_*-\nu)\right)>0
 \qquad
 \bigl(\varphi\in[a,b],\quad
       \nu\in\operatorname{Spec}(\rho_V(\Lambda))\setminus\{\nu_*\}\bigr),
\end{equation}
there is a unique solution $y_*(z;u)$ of the FFH equation associated with \eqref{eq:ffh} with
\begin{equation}\label{eq:FFH-subdominant-large-sector}
 y_*(z;u)
 =\rho_V\!\left(z^{\rho^\vee/h}\right)
    \bigl(\xi_*+O(|z|^{-\tau})\bigr)
       e^{-\nu_*\mathcal P_k(z;u)}
\end{equation}
as $|z|\to\infty$ on closed subsectors of
\[
 h(a-\pi)<\arg z<h(b+\pi).
\]
The solution is jointly holomorphic in $z,u$ near infinity and entire in $u$ at every fixed ordinary point of the differential equation reached by continuation.  The estimate is locally uniform on the product of the spectral plane and the fixed FFH parameter family.
\end{lemma}
\begin{proof}
In the same reduced equation, condition~\eqref{eq:FFH-subdominant-interval} is the hypothesis of Proposition~\ref{prop:sectorial-reduced-system-input}(b), with $D=\rho_V(\Lambda)$, $\mathcal R=-\rho_V(R)$ and $\alpha=1/h$.  Its unique solution has the stated larger sector and leading vector $\xi_*$.  Multiplication by $\rho_V(T_\infty)^{-1}$, using~\eqref{eq:FFH-reduction-estimates}, proves~\eqref{eq:FFH-subdominant-large-sector}.  The parameter assertion and continuation of the original equation give the remaining conclusions.
\end{proof}

The column in Lemma~\ref{lem:FFH-subdominant-solution} is called \emph{subdominant}: along rays with $\arg z/h\in[a,b]$, its exponential decays fastest relative to the others.  Indeed,
\[
 \left|e^{-\nu\mathcal P_k(z;u)}\right|
 =\exp\!\left(-h|z|^{1/h}\operatorname{Re}(e^{i\varphi}\nu)
                 +o(|z|^{1/h})\right),\qquad \arg z=h\varphi.
\]
For the distinguished $\nu_i$, this real part is positive on $-\pi/h<\varphi<0$.  The strict inequality for the real parts selects the normalized column on its sector of subdominance; its leading asymptotic estimate extends to the larger sector stated in the lemma.

\subsubsection{Hypothesis \ref{hyp:symanzik}: Covariance under shifted monodromy}\label{sec:FFH-hyp-covariance}

\begin{proposition}\label{prop:FFH-integral-covariance}
On the generic locus of Definition~\ref{def:generic-FFH-data}, the FFH equation associated with \eqref{eq:ffh} satisfies Hypothesis~\ref{hyp:symanzik} with $T_0=g_0$ and $\Gamma_{T_0,r}=1$ for every $r\in\mathbb Z$.
\end{proposition}
\begin{proof}
The coefficient of the operator \eqref{eq:ffh} is invariant under $\mathsf S_r$: its finite-pole terms are unchanged, and
\[
 (e^{2\pi ir}z)^kq^ru=z^ku.
\]
For each $r$, the convergent Frobenius prefactors in $(z,\lambda_{\mathrm{aff}})$ give the normalization of Lemma~\ref{lem:FFH-Frobenius-basis} on both $D_r$ and $\mathsf S_r(D_r)$ for any bounded nonempty open set $U_r\subset U$ and a sufficiently small sector $\Sigma_r$.

The prescribed series retain the separate indices $(m,n)$, even when two exponents $m+kn$ coincide.  Their monomials satisfy
\[
 (e^{2\pi ir}z)^{m+kn}(q^ru)^n=z^{m+kn}u^n,
\]
whereas the leading power is multiplied by $e^{-2\pi ir\omega(\ell)}$.  Thus pullback preserves each class of Frobenius expansions before normalization of its leading coefficient.  Finally, $g_0$ is constant, so $\Gamma_{T_0,r}=1$, and the limit in Hypothesis~\ref{hyp:symanzik} holds identically.
\end{proof}

The Frobenius covariance identity~\eqref{eq:Frobenius-Symanzik-weight} holds on $D_r$ for the columns prescribed on $D_r$ and $\mathsf S_r(D_r)$.

\subsubsection{Hypothesis \ref{hyp:common-orbit}: Phase and gauge under shifted monodromy}\label{sec:FFH-hyp-matching}\label{sec:FFH-log-normalization}

Here we verify Hypothesis~\ref{hyp:common-orbit} by comparing $\mathcal P_k$ and $T_\infty$ before and after the shift.  We denote by $\lambda_{\log}(u)$ the coefficient of $\log z$ in $\mathcal P_k$.  A turn adds $2\pi i$ to the logarithm, so this term contributes an additive constant to the shifted phase.  The fixed primitive~\eqref{eq:FFH-exact-phase-primitive} gives
\begin{equation}\label{eq:FFH-log-coefficient}
 \lambda_{\log}(u):=\begin{cases}
 \displaystyle\binom{1/h}{m}u^m,
       &\displaystyle m=\frac1{h(1-k)}\in\mathbb Z_{>0},\\[4pt]
 0,&\displaystyle\frac1{h(1-k)}\notin\mathbb Z_{>0}.
 \end{cases}
\end{equation}
When the phase contains a logarithmic term, $k=1-1/(hm)$ and $q^m=\gamma_h$.  Consequently
\[
 \lambda_{\log}(qu)=\gamma_h\lambda_{\log}(u)
\]
for every $0<k<1$, including values for which $\lambda_{\log}=0$.  These logarithms belong to the phase at infinity; the classes of Frobenius expansions at zero are the nonresonant classes already verified above.

\begin{proposition}\label{prop:FFH-log-transport}
The fixed primitive satisfies, for every $n\in\mathbb Z$,
\begin{equation}\label{eq:FFH-phase-integer-rotation}
 \mathcal P_k(e^{2\pi in}z;q^nu)
 =\gamma_h^n\mathcal P_k(z;u)+2\pi in\gamma_h^n\lambda_{\log}(u).
\end{equation}
On the inverse image under $\mathsf S_n$ of the normalization sector, the pulled-back column $\Phi_n^{(i)}$ has leading vector $a_{i,n}(u)\xi_{i,n}$ relative to the scalar phase $e^{-\gamma_h^n\nu_i\mathcal P_k(z;u)}$, where
\begin{equation}\label{eq:FFH-actual-leading-scalar}
 a_{i,n}(u)=\exp(-2\pi in\gamma_h^n\nu_i\lambda_{\log}(u)).
\end{equation}
\end{proposition}
\begin{proof}
For each positive-power term in \eqref{eq:FFH-exact-phase-primitive}, we have
\[
 e^{2\pi in\beta_a}q^{na}=e^{2\pi in(1/h-a)}=\gamma_h^n.
\]
For the exceptional term we use $q^{nm}=\gamma_h^n$, together with the fact that deck translation replaces $\log z$ by $\log z+2\pi in$.  We obtain \eqref{eq:FFH-phase-integer-rotation} for all integer indices, including the negative ones.  The grading factor in \eqref{eq:subdom} transforms according to
\[
 (e^{2\pi in}z)^{\rho^\vee/h}
 =z^{\rho^\vee/h}\gamma_h^{-n\mathbf h},\qquad \mathbf h=-\rho^\vee.
\]
It follows that the leading vector after pullback is $\xi_{i,n}=\gamma_h^{-n\mathbf h}\xi_i$.  The additive constant in the exponential phase then contributes the scalar in \eqref{eq:FFH-actual-leading-scalar}.
\end{proof}

Since $\lambda_{\log}(u)$ is a polynomial, each $a_{i,n}(u)$ in \eqref{eq:FFH-actual-leading-scalar} is entire and nowhere vanishing.  All these scalars are identically one when $\lambda_{\log}\equiv0$.

\begin{corollary}\label{cor:FFH-rotation-hypothesis}
The FFH phase and gauge satisfy Hypothesis~\ref{hyp:common-orbit}, with
\[
 \delta(u)=2\pi i\gamma_h\lambda_{\log}(u).
\]
\end{corollary}
\begin{proof}
The phase assertion is the case $n=1$ of \eqref{eq:FFH-phase-integer-rotation}.  For the gauge limit, \eqref{eq:FFH-reduction-limits} gives
\[
 P_\infty(z,u):=T_\infty(z,u)z^{-\mathbf h/h}\longrightarrow1.
\]
On the overlap of the normalization sector and its preimage under $\mathsf S_1$,
\[
 T_\infty(z,u)T_\infty(e^{2\pi i}z,qu)^{-1}
 =P_\infty(z,u)\gamma_h^{-\mathbf h}P_\infty(e^{2\pi i}z,qu)^{-1}
 \longrightarrow\gamma_h^{-\mathbf h}.
\]
Conjugation by $a_\Lambda$ and \eqref{eq:Kostant-Coxeter-lift} give
\[
 a_\Lambda^{-1}T_\infty(z,u)
 T_\infty(e^{2\pi i}z,qu)^{-1}a_\Lambda
 \longrightarrow\dot c^{-1}.
\]
All limits are locally uniform on compact parameter sets by the gauge reduction and the finite expression for the phase.
\end{proof}

\subsubsection{Hypothesis \ref{hyp:ccm}: Global continuation}\label{sec:FFH-hyp-global}

Here we use the global basis theorem to continue the normalized local frames on $\widetilde X_J$ and establish their entire spectral dependence.  We state the results of \cite[Theorem~3.9 and Propositions~3.14, 3.20]{MR2024} using the current--loop identification of Proposition~\ref{lem:loop-current-comparison}.
\begin{lemma}\label{lem:globalize-FFH-frames}
For every finite-dimensional $G$-module $V$, the global solution module $\mathcal A_V^{\mathrm{cur}}$ is free over $\mathcal O'$ of rank $\dim V$.  On the generic locus of Definition~\ref{def:generic-FFH-data}, the normalized fundamental Frobenius and sectorial frames in the original loop gauge extend uniquely to $\widetilde X_J$.  Their columns are entire in $u$ and, before loop evaluation, form $\mathcal O'$-bases of the corresponding global solution modules.
\end{lemma}

\begin{proposition}\label{prop:tensorial-FFH-frames}\label{prop:FFH-central-existence}
On the generic locus of Definition~\ref{def:generic-FFH-data}, the FFH family satisfies Hypothesis~\ref{hyp:ccm} on $\widetilde X_J$.  At each fixed $z\in\widetilde X_J$, the continued normalized frames $\mathcal Y_{\mathrm F}(z;u)$ and $\mathcal Y_{\mathrm{as}}(z;u)$ are entire $G$-valued functions of $u$.  The central connection $C(u)$ is entire with values in $G$.
\end{proposition}
\begin{proof}
Lemma~\ref{lem:globalize-FFH-frames} gives single-valued analytic continuation of the normalized fundamental frames to $\widetilde X_J$, with entire spectral dependence. Near zero, the continued Frobenius columns are represented on every lifted sector by the same convergent series~\eqref{eq:declared-FFH-Frobenius-class}, with the fixed branch of $\log z$. They therefore restrict to the prescribed frames on $D_n$ and $\mathsf S_n(D_n)$, verifying also the normalization compatibility in Hypothesis~\ref{hyp:ccm}.

Their local defining-tensor identities persist by analytic continuation.  The argument of Section~\ref{sec:general-analytic-consequences}, using Corollary~\ref{cor:finite-tensor-analytic}, therefore gives $G$-valued extensions of the normalized frames.  At each fixed $z$, their matrices in $E$ are entire in $u$ and invertible for every $u$ by Lemma~\ref{lem:globalize-FFH-frames}; thus these extensions are holomorphic in $u$.  The central connection from \eqref{eq:group-frame-central-ratio} is entire because inversion and multiplication in $G$ are regular maps.
\end{proof}

\begin{proof}[Completion of the proof of Theorem~\ref{thm:FFH-verification}]
The five verifications above establish Hypotheses~\ref{hyp:frob}--\ref{hyp:ccm} for the FFH equation associated with \eqref{eq:ffh} on $\widetilde X_J$.  Restricting the entire spectral families to a connected q-stable spectral domain $U$, with $K=\mathcal M(U)$, proves the theorem.
\end{proof}

The verified hypotheses give Steinberg form in the FFH reference frame.  The shifted phase determines its Cartan factor, and the leading scalars in \eqref{eq:FFH-actual-leading-scalar} give the column multipliers for the shifts covered by Corollary~\ref{cor:analytic-rotated-columns}.
\begin{corollary}\label{lem:FFH-common-reference-sector}
For the FFH frame normalized on \eqref{eq:FFH-reference-wide-sector},
\[
 \mathscr M\in U_c(K)h_{\mathrm{as}}(u)\dot c,\qquad
 h_{\mathrm{as}}(u)
 =\exp\!\left(2\pi i\gamma_h\lambda_{\log}(u)\Lambda_H\right).
\]
For every fundamental node $i$,
\[
 \Phi_n^{(i)}(z,u)
 =a_{i,n}(u)\rho_i(\mathcal Y_{\mathrm{as}}(z,u))
       \dot c^{-n}v_{\omega_i}
\]
holds at $n=0,1$ in every ADE type and at $0\le n\le h/2$ whenever $h$ is even, with $a_{i,n}$ given by \eqref{eq:FFH-actual-leading-scalar}.
\end{corollary}
\begin{proof}
Theorem~\ref{prop:analytic-sectorial-Steinberg}, together with Corollary~\ref{cor:FFH-rotation-hypothesis}, gives the displayed Steinberg factorization in the chosen sectorial frame.  Corollary~\ref{cor:analytic-rotated-columns} gives the column identities in the indicated ranges.  To compute their scalars, we use the ordered-product identity~\eqref{eq:finite-infinity-transport-products} and its one-step case in Lemma~\ref{lem:analytic-finite-columns}.  The relations
\[
 \lambda_{\log}(qu)=\gamma_h\lambda_{\log}(u),\qquad
 \operatorname{Ad}_{\dot c}\Lambda_H=\gamma_h\Lambda_H
\]
show that the ordered inverse Coxeter factors act on $v_{\omega_i}$ with scalar
$\exp(-2\pi in\gamma_h^n\nu_i\lambda_{\log}(u))$, which is \eqref{eq:FFH-actual-leading-scalar}.  The identities hold meromorphically after continuation, and analytically wherever the normalized frames are regular.
\end{proof}

\subsubsection{Determinants of consecutive solutions for FFH opers in type $A$}\label{sec:FFH-typeA-basis}
Here we compute the determinant of a consecutive block in the defining representation of type $A_{N-1}$ from its asymptotics at infinity.  We separate the factors coming from the leading eigenvectors and the scalar multipliers of the shifted columns.  This calculation uses the normalization at infinity and applies on the Frobenius resonant locus as well.

In the defining representation chosen in Section~2, $\Lambda^N=aI$.  The distinguished eigenvalue is $\nu_1=1$, so $a=1$, and the cyclic element has the following matrix:
\[
 \Lambda=\sum_{a=1}^{N-1}E_{a+1,a}+E_{1N},\qquad
 \rho^\vee=\operatorname{diag}\left(\frac{N-1}{2},\frac{N-3}{2},\ldots,-\frac{N-1}{2}\right),
 \qquad \gamma=e^{2\pi i/N}.
\]
Here $\Phi_0^{(1)}$ is the distinguished column of the sectorial frame supplied by Lemma~\ref{lem:uniform-tensorial-FFH-reduction}, whose existence does not require Frobenius nonresonance.  The bipartition chosen in Section~2 satisfies $\epsilon_1=0$, and the prescribed leading vector of this column is a $\Lambda$-eigenvector with eigenvalue $1$.  We denote this leading vector by $v_0$ and form the eigenvectors $v_j$ and their determinant by
\begin{equation}\label{eq:typeA-leading-volume}
 v_j=e^{2\pi i j\rho^\vee/N}v_0,\qquad
 \Delta_N=\det(v_0,v_1,\ldots,v_{N-1}).
\end{equation}
These eigenvectors of $\Lambda$ have distinct eigenvalues $1,\gamma,\ldots,\gamma^{N-1}$, so $\Delta_N\ne0$.  The chosen normalization is \eqref{eq:subdom} with eigenvalue $1$ and leading vector $v_0$, and the columns are the iterates in \eqref{eq:abstract-integral-orbit}.

The second factor in the determinant is the product of the first $N$ leading scalars.  Since $\nu_1=1$, Proposition~\ref{prop:FFH-log-transport} gives
\begin{equation}\label{eq:FFH-volume-multiplier}
 d_N(u):=\prod_{j=0}^{N-1}a_{1,j}(u)
 =\exp\!\left(\frac{2\pi iN\lambda_{\log}(u)}{1-\gamma}\right),
 \qquad
 \frac{d_N(qu)}{d_N(u)}
 =e^{-2\pi iN\lambda_{\log}(u)}.
\end{equation}
Here we have used $\sum_{j=0}^{N-1}j\gamma^j=-N/(1-\gamma)$.  It follows that $d_N$ and its reciprocal are entire, and both equal one when $\lambda_{\log}=0$.

\begin{theorem}\label{thm:FFH-typeA-cyclic-basis}
For an untwisted type-$A_{N-1}$ FFH connection with $0<k<1$ and trivial local monodromy at its additional singularities, the determinant of the $N$ consecutive solutions is independent of $z$ and, for every $s\in\mathbb Z$ and $u\in\mathbb C$, equals
\begin{equation}\label{eq:typeA-entire-cyclic-determinant}
 \det(\Phi_s^{(1)},\ldots,\Phi_{s+N-1}^{(1)})
 =\Delta_Nd_N(q^su).
\end{equation}
When $d_N=1$, every consecutive determinant is the same nonzero constant $\Delta_N$.  These columns and the inverse of their solution matrix are entire in $u$.  They form a basis at every complex spectral value and over the ring of entire spectral functions, equivalently over $\mathcal O'$ before loop evaluation.
\end{theorem}

\begin{proof}
We compute the determinant for $s=0$ on a ray where all $N$ shifted columns have a leading asymptotic estimate.  The one-column estimate gives the needed wider sector for $\Phi_0^{(1)}$ once we identify its solution with the reference column.  The eigenvalues of $\Lambda$ in the defining representation satisfy
\[
 \operatorname{Re}\!\left(e^{i\varphi}(1-\gamma^j)\right)>0
 \qquad\left(1\le j<N,\quad |\varphi|<\frac\pi N\right).
\]
Lemma~\ref{lem:FFH-subdominant-solution}, with $\xi_*=v_0$, $\nu_*=1$, and $[a,b]=[-\pi/(2N),\pi/(2N)]$, supplies a solution $y_*$ with
\[
 y_*(z;u)=z^{\rho^\vee/N}(v_0+\eta(z;u))e^{-\mathcal P_k(z;u)},
 \qquad \eta(z;u)\longrightarrow0,
\]
on closed subsectors of $|\arg z|<N\pi+\pi/2$, locally uniformly in $u$.

On a common ray in the reference sector with $-\pi/N<\arg z/N<0$, the full frame of Lemma~\ref{lem:uniform-tensorial-FFH-reduction} has the form
\[
 Y_{\mathrm{ref}}
 =z^{\rho^\vee/N}(\Xi+\varepsilon_{\mathrm{ref}})
       \operatorname{diag}_b(e^{-\nu_b\mathcal P_k}),
 \qquad
 \Xi=(\xi_b)_b,\quad \varepsilon_{\mathrm{ref}}\longrightarrow0.
\]
The label $b_0$ corresponds to $\xi_{b_0}=v_0$ and $\nu_{b_0}=1$.  The coefficient vector $c(u)=Y_{\mathrm{ref}}^{-1}y_*$ is independent of $z$ and satisfies
\[
 c(u)=\operatorname{diag}_b\!\left(e^{(\nu_b-1)\mathcal P_k}\right)
       (\Xi+\varepsilon_{\mathrm{ref}})^{-1}(v_0+\eta)
       \longrightarrow e_{b_0}.
\]
Indeed, $\operatorname{Re}((\nu_b-1)\mathcal P_k)\to-\infty$ for $b\ne b_0$, while $\Xi^{-1}v_0=e_{b_0}$.  Hence $c=e_{b_0}$ and $y_*=\Phi_0^{(1)}$.

The covariance identity $(e^{2\pi ij}z)^kq^ju=z^ku$ shows that all these columns solve the same FFH equation at $u$.  Pulling back the asymptotic estimate for $\Phi_0^{(1)}$ gives the intervals below; on their intersection the estimates hold simultaneously.  For $0\le j<N$, the asymptotic interval of $\Phi_j^{(1)}$ is
\[
 I_j=\left(-N\pi-\frac\pi2-2\pi j,\,
             N\pi+\frac\pi2-2\pi j\right).
\]
Their intersection gives the common sector
\begin{equation}\label{eq:typeA-common-large-sector}
 \arg z\in\bigcap_{j=0}^{N-1}I_j
 \quad\Longleftrightarrow\quad
 \left|\arg z+\pi(N-1)\right|<\frac{3\pi}{2}.
\end{equation}

With $V=(v_0,\ldots,v_{N-1})$ and the scalars of Proposition~\ref{prop:FFH-log-transport}, their solution matrix is
\begin{equation}\label{eq:typeA-large-sector-factorization}
 \begin{aligned}
 Y_{\mathrm{cyc}}(z;u)
 &:=(\Phi_0^{(1)},\ldots,\Phi_{N-1}^{(1)})\\
 &=z^{\rho^\vee/N}(V+\varepsilon(z;u))
   \operatorname{diag}_{0\le j<N}
   \left(a_{1,j}(u)e^{-\gamma^j\mathcal P_k(z;u)}\right),
 \end{aligned}
\end{equation}
where $\varepsilon\to0$ on closed subsectors of \eqref{eq:typeA-common-large-sector}, locally uniformly in $u$.  The determinant factors satisfy
\[
 \det z^{\rho^\vee/N}=1,\qquad
 \prod_{j=0}^{N-1}e^{-\gamma^j\mathcal P_k}=1,\qquad
 \prod_{j=0}^{N-1}a_{1,j}(u)=d_N(u).
\]
Since $\mathscr A(z;u)\in\mathfrak{sl}_N$, the top exterior-power equation gives
\[
 \partial_z\det Y_{\mathrm{cyc}}
 =-\operatorname{tr}(\mathscr A(z;u))\det Y_{\mathrm{cyc}}=0.
\]
Consequently, along $\arg z=-\pi(N-1)$,
\[
 \det Y_{\mathrm{cyc}}
 =\lim_{|z|\to\infty}\det(V+\varepsilon(z;u))d_N(u)
 =\Delta_Nd_N(u)\ne0.
\]
For every $s\in\mathbb Z$, the matrix $Y_s:=(\Phi_s^{(1)},\ldots,\Phi_{s+N-1}^{(1)})$ is $\mathsf S_s^*Y_{\mathrm{cyc}}$, so
\[
 \det Y_s(z;u)=\Delta_Nd_N(q^su),
\]
which proves \eqref{eq:typeA-entire-cyclic-determinant}.

Lemma~\ref{lem:FFH-subdominant-solution} gives entire spectral dependence, and trivial monodromy at the additional singularities permits single-valued continuation to $\widetilde X_J$.  At a fixed $z_0\in\widetilde X_J$,
\[
 Y_s(z_0;u)^{-1}
 =\frac{\operatorname{adj}Y_s(z_0;u)}
        {\Delta_Nd_N(q^su)}
 \in\operatorname{Mat}_N(\mathcal O(\mathbb C)).
\]
Every entire family of solutions $y$ therefore has the unique expansion
\[
 y(z;u)=Y_s(z;u)c_y(u),\qquad
 c_y(u)=Y_s(z_0;u)^{-1}y(z_0;u)\in\mathcal O(\mathbb C)^N.
\]
Proposition~\ref{lem:loop-current-comparison} identifies this basis with an $\mathcal O'$-basis in the current realization.
\end{proof}

\begin{remark}\label{rem:typeA-volume-constant}
If $v_0=c(1,\ldots,1)^t$ in the standard matrix realization above, evaluation of the Vandermonde determinant gives
\[
 \Delta_N=c^N(-i)^{N(N-1)/2}N^{N/2}.
\]
Taking $c=e^{\pi i(N-1)/4}/\sqrt N$ makes $\Delta_N=1$.  When the tensor normalization prescribes a different value, we keep the resulting nonzero $\Delta_N$.  If $d_N$ is nonconstant, it cannot be removed by a constant choice of leading volume.  The volume correction in Section~6 therefore divides by the full determinant $\Delta_Nd_N(u)$.
\end{remark}

\subsection{The FFH tensor identity and its normalization}\label{sec:FFH-Phi-verification}
In this subsection we apply the abstract tensor identity to all untwisted simply-laced FFH opers on the generic locus and compute its scalar multiplier from the phase normalization.  Theorem~\ref{thm:FFH-verification} supplies the hypotheses for this application.  The result is the FFH identity of \cite[Proposition~3.22]{MR2024} and \cite[Theorem~3.6]{MRV2016} in our integer-shift convention.
\begin{corollary}\label{prop:FFH-integral-Phi-system}
Under the assumptions of Theorem~\ref{thm:FFH-verification}, the distinguished columns $\Phi_0^{(i)}$ of the normalized frame \eqref{eq:subdom} have shifted-monodromy iterates satisfying
\begin{equation}\label{eq:FFH-integral-Phi-system}
 m_i\!\left(\Phi_n^{(i)}\wedge\Phi_{n+1}^{(i)}\right)
 =\varphi_i(q^nu)\bigotimes_{j\sim i}\Phi_{n+\epsilon_i}^{(j)},
 \qquad i\in I,\quad n\in\mathbb Z.
\end{equation}
The multiplier is the nowhere-vanishing entire function
\begin{equation}\label{eq:FFH-log-Plucker-factor}
 \varphi_i(u)=\exp\!\left[-2\pi i\nu_i
   \bigl(\gamma_h-\epsilon_i(1+\gamma_h)\bigr)\lambda_{\log}(u)\right].
\end{equation}
It is one when $\lambda_{\log}=0$.
\end{corollary}
\begin{proof}
By Theorem~\ref{thm:FFH-verification}, Proposition~\ref{prop:derived-Phi-system} applies and gives the stated identity with multiplier \eqref{eq:asymptotic-tensor-multiplier}.  Proposition~\ref{prop:FFH-log-transport} gives $a_{i,1}(u)=\exp(-2\pi i\gamma_h\nu_i\lambda_{\log}(u))$.
Substituting these factors into \eqref{eq:asymptotic-tensor-multiplier}, we find
\[
 \varphi_i=\exp\!\left[-2\pi i\gamma_h\lambda_{\log}
           \left(\nu_i-\epsilon_i\sum_{j\sim i}\nu_j\right)\right].
\]
Using the eigenvalue identity \eqref{eq:oriented-cyclic-eigenvalues}, we obtain \eqref{eq:FFH-log-Plucker-factor}.
\end{proof}

\paragraph{Tensor identities with coefficient one.}
Dividing each iterate by its leading scalar $a_{i,n}$ removes the multiplier in the tensor identity.  We compare these rescaled columns with those obtained by successive pullbacks.  For the FFH columns normalized by \eqref{eq:subdom}, the rescaling is
\begin{equation}\label{eq:FFH-unit-normalized-orbit}
 \Phi_{n,\mathrm{unit}}^{(i)}(z;u)
 =\exp\bigl(2\pi i n\gamma_h^n\nu_i\lambda_{\log}(u)\bigr)
   \Phi_n^{(i)}(z;u).
\end{equation}

\begin{proposition}\label{prop:raw-unit-normalization}
These solutions have leading vector $\xi_{i,n}$ with phase
$e^{-\gamma_h^n\nu_i\mathcal P_k}$ and satisfy
\[
 m_i\bigl(\Phi_{n,\mathrm{unit}}^{(i)}\wedge
              \Phi_{n+1,\mathrm{unit}}^{(i)}\bigr)
 =\bigotimes_{j\sim i}\Phi_{n+\epsilon_i,\mathrm{unit}}^{(j)}.
\]
Successive rescaled columns are related by
\begin{equation}\label{eq:FFH-unit-orbit-pullback}
 \Phi_{n+1,\mathrm{unit}}^{(i)}
 =e^{2\pi i\gamma_h^{n+1}\nu_i\lambda_{\log}(u)}
     \mathsf S_1^*\Phi_{n,\mathrm{unit}}^{(i)}.
\end{equation}
The extra factor in this formula is $1$ when $\lambda_{\log}=0$.
\end{proposition}
\begin{proof}
By the phase identity \eqref{eq:FFH-phase-integer-rotation}, the exponential factor in \eqref{eq:FFH-unit-normalized-orbit} cancels the additive logarithmic constant.  In the tensor identity for the rescaled columns the coefficient is $\varphi_i(q^nu)\prod_{j\sim i}a_{j,n+\epsilon_i}(u)/(a_{i,n}(u)a_{i,n+1}(u))$.  The logarithmic exponent of the ratio is
\[
 2\pi i\gamma_h^n\nu_i
       \bigl(\gamma_h-\epsilon_i(1+\gamma_h)\bigr)\lambda_{\log}(u).
\]
In this calculation we used \eqref{eq:oriented-cyclic-eigenvalues}.  Proposition~\ref{prop:derived-Phi-system} identifies the tensor coefficient before rescaling with \eqref{eq:FFH-log-Plucker-factor}, so the displayed exponent cancels that coefficient at $q^nu$.  To obtain the pullback formula, we use $\lambda_{\log}(qu)=\gamma_h\lambda_{\log}(u)$ in the rescaling, which gives \eqref{eq:FFH-unit-orbit-pullback}.
\end{proof}

The next example asks whether tensor coefficients equal to one can be obtained by rescaling the initial column and then taking its shifted iterates.
\begin{example}\label{ex:typeA-logarithmic-volume}
For $h=2$, $k=3/4$, and $q=i$, shifted monodromy gives
\[
 \Phi_n^{(1)}(z;u)=\mathsf S_n^*\Phi_0^{(1)}(z;u)
 =\Phi_0^{(1)}(e^{2\pi in}z;i^nu),\qquad n\in\mathbb Z,
\]
with $z$ on the universal cover.  A nonzero meromorphic rescaling $b(u)$ of the initial solution determines the rescaling of every iterate through
\[
 \mathsf S_n^*\bigl(b(u)\Phi_0^{(1)}\bigr)
 =b(i^nu)\Phi_n^{(1)}.
\]
Here $\lambda_{\log}(u)=-u^2/8$, so \eqref{eq:FFH-log-Plucker-factor} gives $\varphi_1(u)=e^{-\pi i u^2/4}$.  Tensor coefficients equal to one for the rescaled solutions would therefore require
\begin{equation}\label{eq:rank-one-normalization-obstruction}
 b(u)b(iu)=e^{\pi i u^2/4}.
\end{equation}
On a nonempty open $i$-stable domain, the four spectral shifts of this identity imply
\[
 1=\frac{[b(u)b(iu)][b(-u)b(-iu)]}
         {[b(iu)b(-u)][b(-iu)b(u)]}
   =e^{\pi i u^2},
\]
which is impossible.  Thus the columns with tensor coefficients equal to one need not be the shifted-monodromy iterates of a single rescaled solution.
\end{example}

\begin{remark}\label{rem:FFH-normalization-comparison}
Here we compare our central connection with the matrices of \cite[Section~3.3]{MR2024} after matching the ordered eigenbases and sectorial intervals.  At a vertex $i$ of color $\epsilon_i\in\{0,1\}$, the cited convention rotates the connection through $\epsilon_i/2$ of a full turn \cite[Table~3 and (3.19)]{MR2024}.  The rotation factor $e^{\pi i\epsilon_i}=(-1)^{\epsilon_i}$ multiplies the coefficient of $zv_\theta$: after the current-to-loop substitution of Proposition~\ref{lem:loop-current-comparison}, the highest-root term is $(e^{\pi i\epsilon_i}z+z^ku)v_\theta$ \cite[(3.6)]{MR2024}.  This node-dependent rotation is distinct from the substitution $\lambda_{\mathrm{aff}}=z^ku$; we call its convention centered below.

The comparison uses a further pullback $z\mapsto e^{\pi i\epsilon_i}z$ at fixed $\lambda_{\mathrm{aff}}$.  Since $u=z^{-k}\lambda_{\mathrm{aff}}$, this replaces $u$ by $q^{\epsilon_i/2}u$, where $q^{\epsilon_i/2}=e^{-\pi ik\epsilon_i}$.  The highest-root coefficient then becomes
\[
 e^{\pi i\epsilon_i}e^{\pi i\epsilon_i}z
 +(e^{\pi i\epsilon_i}z)^kq^{\epsilon_i/2}u
 =z+z^ku.
\]
The additional-pole terms return to those of \eqref{eq:ffh} as well, since the two rotations together give $\epsilon_i$ full turns.

For this comparison, the centered Frobenius residue eigenvectors are $\rho_i(g_0^{-1})v_{\mu,a}$.  At infinity, the centered ordered eigenbasis is obtained by applying $\rho_i(e^{\pi i\epsilon_i\mathbf h/h})$ to the reference eigenbasis, including within repeated eigenspaces.  Its separating interval is $\mathcal I_i=\mathcal I_{\mathrm{ref}}+\pi\epsilon_i/h$, and its phase primitive is $\mathcal P_k(z;e^{-\pi i\epsilon_i}u)$, with the additive constant fixed by \eqref{eq:FFH-exact-phase-primitive}.  The matrix $Q_{\mathrm{match}}^{(i)}$ expresses the sectorial frame just specified in the centered Frobenius frame, with both frames in the original differential gauge.  The pullback changes the leading terms at zero and infinity by scalar factors.  Restoring the prescribed normalizations gives
\begin{equation}\label{eq:complete-frame-normalized-central-matrix}
 \begin{aligned}
 \rho_i(C(u))
 &=D_{\mathrm F,i}^{-1}Q_{\mathrm{match}}^{(i)}(q^{\epsilon_i/2}u)D_{\infty,i}(u),\\
 D_{\mathrm F,i}
 &=\exp\!\bigl(\pi i\epsilon_i\rho_i(\ell)\bigr),\\
 D_{\infty,i}(u)
 &=\exp\!\bigl(\pi i\epsilon_i\lambda_{\log}(u)\rho_i(\Lambda_H)\bigr).
 \end{aligned}
\end{equation}
Indeed, the pullback multiplies the Frobenius leading term of weight $\mu$ by $e^{-\pi i\epsilon_i\mu(\ell)}$.  Substitution in the fixed phase primitive gives the multiplier $e^{-\pi i\epsilon_i\mu(\Lambda_H)\lambda_{\log}(u)}$ on the sectorial exponential labeled by $\mu$.  The two displayed factors restore the prescribed normalizations, and uniqueness of the normalized frames gives \eqref{eq:complete-frame-normalized-central-matrix}.  When $\lambda_{\log}=0$, the correction at infinity is the identity.

To pass from $Q_{\mathrm{match}}^{(i)}$ to the particular transition matrix in \cite[Section~3.3, (3.47)--(3.49)]{MR2024}, one must also change from the matched sectorial frame to the frame used there.  This further change is the Stokes transition between the two sectorial frames; it is separate from the normalization factors displayed above.
\end{remark}

\section{The central connection, its flags, and functional relations}\label{sec:central-functional}
The sectorial and Frobenius frames select two flags in each spectral fiber. In this section we describe how shifted monodromy acts on these flags and how the central connection relates the highest- and lowest-weight lines in the fundamental representations.  We derive the $(G,q)$-Wronskian equation for the central connection and the full $QQ$-system for its Frobenius coefficients, and prove the nonvanishing of the generalized minors used later.  Throughout the section, the differential family satisfies Hypotheses~\ref{hyp:frob}--\ref{hyp:ccm}, and we work over $K=\mathcal M(U)$.

\subsection{Coxeter normalization of the transition in the sectorial frame}\label{sec:infinity-Coxeter-normalization}
Here we will derive the one-step $(G,q)$-Wronskian equation from the coordinates of $\Phi_1^{(i)}$ in the sectorial frame, namely $\mathscr M(u)^{-1}v_{\omega_i}$. We will also identify the Cartan factor in these coordinates.

Theorem~\ref{prop:analytic-sectorial-Steinberg} gives
\begin{equation}\label{eq:asymptotic-pullback-unipotent}
 \begin{aligned}
 \mathscr M(u)^{-1}&=s_{\mathrm{as}}(u)^{-1}n_+(u),\\
 n_+(u)&:=s_{\mathrm{as}}(u)\mathscr M(u)^{-1}\in U_c(K).
 \end{aligned}
\end{equation}
Since $U_c\subset N_+$ fixes every fundamental highest-weight vector,
\begin{equation}\label{eq:infinity-first-column}
 \mathscr M(u)^{-1}v_{\omega_i}=s_{\mathrm{as}}(u)^{-1}v_{\omega_i}.
\end{equation}
The scalars $a_{i,1}(u)$ in \eqref{eq:infinity-normalized-columns}, corresponding to one application of shifted monodromy, are the inverse fundamental-weight characters of $h_{\mathrm{as}}$:
\begin{equation}\label{eq:intrinsic-Cartan-normalization}
 \omega_i(h_{\mathrm{as}}(u))=a_{i,1}(u)^{-1}\quad(i\in I),\qquad
 s_{\mathrm{as}}(u)=h_{\mathrm{as}}(u)\dot c.
\end{equation}
Thus these scalars express the tensor multiplier \eqref{eq:asymptotic-tensor-multiplier} directly in terms of the Cartan factor:
\begin{equation}\label{eq:Plucker-factor-character-formula}
 \varphi_i(u)=\omega_i(h_{\mathrm{as}}(u))^{-1}
                 \prod_{j\sim i}\omega_j(h_{\mathrm{as}}(u))^{\epsilon_i}.
\end{equation}
The differential tensor identity of Proposition~\ref{prop:general-Phi-from-Kostant} has the following form in sectorial coordinates.
\begin{proposition}\label{prop:sectorial-Plucker-identity}
For every fundamental node $i$,
\begin{equation}\label{eq:Plucker-factor-definition}
 m_i(v_{\omega_i}\wedge s_{\mathrm{as}}(u)^{-1}v_{\omega_i})
 =\varphi_i(u)\bigotimes_{j\sim i}s_{\mathrm{as}}(u)^{-\epsilon_i}v_{\omega_j}.
\end{equation}
\end{proposition}
\begin{proof}
Equivariance of $m_i$ under $a_\Lambda^{-1}$, applied to \eqref{eq:normalized-cyclic-Plucker-vector} at $n=0$, and \eqref{eq:Kostant-Coxeter-orbit-vector} give
\[
 m_i(v_{\omega_i}\wedge\dot c^{-1}v_{\omega_i})
 =\bigotimes_{j\sim i}\dot c^{-\epsilon_i}v_{\omega_j}.
\]
Since $s_{\mathrm{as}}=h_{\mathrm{as}}\dot c$ and $\epsilon_i\in\{0,1\}$,
\[
 s_{\mathrm{as}}^{-1}v_{\omega_i}
 =\omega_i(h_{\mathrm{as}})^{-1}\dot c^{-1}v_{\omega_i},
 \qquad
 s_{\mathrm{as}}^{-\epsilon_i}v_{\omega_j}
 =\omega_j(h_{\mathrm{as}})^{-\epsilon_i}
      \dot c^{-\epsilon_i}v_{\omega_j}.
\]
Substitution into the preceding identity gives \eqref{eq:Plucker-factor-definition} with the multiplier \eqref{eq:Plucker-factor-character-formula}.
\end{proof}

For FFH opers, the phase computation of Corollary~\ref{lem:FFH-common-reference-sector} gives the Cartan factor explicitly and hence determines the ordered product below.  The related rescaling of the pulled-back solutions that gives tensor coefficients equal to one is described in Proposition~\ref{prop:raw-unit-normalization}.

\begin{proposition}\label{prop:FFH-phase-transport}
For the FFH frame normalized on \eqref{eq:FFH-reference-wide-sector} with leading scalars \eqref{eq:FFH-actual-leading-scalar}, the normalization factors are
\begin{equation}\label{eq:FFH-log-Cartan-factor}
 h_{\mathrm{as}}(u)=\exp(2\pi i\gamma_h\lambda_{\log}(u)\Lambda_H),\qquad
 s_{\mathrm{as}}(u)=h_{\mathrm{as}}(u)\dot c.
\end{equation}
Here $\Lambda_H=\operatorname{Ad}_{a_\Lambda^{-1}}\Lambda$.  For every $n\ge0$,
\begin{equation}\label{eq:FFH-log-full-leading-orbit}
 s_{\mathrm{as}}(u)^{-1}\cdots s_{\mathrm{as}}(q^{n-1}u)^{-1}
 =\dot c^{-n}\exp(-2\pi i n\gamma_h^n\lambda_{\log}(u)\Lambda_H).
\end{equation}
\end{proposition}
\begin{proof}
Since $\omega_i(\Lambda_H)=\nu_i$, the scalar formula \eqref{eq:FFH-actual-leading-scalar} gives the values $a_{i,1}^{-1}$ of the fundamental-weight characters on the displayed Cartan factor.  These characters determine the element of $H$, so \eqref{eq:intrinsic-Cartan-normalization} gives \eqref{eq:FFH-log-Cartan-factor}.

The product formula holds for $n=0$.  Its induction step uses $\lambda_{\log}(qu)=\gamma_h\lambda_{\log}(u)$ and $\operatorname{Ad}_{\dot c}\Lambda_H=\gamma_h\Lambda_H$:
\begin{align*}
 &\dot c^{-n}e^{-2\pi i n\gamma_h^n\lambda_{\log}(u)\Lambda_H}
     s_{\mathrm{as}}(q^nu)^{-1}\\
 &\qquad=\dot c^{-n}e^{-2\pi i n\gamma_h^n\lambda_{\log}(u)\Lambda_H}
     \dot c^{-1}e^{-2\pi i\gamma_h^{n+1}\lambda_{\log}(u)\Lambda_H}\\
 &\qquad=\dot c^{-n-1}e^{-2\pi i(n+1)\gamma_h^{n+1}\lambda_{\log}(u)\Lambda_H}.
 \qedhere
\end{align*}
\end{proof}
By Lemma~\ref{prop:finite-infinity-transport}, this product applied to $v_{\omega_i}$ equals $\mathscr M^{[n]}(u)^{-1}v_{\omega_i}$ for $n=0,1$ and, when $h$ is even, for $0\le n\le h/2$.

\subsection{Borel reductions and the one-step \texorpdfstring{$(G,q)$}{(G,q)}-Wronskian equation}\label{sec:differential-Wronskian-transport}
Here we will describe the two Borel reductions selected by the sectorial and Frobenius normalizations and explain their behavior under the $q$-connection. We first recall how a Borel reduction of a principal $G$-bundle over the spectral domain is represented in a chosen frame.
In the following definitions and proposition, $U$ is the spectral domain of Section~\ref{sec:general-analytic-setting}, and $B\subset G$ is an arbitrary fixed Borel subgroup.  Bundles, reductions, and sections are holomorphic.

\begin{definition}\label{def:G-flag}
The \emph{flag variety} associated with $B$ is $G/B$, with the left action
\[
 a\cdot(gB)=agB\qquad(a,g\in G).
\]
A \emph{flag of type $G$} is a point of this variety.
\end{definition}

The stabilizer of the flag $gB$ is $gBg^{-1}$.  For $G=SL_N$, the flag variety is the variety of complete flags
\[
 0=F_0\subset F_1\subset\cdots\subset F_N=\mathbb C^N,
 \qquad \dim F_i=i.
\]
To describe such flags over the spectral domain, we use reductions of a principal bundle.

\begin{definition}\label{def:Borel-reduction}
For a principal right $G$-bundle $\pi:P\to U$, the \emph{flag bundle} associated with $B$ is $P/B\to U$, where $\operatorname{pr}_{B}:P\to P/B$ is the quotient by the right action of $B$.  A \emph{$B$-reduction} is a principal $B$-subbundle $Q\subset P$ for which the natural map
\[
 Q\times^{B}G\longrightarrow P,\qquad [p,g]\longmapsto pg,
\]
is an isomorphism.  In the contracted product, $(pb,g)$ is identified with $(p,bg)$ for $b\in B$.
\end{definition}

A Borel reduction determines a section of the flag bundle. The fundamental representations give a second description by line subbundles, using Proposition~\ref{prop:highest-weight-stabilizers}.
\begin{proposition}\label{prop:fundamental-flag-embedding}
\begin{enumerate}[label=(\alph*),topsep=3pt,itemsep=3pt]
\item The $B$-reductions $Q\subset P$ are in bijection with sections $s:U\to P/B$.  Corresponding $Q$ and $s$ satisfy
\[
 \begin{aligned}
 s(u)&=\operatorname{pr}_{B}(p),\qquad p\in Q_u,\\
 Q&=\{p\in P:\operatorname{pr}_{B}(p)=s(\pi(p))\}.
 \end{aligned}
\]
On an open subset $U_0\subset U$ with a trivialization $P|_{U_0}\simeq U_0\times G$, the section has the form $s(u)=(u,F(u))$ for a holomorphic map $F:U_0\to G/B$.
\item There is a $g_0\in G$ with $B=g_0B_+g_0^{-1}$, and the unique $B$-invariant line in $V_i$ is $L_i=\rho_i(g_0)\mathbb Cv_{\omega_i}$.  A $B$-reduction $Q$ determines the line subbundles
\[
 Q\times^{B}L_i\subset P\times^G V_i,
 \qquad i\in I.
\]
Conversely, a collection of line subbundles of $P\times^G V_i$ determines a $B$-reduction precisely when, in every trivialized fiber, its tuple of lines belongs to the image of the closed embedding
\[
 G/B\hookrightarrow\prod_{i\in I}\mathbb P(V_i),\qquad
 gB\longmapsto\bigl(\rho_i(g)L_i\bigr)_{i\in I}.
\]
\end{enumerate}
\end{proposition}
\begin{proof}
For a reduction $Q$, the fiber $Q_u$ is one right $B$-orbit, so $\operatorname{pr}_{B}(p)$ is independent of $p\in Q_u$ and defines a section $s$.  Locally a section is $s(u)=(u,g(u)B)$, and its inverse image under $\operatorname{pr}_{B}$ is the principal $B$-subbundle with fibers $g(u)B$.  The two constructions are therefore inverse.

Conjugacy of Borel subgroups gives $g_0$; see \cite{Milne2017}.  The formula for $L_i$ follows from the unique $B_+$-invariant line in $V_i$.  Under the identification $gB\mapsto gg_0B_+$, the displayed embedding is that of Proposition~\ref{prop:highest-weight-stabilizers}.  Compatible line subbundles therefore give a section of the flag bundle, and part~(a) gives the corresponding reduction.  The converse follows from the associated-bundle construction.
\end{proof}
For $B=B_+$, the lines $L_i$ are $\mathbb Cv_{\omega_i}$; for $B=B_-$, they are $\mathbb Cv_{w_0\omega_i}$.  For meromorphic data, these constructions are applied on open subsets where the data are holomorphic.  For $SL_N$, the fiber lines are $\bigwedge^iF_i\subset\bigwedge^i\mathbb C^N$; the embedding in Proposition~\ref{prop:fundamental-flag-embedding} expresses the nesting of the subspaces $F_i$.

For the trivial bundle $P=U\times G$, the correspondence between a reduction $Q$, its section $s_Q$, and its representing map $F$ is
\[
 s_Q(u)=(u,F(u)),\qquad
 Q=\{(u,g)\in U\times G:gB=F(u)\}.
\]
We write $Q_\pm$ for reductions to the fixed opposite Borel subgroups $B_\pm$, and $F_\pm$ for their representing flag maps.

We compare two flags of type $B_+$ by their relative position. For flags $xB_+$ and $yB_+$ in a trivialized fiber, our convention is
\[
 \operatorname{pos}(xB_+,yB_+)=w
 \quad\Longleftrightarrow\quad
 x^{-1}y\in B_+\dot wB_+,
\]
where $\dot w\in N_G(H)$ represents $w\in W$; the condition depends only on the two flags.

To compare a $B_+$-reduction with a $B_-$-reduction, we will use opposition. Flags $xB_+\in G/B_+$ and $yB_-\in G/B_-$ are \emph{opposite} if the pair is a simultaneous $G$-translate of $(B_+,B_-)$, equivalently if $x^{-1}y\in B_+B_-$.  Two reductions $Q_+$ and $Q_-$ are opposite at $u$ if their associated flags are opposite in the fiber over $u$.  For $SL_N$, this means
\[
 F_j^+(u)\oplus F_{N-j}^-(u)=\mathbb C^N
 \qquad(1\le j<N)
\]
in any trivialization, where $F_j^\pm(u)$ are the subspaces of the two complete flags.  The reductions are generically opposite if this holds on a dense open subset of $U$.

\begin{definition}\label{def:q-connection}
A meromorphic \emph{$q$-connection} $\mathcal A$ on a principal right $G$-bundle $P\to U$ is a meromorphic family of $G$-equivariant isomorphisms $\mathcal A_u:P_u\to P_{qu}$. In a trivialization its matrix $A(u)$ carries $(u,g)$ to $(qu,A(u)g)$. A meromorphic change of fiber coordinates $g'=\eta(u)g$ acts on the connection matrix by the \emph{$q$-gauge transformation} $A(u)\mapsto\eta(qu)A(u)\eta(u)^{-1}$. The connection is \emph{$Z$-twisted} if its matrix in some meromorphic trivialization is the constant $Z\in H$.

For the fixed Coxeter element $c\in W$, a \emph{$(G,q)$-oper} is a triple $(P,\mathcal A,Q_+)$ consisting of this bundle, a $q$-connection, and a meromorphic $B_+$-reduction. If $F_+$ represents the reduction, then at generic $u$ the flag over $qu$ and the image of the flag over $u$ must have relative position
\[
 \operatorname{pos}\bigl(F_+(qu),A(u)F_+(u)\bigr)=c.
\]
\end{definition}

\begin{definition}\label{def:local-twisted-Miura-qoper}
A \emph{Miura $(G,q)$-oper} is a quadruple $(P,\mathcal A,Q_+,Q_-)$, where $(P,\mathcal A,Q_+)$ is a $(G,q)$-oper and $Q_-$ is a meromorphic $B_-$-reduction of $P$ preserved by $\mathcal A$. Thus, wherever the data are regular,
\[
 \mathcal A_u\bigl((Q_-)_u\bigr)=(Q_-)_{qu}.
\]
The reduction $Q_-$ is called the \emph{Miura reduction}. In a trivialization the quadruple is represented by $(A,F_+,F_-)$, where $F_+$ satisfies the oper condition and the preservation condition on $F_-$ reads $A(u)F_-(u)=F_-(qu)$. The Miura $(G,q)$-oper is $Z$-twisted if its underlying $q$-connection is $Z$-twisted.
\end{definition}

Both flag conditions are independent of the trivialization. A change of coordinates $g'=\eta(u)g$ replaces $F_\pm(u)$ by $\eta(u)F_\pm(u)$. Hence each flag being compared in the fiber over $qu$ is multiplied on the left by the same element $\eta(qu)$.

In an associated representation $V$, the induced fiber map sends $a\in V$ to $\rho_V(A(u))a$. A meromorphic section with coordinate column $a(u)$ is \emph{horizontal} if $a(qu)=\rho_V(A(u))a(u)$.

For $G=SL_N$, the oper reduction determines a complete flag of subbundles in the bundle associated with the defining representation. Its rank-one subbundle is called the \emph{oper line}.

\begin{remark}\label{rem:adapted-oper-frames}
A local frame of the principal bundle is a local section $e:U_0\to P$; it identifies $P|_{U_0}$ with $U_0\times G$ by $p=e(u)g$.  The frame is adapted to $Q_+$ when $e(u)\in(Q_+)_u$, so that $Q_+$ has coordinates $U_0\times B_+$.  In frames adapted at both $u$ and $qu$, the flags in the oper condition have coordinates $B_+$ and $A(u)B_+$.  Thus the condition is $A(u)\in B_+\dot cB_+$.  Reversing the ordered pair gives relative position $c^{-1}$.  For $SL_2$, writing the chosen line as $L(u)$, the condition becomes
\[
 A(u)L(u)\ne L(qu).
\]
In an adapted basis this means that the lower-left matrix entry of $A(u)$ is nonzero.
The oper relative-position condition holds at spectral values where the data are regular and the minors specifying this position are nonzero.
\end{remark}

To compare sections over $u$ and $qu$, we use the inverse $q$-connection to bring the value at $qu$ back to $u$. For a section with coordinate column $v(u)$ in an associated representation $\rho_V$, we denote this pullback operator by $\sigma_A$. Its coordinate formula is
\[
 (\sigma_A v)(u)=\rho_V(A(u))^{-1}v(qu).
\]
It is semilinear over the shift $\sigma_q f(u)=f(qu)$ on $K=\mathcal M(U)$.

The $(G,q)$-Wronskian will compare selected sections of the oper bundle with sections transported by the constant connection $Z$. We therefore introduce a reference bundle carrying this connection.

The reference bundle $P_Z=U\times G$ is a separate copy of the trivial principal bundle. Its $q$-connection is given in the product trivialization by
\[
 (u,g)\longmapsto(qu,Zg),\qquad Z\in H.
\]
In the associated bundle $P_Z\times^G V_i$, the connection carries a vector with column $a$ at $u$ to the vector with column $\rho_i(Z)a$ at $qu$. The corresponding pullback $\sigma_Z$ takes a section with column $a(u)$ to the section with column $\rho_i(Z)^{-1}a(qu)$.

In this trivialization, a meromorphic section of $P_Z$ has the form $u\mapsto(u,\mathscr W(u))$. Using a meromorphic frame of $P$ adapted to the oper reduction, the section determines a meromorphic bundle isomorphism $P\to P_Z$ given by $(u,g)\mapsto(u,\mathscr W(u)g)$. In the associated fundamental bundles, the fixed highest-weight vector $v_{\omega_i}$ maps to the section with column $\rho_i(\mathscr W(u))v_{\omega_i}$, which generates a meromorphic line subbundle of $P_Z\times^G V_i$. A meromorphic Coxeter lift $s\in H(K)\dot c$ gives a second section, the image of $\rho_i(s(u)^{-1})v_{\omega_i}$ under the same map.

\begin{definition}\label{def:generalized-Gq-Wronskian}
A meromorphic section $\mathscr W$ of $P_Z$ is a generalized $(G,q)$-Wronskian for $(Z,s)$ if $\sigma_Z$ takes the section with column $\rho_i(\mathscr W(u))v_{\omega_i}$ to the section with column $\rho_i(\mathscr W(u)s(u)^{-1})v_{\omega_i}$ for every $i\in I$, as in \cite[Definition~4.9]{KZ2021}. In the chosen trivializations, with $\rho_i$ suppressed, this is the equality
\begin{equation}\label{eq:generalized-Wronskian-definition}
 Z^{-1}\mathscr W(qu)v_{\omega_i}
 =\mathscr W(u)s(u)^{-1}v_{\omega_i}\qquad(i\in I).
\qedhere
\end{equation}
\end{definition}

Every $Z$-twisted $(G,q)$-oper admits a generalized $(G,q)$-Wronskian for a lift from its Bruhat factorization. To obtain the $(G,q)$-Wronskian equation from a trivialization to $Z$, we choose an adapted frame in which the inverse connection sends each $v_{\omega_i}$ to $s(u)^{-1}v_{\omega_i}$. The right unipotent factor in a Bruhat factorization can be removed by an adapted frame change.

In a meromorphic frame adapted to the oper reduction, the condition $A\in B_+(K)\dot c B_+(K)$ gives $A=n_L s n_R$ with $n_L,n_R\in N_+(K)$ and $s\in H(K)\dot c$. Changing the frame by right multiplication by $n_R^{-1}$ leaves it adapted and gives the connection matrix
\[
 A'(u)=n_R(qu)A(u)n_R(u)^{-1}=n_R(qu)n_L(u)s(u).
\]
A chosen trivialization to $Z$ has a matrix $g$ in the original frame satisfying $g(qu)A(u)=Zg(u)$. After the adapted change of frame, its matrix is $\mathscr W=g n_R^{-1}$. The identity $Z^{-1}\mathscr W(qu)=\mathscr W(u)A'(u)^{-1}$ then gives
\[
 Z^{-1}\mathscr W(qu)v_{\omega_i}
 =\mathscr W(u)A'(u)^{-1}v_{\omega_i}
 =\mathscr W(u)s(u)^{-1}v_{\omega_i},
\]
because $N_+$ fixes the highest-weight vectors.

For the differential family of Section~\ref{sec:general-analytic-setting} satisfying Hypotheses~\ref{hyp:frob}--\ref{hyp:ccm}, the lift is $s=s_{\mathrm{as}}$. We use the normalized sectorial frame for $P$ and the normalized Frobenius frame for $P_Z$. The $q$-connection matrices in these frames are $\mathscr M$ and $Z$, respectively, and \eqref{eq:central-infinity-transport} relates them through the central connection $C$. Thus $\mathscr M(u)$ carries sectorial coordinates from $u$ to $qu$.
Replacing $\mathcal Y_{\mathrm{as}}$ by $\mathcal Y_{\mathrm{as}}\eta(u)^{-1}$ sends $C$ to $C\eta^{-1}$ and gives
\[
 \mathscr M(u)\longmapsto\eta(qu)\mathscr M(u)\eta(u)^{-1}.
\]

The central connection defines the section $u\mapsto(u,C(u))$ of $P_Z$.  Since $\mathcal Y_{\mathrm{as},V_i}=\mathcal Y_{\mathrm F,V_i}\rho_i(C)$, its highest-weight columns are the Frobenius coefficient sections of the corresponding sectorial solutions.

The two frames also select the reductions. We use the highest-weight lines in the sectorial frame and the lowest-weight lines in the Frobenius frame. Expressed in sectorial coordinates, their flags are given by
\[
 F_+(u)=B_+\in G/B_+,
 \qquad
 F_-(u)=C(u)^{-1}B_-\in G/B_-.
\]
The lowest-weight lines in the Frobenius frame have sectorial coefficient lines $\mathbb C\rho_i(C(u))^{-1}v_{w_0\omega_i}$.  The frame identity \eqref{eq:central-infinity-transport} gives the $(G,q)$-Wronskian relation for $C$ and the $q$-oper structure below.

\begin{proposition}\label{prop:one-step-differential-W}\label{prop:direct-qoper}
For every fundamental representation,
\begin{equation}\label{eq:C-basic-W-relation}
 Z^{-1}C(qu)v_{\omega_i}=C(u)s_{\mathrm{as}}(u)^{-1}v_{\omega_i}.
\end{equation}
Equivalently,
\begin{equation}\label{eq:BC-oper-line}
 \mathscr M(u)s_{\mathrm{as}}(u)^{-1}v_{\omega_i}=v_{\omega_i}.
\end{equation}
Thus $C$ is a generalized $(G,q)$-Wronskian for $(Z,s_{\mathrm{as}})$.  The triple $(\mathscr M,F_+,F_-)$ represents a $Z$-twisted Miura $(G,q)$-oper, with the sectorial oper reduction and the preserved Frobenius reduction.
\end{proposition}
\begin{proof}
The inverse form of \eqref{eq:central-infinity-transport}, combined with \eqref{eq:infinity-first-column}, gives \eqref{eq:C-basic-W-relation}.  Multiplication by $C(qu)^{-1}Z$ gives \eqref{eq:BC-oper-line}.  The $q$-oper condition follows directly from Theorem~\ref{prop:analytic-sectorial-Steinberg}:
\[
 \mathscr M(u)\in U_c(K)s_{\mathrm{as}}(u)\subset B_+(K)\dot cB_+(K).
\]
The Frobenius reduction is preserved because its representing flag map satisfies
\[
 \mathscr M(u)C(u)^{-1}B_-=C(qu)^{-1}Z B_-=C(qu)^{-1}B_-.
\]
Finally, \eqref{eq:central-infinity-transport} gives $C(qu)\mathscr M(u)C(u)^{-1}=Z$, which identifies the twist.
\end{proof}

\subsection{Central minors and the full \texorpdfstring{$QQ$}{QQ}-system}
Here we derive the full $QQ$-system from the differential tensor identity~\eqref{eq:hyp-integral-Phi-system} by expressing the distinguished solutions in the Frobenius frame.

For a $V$-valued solution $y(z;u)$, its coefficient column in the Frobenius frame is
\[
 [y]_{\mathcal Y_{\mathrm F,V}}(u)
 :=\mathcal Y_{\mathrm F,V}(z;u)^{-1}y(z;u).
\]
Its entries are the coefficients of $y$ in the solution basis formed by the columns of $\mathcal Y_{\mathrm F,V}$.  Since the frame and $y$ solve the same differential equation, this column is independent of $z$, and
\[
 y(z;u)=\mathcal Y_{\mathrm F,V}(z;u)
              [y]_{\mathcal Y_{\mathrm F,V}}(u).
\]

By \eqref{eq:group-frame-central-ratio}, the coefficient column of the distinguished solution at spectral value $q^nu$ is
\begin{equation}\label{eq:full-shifted-central-columns}
 [\Phi_0^{(i)}]_{\mathcal Y_{\mathrm F,V_i}}(q^nu)
 =\rho_i(C(q^nu))v_{\omega_i}.
\end{equation}
\begin{proposition}\label{prop:integral-shifted-columns}
For every $i\in I$ and $n\in\mathbb Z$, the Frobenius coefficient column of $\Phi_n^{(i)}$ satisfies
\begin{equation}\label{eq:fixed-fiber-shifted-column}
 [\Phi_n^{(i)}]_{\mathcal Y_{\mathrm F,V_i}}(u)
 =\rho_i(Z)^{-n}\rho_i(C(q^nu))v_{\omega_i}.
\end{equation}
\end{proposition}
\begin{proof}
Equations~\eqref{eq:actual-infinity-frame-column} and~\eqref{eq:group-frame-central-ratio} give
\[
 [\Phi_n^{(i)}]_{\mathcal Y_{\mathrm F,V_i}}(u)
 =\rho_i(C(u))\rho_i(\mathscr M^{[n]}(u)^{-1})v_{\omega_i}.
\]
Using \eqref{eq:central-infinity-iterate}, we obtain
\[
 \rho_i(C(u))\rho_i(\mathscr M^{[n]}(u)^{-1})v_{\omega_i}
 =\rho_i(Z)^{-n}\rho_i(C(q^nu))v_{\omega_i}.\qedhere
\]
\end{proof}

The columns of $\rho_i(\mathcal Y_{\mathrm F})$ are the Frobenius basis
$\{\chi_{\mu,a}\mid \mu\in\operatorname{Wt}(V_i),\ 1\le a\le\dim V_i[\mu]\}$,
normalized by $\rho_i(T_0)\chi_{\mu,a}\sim z^{-\mu(\ell)}v_{\mu,a}$ in their prescribed classes. In this basis the distinguished solution has expansion
\begin{equation}\label{eq:Q-expansion}
 \Phi_0^{(i)}(z;u)
 =\sum_{\mu\in\operatorname{Wt}(V_i)}
  \sum_{a=1}^{\dim V_i[\mu]}
  \mathcal Q_{\mu,a}^{(i)}(u)\chi_{\mu,a}(z;u).
\end{equation}
Here $\operatorname{Wt}(V_i)$ is the set of weights of $V_i$, and the coefficients are meromorphic by Hypothesis~\ref{hyp:ccm}.  Extremal weight spaces have multiplicity one, so we write
\begin{equation}\label{eq:central-minor-notation}
 \mathcal Q_w^{(i)}:=\mathcal Q_{w\omega_i}^{(i)},
 \qquad
 \mathcal Q_{ws_i}^{(i)}=\mathcal Q_{w(\omega_i-\alpha_i)}^{(i)},
 \qquad w\in W.
\end{equation}
Equation~\eqref{eq:full-shifted-central-columns} identifies these coefficients with
\begin{equation}\label{eq:C-minor-Q-u}
 \mathcal Q_w^{(i)}(u)=\Delta_{w\omega_i,\omega_i}(C(u)).
\end{equation}
The notations $\mathcal Q_+^{(i)}$ and $\mathcal Q_-^{(i)}$ refer to $w=1$ and $w=s_i$, respectively.

The shifted Frobenius coefficients acquire phases determined by their weights. For $\ell_W(ws_i)=\ell_W(w)+1$, we write $\tau_{w,i}^{\pm}$ for the two weights evaluated on $\ell$:
\begin{equation}\label{eq:central-minor-Xi}
 \begin{gathered}
 \tau_{w,i}^{+}:=\langle\ell,w\omega_i\rangle,
 \qquad
 \tau_{w,i}^{-}:=\langle\ell,w(\omega_i-\alpha_i)\rangle,\\
 \tau_{w,i}:=\tau_{w,i}^{+}+\tau_{w,i}^{-}
 =\sum_{j\sim i}\langle\ell,w\omega_j\rangle,\\
 \Xi_{w,i}:=e^{\pi i(1-2\epsilon_i)\tau_{w,i}}.
 \end{gathered}
\end{equation}
The relation for the original coefficients contains the constant $\Xi_{w,i}$. The following rescaling removes it. For $\mu\in W\omega_i$, we set
\begin{equation}\label{eq:standard-Q-normalization}
 Q_{\mu}^{(i)}(u)
 :=\exp\!\left[\pi i\left(\epsilon_i-\frac12\right)
                    \langle\ell,\mu\rangle\right]\mathcal Q_{\mu}^{(i)}(u).
\end{equation}
We write $Q_w^{(i)}:=Q_{w\omega_i}^{(i)}$, $Q_+^{(i)}:=Q_1^{(i)}$, and $Q_-^{(i)}:=Q_{s_i}^{(i)}$.

\begin{proposition}\label{prop:central-minor-QQ}\label{prop:standard-QQ-normalization}
For every pair $(w,i)$ with $\ell_W(ws_i)=\ell_W(w)+1$, the central minors satisfy
\begin{equation}\label{eq:central-minor-QQ}
\begin{aligned}
 &e^{\pi i\langle\ell,w(\alpha_i)\rangle}
  \mathcal Q_w^{(i)}(u)\mathcal Q_{ws_i}^{(i)}(qu)
 -e^{-\pi i\langle\ell,w(\alpha_i)\rangle}
  \mathcal Q_w^{(i)}(qu)\mathcal Q_{ws_i}^{(i)}(u)\\
 &\hspace{25mm}=
 \varphi_i(u)\,\Xi_{w,i}
 \prod_{j\sim i}\mathcal Q_w^{(j)}\!\left(q^{\epsilon_i}u\right).
\end{aligned}
\end{equation}
Their rescalings \eqref{eq:standard-Q-normalization} satisfy
\begin{equation}\label{eq:exact-QQ}
\begin{aligned}
 &e^{\pi i\langle\ell,w(\alpha_i)\rangle}
  Q_w^{(i)}(u)Q_{ws_i}^{(i)}(qu)
 -e^{-\pi i\langle\ell,w(\alpha_i)\rangle}
  Q_w^{(i)}(qu)Q_{ws_i}^{(i)}(u)\\
 &\hspace{25mm}=
 \varphi_i(u)
 \prod_{j\sim i}Q_w^{(j)}\!\left(q^{\epsilon_i}u\right).
\end{aligned}
\end{equation}
The two normalizations have the same zero divisors and ratios under spectral shifts.
\end{proposition}

\begin{proof}
At $n=0$, we compare the coefficients in the extremal weight $w(2\omega_i-\alpha_i)$ on both sides of \eqref{eq:hyp-integral-Phi-system}. The normalization \eqref{eq:extremal-Plucker-normalization} sends $v_{w\omega_i}\wedge v_{ws_i\omega_i}$ to $\bigotimes_{j\sim i}v_{w\omega_j}$ with coefficient one. Equation~\eqref{eq:fixed-fiber-shifted-column} gives the shifted coefficients in the resulting relation:
\[
 \begin{aligned}
 &e^{-2\pi i\tau_{w,i}^{-}}\mathcal Q_w^{(i)}(u)
       \mathcal Q_{ws_i}^{(i)}(qu)
 -e^{-2\pi i\tau_{w,i}^{+}}\mathcal Q_w^{(i)}(qu)
       \mathcal Q_{ws_i}^{(i)}(u)\\
 &\hspace{20mm}=\varphi_i(u)e^{-2\pi i\epsilon_i\tau_{w,i}}
       \prod_{j\sim i}\mathcal Q_w^{(j)}(q^{\epsilon_i}u).
 \end{aligned}
\]
Multiplication by $e^{\pi i\tau_{w,i}}$, with
$\tau_{w,i}^{+}-\tau_{w,i}^{-}=\langle\ell,w\alpha_i\rangle$, gives \eqref{eq:central-minor-QQ}.

The rescaling factors on the left and right are respectively
$e^{\pi i(\epsilon_i-1/2)\tau_{w,i}}$ and
$e^{\pi i(1/2-\epsilon_i)\tau_{w,i}}$, because adjacent vertices have opposite colors.  After substituting the rescaled coefficients, the remaining scalar on the right-hand side is $\Xi_{w,i}$ multiplied by the quotient of these two rescaling factors. Its exponent divided by $\pi i$ is
\[
 (1-2\epsilon_i)\tau_{w,i}
 +(\epsilon_i-\tfrac12)\tau_{w,i}
 -(\tfrac12-\epsilon_i)\tau_{w,i}=0,
\]
which proves \eqref{eq:exact-QQ}.  The rescaling factors are nonzero constants in $u$, so zero divisors and shifted ratios are unchanged.
\end{proof}

Equations~\eqref{eq:exact-QQ} for all pairs $(w,i)$ with $\ell_W(ws_i)=\ell_W(w)+1$ constitute the full $QQ$-system \cite[Definition~3.18]{KZ2021} in our normalization.  The rescaling is node-dependent and need not be induced by one Cartan element.  We retain $C$ for the connection between the prescribed differential frames; $Q$ denotes its rescaled coefficients, and the tensor multiplier $\varphi_i$ is unchanged.

With $\mathfrak z_i:=e^{\pi i\langle\ell,\alpha_i\rangle}$, the case $w=1$ is
\begin{equation}\label{eq:oriented-standard-QQ}
 \mathfrak z_iQ_+^{(i)}(u)Q_-^{(i)}(qu)
 -\mathfrak z_i^{-1}Q_+^{(i)}(qu)Q_-^{(i)}(u)
 =\varphi_i(u)
 \prod_{j\sim i}Q_+^{(j)}\!\left(q^{\epsilon_i}u\right).
\end{equation}
The colors record the Coxeter order: neighboring factors preceding $i$ occur at $u$, and those following $i$ at $qu$.  The coefficient $\varphi_i$ is given by \eqref{eq:Plucker-factor-character-formula} and equals one when $h_{\mathrm{as}}=1$.

The fundamental $QQ$-relation gives the Bethe equation at a zero of $Q_+^{(i)}$.
\begin{corollary}\label{cor:Bethe-equations}
For $q\ne0$, let $u_*$ be a zero of $Q_+^{(i)}$.  We assume that all factors in
\eqref{eq:oriented-standard-QQ} are holomorphic at $u_*$ and $q^{-1}u_*$ and that
\[
 Q_-^{(i)}(u_*),\quad Q_+^{(i)}(q^{-1}u_*),\quad
 \varphi_i(q^{-1}u_*),\quad
 Q_+^{(j)}(q^{\epsilon_i-1}u_*)\quad(j\sim i)
\]
are nonzero.  Then
\begin{equation}\label{eq:explicit-Bethe}
 -\mathfrak z_i^{-2}
 \frac{Q_+^{(i)}(qu_*)}{Q_+^{(i)}(q^{-1}u_*)}
 =\frac{\varphi_i(u_*)}{\varphi_i(q^{-1}u_*)}
   \prod_{j\sim i}
   \frac{Q_+^{(j)}(q^{\epsilon_i}u_*)}
        {Q_+^{(j)}(q^{\epsilon_i-1}u_*)}.
\end{equation}
The same equation follows from the unrescaled central minors.
\end{corollary}
\begin{proof}
Evaluation of \eqref{eq:oriented-standard-QQ} at $u_*$ and $q^{-1}u_*$
annihilates its first and second terms, respectively.  Dividing the resulting
identities eliminates $Q_-^{(i)}(u_*)$ and gives \eqref{eq:explicit-Bethe}.
The constant rescalings of the central minors cancel in all shifted ratios.
\end{proof}

\subsection{Generic position and nonvanishing of minors}
Here we will show that the extremal coefficients $\mathcal Q_w^{(i)}$ are nonzero meromorphic functions. For $w=1$, this will give the Gauss factorization of $C$ used in the Miura transformation.

Nonvanishing below means nonvanishing in $K=\mathcal M(U)$; the resulting factorizations hold where their coefficients are regular and the required minors do not vanish.

A Gauss factorization is a decomposition
\[
 C=vn_+,\qquad v\in B_-(K),\quad n_+\in N_+(K).
\]
For $SL_N$, this writes $C$ as a lower triangular matrix times an upper triangular matrix with diagonal entries equal to $1$.  Proposition~\ref{prop:Gauss-minor-criteria} characterizes its existence by nonvanishing of the principal generalized minors.

The generic-position result below yields the Gauss factorization from the $q$-oper condition and a preserved $B_-$ reduction. Here the field automorphism is the spectral shift $\sigma_q(f)(u)=f(qu)$.  The lemma is \cite[Theorem~2.3]{FKSZ2020}, with the Borels interchanged.  Its argument uses the algebraic Coxeter factorization of \cite[Theorem~2.5]{FKSZ2020} and applies to a field automorphism.

\begin{lemma}\label{lem:difference-field-generic-position}
For $K=\mathcal M(U)$, $\sigma=\sigma_q$, and
\[
 s_{\mathrm{as}}=h_{\mathrm{as}}\dot c,\qquad h_{\mathrm{as}}\in H(K),
\]
the conditions
\[
 A\in N_+(K)s_{\mathrm{as}} N_+(K),
 \qquad
 \sigma(g)A g^{-1}\in B_-(K)
\]
for $g\in G(K)$ imply $g\in B_-(K)N_+(K)$.
The same implication holds over any field extension $K/\mathbb C$ with a $\mathbb C$-automorphism $\sigma$, with $G$ obtained from the fixed based root datum by extension of scalars.
\end{lemma}

The geometric meaning of the conclusion is as follows.  The two Borel reductions have coordinates
\[
 F_+=B_+,\qquad F_-=g^{-1}B_-.
\]
Since $A\in B_+\dot cB_+$, the first flag represents the $q$-oper reduction.  Writing $D=\sigma(g)Ag^{-1}\in B_-$ gives
\[
 A(g^{-1}B_-)=\sigma(g)^{-1}DB_-
             =\sigma(g)^{-1}B_-=\sigma(F_-),
\]
so the reduction represented by the second flag is preserved.  They are opposite precisely when
\[
 g^{-1}\in B_+B_-\quad\Longleftrightarrow\quad
 g\in B_-B_+=B_-N_+.
\]
Thus the oper and Miura reductions are opposite on a dense open subset of the spectral domain.  They may fail to be opposite at individual spectral values.

When $\sigma(g)Ag^{-1}$ takes values in $H$, every flag $g^{-1}\widehat wB_-$ is preserved by $A$. The same generic-position argument therefore applies to each corresponding reduction.
\begin{corollary}\label{cor:Cartan-preservation-Weyl-minors}
Under the hypotheses of Lemma~\ref{lem:difference-field-generic-position}, the additional condition $\sigma(g)Ag^{-1}\in H(K)$ implies
\[
 \Delta_{w\omega_i,\omega_i}(g)\in K^\times
 \qquad(w\in W,\ i\in I).
\]
\end{corollary}
\begin{proof}
The minor $\Delta_{w\omega_i,\omega_i}(g)$ is the principal minor of $g_w:=\widehat w^{-1}g$, where $\widehat w$ is the fixed generalized-minor representative. The preservation condition still holds for $g_w$, since $\widehat w$ is constant and normalizes $H$:
\[
 \sigma(g_w)Ag_w^{-1}
 =\widehat w^{-1}\bigl(\sigma(g)Ag^{-1}\bigr)\widehat w
 \in H(K)\subset B_-(K).
\]
Lemma~\ref{lem:difference-field-generic-position} and Proposition~\ref{prop:Gauss-minor-criteria} give a unique factorization
\[
 g_w=n_{-,w}h_wn_{+,w},\qquad
 n_{\pm,w}\in N_\pm(K),\quad h_w\in H(K).
\]
The identities in Proposition~\ref{prop:minor-covariance} then give
\[
 \Delta_{w\omega_i,\omega_i}(g)
 =\Delta_{\omega_i,\omega_i}(g_w)
 =\omega_i(h_w)\in K^\times.
\]
\end{proof}

\begin{proposition}\label{prop:central-principal-nonvanishing}\label{prop:Bruhat-Q-nonvanishing}
For the Miura $(G,q)$-oper of Proposition~\ref{prop:direct-qoper}, every $w\in W$ and $i\in I$ satisfy
\begin{equation}\label{eq:all-Weyl-Q-nonzero}
 \mathcal Q_w^{(i)}
 =\Delta_{w\omega_i,\omega_i}(C)\in K^\times.
\end{equation}
The rescaled functions $Q_w^{(i)}$ are also nonzero in $K$.  In particular,
\begin{equation}\label{eq:principal-central-minors-nonzero}
 \Delta_{\omega_i,\omega_i}(C)\in K^\times
 \qquad(i\in I).
\end{equation}
Thus $C\in B_-(K)N_+(K)$ and has the unique Gauss factorization $C=vn_+$ over $K$.
\end{proposition}

\begin{proof}
The spectral shift is $\sigma_q(f)(u):=f(qu)$.  In the normalized sectorial frame,
Proposition~\ref{prop:one-step-differential-W} gives
$\mathscr M\in N_+(K)s_{\mathrm{as}} N_+(K)$, while
\[
 \sigma_q(C)\mathscr M C^{-1}=C(qu)\mathscr M(u)C(u)^{-1}=Z\in H(K)\subset B_-(K).
\]
Lemma~\ref{lem:difference-field-generic-position}, applied with $A=\mathscr M$ and $g=C$, gives $C\in B_-(K)N_+(K)$.  Proposition~\ref{prop:Gauss-minor-criteria} gives its unique factorization $C=vn_+$.  Corollary~\ref{cor:Cartan-preservation-Weyl-minors} gives \eqref{eq:all-Weyl-Q-nonzero}; taking $w=1$ gives \eqref{eq:principal-central-minors-nonzero}.  The factors in \eqref{eq:standard-Q-normalization} are nonzero constants, so the rescaled functions $Q_w^{(i)}$ are nonzero as well.
\end{proof}

\begin{remark}\label{rem:Miura-structures-fixed-oper}
For a fixed $Z$-twisted $(G,q)$-oper with matrix $A$, a trivialization $g$ satisfying $A(u)=g(qu)^{-1}Zg(u)$ identifies its Miura reductions with meromorphic maps $f:U\to G/B_-$ satisfying
\[
 f(qu)=Zf(u),\qquad F_-(u)=g(u)^{-1}f(u).
\]
Indeed, multiplication by $g(qu)$ transforms $A(u)F_-(u)=F_-(qu)$ into this equation.  The constant solutions form the fixed-flag variety $(G/B_-)^Z$ appearing in \cite[Proposition~3.3]{FKSZ2020}.

The rational case over $\mathbb C(u)$ admits a sharper classification when $g$ is rational and $q$ is not a root of unity.  For $Z=1$, every rational solution $f$ is constant, so the Miura reductions are parametrized by $G/B_-$.  Under the nonresonance condition
\[
 \alpha(Z)\notin q^{\mathbb Z}\qquad\text{for every root }\alpha,
\]
the only rational solutions are $f=\dot wB_-$, $w\in W$, giving exactly $|W|$ Miura reductions.  Both assertions follow by expressing $f$ in root coordinates on the $H$-stable Bruhat cell containing its generic value.  Each coordinate satisfies $r(qu)=\alpha(Z)r(u)$.  A nonzero rational solution of $r(qu)=\lambda r(u)$ has zeros and poles only at $0$ and $\infty$, so it is $r(u)=au^m$ with $\lambda=q^m$.  Thus the coordinates are constant when $Z=1$ and vanish under the nonresonance condition.

Regular semisimplicity excludes $\alpha(Z)=1$ but permits resonances $\alpha(Z)=q^m$ with $m\ne0$.  For example, for $G=SL_2$, $A=Z=\operatorname{diag}(\zeta,\zeta^{-1})$ with $\zeta^2=q\ne1$, and oper line $\mathbb C(1,1)$, the rational lines $F_-^a(u)=\mathbb C(au,1)$, $a\in\mathbb C$, give distinct preserved reductions.  Over a general meromorphic field $\mathcal M(U)$, even the nonresonance condition need not force $f$ to be constant; the reductions are described by the equivariance equation above.
\end{remark}

\section{\texorpdfstring{$(G,q)$}{(G,q)}-opers and universal \texorpdfstring{$(G,q)$}{(G,q)}-Wronskians}
In this section we construct the universal $(G,q)$-Wronskian and use the $Q$-functions to recover the matrix of shifted monodromy. We begin with the one-step Wronskian equation for the central connection. Right multiplication by certain elements of $N_+$ preserves this equation. For the fixed bipartite Coxeter element and even Coxeter number $h$, we show that the equations at the shifts $qu,\ldots,q^{h/2}u$ select a unique representative in each such orbit. We construct this representative by Bruhat factorization of products of inverse connection matrices. We then establish Steinberg reduction for an arbitrary Coxeter element. Returning to the bipartite element of our differential family, we compute the principal minors of the normalized Steinberg representative when $q$ has infinite order. The $q$-Miura transformation expresses these minors through ratios of shifted $Q$-functions and the prescribed Cartan factors, and the minors determine the shifted monodromy matrix.

\paragraph{Coefficient fields and regular loci.}
Throughout, $K=\mathcal M(U)$ and $\sigma_qf(u)=f(qu)$ for a connected $q$-stable spectral domain $U$, and $A(u)\in G(K)$ denotes a $q$-connection matrix.  We retain the Borel conventions of Definition~\ref{def:local-twisted-Miura-qoper}: the oper flag has type $B_+$ and the preserved Miura flag has type $B_-$.  Conditions over $K$ hold where the factors are regular and the required minors are nonzero.  The algebraic arguments also apply over a field extension of $\mathbb C$ with a specified $\mathbb C$-automorphism, and the recursive construction permits every $q\ne0$.

\subsection{Generalized and universal \texorpdfstring{$(G,q)$}{(G,q)}-Wronskians}\label{sec:generalized-universal-Wronskians}
In this subsection we reconstruct the universal $(G,q)$-Wronskian when $h$ is even and $c^{h/2}=w_0$. Starting from a generalized $(G,q)$-Wronskian, we use successive Bruhat factorizations to determine the unique right unipotent factor for which the section satisfies all the shifted equations. For the central connection $C$, the normalization of the sectorial solutions at infinity already puts shifted monodromy in Steinberg form. We show that the right factor is then the identity, so the reconstruction returns $C$ itself. We also describe how the construction behaves under extensions of the coefficient field and unipotent changes of frame.

\subsubsection{Shifted equations and the fixed \texorpdfstring{$(G,q)$}{(G,q)}-Wronskian orbit}
The central connection satisfies the one-step equation of \cite[Section~4]{KZ2021}, as we proved in Proposition~\ref{prop:one-step-differential-W}.  For a generalized $(G,q)$-Wronskian $\mathscr W$ with data $(Z,s_{\mathrm{as}})$, the same equation reads
\begin{equation}\label{eq:W-basic-u}
 Z^{-1}\mathscr W(qu)v_{\omega_i}
 =\mathscr W(u)s_{\mathrm{as}}(u)^{-1}v_{\omega_i},
 \qquad i\in I.
\end{equation}
The group $N_+$ fixes every fundamental highest-weight vector.  Conversely, an element fixing these vectors belongs to $B_+$ by Proposition~\ref{prop:highest-weight-stabilizers}; all fundamental-weight characters take the value one on its Cartan factor, so that factor is the identity.  This gives the following stabilizer description.
\begin{lemma}\label{lem:unipotent-weight-scalar}
The simultaneous stabilizer in $G(K)$ of the vectors $v_{\omega_i}$ in all fundamental representations is $N_+(K)$.
\end{lemma}

For $n\in N_+(K)$, the columns of $\mathscr Wn$ and $\mathscr W$ on the left agree, because $n(qu)$ fixes every $v_{\omega_i}$.  Their right sides agree precisely when $s_{\mathrm{as}}ns_{\mathrm{as}}^{-1}$ fixes all these vectors.  By Lemma~\ref{lem:unipotent-weight-scalar}, the criterion is
\begin{equation}\label{eq:right-N-preserves-Wronskian}
 s_{\mathrm{as}}(u)n(u)s_{\mathrm{as}}(u)^{-1}\in N_+(K).
\end{equation}

\begin{definition}\label{def:Wronskian-right-orbit}
The right-multiplication orbit of a generalized $(G,q)$-Wronskian is
\[
 [\mathscr W]_{\mathrm{Wr}}
 :=\{\mathscr Wn:\ n\in N_+(K)\text{ satisfies }\eqref{eq:right-N-preserves-Wronskian}\}.
\qedhere
\]
\end{definition}

The extra equations compare the values of $\mathscr W(u)$ on successive translates of the fundamental highest-weight vectors.  For the fixed bipartite word with $h$ even, the step $d=h/2$ has $c^{-d}=w_0$ by Corollary~\ref{cor:fixed-datum-half-Coxeter}, so it takes the highest-weight lines to the lowest-weight lines.  The ordered products of the inverse lifts are
\[
 \Pi_0(u):=1,\qquad
 \Pi_m(u):=s_{\mathrm{as}}(u)^{-1}\cdots s_{\mathrm{as}}(q^{m-1}u)^{-1}.
\]
The \emph{shifted $(G,q)$-Wronskian equations} are
\begin{equation}\label{eq:W-41-u}
 Z^{-m}\mathscr W(q^mu)v_{\omega_i}
 =\mathscr W(u)\Pi_m(u)v_{\omega_i},
 \qquad i\in I,\quad 1\le m\le d.
\end{equation}
The first is \eqref{eq:W-basic-u}.  The additional equations evaluate $\mathscr W(u)$ on the further Weyl translates $\Pi_m(u)v_{\omega_i}$ and remove the freedom in Definition~\ref{def:Wronskian-right-orbit}.

\begin{lemma}\label{lem:finite-shift-uniqueness}
For even $h$, with $d=h/2$ and $c^d=w_0$, if $\mathscr W_1,\mathscr W_2\in G(K)$ satisfy the same shifted equations and
\[
 \mathscr W_2=\mathscr W_1n,
 \qquad n\in N_+(K),
\]
then $n=1$.
\end{lemma}

\begin{proof}
Since $n(q^mu)$ fixes $v_{\omega_i}$, the left sides of the two instances of \eqref{eq:W-41-u} coincide.  Cancelling $\mathscr W_1(u)$ from their right sides gives
\[
 n(u)\Pi_m(u)v_{\omega_i}=\Pi_m(u)v_{\omega_i}
 \qquad(i\in I),
\]
so Lemma~\ref{lem:unipotent-weight-scalar} gives $n(u)\in N_+\cap \Pi_m(u)N_+\Pi_m(u)^{-1}$.  At the final step, the image of $\Pi_d\in N_G(H)(K)$ in $W=N_G(H)/H$, which we call its Weyl projection, is $c^{-d}=w_0$.  Since conjugation by a Cartan factor preserves each root subgroup, we have $\Pi_dN_+\Pi_d^{-1}=N_-$.  Consequently,
\[
 n(u)\in N_+\cap N_-=\{1\}.
\]
\end{proof}

\begin{proposition}\label{prop:shifted-Wronskian-longest-cells}
Every $\mathscr W\in G(K)$ satisfying the shifted equations at $m=1$ and $m=d$ satisfies
\begin{equation}\label{eq:W-universal-u}
 \mathscr W\in B_-w_0B_-\cap B_+w_0B_+.
\end{equation}
Moreover, $\Delta_{w\omega_i,\omega_i}(\mathscr W)\in K^\times$ for every $w\in W$ and $i\in I$.
\end{proposition}
\begin{proof}
The one-step equation puts the reconstructed connection in $N_+(K)s_{\mathrm{as}}$, allowing us to apply Corollary~\ref{cor:Cartan-preservation-Weyl-minors} to $\mathscr W$.  The equation at $m=d$ then gives the minors for the other longest cell.

For $A_{\mathscr W}(u):=\mathscr W(qu)^{-1}Z\mathscr W(u)$, the one-step equation gives
\[
 A_{\mathscr W}(u)s_{\mathrm{as}}(u)^{-1}v_{\omega_i}=v_{\omega_i}.
\]
Lemma~\ref{lem:unipotent-weight-scalar} implies $A_{\mathscr W}s_{\mathrm{as}}^{-1}\in N_+(K)$.  Since
\[
 \sigma_q(\mathscr W)A_{\mathscr W}\mathscr W^{-1}=Z\in H,
\]
Corollary~\ref{cor:Cartan-preservation-Weyl-minors} gives $\Delta_{w\omega_i,\omega_i}(\mathscr W)\in K^\times$ for every $w,i$.  The case $w=w_0$ and Proposition~\ref{prop:Gauss-minor-criteria} give $\mathscr W\in B_+w_0B_+$.

For the other longest cell, the image of $\Pi_d$ in $W$ is $w_0$, so
\begin{equation}\label{eq:Pi-d-longest-weight-vector}
 \Pi_d(u)v_{\omega_i}=\kappa_i(u)v_{w_0\omega_i},
 \qquad \kappa_i\in K^\times.
\end{equation}
Pairing the equation at $m=d$ with the highest-weight covector therefore gives
\begin{equation}\label{eq:terminal-opposite-minor}
 \kappa_i(u)\Delta_{\omega_i,w_0\omega_i}(\mathscr W(u))
 =\omega_i(Z^{-d})\Delta_{\omega_i,\omega_i}(\mathscr W(q^du)).
\end{equation}
The right side is nonzero by the preceding nonvanishing result and the automorphism $\sigma_q^d$.  The opposite-cell criterion of Proposition~\ref{prop:Gauss-minor-criteria} gives $\mathscr W\in B_-w_0B_-$.
\end{proof}

\begin{definition}\label{def:universal-Wronskian-notation}
For an orbit $\mathcal C=[\mathscr W]_{\mathrm{Wr}}$, its \emph{universal meromorphic $(G,q)$-Wronskian} is a representative
\[
 \mathscr W_{\mathrm{univ}}(\mathcal C)\in\mathcal C
\]
satisfying all shifted equations~\eqref{eq:W-41-u}.
\end{definition}
Lemma~\ref{lem:finite-shift-uniqueness} gives uniqueness, and Proposition~\ref{prop:shifted-Wronskian-longest-cells} gives both longest-cell conditions.  We now construct such a representative.  For the orbit of the differential central connection we abbreviate $\mathscr W_{\mathrm{univ}}([C]_{\mathrm{Wr}})$ to $\mathscr W_{\mathrm{univ}}$.

\begin{definition}\label{def:qconnection-trivialization}
A meromorphic map $g:U\to G$ trivializes a $Z$-twisted $q$-connection $A$ if
\begin{equation}\label{eq:shifted-cocycle-base}
 Z^{-1}g(qu)=g(u)A(u)^{-1}.
\qedhere
\end{equation}
\end{definition}

We first construct a section satisfying the shifted equations from a trivializing pair $(A,g)$ with
\[
 A\in N_+(K)s_{\mathrm{as}}N_+(K).
\]
To obtain a representative of the orbit of $g$, the initial $g$ must satisfy the one-step equation.  Substitution of \eqref{eq:shifted-cocycle-base} into \eqref{eq:W-basic-u} gives the equivalent condition $As_{\mathrm{as}}^{-1}\in N_+(K)$, because this product must fix every $v_{\omega_i}$.  For the differential pair $(\mathscr M,C)$ the condition follows from \eqref{eq:BC-oper-line}.

\begin{remark}\label{rem:twist-regularity}\label{rem:universal-framing}
The fixed data include $Z$, $s_{\mathrm{as}}$, and the Miura reduction represented by $g^{-1}B_-$, with normalization given by the chosen trivializing map $g$.  Regularity of $Z\in H$ means $\alpha(Z)\ne1$ for every root, or equivalently that its centralizer is $H$.  It is relevant when recovering Miura reductions from an underlying $Z$-twisted $(G,q)$-oper.  Lemma~\ref{lem:finite-shift-uniqueness} requires no such condition.

For a meromorphic map $h_0:U\to G$ satisfying $h_0(qu)=h_0(u)$ and $h_0(u)Z=Zh_0(u)$, the product $h_0\mathscr W$ determines the same $q$-connection as $\mathscr W$, possibly in a different right-multiplication orbit. For example, any constant $t\in H$ gives a section $tC$ satisfying the same shifted equations as $C$, since
\[
 Z^{-m}(tC)(q^mu)v_{\omega_i}
 =tZ^{-m}C(q^mu)v_{\omega_i}
 =(tC)(u)\Pi_m(u)v_{\omega_i}.
\]
When $h$ is even and $c^{h/2}=w_0$, Lemma~\ref{lem:finite-shift-uniqueness} implies that $tC$ belongs to $[C]_{\mathrm{Wr}}$ only for $t=1$. The normalized differential frames prescribe the columns $C(u)v_{\omega_i}$ and thereby select the orbit $[C]_{\mathrm{Wr}}$.
\end{remark}

\subsubsection{Iterating the trivialization identity and successive Bruhat factorization}
Iterating the trivialization identity expresses each shifted value of $g$ through a product of inverse connection matrices. We identify the Bruhat cell of this product and factor it with the prescribed lift. The resulting unipotent factors give the successive corrections to $g$.

We use $P_0(u):=1$ for the empty product.
\begin{lemma}\label{lem:shifted-cocycle-iteration}
If $A$ and $g$ satisfy \eqref{eq:shifted-cocycle-base}, then for every $m\ge1$,
\begin{equation}\label{eq:shifted-cocycle-corrected}
 Z^{-m}g(q^mu)
 =g(u)P_m(u),
 \qquad
 P_m(u):=A(u)^{-1}\cdots A(q^{m-1}u)^{-1}.
\end{equation}
\end{lemma}

\begin{proof}
We iterate \eqref{eq:shifted-cocycle-base} by induction on $m$.  At each step it appends the inverse connection $A(q^mu)^{-1}$ on the right and adds one factor $Z^{-1}$ on the left.  This gives the product with the factors in the stated order.
\end{proof}

\begin{lemma}\label{lem:corrected-Bruhat-count}
If
\[
 A^{-1}\in B_+(K)c^{-1}B_+(K),
\]
then, for $1\le m\le d=h/2$,
\begin{equation}\label{eq:corrected-Bruhat-count}
 P_m\in B_+(K)c^{-m}B_+(K).
\end{equation}
In particular, $P_d$ lies in $B_+(K)w_0B_+(K)$.
\end{lemma}

\begin{proof}
Reversing the fixed reduced Coxeter word gives a reduced word for $c^{-1}$.  We have $c^{-d}=w_0$ and $\ell_W(w_0)=d\ell_W(c)$, so the concatenation of $d$ copies of this word is reduced.  The first $m$ copies thus form a reduced word of length $m\ell_W(c)$.  Proposition~\ref{prop:Bruhat-root-factorization}(a) gives length-additive multiplication of Bruhat cells.  Applying this rule successively to the copies of $c^{-1}$ gives \eqref{eq:corrected-Bruhat-count} by induction.
\end{proof}

For $w\in W$ and a prescribed lift $s\in H(K)\dot w$, the middle factor will be $s^{-1}$.  The left factor only needs the positive root subgroups that $s$ sends to negative ones.  These form the subgroup
\[
 U_s:=N_+\cap s^{-1}N_-s.
\]
Its roots are the positive roots $\beta$ for which $w\beta$ is negative.
\begin{lemma}\label{lem:prescribed-Cartan-Bruhat-factorization}
For a field $K$ over $\mathbb C$, multiplication is an isomorphism of $K$-varieties, whose map on $K$-points is
\begin{equation}\label{eq:prescribed-Cartan-Bruhat-isomorphism}
 U_s(K)\times N_+(K)\xrightarrow{\sim}N_+(K)s^{-1}N_+(K),
 \qquad (x,\widetilde n)\longmapsto xs^{-1}\widetilde n.
\end{equation}
The factorization is unique and commutes with extension of $K$.
\end{lemma}

\begin{proof}
The complementary subgroup is $V_s:=N_+\cap s^{-1}N_+s$.  Proposition~\ref{prop:Bruhat-root-factorization}(b) gives the following root coordinates.  Using the inversion roots of a reduced word
for $w^{-1}$, namely the positive roots $\beta$ for which $w\beta$ is negative,
we obtain coordinates on $U_s$; the remaining positive roots
give coordinates on $V_s$.  Multiplication in this root-subgroup order
therefore gives $N_+=U_sV_s$, with unique factors.  We may write any
left unipotent factor as $xv$, where $x\in U_s(K)$ and $v\in V_s(K)$.  Then
\[
 xv s^{-1}n=x s^{-1}(s v s^{-1})n,
 \qquad s v s^{-1}\in N_+,
\]
which proves existence of the factorization in \eqref{eq:prescribed-Cartan-Bruhat-isomorphism}.
For uniqueness, the equality $xs^{-1}n=x's^{-1}n'$ implies
$x'^{-1}x\in U_s\cap s^{-1}N_+s=\{1\}$, and it follows that $n=n'$.
When $s=h\dot w$, conjugation by the torus merely rescales the root parameters
by characters of $h$.  Thus the factorization and its inverse are both
defined over $\mathbb C$, uniformly in the Cartan parameter.
Specializing this parameter to $h\in H(K)$ gives the compatibility with extension of the field.
\end{proof}

For the pair $(A,g)$ with $A\in N_+s_{\mathrm{as}}N_+$, the prescribed lift at $q^ju$ is
\[
 \sigma_j(u):=s_{\mathrm{as}}(q^ju)=h_{\mathrm{as}}(q^ju)\dot c.
\]
The ordered products $\Pi_m$ are those in the shifted equations \eqref{eq:W-41-u}.  The index $j$ records a shift of the spectral parameter, and each $\sigma_j$ has Weyl projection $c$.
At shift $m$, we factor the preceding right unipotent factor times the next inverse connection.  Starting with $N_0=\widetilde N_0=1$, Lemma~\ref{lem:prescribed-Cartan-Bruhat-factorization} gives unique factors
\begin{equation}\label{eq:chosen-trivialization-recursive-factorization}
 \widetilde N_{m-1}(u)A(q^{m-1}u)^{-1}
 =x_m(u)\sigma_{m-1}(u)^{-1}\widetilde N_m(u),
 \qquad 1\le m\le d,
\end{equation}
with
\begin{equation}\label{eq:chosen-trivialization-one-step-layer}
 x_m(u)\in U_m(u):=N_+\cap \sigma_{m-1}(u)^{-1}N_-\sigma_{m-1}(u),
 \qquad \widetilde N_m(u)\in N_+.
\end{equation}
In the product of the first $m$ inverse connections, the new left factor $x_m$ stands after $\Pi_{m-1}$.  Moving it to the left conjugates it by $\Pi_{m-1}$ and gives the next factor of $N_m$.  The right factor $\widetilde N_m$ is used at the next shift.  Thus the successive correction factors are
\begin{equation}\label{eq:chosen-trivialization-new-layer}
 n_+^m(u):=\Pi_{m-1}(u)x_m(u)\Pi_{m-1}(u)^{-1},
 \qquad N_m(u):=n_+^1(u)\cdots n_+^m(u).
\end{equation}
The sections $W^m(u):=g(u)N_m(u)$, for $0\le m\le d$, have terminal value
\begin{equation}\label{eq:chosen-trivialization-terminal-section}
 \mathscr W_{\mathrm{rec}}(u):=W^d(u)=g(u)N_d(u),
 \qquad
 N_d=n_+^1\cdots n_+^d.
\end{equation}

\begin{proposition}\label{prop:successive-Wronskian-factorization}
For $1\le m\le d$, the factors satisfy
\begin{equation}\label{eq:chosen-trivialization-product-factorization}
 P_m(u)=N_m(u)\Pi_m(u)\widetilde N_m(u).
\end{equation}
Moreover,
\begin{equation}\label{eq:chosen-trivialization-layer-positivity}
 \Pi_j^{-1}n_+^m\Pi_j\in N_+\quad(0\le j<m),
 \qquad
 \Pi_m^{-1}n_+^m\Pi_m\in N_-.
\end{equation}
The section $W^m$ satisfies the first $m$ shifted $(G,q)$-Wronskian equations.
If also $As_{\mathrm{as}}^{-1}\in N_+(K)$, then $g$ is an initial generalized $(G,q)$-Wronskian and $\mathscr W_{\mathrm{rec}}\in[g]_{\mathrm{Wr}}$.
\end{proposition}

\begin{proof}
The assumption $A\in N_+s_{\mathrm{as}}N_+$ implies $A(q^{m-1}u)^{-1}\in N_+\sigma_{m-1}^{-1}N_+$.  If the previous step has produced $\widetilde N_{m-1}\in N_+$, we therefore have
\[
 \widetilde N_{m-1}A(q^{m-1}u)^{-1}
 \in N_+\sigma_{m-1}^{-1}N_+.
\]
We start with $\widetilde N_0=1$ and apply Lemma~\ref{lem:prescribed-Cartan-Bruhat-factorization} at each step.  It gives the required factors in \eqref{eq:chosen-trivialization-recursive-factorization}, and also their uniqueness.

The induction hypothesis for the product formula is $P_{m-1}=N_{m-1}\Pi_{m-1}\widetilde N_{m-1}$.  We multiply this identity by $A(q^{m-1}u)^{-1}$ and substitute \eqref{eq:chosen-trivialization-recursive-factorization}.  This gives
\[
 P_m
 =N_{m-1}\Pi_{m-1}x_m\sigma_{m-1}^{-1}\widetilde N_m
 =N_m\Pi_m\widetilde N_m,
\]
which proves \eqref{eq:chosen-trivialization-product-factorization}.

The later corrections must preserve the equations already obtained.  For this, conjugating the $m$-th correction back by $\Pi_j$ must give an element of $N_+$ for $j<m$.  We verify this using the root subgroups described in Proposition~\ref{prop:Bruhat-root-factorization}(b).  The concatenation
of $d$ copies of the reduced word for $c^{-1}$ is reduced and represents
$w_0$.  Its successive root blocks are described by
\[
 \mathcal R:=\{\beta\in\Delta_+:c\beta\in\Delta_-\}.
\]
These are precisely the roots of the subgroup $U_m$ containing $x_m$.  After conjugation, the subgroup $\Pi_{m-1}U_m\Pi_{m-1}^{-1}$ containing $n_+^m$ has root set $c^{-(m-1)}\mathcal R$, which forms the $m$-th block of the reduced word.  For $0\le j<m$, conjugation by
$\Pi_j^{-1}$ sends this set to $c^{-(m-j-1)}\mathcal R\subset\Delta_+$.
The inclusion holds because we have again obtained one of the positive
blocks of the same reduced concatenation.  For $j=m$, the resulting set
is $c\mathcal R\subset\Delta_-$.  The Cartan factors rescale the root
parameters but do not change the root sets.  This proves the two inclusions in \eqref{eq:chosen-trivialization-layer-positivity}.

It remains to verify the shifted equations.  For $1\le j\le m$, the remaining correction factors form the product $R_{j,m}:=n_+^{j+1}\cdots n_+^m$, with the empty product equal to $1$.  The factorization $N_m=N_jR_{j,m}$ separates the first correction factors from the remaining ones.  The first assertion of \eqref{eq:chosen-trivialization-layer-positivity} implies $\Pi_j^{-1}R_{j,m}\Pi_j\in N_+$.  Both $N_m(q^ju)$ and $\widetilde N_j(u)$ fix every $v_{\omega_i}$, so formulas \eqref{eq:shifted-cocycle-corrected} and \eqref{eq:chosen-trivialization-product-factorization} yield
\[
 Z^{-j}W^m(q^ju)v_{\omega_i}
 =g(u)P_j(u)v_{\omega_i}
 =g(u)N_j(u)\Pi_j(u)v_{\omega_i}.
\]
On the other hand,
\[
 W^m(u)\Pi_j(u)v_{\omega_i}
 =g(u)N_j(u)\Pi_j(u)
  \bigl(\Pi_j(u)^{-1}R_{j,m}(u)\Pi_j(u)\bigr)v_{\omega_i},
\]
where the last factor fixes $v_{\omega_i}$, so the two expressions agree.  Hence $W^m$ satisfies the equations with shifts $1,\ldots,m$; taking $m=d$ gives the full system.  Finally, if $As_{\mathrm{as}}^{-1}\in N_+$, the initial section $g$ also satisfies the one-step equation.  Comparing that equation with the one-step equation for $gN_d$, we find that $s_{\mathrm{as}}N_ds_{\mathrm{as}}^{-1}$ fixes every $v_{\omega_i}$.  It therefore belongs to $N_+$.  This is the conjugation condition of Definition~\ref{def:Wronskian-right-orbit}, so the initial section and the section obtained at the terminal step are in the same $(G,q)$-Wronskian orbit.
\end{proof}

\subsubsection{The universal representative in the two longest Bruhat cells}
At $m=d$, the recursive section satisfies the full shifted system.  When $g$ is a generalized $(G,q)$-Wronskian, the terminal section also belongs to $[g]_{\mathrm{Wr}}$, so the uniqueness statement of Lemma~\ref{lem:finite-shift-uniqueness} applies.

\begin{theorem}\label{thm:universal-Wronskian-both-cells}
We assume that $h$ is even, $c^{h/2}=w_0$, and $A,g\in G(K)$ satisfy
\[
 Z^{-1}g(qu)=g(u)A(u)^{-1},\qquad
 As_{\mathrm{as}}^{-1}\in N_+(K).
\]
These data define a $Z$-twisted Miura $(G,q)$-oper with oper flag $F_+=B_+$ and preserved flag $F_-(u)=g(u)^{-1}B_-$.  The recursive factors produce
\begin{equation}\label{eq:universal-W-from-g}
 \mathscr W_{\mathrm{rec}}=gN_d,\qquad N_d\in N_+(K),
\end{equation}
which is the unique representative of $[g]_{\mathrm{Wr}}$ satisfying all shifted equations:
\begin{equation}\label{eq:recursive-universal-identification}
 \mathscr W_{\mathrm{rec}}=\mathscr W_{\mathrm{univ}}([g]_{\mathrm{Wr}}).
\end{equation}
It also satisfies
\begin{equation}\label{eq:universal-longest-cells}
 \mathscr W_{\mathrm{rec}}\in B_-w_0B_-\cap B_+w_0B_+.
\end{equation}
The reconstructed connection $A_{\mathscr W_{\mathrm{rec}}}(u)=\mathscr W_{\mathrm{rec}}(qu)^{-1}Z\mathscr W_{\mathrm{rec}}(u)$ and preserved flag $\mathscr W_{\mathrm{rec}}^{-1}B_-$ represent the same Miura $(G,q)$-oper as $(A,B_+,F_-)$, through the $q$-gauge transformation $N_d^{-1}$.
\end{theorem}
\begin{proof}
Proposition~\ref{prop:successive-Wronskian-factorization} gives $\mathscr W_{\mathrm{rec}}\in[g]_{\mathrm{Wr}}$ and all shifted equations.  Lemma~\ref{lem:finite-shift-uniqueness} gives uniqueness, so Definition~\ref{def:universal-Wronskian-notation} gives \eqref{eq:recursive-universal-identification}.  Proposition~\ref{prop:shifted-Wronskian-longest-cells} then gives \eqref{eq:universal-longest-cells}.  In particular, since $N_d$ fixes the fundamental highest-weight vectors,
\begin{equation}\label{eq:longest-cell-minor-locus}
 \mathscr W_{\mathrm{rec}}v_{\omega_i}=gv_{\omega_i},\qquad
 \Delta_{w\omega_i,\omega_i}(\mathscr W_{\mathrm{rec}})
 =\Delta_{w\omega_i,\omega_i}(g)\in K^\times.
\end{equation}
The original preserved flag satisfies
\[
 A(u)g(u)^{-1}B_-=g(qu)^{-1}Z B_-=g(qu)^{-1}B_-.
\]
The equality $\mathscr W_{\mathrm{rec}}=gN_d$ gives the change of frame for both reconstructed objects.  For the connection,
\begin{equation}\label{eq:recursive-same-qoper-class}
 \mathscr W_{\mathrm{rec}}(qu)^{-1}Z\mathscr W_{\mathrm{rec}}(u)
 =N_d(qu)^{-1}A(u)N_d(u).
\end{equation}
The gauge $N_d^{-1}$ preserves the oper flag $B_+$ and sends $F_-=g^{-1}B_-$ to $N_d^{-1}g^{-1}B_-=\mathscr W_{\mathrm{rec}}^{-1}B_-$.  Thus it identifies both the connection and the prescribed preserved flag.
\end{proof}

\begin{remark}\label{rem:published-indexing}
Our index $m$ counts the inverse connections in \eqref{eq:shifted-cocycle-corrected}, so $m=k+1$ gives twist $Z^{-(k+1)}$, connection indices $0,\ldots,k$, and Bruhat component $c^{-(k+1)}$.  The terminal correction is therefore $n_+^d$.  This convention corrects the one-step indexing discrepancy in \cite[Eqs.~(4.36)--(4.37), Lemma~4.19, and the terminal display in the proof of Proposition~4.20]{KZ2021}.
\end{remark}

The next corollary identifies the normal form produced by this reconstruction.
\begin{corollary}\label{cor:universal-Wronskian-Steinberg}
The hypotheses are that $h$ is even, $d=h/2$, $c^d=w_0$, and $A,g\in G(K)$ satisfy
\[
 Z^{-1}g(qu)=g(u)A(u)^{-1},\qquad As_{\mathrm{as}}^{-1}\in N_+(K).
\]
The universal $(G,q)$-Wronskian $\mathscr W_{\mathrm{univ}}$ in $[g]_{\mathrm{Wr}}$ has the form $gN_d$ with $N_d\in N_+(K)$. The connection reconstructed from $\mathscr W_{\mathrm{univ}}$ is the $q$-gauge transform of $Z$ by $\mathscr W_{\mathrm{univ}}^{-1}$ and has Steinberg form
\begin{equation}\label{eq:universal-Wronskian-Steinberg}
 \begin{aligned}
 A_{\mathscr W_{\mathrm{univ}}}(u)
 &:=\mathscr W_{\mathrm{univ}}(qu)^{-1}Z\mathscr W_{\mathrm{univ}}(u)\\
 &=N_d(qu)^{-1}A(u)N_d(u)
 \in U_c(K)s_{\mathrm{as}}(u).
 \end{aligned}
\end{equation}
The following three statements are equivalent.
\begin{enumerate}[label=(\roman*),topsep=3pt,itemsep=3pt]
\item $g$ satisfies all shifted $(G,q)$-Wronskian equations~\eqref{eq:W-41-u}.
\item $A$ already lies in the Steinberg section $U_c(K)s_{\mathrm{as}}$.
\item $N_d=1$, so that $g=\mathscr W_{\mathrm{univ}}$.
\end{enumerate}
\end{corollary}
\begin{proof}
To prove that (i) implies (ii), we use $n:=As_{\mathrm{as}}^{-1}\in N_+(K)$.  Substitution of the equation with index $d-1$ at $qu$ into the left side of the equation with index $d$ at $u$ gives the factor $g(u)A(u)^{-1}$ by \eqref{eq:shifted-cocycle-base}.  The right side has the factor $g(u)s_{\mathrm{as}}(u)^{-1}$.  Since $A^{-1}=s_{\mathrm{as}}^{-1}n^{-1}$, cancellation gives
\[
 n(u)\Pi_{d-1}(qu)v_{\omega_i}
 =\Pi_{d-1}(qu)v_{\omega_i}
 \qquad(i\in I).
\]
For $d=1$, the equation with index zero is the identity. Since the Weyl image of $\Pi_{d-1}$ is $c^{-(d-1)}=cw_0$, the joint stabilizer of these vectors is $\dot cN_-\dot c^{-1}$. Thus
\[
 n\in N_+\cap\dot cN_-\dot c^{-1}=U_c.
\]
Conversely, (ii) gives $n\in U_c(K)$. The word obtained by repeating a reduced word for $c$ exactly $d$ times is reduced, since $c^d=w_0$ and $dr=\ell_W(w_0)$. Hence $c^j$ sends every root of $U_c$ to a positive root for $0\le j<d$. Conjugation by $\Pi_j(qu)^{-1}$ therefore sends $n(u)$ into $N_+$, so $n(u)$ fixes every $\Pi_j(qu)v_{\omega_i}$. In $P_m$, each factor $n(q^ju)^{-1}$ can be moved to the right through the subsequent inverse lifts.  The root inclusion just proved keeps the conjugated factors in $N_+$, where they fix $v_{\omega_i}$.  Thus $P_mv_{\omega_i}=\Pi_mv_{\omega_i}$, giving (i) by the iterated trivialization identity \eqref{eq:shifted-cocycle-corrected}.

Theorem~\ref{thm:universal-Wronskian-both-cells} supplies $\mathscr W_{\mathrm{univ}}=gN_d$ satisfying the full shifted system. If (i) holds, uniqueness in Lemma~\ref{lem:finite-shift-uniqueness} gives $g=\mathscr W_{\mathrm{univ}}$, so $N_d=1$. Conversely, if $N_d=1$, then $g=\mathscr W_{\mathrm{univ}}$ and (i) holds.

The first shifted equation for $\mathscr W_{\mathrm{univ}}$ gives $A_{\mathscr W_{\mathrm{univ}}}s_{\mathrm{as}}^{-1}\in N_+(K)$. Applying (i) $\Rightarrow$ (ii) to $(A_{\mathscr W_{\mathrm{univ}}},\mathscr W_{\mathrm{univ}})$ gives the Steinberg condition in \eqref{eq:universal-Wronskian-Steinberg}. The equality there is \eqref{eq:recursive-same-qoper-class}.
\end{proof}

\begin{corollary}\label{cor:differential-universal-Wronskian}
Suppose the ADE differential family satisfies Hypotheses~\ref{hyp:frob}--\ref{hyp:ccm} and its Coxeter number $h$ is even. Then the central connection matrix $C=\mathcal Y_{\mathrm F}^{-1}\mathcal Y_{\mathrm{as}}$ in the normalized frames is the universal $(G,q)$-Wronskian in its orbit. Over $K=\mathcal M(U)$ we have
\[
 C=\mathscr W_{\mathrm{univ}}([C]_{\mathrm{Wr}})
   \in B_-w_0B_-\cap B_+w_0B_+ .
\]
\end{corollary}
\begin{proof}
For $d=h/2$, the fixed bipartite Coxeter element satisfies $c^d=w_0$ by Corollary~\ref{cor:fixed-datum-half-Coxeter}. The frame identity~\eqref{eq:central-infinity-transport} gives $\mathscr M(u)=C(qu)^{-1}ZC(u)$, and Theorem~\ref{prop:analytic-sectorial-Steinberg} gives $\mathscr M(u)\in U_c(K)s_{\mathrm{as}}(u)$. Applying Corollary~\ref{cor:universal-Wronskian-Steinberg} with $(A,g)=(\mathscr M,C)$ therefore gives $N_d=1$ and all the shifted equations for $C$. The two Bruhat cell inclusions follow from Theorem~\ref{thm:universal-Wronskian-both-cells}.
\end{proof}

\subsubsection{Coefficient fields and gauge covariance}
The unique factorization in Lemma~\ref{lem:prescribed-Cartan-Bruhat-factorization} determines how the recursion behaves under extension of the coefficient field and under unipotent changes of frame.

\begin{proposition}\label{prop:Wronskian-base-change}
For an inclusion $K_0\subset K$ of $\sigma_q$-stable fields over $\mathbb C$, we assume the data $A,g,Z,s_{\mathrm{as}}$ of Proposition~\ref{prop:successive-Wronskian-factorization} are defined over $K_0$.  All recursive factors, the terminal section, and its shifted equations are defined over $K_0$ and unchanged after extension to $K$.  The same holds for the universal representative whenever the initial section is a generalized $(G,q)$-Wronskian.  A generalized minor over $K_0$ is nonzero after extension precisely when it was nonzero over $K_0$.
\end{proposition}
\begin{proof}
The factorization and its inverse in Lemma~\ref{lem:prescribed-Cartan-Bruhat-factorization} are defined over the coefficient field and are unique.  Subsequent operations are multiplication, inversion, conjugation, and $\sigma_q$.  Thus the entire construction commutes with extension of the field, and injectivity of $K_0\hookrightarrow K$ preserves nonvanishing.
\end{proof}

\begin{remark}
For rational data $A,g,Z,s_{\mathrm{as}}$, Proposition~\ref{prop:Wronskian-base-change} shows that all factors and the resulting section are rational in $u$.
\end{remark}

In the shifted equations, the vectors representing the shifted lines transform with the frame.  In the terminal factorization, we compare the two input pairs using the same prescribed lifts $s_{\mathrm{as}}$.  For $\eta\in N_+(K)$, the changes of coordinates are
\[
 A^\eta(u)=\eta(qu)A(u)\eta(u)^{-1},\qquad
 g^\eta=g\eta^{-1},\qquad \mathscr W^\eta=\mathscr W\eta^{-1}.
\]
The prescribed shifted line is spanned by $\sigma_{i,m}=\Pi_mv_{\omega_i}$. Under the change of frame this vector becomes $\sigma_{i,m}^\eta=\eta\sigma_{i,m}$.
\begin{proposition}\label{prop:shifted-equations-gauge-covariance}\label{prop:terminal-recursion-gauge-invariance}
The shifted identity
\begin{equation}\label{eq:intrinsic-shifted-equations}
 Z^{-m}\mathscr W(q^mu)v_{\omega_i}
 =\mathscr W(u)\sigma_{i,m}(u)
\end{equation}
holds if and only if it holds for $\mathscr W^\eta$ and $\sigma_{i,m}^\eta$.

For a pair $(A,g)$ as in Proposition~\ref{prop:successive-Wronskian-factorization}, the recursive construction applied to $(A^\eta,g^\eta)$ with the factors $s_{\mathrm{as}}$ unchanged gives
\begin{equation}\label{eq:terminal-gauge-factor-law}
 N_d^\eta=\eta N_d,\qquad
 \widetilde N_d^\eta=\widetilde N_d\eta(q^du)^{-1},
\end{equation}
and
\begin{equation}\label{eq:terminal-W-gauge-invariance}
 \mathscr W_{\mathrm{rec}}(A^\eta,g^\eta)
 =\mathscr W_{\mathrm{rec}}(A,g).
\end{equation}
Consequently the reconstructed $q$-connection matrix is unchanged.
\end{proposition}
\begin{proof}
The left side of \eqref{eq:intrinsic-shifted-equations} is unchanged because $\eta(q^mu)$ fixes $v_{\omega_i}$; its right side is unchanged because $\mathscr W^\eta\sigma_{i,m}^\eta=\mathscr W\sigma_{i,m}$.

For the recursive construction, $A^\eta\in N_+s_{\mathrm{as}}N_+$ and
\[
 P_d^\eta=\eta P_d\eta(q^du)^{-1}
 =(\eta N_d)\Pi_d\bigl(\widetilde N_d\eta(q^du)^{-1}\bigr).
\]
Since $\Pi_d$ has Weyl projection $w_0$, Lemma~\ref{lem:prescribed-Cartan-Bruhat-factorization} gives uniqueness of this terminal factorization with arbitrary left and right $N_+$ factors.  This proves \eqref{eq:terminal-gauge-factor-law}, and $g^\eta N_d^\eta=gN_d$ proves \eqref{eq:terminal-W-gauge-invariance}.
\end{proof}

If the initial section $g$ is a generalized $(G,q)$-Wronskian for $(Z,s_{\mathrm{as}})$, Theorem~\ref{thm:universal-Wronskian-both-cells} identifies the common value in \eqref{eq:terminal-W-gauge-invariance} with the universal $(G,q)$-Wronskian $\mathscr W_{\mathrm{univ}}([g]_{\mathrm{Wr}})$.

\subsection{Shifted monodromy in terms of \texorpdfstring{$Q$}{Q}-functions}\label{sec:difference-field-Steinberg}\label{sec:Miura-Steinberg}\label{sec:consequences}

Our aim in this subsection is to express the shifted monodromy matrix $\mathscr M(u)=C(qu)^{-1}ZC(u)$ in terms of the fundamental $Q$-functions and the prescribed Cartan factors. In the sectorial frame it suffices to compute the principal minors of $h_{\mathrm{as}}(u)^{-1}\mathscr M(u)$, since these are the coordinates of an element in the Steinberg section. The preserved Frobenius flag allows us to write the same $q$-connection in Miura form. We will express its Miura variables through principal minors of the central connection $C$ and hence through ratios of shifted $Q$-functions.

To compute the Steinberg coordinates from the Miura variables, we follow~\cite{ZeitlinFRWronskians} and first reduce the calculation to the formal lattice $\mathbb Z$. The reason is that the matrix of the connection at $q^ku$ depends on the frames at the two points $q^ku$ and $q^{k+1}u$. We may therefore first work with a $G$-equivariant isomorphism between principal right $G$-fibers at each pair of consecutive integers $k,k+1$, allowing a separate frame at each integer. In the discrete $q$-Miura transformation we treat the Miura coordinates at different integers as independent indeterminates. For each vertex $i$, the Steinberg coordinate at $k=0$ is then a Laurent polynomial in finitely many of them. When $q$ has infinite order, this polynomial is the fundamental $q$-character at $i$.

To evaluate this polynomial for our oper, we keep the base point $u$ variable in the spectral domain $U$. The principal right $G$-bundle $P\to U$ of the oper has fiber $P_u=\{u\}\times G$ in the sectorial trivialization. Over the formal lattice, the map $k\mapsto q^ku$ assigns the fiber $P_{q^ku}$ to $k$. The isomorphism to the next fiber is
\[
 \mathcal A_{q^ku}:P_{q^ku}\longrightarrow P_{q^{k+1}u},\qquad
 (q^ku,g)\longmapsto\bigl(q^{k+1}u,\mathscr M(q^ku)g\bigr).
\]
After normalizing the frames at these fibers, we substitute the Miura coordinate of the map from $P_{q^ku}$ to $P_{q^{k+1}u}$ for the formal variable at the same vertex and integer $k$. For each $k$ these coordinates are meromorphic functions of the base point $u$. Evaluating the character therefore gives the corresponding principal minor of $h_{\mathrm{as}}(u)^{-1}\mathscr M(u)$ as a meromorphic function on $U$.

To describe the changes of frame that put the connection in Steinberg form, we fix a Coxeter word $c=s_{i_1}\cdots s_{i_r}$ with its pinned lift $\dot c=\widehat s_{i_1}^{-1}\cdots\widehat s_{i_r}^{-1}$. The group $U_c=N_+\cap\dot cN_-\dot c^{-1}$ depends only on $c$.

\begin{definition}\label{def:discrete-Steinberg-frames}
For a commutative $\mathbb C$-algebra $R$ and each integer $k$, a trivialization identifies the $R$-points of the fiber with $G(R)$ acting on itself on the right. An equivariant isomorphism from the fiber at $k$ to the fiber at $k+1$ is therefore left multiplication by a unique $A_k\in G(R)$. The family of isomorphisms is represented by the sequence $(A_k)_{k\in\mathbb Z}$. A change of coordinates $g\mapsto p_kg$ at each position gives the \emph{discrete gauge transformation}
\begin{equation}\label{eq:discrete-gauge-on-fibers}
                         A'_k=p_{k+1}A_kp_k^{-1}.
\end{equation}
For a lattice connection with $A_k\in N_+(R)\dot c N_+(R)$, a \emph{Steinberg representative} is a sequence $(S_k)_{k\in\mathbb Z}$ obtained by this formula with $p_k\in N_+(R)$ and $S_k\in U_c(R)\dot c$ for every $k$.
\end{definition}

For a meromorphic $q$-connection with matrix $A(u)$, the map from the fiber at $k$ to the fiber at $k+1$ is represented by $A_k=A(q^ku)$. A meromorphic change of frame $p(u)$ gives $p_k=p(q^ku)$ in \eqref{eq:discrete-gauge-on-fibers}. We may also choose a frame separately at each index $k$. We use this freedom to remove the prescribed Cartan factor at each shift.

Frenkel, Reshetikhin, and Semenov-Tian-Shansky proved the cross-section theorem for formal loops in $SL_n$~\cite[Lemma~1]{FrenkelReshetikhinSTS1998}. Semenov-Tian-Shansky and Sevostyanov formulated the formal-loop and finite cyclic versions for semisimple groups~\cite[Theorems~1.1 and~2.9]{SS1998}. Their proof covers the types other than $E_6$~\cite[\S3, footnote~2]{SS1998}. Francone and Leclerc proved the theorem for connection matrices indexed by $\mathbb Z$ in all ADE types~\cite[Theorem~33 and Corollary~34]{FL2025}. In the convention of \eqref{eq:discrete-gauge-on-fibers}, their result takes the following form~\cite[Theorem~4.10]{ZeitlinFRWronskians}.

\begin{theorem}\label{thm:Steinberg-discrete-frames}
For $G$ of ADE type, any Coxeter element $c$, and any commutative $\mathbb C$-algebra $R$, every sequence $A_k\in N_+(R)\dot cN_+(R)$, $k\in\mathbb Z$, admits a unique sequence $p_k\in N_+(R)$ such that
\[
                      p_{k+1}A_kp_k^{-1}\in U_c(R)\dot c
                      \qquad(k\in\mathbb Z).
\]
\end{theorem}

For a $q$-connection the elements at successive shifts come from a single element over a difference field $(K,\tau)$. Here $K$ is a field containing $\mathbb C$ and $\tau$ is an automorphism that fixes $\mathbb C$. For the differential family $K=\mathcal M(U)$ and $(\tau f)(u)=f(qu)$. The sectorial frame can contribute a Cartan factor, so we need the cross-section theorem for a representative $h_*\dot c$ of $c$, with $h_*\in H(K)$. Every representative of $c$ over $K$ has this form.

\begin{theorem}\label{thm:difference-field-Steinberg}
For $G$ of ADE type, a Coxeter element $c$ with its pinned lift $\dot c$, and $h_*\in H(K)$, every $A\in N_+(K)h_*\dot cN_+(K)$ admits a unique $\gamma\in N_+(K)$ such that
\[
                     \tau(\gamma)A\gamma^{-1}\in U_c(K)h_*\dot c.
\]
\end{theorem}
\begin{proof}
The successive images of $A$ and its Cartan factor are $A_k=\tau^k(A)$ and $h_k=\tau^k(h_*)$. The Cartan factor in $t_{k+1}A_kt_k^{-1}$ is $t_{k+1}h_kc(t_k)^{-1}$, where $c(t)=\dot c t\dot c^{-1}$ for $t\in H$. It equals one if the torus elements $t_k$ are determined in both directions by
\[
                     t_0=1,\qquad t_{k+1}=c(t_k)h_k^{-1}.
\]
With this choice, $t_{k+1}A_kt_k^{-1}\in N_+(K)\dot cN_+(K)$ for every $k$. Theorem~\ref{thm:Steinberg-discrete-frames} gives a unique sequence $p_k$ taking these elements to the Steinberg section. In the original coordinates, $\gamma_k=t_k^{-1}p_kt_k$ is the unique sequence in $N_+(K)$ such that
\[
             \gamma_{k+1}A_k\gamma_k^{-1}\in U_c(K)h_k\dot c .
\]
Indeed, the torus normalizes $N_+$ and $U_c$, so conjugation by $t_k$ recovers the changes of frame $p_k$ from the $\gamma_k$ and takes the displayed inclusion back to the Steinberg condition for $t_{k+1}A_kt_k^{-1}$. Applying $\tau$ to the displayed inclusion gives the same reduction problem with $k$ replaced by $k+1$. Uniqueness therefore implies $\tau(\gamma_k)=\gamma_{k+1}$. It follows that $\gamma_0$ is a gauge over $K$ with the required property. Conversely, for any gauge $g$ over $K$ with this property, the changes of frame $\tau^k(g)$ satisfy the displayed inclusion in place of $\gamma_k$. Uniqueness at $k=0$ gives $g=\gamma_0$.
\end{proof}

In the sectorial frame $\mathscr M(u)\in U_c h_{\mathrm{as}}(u)\dot c$ by Theorem~\ref{prop:analytic-sectorial-Steinberg}. Since the torus normalizes $U_c$, removing the Cartan factor gives $S(u)=h_{\mathrm{as}}(u)^{-1}\mathscr M(u)\in U_c\dot c$. With the opposite choice of Borels, the ordered root factors of the Steinberg section give the canonical coordinates for rational $q$-opers defined in~\cite[Proposition~8.1 and Theorem~8.2]{FKSZ2020}. In our normalization, these coordinates are the fundamental principal minors of $S$.

For $g\in G$, the principal minor $\Delta_i(g)=\Delta_{\omega_i,\omega_i}(g)$ is the coefficient of the highest-weight vector $v_{\omega_i}$ in $gv_{\omega_i}$. The roots associated with the word $c$ are $\beta_a=s_{i_1}\cdots s_{i_{a-1}}(\alpha_{i_a})$. We normalize the parametrization $X_a:\mathbb C\to U_{\beta_a}$ by $\Delta_{i_a}(X_a(x)\dot c)=x$.

The principal minors also give coordinates on the elements of $N_+\dot cN_+$ that lie in $B_-$. We record the two parametrizations of~\cite[Proposition~2.6]{ZeitlinFRWronskians} for our choice of Coxeter word and lift.

\begin{proposition}\label{prop:Cartan-factor-removal}
\begin{enumerate}[label=(\roman*),topsep=3pt,itemsep=4pt]
\item The principal minors give an isomorphism $U_c\dot c\longrightarrow\mathbb A^r$. In particular, the functions $T_i(u)=\Delta_i(S(u))$ determine shifted monodromy by
\begin{equation}\label{eq:theta-as-BC-component}
 \mathscr M(u)=h_{\mathrm{as}}(u)\,
 X_1\bigl(T_{i_1}(u)\bigr)\cdots X_r\bigl(T_{i_r}(u)\bigr)\dot c .
\end{equation}
\item On $B_-\cap N_+\dot cN_+$, the map $b\mapsto (\Delta_j(b))_{j\in I}$ is an isomorphism onto $(\mathbb C^\times)^I$. Its inverse is the ordered product
\[
 (x_j)_{j\in I}\longmapsto P_{i_1}(x_{i_1})\cdots P_{i_r}(x_{i_r}),
 \qquad P_j(x)=\exp(x^{-1}f_j)\,\check\alpha_j(x).
\]
Here $f_j$ is the pinned generator for $-\alpha_j$ and $\check\alpha_j:\mathbb C^\times\to H$ is the coroot of $\alpha_j$.
\end{enumerate}
\end{proposition}
\begin{proof}
For (i), moving the simple-reflection lifts to the right in $\prod_{a=1}^r\exp(x_{i_a}e_{i_a})\widehat s_{i_a}^{-1}$ gives the ordered root factors for $\beta_1,\ldots,\beta_r$ followed by $\dot c$. By Proposition~\ref{prop:Bruhat-root-factorization}(b), this product parametrizes $U_c\dot c$. In $V_j$, every factor with $i_a\ne j$ fixes both the highest-weight vector and its dual covector, since $\langle\omega_j,\check\alpha_{i_a}\rangle=0$. The factor at $j$ acts on the two-dimensional simple-root string through $v_{\omega_j}$ by $\left(\begin{smallmatrix}x_j&-1\\1&0\end{smallmatrix}\right)$. Hence the $j$th principal minor of the product is $x_j$. Applying the inverse parametrization to $S(u)=h_{\mathrm{as}}(u)^{-1}\mathscr M(u)$ gives~\eqref{eq:theta-as-BC-component}.

For (ii), in the pinned rank-one subgroup the matrix of $P_j(x)$ is $\left(\begin{smallmatrix}x&0\\1&x^{-1}\end{smallmatrix}\right)$. Hence
\[
 \exp(xe_j)\widehat s_j^{-1}=P_j(x)\exp(-x^{-1}e_j).
\]
If $j\ne\ell$, then $[e_j,f_\ell]=0$. Conjugation by $\check\alpha_\ell(y)$ gives the commutation formula
\[
 \exp(t e_j)P_\ell(y)
 =P_\ell(y)\exp\!\bigl(y^{-\langle\alpha_j,\check\alpha_\ell\rangle}t e_j\bigr).
\]
Each vertex occurs once in the Coxeter word, so each exponential on the right of a rank-one factor can be moved past all subsequent $P_\ell$. Consequently $P_{i_1}(x_{i_1})\cdots P_{i_r}(x_{i_r})$ differs from the Steinberg representative in (i) by right multiplication by an element of $N_+$. The ordered product therefore belongs to $N_+\dot cN_+$. Each $P_j(x_j)$ also belongs to $B_-$ and has Cartan component $\check\alpha_j(x_j)$. Thus the product belongs to $B_-$ and its $j$th principal minor is $x_j$.

Conversely, every $b\in B_-\cap N_+\dot cN_+$ has nonzero principal minors, since they are the fundamental weights evaluated on its Cartan component. By Proposition~\ref{prop:Bruhat-root-factorization}(b), $b$ has the form $S n_+$ with $S\in U_c\dot c$ and $n_+\in N_+$. Since $n_+$ fixes each fundamental highest-weight vector, $\Delta_j(b)=\Delta_j(S)$ for every $j$. If $b_1$ and $b_2$ have the same minors, their Steinberg factors $S$ coincide by (i). Thus $b_1^{-1}b_2\in B_-\cap N_+=\{1\}$, which proves uniqueness. Both maps are regular, since the minors are regular on $G$ and the ordered product is regular for $x_j\ne0$.
\end{proof}

To compute these Steinberg coordinates from the $Q$-functions, we use the preserved Frobenius flag. In the sectorial frame the oper and Miura flags are $B_+$ and $C(u)^{-1}B_-$, respectively. Proposition~\ref{prop:central-principal-nonvanishing} gives the Gauss factorization $C=v_-n$ over $K=\mathcal M(U)$, with $v_-\in B_-(K)$ and $n\in N_+(K)$. Since $n\in N_+$, the change of coordinates by $n$ fixes $B_+$. It takes the preserved flag to $B_-$ because $nC^{-1}=v_-^{-1}\in B_-$. In this trivialization, the fiber isomorphism $P_u\to P_{qu}$ is left multiplication by the element
\[
 A_-(u)=n(qu)\mathscr M(u)n(u)^{-1}
       =v_-(qu)^{-1}Zv_-(u)\in B_-(K).
\]

In these coordinates the isomorphism from $P_{q^ku}$ to $P_{q^{k+1}u}$ is represented by $A_-(q^ku)$. Since $H$ preserves both standard flags, we can still change coordinates by a torus element at each index $k$. As in the proof of Theorem~\ref{thm:difference-field-Steinberg}, this freedom removes the factor $h_{\mathrm{as}}(q^ku)$ from the lift $s_{\mathrm{as}}(q^ku)=h_{\mathrm{as}}(q^ku)\dot c$. The elements $t_k=t_k(u)\in H(K)$ are determined in both directions by
\begin{equation}\label{eq:monodromy-torus-frames}
 t_0=1,\qquad t_{k+1}=c(t_k)h_{\mathrm{as}}(q^ku)^{-1}.
\end{equation}
The matrix of the fiber map in these frames is $A_k^{\mathrm M}=t_{k+1}A_-(q^ku)t_k^{-1}$.

\begin{definition}\label{def:Miura-form-shifted}
On the locus where the two reductions are opposite, a \emph{Miura frame} is a local trivialization in which the oper reduction $Q_+$ and the preserved reduction $Q_-$ are represented by $1B_+$ and $1B_-$, respectively. The \emph{Miura form} of a Miura $(G,q)$-oper is its $q$-connection expressed in such a frame. For the lattice frames normalized by~\eqref{eq:monodromy-torus-frames}, the highest-weight \emph{Miura variables} are the principal minors
\[
 \mathsf q_{j,k}(u):=\Delta_j(A_k^{\mathrm M}),\qquad j\in I,\quad k\in\mathbb Z .
 \qedhere
\]
\end{definition}

\begin{proposition}\label{prop:Miura-central-minors}
\begin{enumerate}[label=(\roman*),topsep=3pt,itemsep=4pt]
\item In a Miura frame the connection matrix $A_-(u)$ belongs to $B_-(K)$, and two Miura frames differ by an $H$-valued change of coordinates. In the normalized lattice frames the variables are $\mathsf q_{j,k}(u)=\omega_j([A_k^{\mathrm M}]_0)$, where $[b]_0$ denotes the image of $b$ in $B_-/N_-=H$. They determine the connection by
\begin{equation}\label{eq:Miura-form-ordered-product}
 A_k^{\mathrm M}=P_{i_1}(\mathsf q_{i_1,k}(u))\cdots P_{i_r}(\mathsf q_{i_r,k}(u)).
\end{equation}
At the base point $u$ one has $A_-(u)=h_{\mathrm{as}}(u)A_0^{\mathrm M}$.
\item The Cartan component of $A_-(u)$ is determined by the principal minors of the central connection matrix $C$ through the identities
\[
 \omega_j([A_-(u)]_0)
 =\omega_j(Z)\frac{\Delta_j(C(u))}{\Delta_j(C(qu))},\qquad j\in I.
\]
The Miura variables in~\eqref{eq:Miura-form-ordered-product} are
\begin{equation}\label{eq:differential-YQ-ratio}
 \mathsf q_{j,k}(u)
 =\omega_j(t_{k+1}Zt_k^{-1})
   \frac{\Delta_j(C(q^ku))}{\Delta_j(C(q^{k+1}u))}
 =\omega_j(t_{k+1}Zt_k^{-1})
   \frac{Q_+^{(j)}(q^ku)}{Q_+^{(j)}(q^{k+1}u)} .
\end{equation}
\item For every $j\in I$ and $k\in\mathbb Z$, the Miura variable $\mathsf q_{j,k}$ belongs to $\mathcal M(U)^\times$. Its relation to the variable obtained by normalizing the frames at $q^ku$ is
\[
 \mathsf q_{j,k}(u)
 =\omega_j\!\bigl(c(t_k(u))\,t_k(u)^{-1}\bigr)\,\mathsf q_{j,0}(q^ku).
\]
In particular, if $h_{\mathrm{as}}\equiv1$, then $\mathsf q_{j,k}(u)=\mathsf q_{j,0}(q^ku)$ for every $j$ and $k$.
\end{enumerate}
\end{proposition}
\begin{proof}
For (i), preservation of the Miura flag gives $A_-(u)B_-=B_-$ in a Miura frame. A change between two such frames fixes both identity cosets, whose common stabilizer is $B_+\cap B_-=H$. For $b\in B_-(K)$, the coefficient of the highest-weight vector in $b v_{\omega_j}$ is $\omega_j([b]_0)$. Applying this to $A_k^{\mathrm M}$ gives the Cartan-coordinate formula. Since $A_-$ is obtained from $\mathscr M$ by an $N_+$-valued change of frame, the Coxeter lift at $q^ku$ is $h_{\mathrm{as}}(q^ku)\dot c$. Equation~\eqref{eq:monodromy-torus-frames} gives $t_{k+1}h_{\mathrm{as}}(q^ku)\dot c\,t_k^{-1}=\dot c$. Thus $A_k^{\mathrm M}\in B_-(K)\cap N_+(K)\dot c N_+(K)$, and Proposition~\ref{prop:Cartan-factor-removal}(ii) gives~\eqref{eq:Miura-form-ordered-product}. For $k=0$, the normalization gives $t_0=1$ and $t_1=h_{\mathrm{as}}(u)^{-1}$.

For (ii), in the Gauss factorization $C=v_-n$, the factor $n$ fixes $v_{\omega_j}$. The coefficient along $v_{\omega_j}$ in $Cv_{\omega_j}$ is therefore $\omega_j([v_-]_0)=\Delta_j(C)$. Projecting $A_-(u)=v_-(qu)^{-1}Zv_-(u)$ to $H$ and applying $\omega_j$ now gives the first identity. At the successive points $q^ku$ and $q^{k+1}u$, the changes of coordinates by $t_k$ and $t_{k+1}$ supply the torus factors in~\eqref{eq:differential-YQ-ratio}. By \eqref{eq:C-minor-Q-u} and~\eqref{eq:standard-Q-normalization}, $Q_+^{(j)}$ differs from $\Delta_j(C)$ by a nonzero constant independent of $u$. This constant cancels in the ratio. Proposition~\ref{prop:central-principal-nonvanishing} ensures that the denominators are nonzero as meromorphic functions.

For (iii), the recurrence~\eqref{eq:monodromy-torus-frames} gives $t_k\in H(K)$ for every integer $k$. The formula in (ii) therefore gives $\mathsf q_{j,k}\in K^\times$. At the base point $v=q^ku$, the initial torus factors are $t_0(v)=1$ and $t_1(v)=h_{\mathrm{as}}(v)^{-1}$. Thus $A_0^{\mathrm M}(v)=h_{\mathrm{as}}(v)^{-1}A_-(v)$. With the frames normalized at $u$, the same fiber map is represented by
\[
 A_k^{\mathrm M}(u)
 =t_{k+1}(u)A_-(v)t_k(u)^{-1}
 =c(t_k(u))A_0^{\mathrm M}(v)t_k(u)^{-1}.
\]
Taking Cartan components and applying $\omega_j$ gives the claimed factor. When $h_{\mathrm{as}}\equiv1$, the recurrence with $t_0=1$ gives $t_k=1$ for all $k$.
\end{proof}

Formula~\eqref{eq:differential-YQ-ratio} gives the Miura variables in terms of the $Q$-functions. To recover shifted monodromy, we now pass from the connection in Miura form to the Steinberg representative of the underlying $q$-oper.

\begin{definition}\label{def:q-Miura-transformation}
The \emph{$q$-Miura transformation}
\[
              (P,\mathcal A,Q_+,Q_-)\longmapsto(P,\mathcal A,Q_+)
\]
assigns to each Miura $(G,q)$-oper the $(G,q)$-oper obtained by forgetting the preserved $B_-$-reduction $Q_-$.
\end{definition}

By Lemma~\ref{lem:difference-field-generic-position}, the oper and Miura flags are opposite on a dense open subset of $U$. There the $q$-Miura transformation has the following description in a Miura frame.

\begin{proposition}\label{prop:Miura-to-Steinberg}
In a fixed Miura frame with Cartan factor $h_*\in H(K)$, the $q$-Miura transformation sends the connection matrix $A_-\in B_-(K)\cap N_+(K)h_*\dot c N_+(K)$ to its Steinberg representative by
\[
 A_-(u)\longmapsto\gamma(qu)A_-(u)\gamma(u)^{-1}\in U_c(K)h_*\dot c ,
\]
with a unique change of frame $\gamma\in N_+(K)$.

For the differential family, $h_*=h_{\mathrm{as}}$ and $\gamma=n^{-1}$, where $C=v_-n$ is the Gauss factorization of the central connection. The Steinberg representative is $\mathscr M(u)$. In the lattice frames normalized by~\eqref{eq:monodromy-torus-frames}, the corresponding representatives are
\[
 S_k(u)=t_{k+1}(u)\mathscr M(q^ku)t_k(u)^{-1}\in U_c(K)\dot c,
 \qquad S_0(u)=h_{\mathrm{as}}(u)^{-1}\mathscr M(u).
\]
\end{proposition}
\begin{proof}
The Miura frame represents the two reductions by the identity cosets $1B_+$ and $1B_-$. A change of frame by $N_+(K)$ fixes $1B_+$, so it preserves the oper reduction retained in Definition~\ref{def:q-Miura-transformation}. Theorem~\ref{thm:difference-field-Steinberg} gives the unique such change of frame that puts $A_-$ in the stated Steinberg section.

For the differential family, $A_-(u)=n(qu)\mathscr M(u)n(u)^{-1}$ with $n\in N_+(K)$. Since $\mathscr M\in U_c(K)h_{\mathrm{as}}\dot c$ by Theorem~\ref{prop:analytic-sectorial-Steinberg}, uniqueness gives $\gamma=n^{-1}$. In the lattice frames the corresponding change of coordinates at $k$ is $p_k=t_k n(q^ku)^{-1}t_k^{-1}\in N_+(K)$. Substitution of $A_k^{\mathrm M}=t_{k+1}A_-(q^ku)t_k^{-1}$ gives
\[
 p_{k+1}A_k^{\mathrm M}p_k^{-1}=t_{k+1}\mathscr M(q^ku)t_k^{-1}.
\]
The torus normalizes $U_c$, and~\eqref{eq:monodromy-torus-frames} gives $t_{k+1}h_{\mathrm{as}}(q^ku)\dot c\,t_k^{-1}=\dot c$. Thus the right side lies in $U_c(K)\dot c$. At $k=0$ the formula uses $t_0=1$ and $t_1=h_{\mathrm{as}}(u)^{-1}$.
\end{proof}

It remains to express the Steinberg coordinates $T_i(u)=\Delta_i(S_0(u))$ in terms of the Miura variables $\mathsf q_{j,k}(u)$. For $q$ of infinite order and $\hbar^{-2}=q$, the Laurent polynomial giving $T_i(u)$ is the $q$-character of the fundamental representation at vertex $i$ of the quantum loop algebra $U_\hbar(L\g)$. We recall how its variables specify the Miura coordinates to be substituted at the different lattice positions.

The ordinary character of a finite-dimensional type~$1$ module $V$ records the simultaneous eigenspaces of the Cartan generators. These generators act diagonally, with eigenvalues in $\hbar^{\mathbb Z}$. In Drinfeld's presentation of the quantum loop algebra, each Cartan generator and its inverse are the constant terms of two series, the Cartan currents at that vertex. The coefficients of all these series commute. Their simultaneous generalized eigenspaces refine the ordinary weight spaces, and the eigenvalues on each such summand constitute its \emph{$\ell$-weight}.

An $\ell$-weight is encoded by a Laurent monomial $m=\prod_{j,b}Y_{j,b}^{n_{j,b}}$, where $j$ runs over the Dynkin vertices and $b\in\mathbb C^\times$. Only finitely many exponents are nonzero. At vertex $j$ the two current eigenvalues are the expansions at zero and infinity of the rational function
\[
 \prod_b\left(\hbar\,\frac{1-b\hbar^{-1}x}{1-b\hbar x}\right)^{n_{j,b}}
\]
in the formal current variable $x$. Its value at zero is $\hbar^{\sum_b n_{j,b}}$. Thus $m$ has ordinary weight $\sum_{j,b}n_{j,b}\omega_j$, while its individual factors retain the spectral labels $b$. Writing $V_m$ for the summand of $\ell$-weight $m$, the \emph{$q$-character} of $V$ is
\[
                  \chi_\hbar(V)=\sum_{m:\,V_m\ne0}(\dim V_m)\,m .
\]
Replacing every $Y_{j,b}$ by the ordinary weight symbol $e^{\omega_j}$ recovers the ordinary character. The $q$-character is additive on short exact sequences and multiplicative on tensor products. It therefore defines an injective homomorphism from the Grothendieck ring of finite-dimensional type~$1$ modules to the ring of Laurent polynomials in the $Y_{j,b}$ with integer coefficients~\cite{FR1998}.

A finite-dimensional simple module is generated by a simultaneous eigenvector of the Cartan currents annihilated by all positive-root current generators. Its $\ell$-weight is the \emph{highest $\ell$-weight} of the module. The exponents in its monomial are nonnegative and determine the Drinfeld polynomial $\prod_b(1-bx)^{n_{j,b}}$ at each vertex $j$. For the fundamental module $L(Y_{i,b})$, the polynomial is $1-bx$ at $i$ and $1$ at every other vertex. More generally, a string of factors $Y_{i,b}Y_{i,b\hbar^2}\cdots Y_{i,b\hbar^{2(p-1)}}$ is the highest monomial of the Kirillov--Reshetikhin module of length $p$ at $i$. Here the Steinberg coordinates are fundamental principal minors, and the characters we use are those of the length-one modules.

To evaluate these characters we use the bipartition of the Coxeter word $c=c_1c_0$. At vertex $j$ of color $\epsilon_j$, the formal variable for the map from the fiber at $k$ to the fiber at $k+1$ has spectral label
\[
                 b_{j,k}=\hbar^{2k+\epsilon_j},\qquad k\in\mathbb Z .
\]
Increasing $k$ by one multiplies the character label by $\hbar^2=q^{-1}$ and moves the source of the fiber map from $q^ku$ to $q^{k+1}u$. Thus $b_{j,k}$ specifies which Miura coordinate to substitute, while $u$ remains the variable of the resulting meromorphic function. The following identity first makes this substitution for a lattice connection over a field.

\begin{theorem}[{\cite[Theorem~6.5]{ZeitlinFRWronskians}}]\label{thm:Miura-character-sequences}
Over a field $F$ containing $\mathbb C$, suppose the isomorphism from the trivialized fiber at $k$ to that at $k+1$ is represented by $D_k\in B_-(F)\cap N_+(F)\dot c N_+(F)$ for every $k\in\mathbb Z$. The discrete $q$-Miura transformation sends these fiber maps in Miura form to their Steinberg representatives $S_k=p_{k+1}D_kp_k^{-1}\in U_c(F)\dot c$, where the changes of frame $p_k\in N_+(F)$ are uniquely determined by Theorem~\ref{thm:Steinberg-discrete-frames}. The principal minors of $S_0$ are
\[
 \Delta_i(S_0)=
 \left.\chi_\hbar\bigl(L(Y_{i,b_{i,0}})\bigr)
       \right|_{Y_{j,b_{j,k}}=\Delta_j(D_k)}
       \qquad(i\in I).
\]
Each variable in this character has a label $b_{j,k}$, so the substitution uses only finitely many of the elements $D_k$.
\end{theorem}

Francone and Leclerc first proved this character formula for every Coxeter element in types $A$, $D$, and $E$. Their proof uses generalized Baxter relations and the cluster algebra of $(G,c)$-bands, the discrete form of the universal $(G,q)$-Wronskians~\cite{FLBands2025,FL2025}. The direct proof in~\cite{ZeitlinFRWronskians} also applies beyond ADE. For the untwisted quantum affine algebras of types $B_n$, $C_n$, $F_4$, and $G_2$, it recovers the ordinary fundamental $q$-characters by combining the spectral shift with a Dynkin diagram automorphism and choosing an order of representatives of its vertex orbits.

The argument starts from the geometric meaning of the Miura transformation. Varying the preserved flag changes the Miura variables while fixing the underlying oper, so every Steinberg coordinate is constant under these variations. In the formal lattice variables their derivatives give the quantum screening operators, one for each simple root in ADE. The Laurent polynomial for a Steinberg coordinate therefore belongs to their common kernel. By the Frenkel--Mukhin characterization, an element of this kernel is determined by the coefficients of its dominant monomials, those with all exponents nonnegative~\cite{FrenkelMukhin2001}. For the $i$th coordinate at position zero, the calculation in~\cite{ZeitlinFRWronskians} gives the single dominant monomial $Y_{i,b_{i,0}}$ with coefficient one. This identifies the polynomial with $\chi_\hbar(L(Y_{i,b_{i,0}}))$ before evaluation on our meromorphic $q$-connection.

For shifted monodromy we write $\chi_i=\chi_\hbar(L(Y_{i,b_{i,0}}))$ for the fundamental character at vertex $i$. Proposition~\ref{prop:Miura-central-minors} expresses each Miura variable in terms of $C$, so evaluation of the Laurent variables gives meromorphic functions on $U$ by
\begin{equation}\label{eq:differential-Y-evaluation}
 \operatorname{ev}_C(Y_{j,b_{j,k}})(u)=\mathsf q_{j,k}(u).
\end{equation}
\begin{theorem}\label{thm:differential-qcharacter}
Under Hypotheses~\ref{hyp:frob}--\ref{hyp:ccm}, for $q$ of infinite order the Steinberg coordinates of shifted monodromy are
\begin{equation}\label{eq:differential-TQ-qcharacter}
 T_i(u)=\Delta_i\!\left(h_{\mathrm{as}}(u)^{-1}\mathscr M(u)\right)
       =\operatorname{ev}_C(\chi_i)(u),\qquad i\in I .
\end{equation}
The substitution~\eqref{eq:differential-YQ-ratio} expresses these coordinates through the principal minors of $C$ or, equivalently, through the normalized $Q$-functions. They determine $\mathscr M(u)$ by \eqref{eq:theta-as-BC-component}.
\end{theorem}
\begin{proof}
Proposition~\ref{prop:Miura-to-Steinberg} identifies $S_0(u)=h_{\mathrm{as}}(u)^{-1}\mathscr M(u)$ as the Steinberg representative at $k=0$ for the fiber maps $D_k=A_k^{\mathrm M}$. Their principal minors are $\Delta_j(D_k)=\mathsf q_{j,k}(u)$ by Definition~\ref{def:Miura-form-shifted}. With this substitution Theorem~\ref{thm:Miura-character-sequences} gives~\eqref{eq:differential-TQ-qcharacter}.
\end{proof}

For the FFH family of Section~\ref{sec:FFH-verification}, the infinite-order assumption gives $h_{\mathrm{as}}=1$. Indeed, the logarithmic phase in \eqref{eq:FFH-log-coefficient} can occur only when $q^m=\gamma_h=e^{2\pi i/h}$ for a positive integer $m$. Thus at infinite order the phase increment vanishes by Corollary~\ref{cor:FFH-rotation-hypothesis}, and $\mathsf q_{j,k}(u)=\mathsf q_{j,0}(q^ku)$. At the exceptional FFH parameters Corollary~\ref{lem:FFH-common-reference-sector} retains the Cartan factor contributed by the logarithm.

\begin{example}\label{sec:rank-one-T}\label{rem:A1-differential-YQ}
For $G=SL_2$, we write $Z=\operatorname{diag}(\zeta,\zeta^{-1})$, $\phi=\varphi_1$ and $Q=Q_+^{(1)}$. In rank one $h_{\mathrm{as}}(u)=\operatorname{diag}(\phi(u)^{-1},\phi(u))$ and $c$ acts by inversion. Thus $\omega_1(t_1)=\phi(u)$ and $\omega_1(t_2)=\phi(qu)/\phi(u)$. The dependence on the initial point in~\eqref{eq:differential-YQ-ratio} is visible in the first two variables, which satisfy
\[
 \mathsf q_{1,0}(u)=\zeta\phi(u)\frac{Q(u)}{Q(qu)},\qquad
 \mathsf q_{1,1}(u)=\phi(u)^{-2}\mathsf q_{1,0}(qu).
\]
The vertex has color $0$, so the character at $b_{1,0}=1$ is $Y_{1,1}+Y_{1,\hbar^2}^{-1}$. Its two monomials give
\[
 T_1(u)=\mathsf q_{1,0}(u)+\mathsf q_{1,1}(u)^{-1}
 =\zeta\phi(u)\frac{Q(u)}{Q(qu)}
  +\zeta^{-1}\frac{\phi(u)^2}{\phi(qu)}\frac{Q(q^2u)}{Q(qu)}.
\]
The same expression follows directly from $C(qu)\mathscr M(u)=ZC(u)$. The Steinberg form of $\mathscr M(u)$ is
\[
 \begin{pmatrix}\phi(u)^{-1}T_1(u)&-\phi(u)^{-1}\\ \phi(u)&0\end{pmatrix}.
\]
If $v(u)$ is the first column of $C(u)$, comparison of the second columns of $C(qu)\mathscr M(u)$ and $ZC(u)$ gives $-\phi(u)^{-1}Z^{-1}v(qu)$ for the second column of $C(u)$. This also gives the second column of $C(qu)$ after shifting $u$ to $qu$. Comparing the first columns and taking their highest-weight coefficient now gives the same formula for $T_1(u)$.
\end{example}

\section{The ADE central connection and its Steinberg frames}\label{sec:ADE-central-connection}
In this section we first show how, in type $A$, a section generating the oper line bundle and its prescribed image in the reference bundle determine the trivialization to the constant connection $Z$.  Comparing their successive pullbacks gives an explicit formula for the corresponding $(G,q)$-Wronskian.  We apply this construction to the differential family and combine it with the recursive construction of Section~\ref{sec:generalized-universal-Wronskians} to identify the central connection with the corresponding $(G,q)$-Wronskian in every ADE type.  We summarize these identifications and the Steinberg forms in one theorem, including the $q$-character parametrization of shifted monodromy when $q$ is not a root of unity.  We conclude by discussing the choice of Coxeter element.

\subsection{Algebraic reconstruction in type \texorpdfstring{$A$}{A}}\label{sec:typeA-direct}
Here we will recover a chosen trivialization to the constant connection from a generator of the oper line and its image. Comparing their successive pullbacks in the two bundles will also give a new frame in which to study the connection matrix.

\subsubsection{Type-$A$ \texorpdfstring{$(G,q)$}{(G,q)}-Wronskians}
In this subsection we use the defining representation $V=\mathbb C^N$ of $SL_N$ and omit $\rho_1$ in formulas for matrices acting on $V$. We work on a connected domain $U$ invariant under $u\mapsto q^{\pm1}u$, and write $K=\mathcal M(U)$ and $\sigma_q f(u)=f(qu)$.

We introduce the reference bundle $\mathcal P_Z=U\times SL_N$ with constant diagonal $q$-connection $Z\in H(\mathbb C)$, as in Definition~\ref{def:generalized-Gq-Wronskian}. We compare it with a $Z$-twisted $(SL_N,q)$-oper on a principal bundle $\mathcal P_A$. Its oper condition compares the transported flag at $qu$ with the chosen flag there, in the order specified in Section~\ref{sec:differential-Wronskian-transport}. We allow any Coxeter element $c\in W$ for this relative position and choose a representative $\dot c\in N_G(H)$.

The associated bundles in the defining representation are
\[
 E_A=\mathcal P_A\times^{SL_N}\mathbb C^N,
 \qquad E_Z=\mathcal P_Z\times^{SL_N}\mathbb C^N.
\]
The standard determinant on $\mathbb C^N$ induces a trivialization $\operatorname{vol}_X:\bigwedge^N E_X\xrightarrow{\sim}\mathcal O_U$ for $X=A,Z$.  At each $u$ it evaluates the top exterior product of $N$ vectors in the fiber to a scalar.  For an ordered tuple of meromorphic sections $P=(s_1,\ldots,s_N)$, we write $\operatorname{vol}_X(P)=\operatorname{vol}_X(s_1\wedge\cdots\wedge s_N)$ for this meromorphic function.  A frame has volume one when this function is identically one.

The oper reduction gives a complete flag of subbundles $\mathcal F_+$ in $E_A$ with
\[
 0=\mathcal F_{+,0}\subset\mathcal F_{+,1}\subset\cdots
   \subset\mathcal F_{+,N}=E_A,
 \qquad \operatorname{rank}\mathcal F_{+,a}=a.
\]
We say that a meromorphic frame $(e_1,\ldots,e_N)$ of $E_A$ is adapted to the oper reduction if, at every $u$ where the frame and flag are defined, the vectors $e_1(u),\ldots,e_a(u)$ form a basis over $\mathbb C$ of the fiber $\mathcal F_{+,a}(u)$ for each $a$.  We choose such a frame of volume one.  In $E_Z$, the standard volume-one frame $(f_1,\ldots,f_N)$ gives the lower flag $\mathcal F_-^Z$. Its rank-$a$ member has basis $f_N(u),\ldots,f_{N-a+1}(u)$ in the fiber at $u$.  The diagonal connection $Z$ preserves this flag.

These frames identify the meromorphic sections of $E_A$ and $E_Z$ with columns in $K^N$.  The two $q$-connections have matrices $A(u)\in SL_N(K)$ and $Z$, representing fiber maps from $u$ to $qu$.

In the chosen frames of $E_A$ and $E_Z$, the pullbacks of Section~\ref{sec:differential-Wronskian-transport} act by
\begin{equation}\label{eq:typeA-semilinear-shifts}
 (\sigma_A\mathbf v)(u)=A(u)^{-1}\mathbf v(qu),
 \qquad
 (\sigma_Z\mathbf w)(u)=Z^{-1}\mathbf w(qu).
\end{equation}
We denote the oper line by $\mathcal L_A:=\mathcal F_{+,1}\subset E_A$. We choose $p_0^A=e_1$, the first section of the adapted frame, as its meromorphic generator. At each $u$ where the frame is defined, $p_0^A(u)$ spans the fiber $(\mathcal L_A)_u$. We write $\mathbf e_i$ for the standard unit columns of $K^N$, so the coordinate column of $p_0^A$ in this frame is the constant column $\mathbf e_1=(1,0,\ldots,0)^{\mathsf T}$.

The successive inverse pullbacks of this generator form a frame, as the following proposition shows.

\begin{proposition}\label{prop:typeA-cyclicity-from-reference-matching}
For an $(SL_N,q)$-oper of Coxeter relative position $c$ with $q$-connection $A$ and a meromorphic generator $p$ of its oper line, the ordered tuple $(p,\sigma_Ap,\ldots,\sigma_A^{N-1}p)$ is a meromorphic frame of the associated bundle in the defining representation.
\end{proposition}
\begin{proof}
We show that the smallest subspace containing the first section and stable under the inverse pullback is the whole space. The Coxeter factorization allows us to enlarge this subspace along the standard flag.

We use a meromorphic frame of volume one adapted to the oper flag with first section $p$, and denote the connection matrix in that frame by $A(u)$. In this frame $p$ has coordinate column $\mathbf e_1$, where $\mathbf e_i$ are the unit columns of $K^N$.

The oper condition gives a factorization $A=b_1\dot c b_2$ with $b_1,b_2\in B_+(K)$. We write $E_k=\operatorname{span}_K(\mathbf e_1,\ldots,\mathbf e_k)$.  The span $L$ of the nonnegative iterates of $\mathbf e_1$ under $\sigma_A$ satisfies $\sigma_A(L)\subset L$.  Since $\sigma_A$ is invertible and the shift $\sigma_q$ is an automorphism of $K$, we have $\dim_K\sigma_A(L)=\dim_KL$.  Thus $L$ is stable under both $\sigma_A$ and $\sigma_A^{-1}$.

We identify the Weyl group of $SL_N$ with $S_N$ through its action on the coordinate lines, with $s_i=(i\,\,i+1)$.  For each $k>1$, either $c(k)<k$ or $c^{-1}(k)<k$.  For $k<N$, one of the two words $c,c^{-1}$ applies $s_{k-1}$ before $s_k$.  This is the first reflection to move $k$, and it moves it into $\{1,\ldots,k-1\}$.  Only $s_{k-1}$ crosses the boundary of this set, and it occurs once.  For $k=N$, the same argument uses $s_{N-1}$.

To pass from $E_{k-1}\subset L$ to $E_k\subset L$, we use the two cases above.  If $c(k)<k$, then $\sigma_A(E_{k-1})=b_2^{-1}\dot c^{-1}E_{k-1}$ contains $b_2^{-1}\mathbf e_k$.  If $c^{-1}(k)<k$, then $\sigma_A^{-1}(E_{k-1})=\sigma_q^{-1}(b_1)\dot c E_{k-1}$ contains $\sigma_q^{-1}(b_1)\mathbf e_k$.  In either case the vector lies in $L$ and, by upper triangularity, is a nonzero multiple of $\mathbf e_k$ modulo $E_{k-1}$.  Starting from $E_1\subset L$, induction gives $L=K^N$.  Any dependence among the first $N$ iterates would make the span preceding the first dependence $\sigma_A$-stable, contradicting $L=K^N$.
\end{proof}

The $Z$-twist gives an isomorphism to the reference connection. We choose such an isomorphism to specify the image of the chosen oper-line generator. The cyclic frame will show that this image and the two connections already determine the isomorphism.

\begin{definition}\label{def:typeA-universal-Wronskian}
For the $(SL_N,q)$-oper, reference connection, and frames above, a \emph{type-$A$ $(G,q)$-Wronskian datum} is a meromorphic volume-preserving bundle isomorphism $g:E_A\to E_Z$ intertwining the two $q$-connections. Its coordinate matrix, also denoted by $g$, belongs to $SL_N(K)$ and satisfies
\[
 g(qu)A(u)=Zg(u).
\]
\end{definition}
The chosen $g$ determines the section, line subbundle, and complete flag
\[
 p_0^Z:=g p_0^A,\qquad \mathcal L_Z:=g\mathcal L_A,
 \qquad \mathcal F_-:=g^{-1}\mathcal F_-^Z.
\]
In the chosen frames, $p_0^Z(u)=g(u)\mathbf e_1$, including its normalization. Since $Z$ preserves $\mathcal F_-^Z$, the connection preserves the induced flag,
\[
 A(u)\mathcal F_-(u)=g(qu)^{-1}Z\mathcal F_-^Z
                  =\mathcal F_-(qu).
\]
Thus $\mathcal F_-$ supplies a Miura reduction.

Applying the two pullback operators gives sections $p_j^A$ of $E_A$ and $p_j^Z$ of $E_Z$. We collect the first $N$ iterates of each initial section in the ordered tuples $P^A$ and $P^Z$.
\begin{equation}\label{eq:typeA-cyclic-orbits}
 p_j^A:=\sigma_A^{\,j}p_0^A,
 \qquad p_j^{Z}:=\sigma_Z^{\,j}p_0^Z,
 \qquad
 P^A:=(p_0^A,\ldots,p_{N-1}^{A}),
 \qquad P^Z:=(p_0^Z,\ldots,p_{N-1}^Z).
\end{equation}
Since $A\in B_+(K)\dot c B_+(K)$ in the adapted frame, Proposition~\ref{prop:typeA-cyclicity-from-reference-matching} shows that the first $N$ iterates of $p_0^A$ form a meromorphic frame for every Coxeter element $c$.  An ordered frame of successive iterates is called a \emph{cyclic frame}.  Thus $P^A$ is the cyclic frame generated by $p_0^A$.  The defining connection identity for $g$ gives $g\sigma_A=\sigma_Zg$ on meromorphic sections, and hence
\[
 p_j^Z=g p_j^A\quad(j\ge0),\qquad P^Z=gP^A.
\]
Thus $P^Z$ is the cyclic frame of $E_Z$ generated by $p_0^Z$.

We use $P^A$ and $P^Z$ also for the coordinate matrices of the two tuples.  In the chosen frames their coordinate columns are
\[
 \begin{aligned}
 p_j^A(u)&=A(u)^{-1}A(qu)^{-1}\cdots A(q^{j-1}u)^{-1}\mathbf e_1,\\
 p_j^Z(u)&=Z^{-j}p_0^Z(q^ju),\qquad j\ge0.
 \end{aligned}
\]
The matrix product equals the identity for $j=0$.

Following the principal-bundle interpretation of \cite[Definition~4.9]{KZ2021}, a meromorphic section of $\mathcal P_Z$ is a volume-one frame of $E_Z$.  Sending the chosen frame $(e_1,\ldots,e_N)$ to this frame defines a meromorphic volume-preserving bundle isomorphism $E_A\to E_Z$.

\begin{definition}\label{def:typeA-Wronskian-section}
A \emph{$(G,q)$-Wronskian section} for the given data is a meromorphic section of $\mathcal P_Z$ whose induced isomorphism $E_A\to E_Z$ intertwines the two $q$-connections and sends $p_0^A$ to $p_0^Z$.
\end{definition}

In coordinates the two requirements of Definition~\ref{def:typeA-Wronskian-section} are
\begin{equation}\label{eq:typeA-Wronskian-section}
 A(u)=\mathscr W(qu)^{-1}Z\mathscr W(u),
 \qquad \mathscr W(u)\mathbf e_1=p_0^Z(u).
\end{equation}
For every $Z$-twisted $(SL_N,q)$-oper, the chosen trivialization $g$ supplies such a section with matrix $g$. We recover that matrix from $p_0^Z$ by comparing the iterates of the two generators.

The same comparison identifies the recurrence governing the iterates in both bundles. Since $P^A$ and $P^Z$ are frames, the next iterate $p_N^X$ has unique coefficients in $P^X$ for $X=A,Z$. We write
\begin{equation}\label{eq:typeA-cyclic-recurrences}
 p_N^X=\sum_{j=0}^{N-1}\tau_{X,j}\,p_j^X,
 \qquad
 \boldsymbol\tau_X=(\tau_{X,0},\ldots,\tau_{X,N-1})^{\mathsf T}.
\end{equation}
Thus $\boldsymbol\tau_X$ is the coordinate column of $p_N^X$ in $P^X$.

\begin{proposition}\label{prop:typeA-cyclic-reconstruction}
For the data of Definition~\ref{def:typeA-universal-Wronskian}, the unique $(G,q)$-Wronskian section, denoted $\mathscr W_{\mathrm{univ}}^A$, has coordinate matrix
\begin{equation}\label{eq:typeA-abstract-universal-formula}
 \mathscr W_{\mathrm{univ}}^A=P^Z\bigl(P^A\bigr)^{-1}=g.
\end{equation}
The two frames have the same recurrence coefficients and volume, and the induced Miura flags have the same coordinates in these frames,
\begin{equation}\label{eq:typeA-frame-identities}
 \boldsymbol\tau_A=\boldsymbol\tau_Z,
 \qquad \operatorname{vol}_Z(P^Z)=\operatorname{vol}_A(P^A),
 \qquad (P^A)^{-1}\mathcal F_-=(P^Z)^{-1}\mathcal F_-^Z.
\end{equation}
In the frames $P^A$ and $P^Z$, both pullbacks act on coordinate columns by $\mathbf x(u)\mapsto R_{\mathrm{cyc}}(u)\mathbf x(qu)$, where
\[
 R_{\mathrm{cyc}}=(\mathbf e_2,\ldots,\mathbf e_N,\boldsymbol\tau_A).
\]
The matrix $R_{\mathrm{cyc}}$ is invertible, and
\begin{equation}\label{eq:typeA-abstract-universal-transport}
 A(u)^{-1}P^A(qu)=P^A(u)R_{\mathrm{cyc}}(u),
 \qquad Z^{-1}P^Z(qu)=P^Z(u)R_{\mathrm{cyc}}(u).
\end{equation}
Both $q$-connections from $u$ to $qu$ have matrix $R_{\mathrm{cyc}}^{-1}$ in these frames.
\end{proposition}
\begin{proof}
Any $(G,q)$-Wronskian section with matrix $\mathscr W$ satisfies $\mathscr W\sigma_A=\sigma_Z\mathscr W$. Its prescribed image of $p_0^A$ therefore gives $\mathscr Wp_j^A=p_j^Z$ for every $j\ge0$. The first $N$ iterates form a frame, so $\mathscr W=P^Z(P^A)^{-1}=g$. This proves uniqueness and the formula for the section.

Applying $g$ to the expansion of $p_N^A$ in \eqref{eq:typeA-cyclic-recurrences} gives the expansion of $p_N^Z$ with the same coefficients. Volume preservation gives the second identity in \eqref{eq:typeA-frame-identities}. Its last identity follows from $P^Z=gP^A$ and $\mathcal F_-=g^{-1}\mathcal F_-^Z$.

For $X=A,Z$, the pullback sends $p_j^X$ to $p_{j+1}^X$. Its coordinate columns in $P^X$ are consequently $\mathbf e_2,\ldots,\mathbf e_N,\boldsymbol\tau_X$. This proves \eqref{eq:typeA-abstract-universal-transport}. The identity $R_{\mathrm{cyc}}=(P^A)^{-1}A^{-1}P^A(qu)$ gives invertibility, and inverting the identities gives the two $q$-connection matrices.
\end{proof}

The matrix $R_{\mathrm{cyc}}$ is the companion matrix of the common recurrence: its first $N-1$ columns record the successive iterates, and the last records their recurrence coefficients.

In the original volume-one frames of $E_A$ and $E_Z$, the equality of volumes in \eqref{eq:typeA-frame-identities} reads $\det P^Z=\det P^A$. These equal volumes need not be one. The formula \eqref{eq:typeA-abstract-universal-formula} represents the $(G,q)$-Wronskian section in the original adapted frame of $E_A$. To express the same section in $P^A$, we next record the change of coordinates for a frame of arbitrary volume.

\begin{proposition}
\label{prop:typeA-Wronskian-frame-covariance}
For a $(G,q)$-Wronskian section $\mathscr G$ for the given data, we write $e=(e_1,\ldots,e_N)$ for the chosen adapted frame of $E_A$. We express the induced bundle map $E_A\to E_Z$ in the meromorphic frame $e^T=eT$, with $T\in GL_N(K)$, and the fixed frame $(f_1,\ldots,f_N)$ of $E_Z$. The relation $\mathbf v=T\mathbf x$ between the old and new coordinates in $E_A$ gives the following matrices, generator, and flags
\begin{equation}\label{eq:typeA-Wronskian-frame-covariance}
 \begin{gathered}
 A^T(u)=T(qu)^{-1}A(u)T(u),\qquad
 \mathscr G^T(u)=\mathscr G(u)T(u),\\
 p_0^{A,T}=T^{-1}\mathbf e_1,\qquad
 \mathcal F_\pm^T=T^{-1}\mathcal F_\pm.
 \end{gathered}
\end{equation}
The volume of the new frame is
\[
 \operatorname{vol}_A(e_1^T\wedge\cdots\wedge e_N^T)=\det T.
\]
Consequently,
\[
 \det A^T(u)=\frac{\det T(u)}{\det T(qu)},
 \qquad
 \det\mathscr G^T(u)=\det T(u).
\]
Both matrices are $SL_N(K)$-valued when $\det T=1$.  The iterates $p_j^A$ and the frame $P^A$ have new coordinate columns
\begin{equation}\label{eq:typeA-cyclic-columns-frame-covariance}
 p_j^{A,T}=T^{-1}p_j^A,\qquad
 P^{A,T}=T^{-1}P^A.
\end{equation}
\end{proposition}
\begin{proof}
The substitution $\mathbf v=T(u)\mathbf x$ gives the connection, generator, and flag formulas.  Taking the top exterior product of $e^T=eT$ gives the volume factor $\det T$.  The pullback on the new coordinate columns is $T^{-1}\sigma_A T$, so its iterates give \eqref{eq:typeA-cyclic-columns-frame-covariance}.  The inverse coordinate change recovers the defining conditions in the original frame.
\end{proof}

These formulas express the same bundle map in different coordinates: $\mathscr G^T$ is its matrix in the new frame of $E_A$. The fiberwise volume forms used in the defining conditions remain unchanged.

The cyclic reconstruction starts from the oper line in the defining representation. To compare the reconstructed $(G,q)$-Wronskian section with the generalized equations of Definition~\ref{def:generalized-Gq-Wronskian}, we use the remaining members of the oper flag through their top exterior powers. For $1\le a<N$, the fundamental representation is $\rho_a(g)=\bigwedge^a g$ on $\bigwedge^a\mathbb C^N$, and its associated bundle is $\bigwedge^a E_A$. The section $e_1\wedge\cdots\wedge e_a$ generates the line subbundle $\bigwedge^a\mathcal F_{+,a}$ and has coordinate column $v_{\omega_a}$.

\begin{proposition}\label{prop:typeA-exterior-Wronskian}
Suppose that, in the chosen oper frame, $A s_A^{-1}\in N_+(K)$ for a lift $s_A\in H(K)\dot c$. Then $\mathscr W_{\mathrm{univ}}^A$ satisfies the generalized $(G,q)$-Wronskian equations
\begin{equation}\label{eq:typeA-fundamental-Wronskian}
 \rho_a(Z)^{-1}\rho_a(\mathscr W_{\mathrm{univ}}^A(qu))v_{\omega_a}
 =\rho_a(\mathscr W_{\mathrm{univ}}^A(u))
   \rho_a(s_A(u))^{-1}v_{\omega_a},
 \qquad 1\le a<N.
\end{equation}
For a change of frame $e^T=eT$ as in Proposition~\ref{prop:typeA-Wronskian-frame-covariance}, its matrix $\mathscr W^T=\mathscr W_{\mathrm{univ}}^A T$ satisfies
\begin{equation}\label{eq:typeA-Wronskian-lines-arbitrary-frame}
 \begin{aligned}
 &\rho_a(Z)^{-1}(\textstyle\bigwedge^a\mathscr W^T(qu))
             (\textstyle\bigwedge^aT(qu))^{-1}v_{\omega_a}\\
 &\qquad=(\textstyle\bigwedge^a\mathscr W^T(u))
             (\textstyle\bigwedge^aT(u))^{-1}
                    \rho_a(s_A(u))^{-1}v_{\omega_a}.
 \end{aligned}
\end{equation}
\end{proposition}
\begin{proof}
Since $N_+$ fixes $v_{\omega_a}$, the hypothesis gives
\begin{equation}\label{eq:typeA-oper-line}
 \rho_a(A(u))\rho_a(s_A(u))^{-1}v_{\omega_a}=v_{\omega_a},
 \qquad 1\le a<N.
\end{equation}
Thus $\rho_a(s_A(u))^{-1}v_{\omega_a}$ is the coordinate column of the inverse pullback of $e_1\wedge\cdots\wedge e_a$. Applying the exterior power of the connection identity \eqref{eq:typeA-Wronskian-section} to $v_{\omega_a}$ gives \eqref{eq:typeA-fundamental-Wronskian}, in the conventions of \cite[Definition~4.9, Eq.~(4.9)]{KZ2021}. Substitution of $\mathscr W^T=\mathscr W_{\mathrm{univ}}^A T$ cancels $T$ and gives the formula in the new frame. The factors $(\bigwedge^aT)^{-1}$ express the specified vectors in that frame.
\end{proof}

\begin{remark}\label{rem:typeA-fixed-data}
If the data are defined over a subfield $K_0\subset K$ stable under the spectral shift and its inverse, the matrices $P^A$, $P^Z$ and the recurrence coefficients are defined over $K_0$. The identities \eqref{eq:typeA-frame-identities} and $P^Z(P^A)^{-1}=g$ hold after extension of the coefficient field, giving the $(G,q)$-Wronskian section of the extended data.
\end{remark}

\subsubsection{Steinberg form in the cyclic frame}\label{sec:typeA-cyclic-algebra}

Here we will express the connection in the cyclic frame $P^A=(p_0^A,\ldots,p_{N-1}^A)$ and compare the original oper flag with the new flag defined by the spans of the first vectors of this frame. The original flag retains its relative position $c$, because its value $\mathcal F_+(qu)$ and its transported value $A(u)\mathcal F_+(u)$ are both transformed by $(P^A(qu))^{-1}$. The new flag will give the linear Coxeter position.

\begin{definition}\label{def:typeA-cyclic-flag}
The \emph{cyclic oper flag} $\mathcal F_+^{\mathrm{cyc}}$ associated with the cyclic frame $P^A$ is the complete flag in $E_A$ whose rank-$a$ subspaces are
\[
 \mathcal F_{+,a}^{\mathrm{cyc}}(u)
 =\operatorname{span}_{\mathbb C}\{p_0^A(u),\ldots,p_{a-1}^A(u)\},
 \qquad 1\le a\le N,
\]
at every $u$ where the frame is defined.
\end{definition}
In the frame $P^A$, the rank-$a$ subspace is spanned by the first $a$ coordinate vectors, so $\mathcal F_+^{\mathrm{cyc}}$ is represented by the standard complete flag.  We show that the connection has relative position $c_{\mathrm{lin}}=s_{N-1}\cdots s_1$ with respect to $\mathcal F_+^{\mathrm{cyc}}$.  This is the \emph{linear Coxeter element}, whose simple reflections are ordered consecutively along the type-$A_{N-1}$ Dynkin diagram.

We write $P=P^A$ for the coordinate matrix of the cyclic frame and $\tau_j=\tau_{A,j}$ for its recurrence coefficients.  Its volume is $\Delta_A(u)=\det P(u)\in K^\times$.  We normalize this frame to volume one and show that the $q$-connection matrix in the resulting frame is in Steinberg form, with factors given by an upper unipotent matrix, a diagonal matrix, and a fixed lift of $c_{\mathrm{lin}}$.

We write $R=R_{\mathrm{cyc}}$ for the pullback matrix of Proposition~\ref{prop:typeA-cyclic-reconstruction}. In the frame $P=P^A$, its companion form is
\begin{equation}\label{eq:typeA-companion}
 R(u)=
 \begin{pmatrix}
 0&0&\cdots&0&\tau_0(u)\\
 1&0&\cdots&0&\tau_1(u)\\
 0&1&\cdots&0&\tau_2(u)\\
 \vdots&\ddots&\ddots&\vdots&\vdots\\
 0&0&\cdots&1&\tau_{N-1}(u)
 \end{pmatrix}.
\end{equation}

The determinant of this matrix measures the change of volume under pullback.

\begin{theorem}\label{thm:typeA-Wronskian}\label{lem:typeA-matrix-shift}
The inverse pullback $\sigma_A$ acts on coordinate columns in the cyclic frame $P=(p_0^A,\ldots,p_{N-1}^A)$ by $\mathbf x(u)\mapsto R(u)\mathbf x(qu)$. Equivalently,
\begin{equation}\label{eq:typeA-matrix-shift}
 A(u)^{-1}P(qu)=P(u)R(u).
\end{equation}
The determinant of $R$ and the recurrence coefficient $\tau_0$ satisfy
\begin{equation}\label{eq:typeA-tau0}
 \det R(u)=(-1)^{N-1}\tau_0(u)=\frac{\Delta_A(qu)}{\Delta_A(u)}.
\end{equation}
The $q$-connection from $u$ to $qu$ has matrix $R(u)^{-1}$ in this frame:
\begin{equation}\label{eq:typeA-qoper}
 P(qu)^{-1}A(u)P(u)=R(u)^{-1}.
\end{equation}
\end{theorem}
\begin{proof}
The first identity is \eqref{eq:typeA-abstract-universal-transport} with $P=P^A$ and $R=R_{\mathrm{cyc}}$. Taking determinants gives $\Delta_A(qu)=\Delta_A(u)\det R(u)$, since $\det A=1$. Expansion of the companion determinant gives $\det R=(-1)^{N-1}\tau_0$, and inversion gives \eqref{eq:typeA-qoper}.
\end{proof}

The determinant formula shows how the volume of the frame changes under pullback. Rescaling the last section gives a frame of volume one and keeps the oper-line generator fixed
\begin{equation}\label{eq:typeA-volume-correction}
 D_A(u)=\operatorname{diag}(1,\ldots,1,\Delta_A(u)^{-1}),
 \qquad P^\circ=P D_A,\qquad \det P^\circ=1.
\end{equation}
\begin{corollary}\label{cor:typeA-companion-volume-normalization}
We write $R^\circ(u)=D_A(u)^{-1}R(u)D_A(qu)$ for the pullback matrix in $P^\circ$. The $q$-connection in this frame has matrix
\begin{equation}\label{eq:typeA-volume-qtransport}
 (R^\circ(u))^{-1}
 =D_A(qu)^{-1}R(u)^{-1}D_A(u)
 =(P^\circ(qu))^{-1}A(u)P^\circ(u)\in SL_N(K).
\end{equation}
The original oper and Miura flags have coordinates $(P^\circ)^{-1}\mathcal F_+$ and $(P^\circ)^{-1}\mathcal F_-$.
\end{corollary}
\begin{proof}
Changing the frame of $E_A$ by $P^\circ=P D_A$ gives the displayed matrix and the two flags.  By \eqref{eq:typeA-tau0},
\[
 \det R^\circ(u)
 =\Delta_A(u)\frac{\Delta_A(qu)}{\Delta_A(u)}\Delta_A(qu)^{-1}=1.
\]
\end{proof}

For the cyclic flag, the following factorization of $(R^\circ(u))^{-1}$ separates the permutation of the basis vectors, the factors from the volume normalization, and the unipotent entries coming from the recurrence. For the linear word $c_{\mathrm{lin}}=s_{N-1}\cdots s_1$ the pinned lift $\dot c_{\mathrm{lin}}=\widehat s_{N-1}^{-1}\cdots\widehat s_1^{-1}$ acts by
\begin{equation}\label{eq:typeA-linear-Coxeter-lift}
 \dot c_{\mathrm{lin}}\mathbf e_1=\mathbf e_N,\qquad
 \dot c_{\mathrm{lin}}\mathbf e_j=-\mathbf e_{j-1}\quad(2\le j\le N).
\end{equation}
The matrices $E_{ij}$ are the standard matrix units. With this lift the unipotent subgroup and the Cartan factor in the volume-one cyclic frame are
\[
 U_{\mathrm{lin}}:=N_+\cap
       \dot c_{\mathrm{lin}}N_-\dot c_{\mathrm{lin}}^{-1},
 \qquad
 h_{\mathrm{cyc}}(u):=
       \operatorname{diag}\!\bigl(\underbrace{-1,\ldots,-1}_{N-2},
                      -\Delta_A(u)^{-1},(-1)^{N-1}\Delta_A(u)\bigr).
\]

\begin{proposition}\label{prop:typeA-cyclic-Steinberg}
The $q$-connection matrix $(R^\circ(u))^{-1}$ in the volume-one cyclic frame $P^\circ$ factors as
\begin{equation}\label{eq:typeA-cyclic-Steinberg}
 (R^\circ(u))^{-1}
 =\left(1-\sum_{j=1}^{N-1}
        \frac{\tau_j(u)}{\Delta_A(qu)}E_{jN}\right)
       h_{\mathrm{cyc}}(u)\dot c_{\mathrm{lin}}
 \in U_{\mathrm{lin}}(K)h_{\mathrm{cyc}}(u)\dot c_{\mathrm{lin}}.
\end{equation}
\end{proposition}

\begin{proof}
The positive roots taken to negative roots by $c_{\mathrm{lin}}^{-1}$ are exactly the roots of the matrix units $E_{jN}$, $1\le j<N$.  Hence
\[
 U_{\mathrm{lin}}
   =\left\{1+\sum_{j=1}^{N-1}x_jE_{jN}\right\}.
\]
For the factorization, we use $\varepsilon=(-1)^{N-1}$.
Its right-hand side has first column
\[
 \varepsilon \Delta_A(u)\mathbf e_N
   -\frac{\varepsilon \Delta_A(u)}{\Delta_A(qu)}
                \sum_{j=1}^{N-1}\tau_j(u)\mathbf e_j,
\]
columns $2,\ldots,N-1$ equal to $\mathbf e_1,\ldots,\mathbf e_{N-2}$, and last column
$\Delta_A(u)^{-1}\mathbf e_{N-1}$.  These are the columns of
$D_A(qu)^{-1}R(u)^{-1}D_A(u)$, since
$\tau_0(u)=\varepsilon \Delta_A(qu)/\Delta_A(u)$ by
\eqref{eq:typeA-tau0}.  Equation~\eqref{eq:typeA-volume-qtransport} identifies it with $(P^\circ(qu))^{-1}A(u)P^\circ(u)$.  The calculation also covers $N=2$, when there are no intermediate columns.
\end{proof}

Rescaling the last section to obtain $P^\circ$ preserves every subspace of $\mathcal F_+^{\mathrm{cyc}}$.  Hence this flag is also standard in the frame $P^\circ$.  By \eqref{eq:typeA-cyclic-Steinberg}, it gives a $q$-oper reduction with relative position $c_{\mathrm{lin}}$.

The Steinberg calculation uses the volume-one frame $P^\circ$. For the comparison with the column convention of \cite[\S4.1--4.2]{KSZ2021}, we use the cyclic frame $P^A$ before the volume correction. In this frame the columns of the $(G,q)$-Wronskian section are precisely the iterates of $p_0^Z$ in the fixed frame of $E_Z$.
\begin{proposition}
\label{prop:typeA-standard-Wronskian-comparison}
For the data of Definition~\ref{def:typeA-universal-Wronskian}, the $(G,q)$-Wronskian section has matrix $g$ in the chosen frames. Using $P^A$ on $E_A$ and the fixed frame on $E_Z$ gives the matrix $gP^A=P^Z$ by Proposition~\ref{prop:typeA-Wronskian-frame-covariance}, with columns
\begin{equation}\label{eq:typeA-standard-column-Wronskian}
 P^Z(u)
  =\bigl(p_0^Z(u),Z^{-1}p_0^Z(qu),\ldots,
                 Z^{-(N-1)}p_0^Z(q^{N-1}u)\bigr).
\end{equation}
The $q$-connection matrix in this frame is $R^{-1}$, with $R=R_{\mathrm{cyc}}$ from \eqref{eq:typeA-abstract-universal-transport}.
\end{proposition}

The columns in \eqref{eq:typeA-standard-column-Wronskian} give the unnormalized type-$A$ quantum Wronskian of \cite[Eqs.~(4.2), (4.6)--(4.8)]{KSZ2021}, where the twist parameter is $Z^{-1}$.  For nonzero meromorphic functions $d_0,\ldots,d_{N-1}$, the rescaled columns $\widetilde p_j^Z=d_jp_j^Z$ satisfy
\[
 (\sigma_Z\widetilde p_j^Z)(u)
 =\frac{d_j(qu)}{d_{j+1}(u)}\widetilde p_{j+1}^Z(u),
 \qquad 0\le j<N-1.
\]
These ratios are the diagonal multipliers in \cite[\S4.5, Eq.~(4.43)]{KZ2021}.  For data over $\mathbb C(u)$, the column construction is unchanged upon extension to $\mathcal M(U)$ by Remark~\ref{rem:typeA-fixed-data}.

\begin{remark}
Reversing the order of the cyclic basis conjugates $R^{-1}$ by $J_N=(\delta_{i,N+1-j})_{1\le i,j\le N}$ and puts it in subdiagonal companion form. The diagonal sign matrix $D_N=\operatorname{diag}((-1)^{N-1},\ldots,-1,1)$ gives the signed longest-Weyl lift
\begin{equation}\label{eq:typeA-signed-longest-lift}
 \dot w_{0,\mathrm{cyc}}:=D_NJ_N\in SL_N,
\end{equation}
obtained from simple-reflection lifts with block $\left(\begin{smallmatrix}0&-1\\1&0\end{smallmatrix}\right)$. Direct inversion of \eqref{eq:typeA-companion} gives
\begin{equation}\label{eq:typeA-standard-companion}
 J_NR(u)^{-1}J_N^{-1}
 =\left(\mathbf e_2,\ldots,\mathbf e_N,
 \frac{\mathbf e_1-\sum_{i=2}^N\tau_{N+1-i}(u)\mathbf e_i}{\tau_0(u)}\right).
\end{equation}
The signed lift instead gives
\begin{equation}\label{eq:typeA-signed-companion}
 \dot w_{0,\mathrm{cyc}}R(u)^{-1}\dot w_{0,\mathrm{cyc}}^{-1}
 =D_N\bigl(J_NR(u)^{-1}J_N^{-1}\bigr)D_N^{-1},
\end{equation}
with subdiagonal entries $-1$ and last-column entries multiplied by $(-1)^{N-i}$ in row $i$. Both matrices have determinant $\Delta_A(u)/\Delta_A(qu)$ by \eqref{eq:typeA-tau0}.
\end{remark}

\subsubsection{Relation to the recursive construction}

Here we will compare the two constructions for even $N$, when $A$ is in the bipartite Steinberg form and both apply to the same trivialization. We will compare the sections in the original frame and then express the connection in the volume-normalized cyclic frame.

\begin{corollary}\label{cor:typeA-even-Coxeter-comparison}
Suppose that $N=h=2d$ and $(A,g)$ is a datum of Definition~\ref{def:typeA-universal-Wronskian}, with
\[
 A\in U_c(K)s_{\mathrm{as}},\qquad c^d=w_0,
\]
for the prescribed factor $s_{\mathrm{as}}=h_{\mathrm{as}}\dot c$ of the bipartite word. The recursive and type-$A$ constructions give the same $(G,q)$-Wronskian. In the defining representation, their matrices satisfy
\begin{equation}\label{eq:typeA-even-Coxeter-exact-comparison}
 \mathscr W_{\mathrm{univ}}
 =\mathscr W_{\mathrm{univ}}^A
 =P^Z(P^A)^{-1}=g.
\end{equation}
In the volume-normalized cyclic frame $P^\circ=P^AD_A$, we use $P^{Z,\circ}:=P^ZD_A$ for the coordinate matrix of the corresponding volume-one frame of $E_Z$. The section and connection matrices satisfy
\begin{equation}\label{eq:typeA-finite-frame-qgauge}
 \begin{aligned}
 P^{Z,\circ}&=\mathscr W_{\mathrm{univ}}P^\circ,\\
 A_{\mathrm{cyc}}(u)
 &:=(P^\circ(qu))^{-1}A(u)P^\circ(u)\\
 &=(R^\circ(u))^{-1}
 =P^{Z,\circ}(qu)^{-1}ZP^{Z,\circ}(u).
 \end{aligned}
\end{equation}
The matrix $A$ has Steinberg form for $c$ with the original oper flag, while $A_{\mathrm{cyc}}$ has Steinberg form for $c_{\mathrm{lin}}$ with the cyclic oper flag. In the volume-normalized cyclic frame, the original oper flag is represented by $(P^\circ)^{-1}B_+$ and retains relative position $c$. The cyclic oper flag is represented by $B_+$ and has relative position $c_{\mathrm{lin}}$.
\end{corollary}

\begin{example}\label{ex:typeA-distinct-flag-orders}
In $SL_4$ the original and cyclic oper flags can differ. With the original flag standard, we use the bipartite word $w=s_1s_3s_2=c^{-1}$ and its pinned product $\dot w=\widehat s_1^{-1}\widehat s_3^{-1}\widehat s_2^{-1}$. A Steinberg matrix for this word and its cyclic frame have the form
\begin{equation}\label{eq:SL4-explicit-bipartite-cyclic}
 A(u)=\begin{pmatrix}
 a&-b&1&0\\1&0&0&0\\0&c_3&0&-1\\0&1&0&0
 \end{pmatrix},\qquad
 P(u)=\begin{pmatrix}
 1&0&0&0\\0&0&0&-1\\0&1&0&-b\\0&0&-1&-c_3
 \end{pmatrix},
\end{equation}
where $a=a(u)$, $b=b(u)$, and $c_3=c_3(u)$ are arbitrary meromorphic functions.  The second matrix gives the coordinates of the successive iterates of $p_0^A$ under $\sigma_A$, as in \eqref{eq:typeA-cyclic-orbits}.  Both matrices have determinant one.  Multiplication gives
\begin{equation}\label{eq:SL4-explicit-companion-gauge}
 P(qu)^{-1}A(u)P(u)
 =\begin{pmatrix}
 a(u)&1&0&0\\-b(qu)&0&1&0\\c_3(qu)&0&0&1\\-1&0&0&0
 \end{pmatrix}.
\end{equation}
The second cyclic subspace is $\operatorname{span}_{\mathbb C}(\mathbf e_1,\mathbf e_3)$ in the original coordinates, whereas the original second subspace is $\operatorname{span}_{\mathbb C}(\mathbf e_1,\mathbf e_2)$. The cyclic frame depends on the coefficients of $A$.  Its third column $P\mathbf e_3=-\mathbf e_4$ gives the coordinates of $p_2^A$, the second inverse-connection pullback of $p_0^A$.  This is also the image of $\mathbf e_1$ under the inverse square of the chosen lift, $\dot w^{-2}\mathbf e_1=-\mathbf e_4$.

In permutation notation, the linear Coxeter element $c_{\mathrm{lin}}=(1\,4\,3\,2)$ and the inverse bipartite element $c^{-1}=(1\,2\,4\,3)$ have different squares relative to the same standard flag.  The square of the first takes $\mathbf e_1$ to $\mathbf e_3$, whereas the square of the second takes it to $\mathbf e_4$.  The matrix $P$ relates the coordinates of the two oper flags.
\end{example}

\subsection{ADE differential reconstruction and monodromy parametrization}\label{sec:ADE-differential-reconstruction}\label{sec:finite-shift-realization}
In this subsection we apply the algebraic constructions to the differential pair $(\mathscr M,C)$. We first identify the type-$A$ frames of successive pullbacks with frames of differential solutions and derive their scalar Baxter equation. We then collect the ADE conclusions in one theorem. The recursive construction applies whenever $h$ is even, and the type-$A$ construction applies for every $N$.

\subsubsection{Differential realization in type \texorpdfstring{$A$}{A}}\label{sec:typeA-differential-realization}

Here we will apply the algebraic reconstruction to the differential pair $A=\mathscr M$, $g=C$ by expressing the consecutive solutions in the sectorial and Frobenius frames. Under Hypotheses~\ref{hyp:frob}--\ref{hyp:ccm}, these frames supply the two trivializations, and \eqref{eq:central-infinity-transport} gives the datum of Definition~\ref{def:typeA-universal-Wronskian}. For a solution $y$, the notation $[y]_{\mathcal Y}(u)$ denotes its coefficient column in $\mathcal Y$. Thus $\Phi_0^{(1)}$ has sectorial column $\mathbf e_1$, and its Frobenius column is the prescribed image of the oper-line generator
\begin{equation}\label{eq:typeA-first-coefficient-column}
 p_0^Z(u)=[\Phi_0^{(1)}]_{\mathcal Y_{\mathrm F,V_1}}(u)=C(u)\mathbf e_1.
\end{equation}
Its iterates have Frobenius coordinate columns
\begin{equation}\label{eq:typeA-column}
 \mathbf w_j(u):=Z^{-j}C(q^ju)\mathbf e_1,\qquad j\in\mathbb Z.
\end{equation}

\begin{theorem}\label{thm:typeA-universal-realization}
Under Hypotheses~\ref{hyp:frob}--\ref{hyp:ccm} in type $A$, the solutions $\Phi_s^{(1)},\ldots,\Phi_{s+N-1}^{(1)}$ form a meromorphic frame for every $s\in\mathbb Z$. The frame $\mathcal Y_{\mathrm{cyc}}=(\Phi_0^{(1)},\ldots,\Phi_{N-1}^{(1)})$ has sectorial coefficient matrix $P=P^A$ and Frobenius coefficient matrix $P^Z=(\mathbf w_0,\ldots,\mathbf w_{N-1})$, with
\begin{equation}\label{eq:typeA-analytic-frame-transition}
 \mathcal Y_{\mathrm{cyc}}(z;u)
 =\mathcal Y_{\mathrm{as}}(z;u)P(u)
 =\mathcal Y_{\mathrm F}(z;u)P^Z(u),
 \qquad P^Z=CP.
\end{equation}
The matrix $P\in GL_N(K)$ is independent of $z$. With $D_A$ from \eqref{eq:typeA-volume-correction}, the matrices $P^\circ=PD_A$ and $P^{Z,\circ}=P^ZD_A$ have determinant one, and the unique $(G,q)$-Wronskian section for the differential datum above is represented by
\begin{equation}\label{eq:typeA-universal-reference-comparison}
 \rho_1(\mathscr W_{\mathrm{univ}}^A)
 =P^ZP^{-1}
 =P^{Z,\circ}(P^\circ)^{-1}
 =C.
\end{equation}
It satisfies \eqref{eq:typeA-fundamental-Wronskian} in every exterior power with $s_A=s_{\mathrm{as}}$. When $N=h$ is even, the recursive construction gives the same section,
\[
 \mathscr W_{\mathrm{univ}}^A=\mathscr W_{\mathrm{univ}},
 \qquad \rho_1(\mathscr W_{\mathrm{univ}})=C.
\]

The common recurrence coefficients of Section~\ref{sec:typeA-cyclic-algebra} satisfy
\begin{equation}\label{eq:typeA-recurrence}
 \mathbf w_N(u)=\sum_{a=0}^{N-1}\tau_a(u)\mathbf w_a(u).
\end{equation}
The connection matrices in $\mathcal Y_{\mathrm{cyc}}$ and its volume-one normalization $\mathcal Y_{\mathrm{cyc}}^\circ:=\mathcal Y_{\mathrm{cyc}}D_A$ are $R^{-1}$ and $(R^\circ)^{-1}$, respectively. Their comparison with the sectorial realization is
\begin{equation}\label{eq:typeA-differential-Steinberg-comparison}
 \begin{aligned}
 P^{Z,\circ}&=CP^\circ,\\
 (R^\circ(u))^{-1}
 &=(P^\circ(qu))^{-1}\mathscr M(u)P^\circ(u)\\
 &=P^{Z,\circ}(qu)^{-1}ZP^{Z,\circ}(u).
 \end{aligned}
\end{equation}
The matrices $\mathscr M$ and $(R^\circ)^{-1}$ express the same shifted monodromy in the sectorial and cyclic frames. They have Steinberg forms for the bipartite and linear Coxeter elements with respect to the sectorial and cyclic oper flags, respectively.
\end{theorem}
\begin{proof}
Expanding each $\Phi_j^{(1)}$ in the sectorial and Frobenius frames gives two columns independent of $z$. Equations~\eqref{eq:actual-infinity-frame-column}--\eqref{eq:central-infinity-iterate} identify them with the iterates of $p_0^A$ and $p_0^Z$, respectively. In particular,
\[
 P(u)=\bigl(\mathbf e_1,\mathscr M^{[1]}(u)^{-1}\mathbf e_1,\ldots,
                  \mathscr M^{[N-1]}(u)^{-1}\mathbf e_1\bigr).
\]
This gives the frame identities and independence of $z$. The $q$-oper condition of Theorem~\ref{prop:analytic-sectorial-Steinberg} and Proposition~\ref{prop:typeA-cyclicity-from-reference-matching} give invertibility. Applying the invertible semilinear pullback $\sigma_A^s$ gives independence of every shifted block.

Proposition~\ref{prop:typeA-cyclic-reconstruction} gives the unique section, the common recurrence and volume, and the inverse-pullback identities \eqref{eq:typeA-abstract-universal-transport}. The common right factor $D_A$ cancels in the quotient, and \eqref{eq:typeA-volume-qtransport} gives the connection matrix in the volume-one cyclic frame.

The sectorial Steinberg form of Theorem~\ref{prop:analytic-sectorial-Steinberg} gives $\mathscr M s_{\mathrm{as}}^{-1}\in U_c(K)\subset N_+(K)$. This supplies the hypothesis for the exterior-power equations in Proposition~\ref{prop:typeA-exterior-Wronskian}. For even $N$, the same Steinberg form also allows the comparison with the recursive construction. Corollary~\ref{cor:typeA-even-Coxeter-comparison}, applied to $(A,g)=(\mathscr M,C)$, identifies the two universal $(G,q)$-Wronskians. Equation~\eqref{eq:typeA-differential-Steinberg-comparison} follows from the cyclic change of frame, and Proposition~\ref{prop:typeA-cyclic-Steinberg} gives the Steinberg form for the linear Coxeter element.
\end{proof}

The cyclic construction applies for every $N$ and agrees with the recursive construction when $N$ is even. The common recurrence can be computed by determinants in the sectorial or Frobenius frame. Cramer's rule applied to \eqref{eq:typeA-cyclic-recurrences} gives
\begin{equation}\label{eq:typeA-tau-minor}
 \tau_a(u)
 =\frac{
  \det\bigl(p^A_0,\ldots,p^A_{a-1},
             p^A_N,p^A_{a+1},\ldots,p^A_{N-1}\bigr)(u)
 }{\Delta_A(u)},
 \qquad 0\le a<N.
\end{equation}
The numerator is the determinant of $P(u)$ with the column $p_a^A(u)$ replaced by $p_N^A(u)$.
Since $\mathbf w_j=Cp_j^A$ and $\det C=1$, replacing the columns $p_j^A$ by $\mathbf w_j$ gives the same $\tau_a$. Taking the component at a Frobenius weight gives the scalar Baxter equation.

\begin{proposition}\label{prop:typeA-TQ}
For a weight vector $e_\mu$ in the fixed basis used to normalize the Frobenius solutions in the defining representation, with
$Z e_\mu=\zeta_\mu e_\mu$ and dual covector $e_\mu^*$, the corresponding central connection coefficient is
\begin{equation}\label{eq:typeA-Qmu}
 \mathcal Q_\mu(u):=\bigl\langle e_\mu^*,\rho_1(C(u))\mathbf e_1\bigr\rangle.
\end{equation}
The coefficient satisfies the scalar Baxter equation
\begin{equation}\label{eq:typeA-TQ-general}
 \zeta_\mu^{-N}\mathcal Q_\mu(q^Nu)
 =\sum_{a=0}^{N-1}\tau_a(u)\,\zeta_\mu^{-a}\mathcal Q_\mu(q^au).
\end{equation}
The same equation holds for
the rescaled function $Q_\mu^{(1)}$ defined by \eqref{eq:standard-Q-normalization}, and is the scalar form of
$P^Z(qu)^{-1}ZP^Z(u)=R^{-1}$.
\end{proposition}

\begin{proof}
We first compute the components of $\mathbf w_a$ in the Frobenius basis.  Equation~\eqref{eq:typeA-column} gives
\begin{equation}\label{eq:typeA-row-column}
 \bigl\langle e_\mu^*,\mathbf w_a(u)\bigr\rangle
 =\zeta_\mu^{-a}\mathcal Q_\mu(q^au).
\end{equation}
Applying $e_\mu^*$ to the column relation
$\mathbf w_N=\sum_{a=0}^{N-1}\tau_a\mathbf w_a$, we obtain the asserted scalar equation.  Conversely, since the Frobenius covectors form a basis, the scalar equations determine all components of this vector relation.  They therefore recover the companion recurrence itself.
\end{proof}

The inverse-pullback identities \eqref{eq:typeA-abstract-universal-transport} and the recurrence \eqref{eq:typeA-recurrence} also hold for any type-$A$ differential family with a global central connection and independent consecutive solution columns, by the same comparison of iterates. The determinant expressions have the same form as those used for the classical scalar Stokes relations in \cite[Sections~3--4]{DDT2000}.

For FFH connections, the determinant calculation in Section~\ref{sec:FFH-typeA-basis} gives a consecutive frame that is invertible at every spectral value. This also controls the coefficients of the recurrence.
\begin{corollary}\label{cor:FFH-typeA-entire-companion}
For a type-$A$ FFH connection with a normalized entire Frobenius frame, $P^Z$ and $(P^Z)^{-1}$ are entire in $u$. All coefficients $\tau_a$, the inverse companion $R^{-1}$, and the volume-normalized inverse companion $(R^\circ)^{-1}$ are entire as well. Moreover,
\begin{equation}\label{eq:typeA-log-tau0}
 \tau_0(u)=(-1)^{N-1}\frac{d_N(qu)}{d_N(u)},
\end{equation}
so $\tau_0$ has no zeros and equals $(-1)^{N-1}$ when $d_N=1$. The transition $P$ and its inverse are entire whenever $C$ is entire. These assertions remain valid at zeros of individual $Q$-functions.
\end{corollary}
\begin{proof}
Theorem~\ref{thm:FFH-typeA-cyclic-basis} gives an entire consecutive frame with nowhere-zero determinant $\Delta_Nd_N(u)$. Since the normalized Frobenius frame is entire with determinant one, $P^Z=\mathcal Y_{\mathrm F}^{-1}\mathcal Y_{\mathrm{cyc}}$ and its inverse are entire. Hence the recurrence vector $(P^Z)^{-1}\mathbf w_N$ is entire. The equality $\Delta_A=\det P=\det P^Z=\Delta_Nd_N$ and \eqref{eq:typeA-tau0} give \eqref{eq:typeA-log-tau0}. The nonvanishing of $\Delta_A$ and $\tau_0$, together with \eqref{eq:typeA-volume-qtransport}, gives the asserted entire inverses. Finally, $P=C^{-1}P^Z$ gives the last claim.
\end{proof}

The determinant of an entire consecutive frame can still vary with $u$. In the logarithmic rank-one example, the volume correction can be seen explicitly. For $k=3/4$ and $q=i$, Example~\ref{ex:typeA-logarithmic-volume} gives the phase and \eqref{eq:FFH-volume-multiplier} gives $d_2(u)=e^{-\pi i u^2/4}$. Using the standard realization $f_1=E_{21}$ of Section~\ref{sec:FFH-typeA-basis}, the normalization $m_1(f_1\mathbf e_1\wedge\mathbf e_1)=1$ gives $m_1=-\det$ and, by \eqref{eq:normalized-cyclic-Plucker-vector}, $\Delta_2=-1$. Thus \eqref{eq:typeA-entire-cyclic-determinant}, \eqref{eq:typeA-tau0}, and \eqref{eq:typeA-volume-correction} give
\[
 \Delta_A=-e^{-\pi i u^2/4},\qquad
 \tau_0=-e^{\pi i u^2/2},\qquad
 D_A=\operatorname{diag}(1,-e^{\pi i u^2/4}).
\]

\subsubsection{The ADE reconstruction theorem}

Here we will summarize the ADE results in one theorem. The sectorial frame gives the connection and the full $QQ$-system in every ADE type. The two reconstruction results identify the same central connection in their respective ranges, and the principal minors of the normalized sectorial matrix give the $q$-character parametrization.

\begin{theorem}\label{thm:general-final}\label{thm:general-orbit}\label{cor:typeA-even-finite}\label{cor:FFH-typeA-even-finite}
For an ADE differential family $\mathcal L_u$ as in Section~\ref{sec:general-hypotheses} satisfying Hypotheses~\ref{hyp:frob}--\ref{hyp:ccm}, the following statements hold over $K=\mathcal M(U)$.
\begin{enumerate}[label=(\roman*),topsep=3pt,itemsep=5pt]
\item In the normalized sectorial frame, shifted monodromy has the form
\begin{equation}\label{eq:unified-ADE-central-Wronskian}
 \mathscr M(u)=C(qu)^{-1}ZC(u)
 \in U_c(K)h_{\mathrm{as}}(u)\dot c.
\end{equation}
Together with the sectorial oper reduction and the preserved Frobenius reduction, it defines a $Z$-twisted Miura $(G,q)$-oper. The prescribed rescalings of the fundamental highest-weight minors of $C$ satisfy the full $QQ$-system of Proposition~\ref{prop:central-minor-QQ}.

\item If $h$ is even, the fixed bipartite Coxeter element satisfies $c^d=w_0$ for $d=h/2$. The central connection $C$ is the unique representative in $[C]_{\mathrm{Wr}}$ satisfying all shifted $(G,q)$-Wronskian equations~\eqref{eq:W-41-u}, and
\begin{equation}\label{eq:general-orbit-identification}
 C=\mathscr W_{\mathrm{rec}}
  =\mathscr W_{\mathrm{univ}}([C]_{\mathrm{Wr}}),
 \qquad N_d=1.
\end{equation}
\makeatletter
\begingroup
\edef\@currentlabel{\theequation}
\label{eq:C-equals-Wrec}
\endgroup
\makeatother
Moreover,
\begin{equation}\label{eq:C-two-longest-cells}
 C\in B_-w_0B_-\cap B_+w_0B_+.
\end{equation}
Pointwise, the cell statement holds wherever $C$ is regular and the required minors do not vanish.

\item In type $A_{N-1}$, for every $N$, the datum $A=\mathscr M$, $g=C$ has a unique type-$A$ $(G,q)$-Wronskian section $\mathscr W_{\mathrm{univ}}^A$, represented by $C$ in the normalized sectorial and Frobenius frames. Its expression through the two frames of successive pullbacks is \eqref{eq:typeA-universal-reference-comparison}. The connection matrix $(R^\circ)^{-1}$ in the volume-normalized cyclic frame has Steinberg form for $c_{\mathrm{lin}}$ with the cyclic oper flag.

\item If $q$ is not a root of unity, the principal minors of $S=h_{\mathrm{as}}^{-1}\mathscr M$ are
\[
 T_i=\operatorname{ev}_C(\chi_i),
\]
where $\chi_i$ is the fundamental $q$-character with the spectral normalization of Section~\ref{sec:consequences}, and $\operatorname{ev}_C$ is the shifted $Q$-ratio substitution of \eqref{eq:differential-Y-evaluation}. These minors determine $S$ uniquely and hence parametrize shifted monodromy by $\mathscr M=h_{\mathrm{as}}S$.
\end{enumerate}
For an untwisted ADE FFH oper with $0<k<1$, generic in the sense of Definition~\ref{def:generic-FFH-data}, the same statements apply to the $q$-connection associated with current differential monodromy, with $q=e^{-2\pi ik}$.
\end{theorem}

The sectorial form in (i) is Theorem~\ref{prop:analytic-sectorial-Steinberg}, and the statements about the Miura structure and the full $QQ$-system are Propositions~\ref{prop:one-step-differential-W} and~\ref{prop:central-minor-QQ}.

The identification of $C$ in (ii), including the two Bruhat cell inclusions, is Corollary~\ref{cor:differential-universal-Wronskian}. Part (iii) is Theorem~\ref{thm:typeA-universal-realization}.

Part (iv) follows from Proposition~\ref{prop:Cartan-factor-removal} and Theorem~\ref{thm:differential-qcharacter}. The FFH specialization uses Theorem~\ref{thm:FFH-verification} and Corollary~\ref{cor:FFH-current-monodromy}.

\subsection{The choice of Coxeter element}\label{sec:final-remarks}

The sectorial normalization singles out a Coxeter element. Varying the direction at infinity gives the two bipartite elements $c_1c_0$ and $c_0c_1$, which coincide in rank one. To realize the other Coxeter elements, we keep the underlying $q$-connection and its preserved flag fixed and change the oper reduction.

\begin{proposition}\label{prop:phase-Coxeter-orders}\label{rem:choice-Coxeter-element}
For the Cartan element $\Lambda_H$ of \eqref{eq:cyclic-Cartan-PF-vector} and the bipartite element $c=c_1c_0$, let
\[
 H_\theta=\operatorname{Re}(e^{i\theta}\Lambda_H),
\]
and let $\mathcal C_+$ denote the dominant Weyl chamber. If $H_\theta$ is regular and belongs to $w\mathcal C_+$, then
\[
 w^{-1}cw\in\{c,c^{-1}\}.
\]
Thus the choice of Weyl chamber used in the proof of Theorem~\ref{prop:analytic-sectorial-Steinberg} yields only the bipartite elements $c_1c_0$ and $c_0c_1$ as the regular phase direction varies; they coincide in rank one.
\end{proposition}
\begin{proof}
In rank one, $c=c^{-1}$ is the unique Coxeter element. Suppose the rank is greater than one. Both color classes are then nonempty, and the simple-root values in the proof of Theorem~\ref{prop:analytic-sectorial-Steinberg} give $H_\theta\in\mathcal C_+$ for $-\pi/h<\theta<0$. The first identity below follows by applying $c_1$ to the expression \eqref{eq:cyclic-Cartan-PF-vector} for $\Lambda_H$ and using the eigenvector equation for the Dynkin adjacency matrix from Lemma~\ref{lem:bipartite-Coxeter-identification}.  The second is the Coxeter eigenvalue relation in that lemma:
\[
 c_1\Lambda_H=\overline{\Lambda_H},
 \qquad c\Lambda_H=e^{2\pi i/h}\Lambda_H,
\]
and hence
\[
 H_{-\theta}=c_1H_\theta,\qquad
 H_{\theta+2\pi/h}=cH_\theta.
\]
It follows that $H_\theta\in c_1\mathcal C_+$ for $0<\theta<\pi/h$. These two intervals and their translates by $2\pi/h$ contain every regular phase direction: at their endpoints, a simple root vanishes after the corresponding Weyl transformation. The chambers encountered are therefore $c^k\mathcal C_+$ and $c^kc_1\mathcal C_+$. Their Weyl elements conjugate $c$ to $c$ and $c_1cc_1=c^{-1}$, respectively.
\end{proof}

The next proposition changes the $q$-oper reduction while retaining the underlying $q$-connection and its preserved $B_-$-reduction. It applies to every Coxeter element, without the half-power condition used in the recursive construction.

\begin{proposition}\label{prop:change-Coxeter-oper-reduction}
Suppose $G$ has ADE type and rank $r>1$, and $c=s_{i_1}\cdots s_{i_r}$ is a Coxeter element. Assume that
\[
 s(u)=h_*(u)\dot c,\qquad
 A(u)=g(qu)^{-1}Zg(u)\in U_c(K)s(u),
\]
where $h_*\in H(K)$, $\dot c$ is the pinned product for this word, and $g\in G(K)$. For $i=i_r$ and any constant representative $k_i\in N_G(H)$ of $s_i$, the element
\[
 s'=k_i^{-1}sk_i
 \quad\text{has Weyl image}\quad
 c'=s_i c s_i=s_i s_{i_1}\cdots s_{i_{r-1}}.
\]
There is a unique $\eta\in N_+(K)$ such that
\begin{equation}\label{eq:change-Coxeter-Steinberg}
 \begin{aligned}
 A'(u)&:=\eta(qu)k_i^{-1}A(u)k_i\eta(u)^{-1}
          \in U_{c'}(K)s'(u),\\
 g'(u)&:=g(u)k_i\eta(u)^{-1},\qquad
 A'(u)=g'(qu)^{-1}Zg'(u).
 \end{aligned}
\end{equation}
The oper flag represented by $B_+$ in the resulting frame is represented by $k_iB_+$ in the original frame. The preserved reduction has flag $g'^{-1}B_-$ in the resulting frame and is the same intrinsic reduction as before. Iterating these changes and commuting adjacent commuting simple reflections gives a Steinberg realization for every Coxeter element of $W$.
\end{proposition}
\begin{proof}
Writing $c=w s_i$, the word $w$ contains every simple reflection except $s_i$ once. Applying $w^{-1}$ to $\alpha_i$ gives a positive root different from $\alpha_i$. Indeed, before each reflection is applied, the coefficient at its node is zero; the reflection replaces it by the sum of the nonnegative coefficients at neighboring nodes. Since the diagram is connected and has rank greater than one, a neighbor of $i$ acquires a positive coefficient. A simple reflection sends every positive root other than its own simple root to a positive root, so
\[
 c^{-1}\alpha_i=s_iw^{-1}\alpha_i\in\Delta_+.
\]
Consequently $\alpha_i$ is absent from the roots of $U_c=N_+\cap\dot cN_-\dot c^{-1}$. Conjugation by $k_i^{-1}$ sends every root subgroup of $U_c$ into $N_+$, and therefore
\[
 k_i^{-1}Ak_i\in N_+(K)s'.
\]
For the pinned product $\dot c'$ of the new word $c'$, the factor $s'(\dot c')^{-1}$ lies in $H(K)$. Theorem~\ref{thm:difference-field-Steinberg} with this Cartan factor gives the unique $\eta$ and the Steinberg condition in \eqref{eq:change-Coxeter-Steinberg}. In particular $k_i=\widehat s_i$ gives $k_i^{-1}\dot c k_i=\dot c'$ by cancellation of the last simple-reflection factor. Substitution of $g'=gk_i\eta^{-1}$ gives its connection identity.

The change of frame is $k_i\eta^{-1}$. Since $\eta$ preserves $B_+$, the new oper flag has original coordinates $k_i\eta^{-1}B_+=k_iB_+$. The preserved flag transforms from $g^{-1}B_-$ to $\eta k_i^{-1}g^{-1}B_-=g'^{-1}B_-$. Finally, all Coxeter elements of a finite ADE Weyl group are conjugate by Proposition~\ref{prop:finite-Weyl-data}, and cyclic changes of the word together with commutations connect conjugate Coxeter elements \cite{ErikssonEriksson2009}. Iteration gives the last assertion.
\end{proof}

For the differential pair $(A,g)=(\mathscr M,C)$, this change gives the solution frame $\mathcal Y_{\mathrm{as}}k_i\eta^{-1}$ with transition matrix $C'=Ck_i\eta^{-1}$ to the original Frobenius frame. The recursive construction applies in these coordinates when the new Coxeter element still has half-power $w_0$. This frame is meromorphic in $u$ and need not have the original sectorial normalization. In rank one there is no change of Coxeter element.

\begin{corollary}\label{cor:changed-Coxeter-universal-Wronskian}
Suppose $h$ is even and $c'$ is a Coxeter element satisfying $(c')^{h/2}=w_0$. Starting from the differential pair $(\mathscr M,C)$ of Theorem~\ref{thm:general-final}, let $(\mathscr M',C')$ be obtained by the changes of Proposition~\ref{prop:change-Coxeter-oper-reduction}, with resulting lift $s'$ and $\mathscr M'\in U_{c'}(K)s'$. In rank one, the pair and lift are unchanged. The central connection $C'$ is the unique section in its right $N_+(K)$-orbit satisfying the shifted $(G,q)$-Wronskian equations~\eqref{eq:W-41-u} for the constant Cartan twist $Z$ and the lift $s'$ of $c'$. We denote this universal section by $\mathscr W_{\mathrm{univ}}^{\prime}$, so
\begin{equation}\label{eq:changed-Coxeter-universal-Wronskian}
 C'=\mathscr W_{\mathrm{univ}}^{\prime}
   =\mathscr W_{\mathrm{univ}}([C']_{\mathrm{Wr}}).
\end{equation}
The connection reconstructed from this section is
\begin{equation}\label{eq:changed-Coxeter-reconstructed-connection}
 \mathscr M'(u)
 =\mathscr W_{\mathrm{univ}}^{\prime}(qu)^{-1}
     Z\mathscr W_{\mathrm{univ}}^{\prime}(u)
 \in U_{c'}(K)s'(u).
\end{equation}
\end{corollary}
\begin{proof}
The transformed pair satisfies
\[
 Z^{-1}C'(qu)=C'(u)\mathscr M'(u)^{-1},
 \qquad \mathscr M'(u)\in U_{c'}(K)s'(u).
\]
Corollary~\ref{cor:universal-Wronskian-Steinberg} therefore identifies $C'$ with the universal representative for $(Z,s')$. Its connection identity gives \eqref{eq:changed-Coxeter-reconstructed-connection}.
\end{proof}

\section{The non-simply-laced correspondence}\label{sec:non-simply-laced-sketch}

In the preceding sections we studied the $q$DE/ODE/IM correspondence for Lie algebras of ADE type. For a generic untwisted Feigin--Frenkel--Hernandez (FFH) oper of one of these types, a full turn in $z$ accompanied by the spectral shift $u\mapsto qu$ gives the $q$-connection of a $Z$-twisted Miura $(G,q)$-oper. For even Coxeter number, the central connection matrix $C(u)$ between the Frobenius and sectorial frames is the universal $(G,q)$-Wronskian. In type~$A$ we also recovered $C(u)$ by comparing the Frobenius and sectorial realizations of a cyclic frame. In every ADE type the normalized minors of $C(u)$ satisfy the $QQ$-relations of the same affine type as the FFH oper. This agreement of types reflects affine Langlands duality, since the affine Cartan matrix in ADE is symmetric.

In this section we review how the correspondence extends beyond ADE. The results will appear with their proofs in~\cite{ZeitlinNonSimplyLaced}. For the untwisted quantum affine types $B_n^{(1)},C_n^{(1)},F_4^{(1)}$ and $G_2^{(1)}$, transposing the affine Cartan matrix gives the twisted affine type of the corresponding FFH oper. If the defining Dynkin diagram automorphism $\sigma$ has order $r$, rotation through $2\pi/r$ accompanied by a spectral shift $u\mapsto\eta u$ permutes the fundamental representations according to $\sigma$. The induced map from the fiber at $u$ to the fiber at $\eta u$ intertwines right multiplication by $g\in G$ with right multiplication by $\sigma(g)$. Together with the sectorial and Frobenius flags, it defines a Miura $(G,\eta,\sigma)$-oper, with full-turn multiplier $q=\eta^r$. The central connection matrix now represents the universal $(G,\eta,\sigma)$-Wronskian. Its minors give the $Q$-functions. From these minors the $q$-character formulas recover the Steinberg coordinates of the matrix representing the fractional rotation in the sectorial frame.

For an untwisted FFH oper of non-simply-laced type, a full turn gives an ordinary Miura $(G,q)$-oper whose Wronskian minors satisfy the folded $QQ$-relations. The self-dual twisted type $A_{2n}^{(2)}$ involves two flags, each of which defines an oper reduction of the full-turn $q$-connection. To see where these differences enter, we begin with the vacuum FFH equation and examine how rotation acts on the coefficient belonging to the affine vertex.

\subsection{The twisted and untwisted vacuum FFH equations}\label{sec:nsl-vacuum-models}

The twisted affine Lie algebras that occur on the differential side have a common description in terms of a Dynkin diagram automorphism. For quantum types $B_n^{(1)}$, $C_n^{(1)}$, $F_4^{(1)}$, and $G_2^{(1)}$, the affine Langlands-dual types are $A_{2n-1}^{(2)}$, $D_{n+1}^{(2)}$, $E_6^{(2)}$, and $D_4^{(3)}$, respectively. In each case we begin with the simply connected group $G$ of the corresponding simply-laced type and its Lie algebra $\g$. An automorphism of its Dynkin diagram is a permutation of the vertices that preserves the Cartan matrix. The pinning determines a lift $\sigma$ to $G$ whose action on $\g$ permutes the generators $e_i,f_i,\alpha_i^\vee$ in the same way. For $A_{2n-1}$ the automorphism reverses the chain, for $D_{n+1}$ it interchanges the two ends of the fork, and for $E_6$ it interchanges the two long arms. These involutions have order $r=2$. For $D_4^{(3)}$ we use the automorphism of order $r=3$ that cycles the three outer vertices of $D_4$. In the two classical families $n\ge2$, with $D_3=A_3$.

For $\varepsilon=e^{2\pi i/r}$, a twisted loop is a Laurent polynomial $X(t)$ with values in $\g$ satisfying $X(\varepsilon t)=\sigma(X(t))$. Its constant coefficient lies in the fixed Lie algebra $\g^\sigma$, of type $C_n$, $B_n$, $F_4$, or $G_2$, respectively. Besides twisted loops in $\g$, we can therefore consider ordinary loops with values in $\g^\sigma$. In either case the central extension together with the degree derivation gives the affine Lie algebra of the FFH oper.

As in Section~\ref{sec:FFH-verification}, a vacuum FFH oper has singularities only at zero and infinity. In rank-one quantum KdV its $Q$-functions give the eigenvalues of the Baxter $Q$-operators on the highest-weight line of the chosen Virasoro module~\cite{BLZ2001}. To describe higher states, the untwisted ADE construction allows additional regular singularities that are required to be apparent for every spectral parameter~\cite{BLZhigher2003,MRhigher2020}. A circuit around such a pole fixes every local solution, so the comparison of the Frobenius and sectorial frames is independent of windings around the added poles. For the four twisted affine types listed above, Masoero and Raimondo proved that with the residues prescribed by the standard FFH form, this condition forces the set of additional poles to be empty~\cite[Theorem~4.5(2)]{MR2024}. We also use the vacuum for the untwisted and self-dual families to compare the actions of rotation with the same two singular points. In the untwisted non-simply-laced case, the folded $QQ$-relations extend to perturbations with additional apparent poles satisfying the analytic hypotheses~\cite{ZeitlinNonSimplyLaced}.

In the current realization of the vacuum oper the differential part is $z\partial_z+k\lambda_{\mathrm{aff}}\partial_{\lambda_{\mathrm{aff}}}$ with $0<k<1$, so $u=z^{-k}\lambda_{\mathrm{aff}}$ is constant along its integral curves. Along $\lambda_{\mathrm{aff}}=z^ku$ the vacuum differential operator in either case has the form
\begin{equation}\label{eq:nsl-vacuum}
 \mathcal L_u=\partial_z+\frac{f+\ell+(z+uz^k)v}{z}.
\end{equation}
Here $f$ is the sum of the pinned negative simple-root generators of $\g^\sigma$, and $\ell$ lies in its Cartan subalgebra $\gh^\sigma$.

The coefficient of $t$ in the loop realization determines where $v$ lies. For the four twisted families above, this coefficient belongs to the $\varepsilon$-eigenspace of $\sigma$ on $\g$. As a representation of $\g^\sigma$, this eigenspace has highest weight equal to the highest short root of $\g^\sigma$. The vector $v$ spans its highest-weight line, and we consider the equation in representations of $G$. For ordinary loops, the coefficient of $t$ lies in $\g^\sigma$. In this case $v$ is a highest-root vector of $\g^\sigma$, and we use representations of the simply connected group with that Lie algebra.

The action of $\sigma$ on $v$ determines the rotation of the twisted equation in current coordinates. On the logarithmic cover of $\mathbb C^\times$, the rotation $R$ increases $\log z$ by $2\pi i/r$ and projects to $z'=\varepsilon z$. Rotating the affine loop coordinate at the same time sends
\[
 (z,\lambda_{\mathrm{aff}})\longmapsto(Rz,\varepsilon\lambda_{\mathrm{aff}}).
\]
Both terms in $z+\lambda_{\mathrm{aff}}$ then acquire the factor $\varepsilon$. Since $\sigma$ fixes $f$ and $\ell$ and sends $v$ to $\varepsilon v$, this rotation acts on the coefficient of $dz/z$ by $\sigma$. The arc from $z$ to $Rz$ projects to one positive circuit around zero in the coordinate $z^r$, and $\sigma$ identifies the solution fibers at the two endpoints of the rotation.

In the loop realization, the same rotation changes the spectral parameter $u=z^{-k}\lambda_{\mathrm{aff}}$ to
\[
 (Rz)^{-k}\varepsilon\lambda_{\mathrm{aff}}
   =\eta u,
 \qquad \eta=e^{2\pi i(1-k)/r}.
\]
Here $(Rz)^{-k}=\exp(-k(\log z+2\pi i/r))$ on the cover. After $r$ rotations, $z$ has made a full positive turn, the current coordinate $\lambda_{\mathrm{aff}}$ has returned to its initial value, and $\sigma^r=1$. The resulting spectral multiplier is $q=\eta^r=e^{-2\pi ik}$. For the untwisted equation we use a full turn at fixed $\lambda_{\mathrm{aff}}$, which gives the same shift $u\mapsto qu$.

The vacuum equation of self-dual twisted type $A_{2n}^{(2)}$ also has the form \eqref{eq:nsl-vacuum}. Here $G=\mathrm{SL}_{2n+1}$ and $\sigma$ reverses the chain of type $A_{2n}$. The coefficients are $f=\sum_{i=1}^{2n}E_{i+1,i}$ and $v=E_{1,2n+1}$, where $E_{ij}$ are matrix units. The pinned involution fixes $f$ and $\ell\in\gh^\sigma$ and sends $v$ to $-v$, giving the same rotation formula with $r=2$. There is one difference in the action on vertices. The middle edge of $A_{2n}$ joins two vertices in the same orbit, whereas in the four preceding twisted types the vertices of each orbit are pairwise nonadjacent.

For each of these equations we normalize the solutions at zero by their Frobenius expansions. The regularity and nonresonance conditions on $\ell$ in the residue $f+\ell$ are
\begin{equation}\label{eq:nsl-nonresonance}
 \alpha(\ell)\ne0,\qquad
 \lambda(\ell)\notin\{a+bk:a,b\in\mathbb Z_{\ge0},\ a+b>0\}.
\end{equation}
Here $\alpha$ ranges over the roots of the differential group of that equation, and $\lambda$ ranges over its roots and the differences of two weights in each of its fundamental representations. The first condition allows a constant lower unipotent change of frame taking the residue $f+\ell$ to $\ell$. The second excludes resonances in the Frobenius expansions in $z$ and $z^k$ at zero.

\subsection{Rotational covariance and the five analytic hypotheses}

The ADE construction of shifted monodromy uses the action of rotation on the frames normalized at zero and at infinity. The vacuum equations show what changes when a fractional rotation acts through $\sigma$. We can make the same comparison for any differential family with this symmetry and with the corresponding analytic properties of the two frames. In what follows $G$ denotes the group in whose representations we consider the differential equation, and $\g$ its Lie algebra. For a twisted family $\sigma$ is its pinned automorphism of order $r$. For an untwisted family $\sigma$ is the identity on $G$ and $r=1$. The rotation $R$ adds $2\pi i/r$ to $\log z$. Its spectral multiplier $\eta$ is the one computed from \eqref{eq:nsl-vacuum} for the vacuum families. In general $\eta$ is specified as part of the family, on a connected spectral domain $U$ satisfying $\eta U=U$.

In a trivialization the differential connection is $\nabla_u=d_z+\mathscr A(z;u)\,dz$, where $\mathscr A(z;u)\in\g$. The rotation law of the vacuum equation takes the intrinsic form $R^*\nabla_{\eta u}=\sigma(\nabla_u)$, which we require of the general family. In coordinates this means
\begin{equation}\label{eq:nsl-coefficient-covariance}
 \varepsilon\,\mathscr A(Rz;\eta u)=\sigma\bigl(\mathscr A(z;u)\bigr).
\end{equation}
The argument $Rz$ denotes the rotated point of the cover. The factor $\varepsilon=e^{2\pi i/r}$ outside $\mathscr A$ is the scalar in $dz'=\varepsilon\,dz$. For a full turn this scalar is $1$ and $R$ still advances $\log z$ by $2\pi i$. We use $\varepsilon^j z$ for the lifted point $R^jz$ in the frame comparisons.

By \eqref{eq:nsl-coefficient-covariance}, a horizontal group-valued frame $Y$ at $\eta u$ pulls back to the horizontal frame $\sigma^{-1}(Y(\varepsilon z;\eta u))$ at $u$. In the other direction, from $u$ to $\eta u$, right multiplication of a frame by $g$ becomes right multiplication of the rotated frame by $\sigma(g)$. This accounts for the change in the equivariance condition of Definition~\ref{def:q-connection}.

An \emph{$(\eta,\sigma)$-connection} on a principal right $G$-bundle $P$ is a meromorphic family of fiber isomorphisms $\mathcal A_u:P_u\to P_{\eta u}$ satisfying $\mathcal A_u(pg)=\mathcal A_u(p)\sigma(g)$. In a trivialization the map is $g\mapsto A(u)\sigma(g)$. A change of coordinates $g\mapsto b(u)g$ replaces $A(u)$ by $b(\eta u)A(u)\sigma(b(u))^{-1}$. The fiber map sends the flag $gB_\pm$ at $u$ to $A(u)\sigma(g)B_\pm$ at $\eta u$, since $\sigma$ preserves both Borel subgroups. The composite of $r$ such maps takes $P_u$ to $P_{qu}$ with $q=\eta^r$. It commutes with the right action of $G$ because $\sigma^r=1$. Thus a full turn gives an ordinary $q$-connection.

\begin{definition}\label{def:nsl-oper-flags}
An ordering of the $\sigma$-orbits of vertices and a choice of one vertex in each orbit determine an element $c$ of the Weyl group of $G$, the product of the corresponding simple reflections in the given order. When $\sigma=1$, this is an ordinary Coxeter element. A $(G,\eta,\sigma)$-oper consists of an $(\eta,\sigma)$-connection as above and a meromorphic $B_+$-reduction $F_+$ satisfying, at generic $u$,
\begin{equation}\label{eq:nsl-oper-flags}
 \operatorname{pos}\bigl(F_+(\eta u),A(u)\sigma(F_+(u))\bigr)=c .
\end{equation}
A Miura $(G,\eta,\sigma)$-oper also carries a meromorphic $B_-$-reduction $F_-$ satisfying $F_-(\eta u)=A(u)\sigma(F_-(u))$. It is $Z_0$-twisted if, in a meromorphic trivialization, the connection acts by $g\mapsto Z_0\sigma(g)$ for a fixed $Z_0\in H$.
\end{definition}

The rotation law identifies solution spaces at $u$ and $\eta u$. To determine its action on the flags we compare the normalized frames, as in ADE. At zero the Frobenius normalization gives the Cartan twist and the preserved flag. At infinity the sectorial asymptotics determine the relative position of the oper flags. Continuation of the two frames to a common domain relates these comparisons by the central connection matrix. We therefore use the five analytic hypotheses of Section~\ref{sec:five-analytic-hypotheses} with the following modifications for the action of $\sigma$.

The gauges $T_0,T_\infty$ and the choices of branches are as in that section. The asymptotic comparisons are uniform on closed subsectors and compact regular spectral sets. The exponent $1/h$ in the leading phase at infinity fixes the sector angle $\arg z/h$. For the four distinct affine Langlands pairs, $h$ is the dual Coxeter number of the finite quantum type, respectively $2n-1,n+1,9,4$. For an untwisted full turn, $h$ is the Coxeter number of the differential group and $\eta=q$. In the self-dual half-turn model $h=2n+1$. In all cases $R$ increases $\arg z/h$ by $2\pi/(rh)$.
\begin{enumerate}[label=(\arabic*),leftmargin=*,topsep=3pt,itemsep=2pt]
\item At zero the exponent $\ell$ belongs to the fixed Cartan subalgebra $\gh^\sigma$. The prescribed expansion class and leading weight vector together determine a unique Frobenius column. The resulting frames preserve the defining tensors, as in Hypothesis~\ref{hyp:frob}.
\item At infinity the sectorial normalization of Hypothesis~\ref{hyp:subdominant} uses a regular semisimple leading element $\Lambda$, a diagonalizer $a_\Lambda$ with $D=\operatorname{Ad}(a_\Lambda^{-1})\Lambda\in\gh$, and the exponential $e^{-S D}$. The phase $S$ has its additive constant fixed and includes any logarithmic term. If $\operatorname{Re}(e^{i\vartheta}D)$ is strictly dominant for $a_0<\vartheta<b_0$, the sector is $a_0-\pi<\arg z/h<b_0$. If its negative is strictly dominant there, the sector is $a_0<\arg z/h<b_0+\pi$.
\item The differential family satisfies \eqref{eq:nsl-coefficient-covariance}. At zero the pinning identifies $V_i$ with $V_{\sigma(i)}$ by the intertwiner carrying the chosen highest-weight vector to the chosen highest-weight vector. For each $j\in\mathbb Z$, the normalization is prescribed on a source sector and its $j$th rotated image. On the rotated sector the leading vector is the image of the original one under the intertwiner for $\sigma^j$. The column with this leading vector pulls back to the original Frobenius expansion class. The comparison of gauges $T_0(z;u)\sigma^{-j}(T_0(\varepsilon^jz;\eta^ju)^{-1})$ tends to $1$ at zero.
\item On the overlap of the sectors at infinity the complete phase satisfies $S(\varepsilon z;\eta u)=\omega S(z;u)+\delta(u)+o(1)$, where $\omega=e^{2\pi i/(rh)}$ and $\delta$ is holomorphic. In the diagonalizing frames the corresponding comparison of gauges is
\[
 a_\Lambda^{-1}T_\infty(z;u)\,\sigma^{-1}\!\bigl(T_\infty(\varepsilon z;\eta u)^{-1}a_\Lambda\bigr)
 \longrightarrow \sigma^{-1}(K_+^{-1}).
\]
Here $K_+\in N_G(H)$ is constant, represents $c$ in $N_G(H)/H$, and satisfies $\operatorname{Ad}(K_+)\sigma(D)=\omega D$.
\item The two normalized frames admit single-valued continuation to a common connected domain in the logarithmic cover, invariant under $R$, and depend meromorphically on $u$. On the prescribed lifts their restrictions agree with the normalized Frobenius germs on the source and rotated sectors at zero and with the normalized sectorial frames on the tails at infinity.
\end{enumerate}

For the four twisted families in distinct affine Langlands pairs, a constant change of frame makes $K_+$ in (4) equal to a prescribed representative of $c$. We use the ordered product $\dot c$ of $\widehat s_i^{-1}$ for the chosen orbit representatives. Right multiplication of the diagonalizer and sectorial frame by $b\in H$ changes $K_+$ to $b^{-1}K_+\sigma(b)$. Such a $b$ can be chosen to make this element equal to $\dot c$. Indeed, the torus homomorphism $b\mapsto b\,c(\sigma(b))^{-1}$ is surjective because $c\sigma$ has no fixed vector on $\gh$~\cite{ZeitlinNonSimplyLaced}. We use the sectorial frame with this normalization in what follows.

\begin{theorem}\label{thm:nsl-differential}
For the four affine Langlands pairs above, assume the five hypotheses hold. The central connection matrix $C(u)=\mathcal Y_{\mathrm F}(z;u)^{-1}\mathcal Y_{\mathrm{as}}(z;u)$ is meromorphic in $u$ and independent of $z$. In Frobenius coordinates the forward comparison from $u$ to $\eta u$ sends $g$ to $Z_0\sigma(g)$, where $Z_0=e^{2\pi i\ell/r}$. In sectorial coordinates it sends $g$ to $\mathscr M(u)\sigma(g)$, where
\begin{equation}\label{eq:nsl-differential-bridge}
 \mathscr M(u)=\mathcal Y_{\mathrm{as}}(\varepsilon z;\eta u)^{-1}\sigma(\mathcal Y_{\mathrm{as}}(z;u))
              =C(\eta u)^{-1}Z_0\sigma(C(u)).
\end{equation}
The matrix $\mathscr M(u)$ belongs to $U_cs(u)$, where $U_c=N_+\cap\dot cN_-\dot c^{-1}$ and $s(u)=e^{\delta(u)D}\dot c$. The sectorial flag $B_+$ is the oper reduction, and the Frobenius flag $B_-$ gives the preserved reduction, represented by $C(u)^{-1}B_-$ in sectorial coordinates. The full-turn $q$-connection has matrix $C(qu)^{-1}e^{2\pi i\ell}C(u)$ in the sectorial frame.
\end{theorem}

For a twisted FFH oper, we can read \eqref{eq:nsl-differential-bridge} directly in the current realization. The sectorial frame in these coordinates is $\Phi(z,\lambda_{\mathrm{aff}})=\mathcal Y_{\mathrm{as}}(z;z^{-k}\lambda_{\mathrm{aff}})$. It satisfies
\[
 \sigma\!\left(\Phi(z,\lambda_{\mathrm{aff}})\right)
   =\Phi(Rz,\varepsilon\lambda_{\mathrm{aff}})\,\mathscr M(u),
 \qquad u=z^{-k}\lambda_{\mathrm{aff}}.
\]
In each fundamental representation, the columns of $\mathscr M(u)$ express the vectors of the frame $\sigma(\Phi(z,\lambda_{\mathrm{aff}}))$ in the sectorial basis at the endpoint. The forward rotation therefore carries a frame with coordinate $g$ at $u$ to one with coordinate $\mathscr M(u)\sigma(g)$ at $\eta u$.

To recover the differential monodromy around zero, we compose the $r$ rotations. Each rotation applies $\sigma$ to the factors contributed by the preceding rotations, so the full-turn forward matrix is
\[
 \mathscr M(\eta^{r-1}u)\,\sigma\!\left(\mathscr M(\eta^{r-2}u)\right)
 \cdots\sigma^{r-1}\!\left(\mathscr M(u)\right)
 =C(qu)^{-1}e^{2\pi i\ell}C(u).
\]
The equality follows from \eqref{eq:nsl-differential-bridge} since $\sigma(Z_0)=Z_0$. With the convention of Corollary~\ref{cor:FFH-current-monodromy}, positive continuation around zero at fixed $\lambda_{\mathrm{aff}}$ acts by pullback. It therefore multiplies the sectorial frame on the right by the inverse matrix $C(u)^{-1}e^{-2\pi i\ell}C(qu)$.

The Cartan element $Z_0$ records the change in the Frobenius frame at zero. On a line of weight $\mu$ the fractional turn multiplies the leading factor $z^{-\mu(\ell)}$ by $e^{-2\pi i\mu(\ell)/r}$, so the forward comparison contributes $Z_0=e^{2\pi i\ell/r}$. At infinity the change of the full phase contributes $e^{\delta(u)D}$ to $s(u)$. To obtain the Wronskian equation from the fractional rotation, we use $\mathscr M(u)s(u)^{-1}\in N_+$. Since $N_+$ fixes each $v_{\omega_i}$, we have $\mathscr M(u)^{-1}v_{\omega_i}=s(u)^{-1}v_{\omega_i}$. Substitution in \eqref{eq:nsl-differential-bridge} gives the equation for each column $C(u)v_{\omega_i}$.

In the trivialization with connection $g\mapsto Z_0\sigma(g)$, a generalized $(G,\eta,\sigma)$-Wronskian is a meromorphic section represented by a $G$-valued function $\mathscr W$ satisfying
\[
 \mathscr W(\eta u)v_{\omega_i}
   =Z_0\sigma(\mathscr W(u))s(u)^{-1}v_{\omega_i}
\]
for every vertex $i$ of $G$. After $j$ rotations each lift has been acted on by $\sigma$ once for every rotation following it. The lifts therefore occur in the ordered product
\[
 s_{[j]}(u)=s(\eta^{j-1}u)\sigma(s(\eta^{j-2}u))\cdots\sigma^{j-1}(s(u)).
\]

\begin{theorem}\label{thm:nsl-universal}
For the four twisted affine families above and the $\sigma$-invariant twist $Z_0=e^{2\pi i\ell/r}$, each generalized Wronskian has a unique representative with the same fundamental highest-weight columns satisfying
\begin{equation}\label{eq:nsl-finite-wronskian}
 \mathscr W(\eta^ju)v_{\omega_i}
   =Z_0^j\sigma^j(\mathscr W(u))s_{[j]}(u)^{-1}v_{\omega_i},
 \qquad 1\le j\le rh/2 .
\end{equation}
These equations are equivalent to $\mathscr W(\eta u)^{-1}Z_0\sigma(\mathscr W(u))\in U_cs(u)$. Under the five analytic hypotheses their unique solution with highest-weight columns $C(u)v_{\omega_i}$ is $C(u)$, the universal $(G,\eta,\sigma)$-Wronskian.
\end{theorem}

At $m=rh/2$ the concatenation of the words for $c,\sigma(c),\ldots,\sigma^{m-1}(c)$ is a reduced expression for the longest Weyl element. The last equation thus relates the highest-weight columns to the lowest-weight lines. The intermediate rotations determine which arguments of the columns occur in the spectral relations.

\subsection{The spectral relations of the twisted differential families}
In each of the four twisted FFH families considered above, a fractional turn permutes the fundamental representations by the diagram automorphism $\sigma$. The vertices of the finite Dynkin diagram on the quantum side correspond to the orbits of this permutation. The edges of the original diagram determine how many neighboring $Q$-functions occur at each vertex of the quantum diagram. More precisely, if $\widetilde A$ is the Cartan matrix of $\g$, $I_\sigma$ is a set of orbit representatives, and $o_j$ is the size of the orbit of $j$, then the finite Cartan matrix on the quantum side is
\begin{equation}\label{eq:nsl-orbit-cartan}
 A_{ij}=\sum_{a=0}^{o_j-1}\widetilde A_{i,\sigma^a(j)},\qquad i,j\in I_\sigma .
\end{equation}
Distinct orbits are adjacent when an edge of the original diagram joins them. For $i\ne j$, the integer $-A_{ij}$ counts the neighbors of $i$ in the orbit of $j$. The fixed Lie algebra $\g^\sigma$, in the same orbit order, has Cartan matrix $A^{\mathsf T}$.

To obtain the $Q$-functions, we expand a distinguished sectorial solution in the Frobenius basis as in~\eqref{eq:Q-expansion}. On a dominant ray this solution is $\mathcal Y_{\mathrm{as},V_i}(z;u)v_{\omega_i}$. Since $\mathcal Y_{\mathrm{as},V_i}=\mathcal Y_{\mathrm F,V_i}\rho_i(C)$, the vector of its Frobenius coefficients is $\rho_i(C(u))v_{\omega_i}$.

The diagram automorphism acts on the Weyl group by $s_j\mapsto s_{\sigma(j)}$. For a Weyl element $w$ fixed by $\sigma$, the extremal weight spaces of weights $w\omega_i$ and $ws_i\omega_i$ are lines spanned respectively by $\rho_i(\widehat w)v_{\omega_i}$ and $\rho_i(\widehat w\widehat s_i)v_{\omega_i}$. We use the Frobenius coefficients along these two vectors. The functional $v_{\omega_i}^*$ takes the value $1$ on $v_{\omega_i}$ and vanishes on the other weight spaces. With the pinned lifts of Section~\ref{sec:lie-theoretic-data}, the two coefficients are
\[
\begin{aligned}
 \mathcal Q_{w,+}^{(i)}(u)
   &=v_{\omega_i}^*\!\left(\rho_i(\widehat w^{-1}C(u))v_{\omega_i}\right)
     =\Delta_{w\omega_i,\omega_i}(C(u)),\\
 \mathcal Q_{w,-}^{(i)}(u)
   &=v_{\omega_i}^*\!\left(\rho_i(\widehat s_i^{-1}\widehat w^{-1}C(u))v_{\omega_i}\right)
     =\Delta_{s_i\omega_i,\omega_i}(\widehat w^{-1}C(u)).
\end{aligned}
\]
Indeed, $\rho_i(\widehat w^{-1})$ carries the first vector to $v_{\omega_i}$ and the second to $\rho_i(\widehat s_i)v_{\omega_i}=-f_iv_{\omega_i}$. In the second formula $\rho_i(\widehat s_i^{-1})$ then carries this vector to $v_{\omega_i}$. This order of the lifts fixes the sign even when $ws_i$ has smaller length than $w$.

After $o_i$ fractional turns, the column that started in $V_i$ belongs to $V_i$ again. We can therefore compare the column at $u$ with the column brought back from $\eta^{o_i}u$ in the same representation. In their exterior product, the coefficient along the wedge of the two extremal vectors above is the determinant of the corresponding Frobenius coefficients. The Pl\"ucker identity expresses this determinant as a product with one factor for each neighbor $\beta$ of $i$ in the Dynkin diagram of $G$. To find the argument of that factor, let $j$ represent the orbit of $\beta$. Bringing the column in $V_j$ back from $\eta^nu$ to $u$ puts it in $V_{\sigma^{-n}(j)}$, so the exponent $n$ must satisfy $\sigma^{-n}(j)=\beta$. Among such integers, the ordered word $c$ selects the largest $n\le o_i$ if $j$ follows $i$ in $c$, and the largest $n<o_i$ if $j$ precedes $i$. We denote the resulting integer by $n_{i,\beta}$.

If the orbits of $i$ and $j$ are adjacent, the $-A_{ij}$ neighbors $\beta$ of $i$ in the orbit of $j$ have consecutive exponents $n_{i,\beta}$. We keep the ADE convention that the determinant at $i$ compares the arguments $u$ and $\eta^{o_i}u$. To express the neighboring factors in this convention, we shift the spectral argument at each vertex by $\eta^{d_i}$ and normalize its two Frobenius coefficients by constant factors determined by $p_i$. For an adjacent pair $i,j$ of orbit representatives, write $\bar n_{ij}$ for the mean of $n_{i,\beta}$ over neighbors $\beta$ of $i$ in the orbit of $j$. With $\epsilon_i\in\{0,1\}$ a coloring that assigns opposite colors to neighbors in the orbit diagram, the integers $d_i$ and the real numbers $p_i$ satisfy
\[
 d_j-d_i=\min_\beta n_{i,\beta}-\epsilon_i\min(o_i,o_j),
 \qquad p_j-p_i=\bar n_{ij}-o_i/2 .
\]
The minimum is over the same neighbors as the mean. The orbit diagram is a tree, so setting $d_i=p_i=0$ at one vertex determines both families. The normalized coefficients are
\[
 Q_{w,\pm}^{(i)}(u)=e^{-2\pi i p_i\mu(\ell)/r}\,
                 \mathcal Q_{w,\pm}^{(i)}(\eta^{d_i}u),
 \qquad
 \mu=w\omega_i\ \text{for }+,\quad \mu=ws_i\omega_i\ \text{for }- .
\]

\Needspace{11\baselineskip}
\begin{theorem}\label{thm:nsl-spectral}
Under the five analytic hypotheses, $Q_{w,+}^{(i)}$ is not identically zero for any $i\in I_\sigma$ and any $w$ fixed by $\sigma$. For each $i$ there is a unique nonzero meromorphic function $\varphi_i$, independent of $w$, such that, with $\mathfrak z_{w,i}=\exp(\pi i o_i(w\alpha_i)(\ell)/r)$, the normalized Frobenius coefficients satisfy
\begin{equation}\label{eq:nsl-qq}
\begin{aligned}
 &\mathfrak z_{w,i}Q_{w,+}^{(i)}(u)Q_{w,-}^{(i)}(\eta^{o_i}u)
 -\mathfrak z_{w,i}^{-1}Q_{w,+}^{(i)}(\eta^{o_i}u)Q_{w,-}^{(i)}(u)\\
 &\qquad=\varphi_i(\eta^{d_i}u)\prod_{j\sim i}\prod_{a=0}^{-A_{ij}-1}
   Q_{w,+}^{(j)}\bigl(\eta^{\epsilon_i\min(o_i,o_j)+a}u\bigr).
\end{aligned}
\end{equation}
Here $j\sim i$ denotes adjacency in the orbit diagram. If $s$ is constant, a constant change of sectorial frame normalizes $s$ to $\dot c$. In this normalization all $\varphi_i=1$, and \eqref{eq:nsl-qq} gives the $QQ$-system of the untwisted quantum affine algebra with finite Cartan matrix $A$.
\end{theorem}

The coefficient $\varphi_i$ in \eqref{eq:nsl-qq} records the torus factor $h_c(u)=s(u)\dot c^{-1}=e^{\delta(u)D}$ of the sectorial lift. We compute it by choosing coordinates in which the lift from $\eta^lu$ to $\eta^{l+1}u$ equals $\dot c$ for every $l\in\mathbb Z$. If the change of coordinates at position $l$ is $g\mapsto t_l(u)g$, with $t_l(u)\in H$, the lift becomes $t_{l+1}(u)s(\eta^lu)\sigma(t_l(u))^{-1}$. Keeping the original coordinates at position zero, we obtain
\begin{equation}\label{eq:nsl-torus-recurrence}
 t_0(u)=1,\qquad t_{l+1}(u)=\dot c\,\sigma(t_l(u))\,\dot c^{-1}h_c(\eta^l u)^{-1}
 \qquad(l\in\mathbb Z).
\end{equation}
Since $\sigma$ and conjugation by $\dot c$ act invertibly on $H$, this recurrence also determines $t_l(u)$ for negative $l$.

In these coordinates $C(\eta^lu)$ is replaced by $C(\eta^lu)t_l(u)^{-1}$, which multiplies its $j$th highest-weight column by $\omega_j(t_l(u))^{-1}$. The Pl\"ucker identity for the pinned lift has coefficient one. Undoing these rescalings in the determinant at positions $0$ and $o_i$ and in the neighboring factors at positions $n_{i,\beta}$ gives
\[
 \varphi_i(u)=\frac{\omega_i(t_{o_i}(u))}
 {\displaystyle\prod_{\substack{\beta\text{ adjacent to }i\\
                                  \text{in the diagram of }G}}
           \omega_{j(\beta)}(t_{n_{i,\beta}}(u))},
\]
where $j(\beta)$ is the representative of the orbit of $\beta$.

For example, in the $A_5$ diagram for the quantum type $B_3$, vertex $3$ has neighbors $2$ and $4$, which lie in the same orbit. With $\epsilon_3=0$, their contribution in \eqref{eq:nsl-qq} is $Q_{w,+}^{(2)}(u)Q_{w,+}^{(2)}(\eta u)$.

For the four twisted vacuum families in \eqref{eq:nsl-vacuum}, the spectral multiplier is $\eta=e^{2\pi i(1-k)/r}$. With the quantum parameter $\qaff=e^{-\pi i(1-k)/r}$, the two arguments in the determinant at $i$ differ by the factor $\eta^{o_i}=\qaff^{-2o_i}$. For these families the verification of the five hypotheses in~\cite{ZeitlinNonSimplyLaced} gives the following corollary.

\begin{corollary}\label{cor:nsl-vacuum-qq}
For each of the four twisted vacuum families, the residue conditions \eqref{eq:nsl-nonresonance} imply that the five hypotheses hold and $C(u)$ is entire in $u$. If $1/[h(1-k)]\notin\mathbb Z_{>0}$, the phase has no logarithm and $\varphi_i=1$ after normalization. At the exceptional exponents \eqref{eq:nsl-qq} retains the coefficient of the sectorial lift.
\end{corollary}

\subsection{Full turns and the folded \texorpdfstring{$QQ$}{QQ}-system}

For the untwisted vacuum equation in \eqref{eq:nsl-vacuum}, $G$ is the simply connected group of type $B_n,C_n,F_4$ or $G_2$ in whose representations we solve the equation. A full turn preserves each fundamental representation. When a neighboring vertex contributes with multiplicity greater than one, all its factors therefore have the same spectral argument. Their product gives the power of a single $Q$-function in the folded relation.

We write $A_{ij}=\alpha_j(\alpha_i^\vee)$ for the Cartan matrix of $G$ and $h$ for its Coxeter number. The vertices in the color class $I_1$ precede those in $I_0$ in the chosen ordering, which determines the meaning of $j<i$ and $j>i$ in the folded relations. Within each color $b$ the reflections in $c_b=\prod_{i\in I_b}s_i$ commute, as do their pinned representatives in $P_b=\prod_{i\in I_b}\widehat s_i^{-1}$. The sector corresponds to $c=c_1c_0$.

\begin{theorem}\label{thm:nsl-folded}
Under the residue conditions \eqref{eq:nsl-nonresonance} for $\g$, the five full-turn hypotheses hold for the vacuum. Its central connection matrix $C(u)$ is entire and represents the universal Wronskian of a $Z$-twisted Miura $(G,q)$-oper, where $Z=e^{2\pi i\ell}$. On the sector $a_0-\pi<\arg z/h<b_0$ of hypothesis~(2), a constant Cartan normalization gives
\[
 C(qu)^{-1}ZC(u)\in U_c\, e^{\delta(u)D}(P_0P_1)^{-1}.
\]
The principal and simple-reflection minors $Q_+^{(i)}(u)=\Delta_{\omega_i,\omega_i}(C(u))$ and $Q_-^{(i)}(u)=\Delta_{s_i\omega_i,\omega_i}(C(u))$ satisfy
\begin{equation}\label{eq:nsl-folded-forward}
\begin{aligned}
 &\left(\zeta_i\prod_{j>i}\zeta_j^{A_{ji}}\right)Q_-^{(i)}(u)Q_+^{(i)}(qu)
 -\left(\zeta_i^{-1}\prod_{j<i}\zeta_j^{-A_{ji}}\right)Q_-^{(i)}(qu)Q_+^{(i)}(u)\\
 &\hspace{2em}=\Lambda_i(u)\prod_{j>i}Q_+^{(j)}(qu)^{-A_{ji}}
                          \prod_{j<i}Q_+^{(j)}(u)^{-A_{ji}},
\end{aligned}
\end{equation}
where $\zeta_i=\omega_i(Z)^{-1}$. The lift gives $\Lambda_i(u)=\exp[-\delta(u)\omega_i(D)]$ for $i\in I_0$ and $\Lambda_i(u)=\exp[\delta(u)(\omega_i-\alpha_i)(D)]$ for $i\in I_1$. When $1/[h(1-k)]\notin\mathbb Z_{>0}$ the phase has no logarithmic term and $\Lambda_i=1$.
\end{theorem}

The order of the factors and the signs in \eqref{eq:nsl-folded-forward} use the same pinning as in ADE. In particular, $\widehat s_i v_{\omega_i}=-f_iv_{\omega_i}$, so the Frobenius coefficient along $f_iv_{\omega_i}$ is $-Q_-^{(i)}$. The component of the exterior product of the columns at $u$ and $qu$ on which we use the Pl\"ucker identity has highest weight $2\omega_i-\alpha_i=\sum_{j\ne i}(-A_{ji})\omega_j$. The coefficient $-A_{ji}$ gives the number of copies of the neighboring highest-weight vector at $j$. The order in $c$ fixes their common argument, $u$ for $j<i$ and $qu$ for $j>i$.

The principal minors in the twisted and folded $QQ$-relations also give the Miura variables of Section~\ref{sec:Miura-Steinberg}, and hence the Steinberg coordinates of shifted monodromy. To use the same character formula in both cases, we express the full-turn sectorial lift relative to $\dot c=P_1P_0$. In the normalization of Theorem~\ref{thm:nsl-folded} its Cartan factor is $h_c(u)=e^{\delta(u)D}(P_0P_1)^{-1}\dot c^{-1}$. The changes of coordinates $t_l(u)$ at the successive points $q^lu$ satisfy \eqref{eq:nsl-torus-recurrence} with this $h_c$, $\sigma=1$, and $\eta=q$. For these choices the recurrence agrees with \eqref{eq:monodromy-torus-frames} in the simply-laced case.

\begin{corollary}\label{cor:nsl-character}
For each of the four twisted vacuum families of Corollary~\ref{cor:nsl-vacuum-qq} and the four untwisted families of Theorem~\ref{thm:nsl-folded}, let $j$ be a vertex whose simple reflection occurs in $c$, and let $l\in\mathbb Z$. In the coordinates given by $t_l(u)$ at the source and $t_{l+1}(u)$ at the target, the Miura variable of the isomorphism from $\eta^lu$ to $\eta^{l+1}u$ is
\begin{equation}\label{eq:nsl-character-minors}
 \mathsf q_{j,l}(u)=
 \frac{\omega_j(t_{l+1}(u)Z_0)}{\omega_{\sigma^{-1}(j)}(t_l(u))}
 \frac{\Delta_{\sigma^{-1}(j)}(C(\eta^lu))}{\Delta_j(C(\eta^{l+1}u))}.
\end{equation}
Here $\Delta_j(g)=\Delta_{\omega_j,\omega_j}(g)$ is the principal minor, and $Z_0=Z$ for a full turn. The minor at the source has index $\sigma^{-1}(j)$ because the automorphism carries the corresponding highest-weight line to the line at $j$ in the target. The Miura variables determine the Steinberg coordinates of $h_c(u)^{-1}\mathscr M(u)$ also when the phase contains a logarithmic term.
\end{corollary}

For the four twisted families the Steinberg coordinates are the fundamental $q$-characters of the untwisted quantum affine algebra evaluated at the Miura variables \eqref{eq:nsl-character-minors}, with indices fixed by the order of the reflections. For untwisted differential type $B_n,C_n,F_4,G_2$, respectively, we use the fundamental $q$-characters of types $D_{n+1}^{(1)},A_{2n-1}^{(1)},E_6^{(1)},D_4^{(1)}$. At each spectral label we identify the variables belonging to vertices of the same diagram orbit. For this evaluation the quantum deformation parameter is chosen to have infinite order independently of $\eta$. The character variables are indexed by integers $l$, and the substitution evaluates the Miura variables at the corresponding shifts $\eta^lu$. Thus the evaluation also applies when $\eta$ has finite order~\cite{ZeitlinNonSimplyLaced}.

\subsection{Two decorated flags in the self-dual type \texorpdfstring{$A_{2n}^{(2)}$}{A2n(2)}}

The twisted affine type $A_{2n}^{(2)}$ is Langlands self-dual up to a relabeling of the affine vertices. Here $G=\mathrm{SL}_N$, with $N=2n+1$ and $n\ge1$, and $\sigma$ reverses the $A_{2n}$ chain. In the four twisted families considered above, no two vertices in the same orbit are adjacent. Here $\sigma$ interchanges the two ends $n$ and $n+1$ of the middle edge. Their simple reflections do not commute, and reversal exchanges their order in every Coxeter word. Thus $\sigma$ fixes no Coxeter element. If a single flag satisfied $F(-u)=\sigma(F(u))$, its generic relative position under the full-turn connection with $\sigma$-fixed twist would be fixed by $\sigma$. This rules out such a flag as an oper reduction. Instead the symmetry exchanges two oper reductions. To keep track of the vectors used to define the $Q$-functions, we use \emph{decorated flags}. The coset $WN_+$ specifies both the flag $WB_+$ and the vector $Wv_{\omega_i}$ on each of its fundamental highest-weight lines.

For the $\mathrm{SL}_N$ vacuum equation of Section~\ref{sec:nsl-vacuum-models}, a half turn sends $u$ to $\eta u$, with $\eta=e^{\pi i(1-k)}$, and acts in the Frobenius frame by $g\mapsto T\sigma(g)$, where $T=e^{\pi i\ell}$. Reversing the chain interchanges its two color classes $I_0,I_1$. We label them so that the half-turn lift represents $c_1$ and use the products $c_b,P_b$ over $I_b$ as above. Composing the half turn with the symmetry $(u,g)\mapsto(-u,\sigma(g))$ gives the ordinary connection $(u,g)\mapsto(\kappa u,Tg)$ with $\kappa=-\eta$. This connection compares one decorated flag at $u$ with the other at $\kappa u$. Its square compares each flag with itself at the full-turn shift $qu$, where $q=\kappa^2$ and the connection matrix is $T^2$.

\begin{theorem}\label{thm:nsl-self-dual}
Under \eqref{eq:nsl-nonresonance} for $\mathfrak{sl}_N$ and $1/[N(1-k)]\notin\mathbb Z_{>0}$, the central connection matrix $C$ is entire. After a constant Cartan normalization it represents the universal $(G,\eta,\sigma)$-Wronskian with half-turn lift $P_1^{-1}$, the inverse of the pinned product over color $1$. The decorated flags represented by $W_0(u)=C(u)$ and $W_1(u)=\sigma(C(-u))$ satisfy, for $b\in\{0,1\}$,
\begin{equation}\label{eq:nsl-two-flags}
 W_{1-b}(\kappa u)\in T W_b(u)N_+P_bN_+,\qquad
 W_{1-b}(-u)N_+=\sigma(W_b(u))N_+.
\end{equation}
Each flag $W_b(u)B_+$ defines an oper reduction of the full-turn $q$-connection with matrix $T^2$. Its relative position is $(c_bc_{1-b})^{-1}$, and the constant flag $B_-$ gives the Miura reduction in both cases.

For each orbit representative $1\le i\le n$, its color $\epsilon_i$ is defined by $i\in I_{\epsilon_i}$. The projections of $W_{\epsilon_i}(\kappa^{\epsilon_i}u)v_{\omega_i}$ onto the $\omega_i$ and $s_i\omega_i$ weight lines in $V_i$ are $Q_+^{(i)}(u)v_{\omega_i}$ and $Q_-^{(i)}(u)f_iv_{\omega_i}$, respectively. With $\mathfrak z_i=\alpha_i(T)$ and $Q_+^{(0)}=1$, these coefficients satisfy
\begin{equation}\label{eq:nsl-self-dual-qq}
\begin{aligned}
 &\mathfrak z_i Q_+^{(i)}(u)Q_-^{(i)}(qu)
   -\mathfrak z_i^{-1}Q_+^{(i)}(qu)Q_-^{(i)}(u)\\
 &\hspace{2em}=\begin{cases}
 Q_+^{(i-1)}(q^{\epsilon_i}u)Q_+^{(i+1)}(q^{\epsilon_i}u),&i<n,\\
 Q_+^{(n-1)}(q^{\epsilon_n}u)Q_+^{(n)}(-\kappa u),&i=n.
 \end{cases}
\end{aligned}
\end{equation}
\end{theorem}

At vertex $n$, the neighbor across the middle edge is $\sigma(n)=n+1$. The second identity in~\eqref{eq:nsl-two-flags} identifies its highest-weight coefficient at $\kappa u$ with that at vertex $n$ at $-\kappa u$. This gives the last factor in~\eqref{eq:nsl-self-dual-qq}. The construction also applies to half-turn families satisfying the five hypotheses on a domain invariant under negation and $\eta^{\pm1}$, provided the twist is $\sigma$-fixed, $\delta=0$, and the half-turn lift represents $c_1$~\cite{ZeitlinNonSimplyLaced}.

\subsection{The correspondence for vacuum affine opers}
The table summarizes the construction from the differential equation to the quantum $QQ$-relation. The first column gives the affine Lie type of the vacuum FFH oper. The second describes the difference oper defined by its shifted monodromy and gives the finite Lie type of its structure group $G$. The third records the $QQ$-system satisfied by the generalized minors of the central connection matrix $C(u)$. The fourth identifies the quantum affine algebra whose representation theory gives the corresponding functional relation. For the representation-theoretic statements we assume that the quantum deformation parameter $\qaff$ has infinite multiplicative order~\cite{ZeitlinNonSimplyLaced}.

\par\addvspace{\medskipamount}
\noindent\begin{minipage}{\linewidth}
\makeatletter\def\@captype{table}\makeatother
\centering
\caption{Difference opers and spectral relations for the vacuum affine opers.}
\label{tab:nsl-correspondences}
\small
\setlength{\tabcolsep}{4pt}
\renewcommand{\arraystretch}{1.12}
\begin{tabular}{@{}llll@{}}
\hline
\shortstack[l]{Differential\\affine type} &
\shortstack[l]{Difference oper\\and group type} &
\shortstack[l]{Spectral\\$QQ$-system} &
\shortstack[l]{Quantum affine\\type}\\
\hline
$X^{(1)}$, $X$ of ADE type & $(G,q)$, $X$ & Simply laced & $X^{(1)}$\\
\hline
$A_{2n-1}^{(2)}$ & $(G,\eta,\sigma)$, $A_{2n-1}$ & Untwisted $B_n$ & $B_n^{(1)}$\\
$D_{n+1}^{(2)}$ & $(G,\eta,\sigma)$, $D_{n+1}$ & Untwisted $C_n$ & $C_n^{(1)}$\\
$E_6^{(2)}$ & $(G,\eta,\sigma)$, $E_6$ & Untwisted $F_4$ & $F_4^{(1)}$\\
$D_4^{(3)}$ & $(G,\eta,\sigma)$, $D_4$ & Untwisted $G_2$ & $G_2^{(1)}$\\
\hline
$A_{2n}^{(2)}$ & Two oper flags, $A_{2n}$ & Twisted $A_{2n}^{(2)}$ & $A_{2n}^{(2)}$\\
\hline
$B_n^{(1)}$ & $(G,q)$, $B_n$ & Folded $B_n$ & $D_{n+1}^{(1)}$\\
$C_n^{(1)}$ & $(G,q)$, $C_n$ & Folded $C_n$ & $A_{2n-1}^{(1)}$\\
$F_4^{(1)}$ & $(G,q)$, $F_4$ & Folded $F_4$ & $E_6^{(1)}$\\
$G_2^{(1)}$ & $(G,q)$, $G_2$ & Folded $G_2$ & $D_4^{(1)}$\\
\hline
\end{tabular}
\end{minipage}\par\addvspace{\medskipamount}

In the first six rows of Table~\ref{tab:nsl-correspondences}, the differential and quantum affine types are Langlands dual. The first row recovers the ADE correspondence of this paper. The affine type $X^{(1)}$ is self-dual, a full turn gives the ordinary $(G,q)$-oper, and the minors of $C(u)$ satisfy the simply-laced $QQ$-system of type $X^{(1)}$. In rows two through five the differential type is twisted and the quantum type is untwisted. A half turn acts on the fundamental representations through the diagram involution $\sigma$ in the first three of these pairs, and a third of a turn acts through triality in type $D_4^{(3)}$. A column in $V_i$ returns to the same representation after $o_i$ fractional turns, where $o_i$ is the length of the orbit of $i$. This gives the argument $\eta^{o_i}u$ on the left of~\eqref{eq:nsl-qq}. On the right, neighbors of $i$ in the same orbit contribute factors at successive spectral arguments. Their product has coefficient one in the quantum identity. In~\eqref{eq:nsl-qq} the sectorial lift gives the coefficient $\varphi_i(\eta^{d_i}u)$, which equals one for the pinned lift $s=\dot c$.

In the sixth row both affine types are $A_{2n}^{(2)}$. The two oper reductions in the second column come from the decorated flags in~\eqref{eq:nsl-two-flags}. The connection with multiplier $\kappa=-\eta$ compares the two flags, while its square compares the flag of each color at $u$ and $qu$. This accounts for the same ratio $q=\kappa^2$ of the two arguments on the left of~\eqref{eq:nsl-self-dual-qq} at every vertex. The factor $Q_+^{(n)}(-\kappa u)$ at the last vertex records the exchange of the ends of the middle edge. Both features occur in Wang's twisted relation~\cite{WangTwistedQQ}.

For each of these six rows, the quantum identity holds in the weight completion of the Grothendieck ring of category $\mathcal O$ for the positive Borel subalgebra of the quantum affine algebra in the fourth column. For example, the differential type $A_{2n-1}^{(2)}$ corresponds to the positive Borel of $U_{\qaff}(B_n^{(1)})$. In the quantum loop convention the positive Borel is generated by the positive Chevalley generators $E_j$ and the invertible Cartan generators $K_j^{\pm1}$ at all vertices of the affine Dynkin diagram, including vertex~$0$. A module in $\mathcal O$ is a direct sum of finite-dimensional Cartan weight spaces. Its weights are bounded above by a finite set in the order defined by the finite simple roots. In the Grothendieck ring, a short exact sequence $0\to U\to V\to W\to0$ gives $[V]=[U]+[W]$, and tensor products give multiplication. The two quantum $Q$-symbols represent classes of auxiliary modules normalized by ordinary characters and one-dimensional weight classes. Completion in the weight order allows the ordinary characters to be treated as formal series. Wherever the corresponding transfer matrices are defined, their simultaneous eigenvalues satisfy the scalar relation that we compare with the identity among Wronskian minors~\cite{FrenkelHernandez2018,WangTwistedQQ}.

The last four rows start with untwisted FFH opers of types $B_n^{(1)}$, $C_n^{(1)}$, $F_4^{(1)}$ and $G_2^{(1)}$. Their full-turn $(G,q)$-opers give the folded relation~\eqref{eq:nsl-folded-forward}. The fourth column now specifies the simply-laced quantum affine algebra whose diagram is folded. Its transfer matrices act on a finite-dimensional simple type~$1$ module whose highest $\ell$-weight is invariant under the diagram automorphism. In terms of the Drinfeld polynomials specifying this highest $\ell$-weight, invariance means that vertices in the same orbit carry the same polynomial.

These transfer matrices depend on a twist $\zeta$ in the Cartan torus of the simply connected group of the indicated simply-laced type. In the trace over the auxiliary module, the contribution from the $\mu$-weight space is multiplied by the scalar $\mu(\zeta)$. For the positive and negative Baxter operators the auxiliary modules are infinite-dimensional, and the trace is first formed as a series in the Cartan weights. We require $\sigma(\zeta)=\zeta$, where $\sigma$ is the diagram automorphism, so that vertices in the same orbit are compared at the same twist. We keep the spectral parameter formal and assume that each coefficient of the Baxter operators and of the normalized transfer matrices entering the folded relation is regular at $\zeta$.

Within a weight space whose Cartan weight is fixed by $\sigma$, the folded relation holds on each joint generalized eigenspace on which the Baxter eigenvalue series agree at vertices in the same orbit. This agreement is coefficient by coefficient at the same spectral argument, separately for the positive and negative operators. For example, the involution of $D_{n+1}$ exchanges the two terminal vertices. Their common Baxter eigenvalue therefore appears squared in the folded $B_n$ relation at the vertex adjacent to both.

In the folded $QQ$-relation at vertex $i$, the product of neighboring Baxter eigenvalues has coefficient one. In~\eqref{eq:nsl-folded-forward}, the corresponding product of principal minors of $C$ carries the factor $\Lambda_i(u)$. This factor records the additive constant in the comparison of the sectorial phases at infinity under a full turn with $u\mapsto qu$. With our normalization the constant vanishes when the phase has no logarithmic term, giving $\Lambda_i(u)=1$. It remains open whether the same untwisted FFH opers also yield spectral functions satisfying the standard twisted $QQ$-system of the affine Langlands dual, in which the diagram automorphism acts on the spectral argument as well as on the vertex~\cite{ZeitlinNonSimplyLaced}.

\appendix
\refstepcounter{section}
\section*{Appendix}\label{sec:appendix-distinguished-solutions}

When the Coxeter number $h$ is even, the distinguished solutions give a direct proof that the central connection $C(u)$ is the universal $(G,q)$-Wronskian. Under Hypotheses~\ref{hyp:frob}--\ref{hyp:ccm}, $C(u)$ expresses the sectorial frame at infinity in the Frobenius frame at zero. We compute the Frobenius coefficients of each rotated distinguished solution in two ways. At infinity its leading asymptotics determine its coefficient vector in the sectorial frame at $u$, including the scalar on the corresponding extremal weight line. At zero the rotation acts through the constant twist $Z$. Comparing the coefficients gives the shifted Wronskian equations, from which we recover the Steinberg form of the $q$-connection from $u$ to $qu$.

For a fundamental representation $\rho_i:G\to GL(V_i)$ of highest weight $\omega_i$, with chosen highest-weight vector $v_{\omega_i}$, the distinguished solution and its pullback under $n$ rotations are
\[
 \Phi_0^{(i)}(z;u)=\mathcal Y_{\mathrm{as},V_i}(z;u)v_{\omega_i},
 \qquad
 \Phi_n^{(i)}(z;u)=\Phi_0^{(i)}(e^{2\pi in}z;q^nu).
\]
Here $e^{2\pi in}z$ denotes the point on the universal cover whose argument has increased by $2\pi n$. Hypothesis~\ref{hyp:symanzik} ensures that the pulled-back column solves the differential equation at $u$. We follow the rotations up to $n=d=h/2$. At this value the bipartite Coxeter element $c$ fixed in Section~2 satisfies $c^d=w_0$, the longest element of the Weyl group (Corollary~\ref{cor:fixed-datum-half-Coxeter}).

The sectorial normalization determines how rotation changes the leading vector at infinity. In the basis determined by $a_\Lambda$, the leading element of the differential equation is $D=\operatorname{Ad}_{a_\Lambda^{-1}}\Lambda\in\gh$. The phase $S(e^{2\pi i}z;qu)$ differs from $\gamma_h S(z;u)$ by a term tending to $\delta(u)$, where $\gamma_h=e^{2\pi i/h}$. This limit includes any additive constant contributed by a logarithmic term in $S$. The corresponding group element is $s_{\mathrm{as}}(u)=h_{\mathrm{as}}(u)\dot c$, with $h_{\mathrm{as}}(u)=\exp(\delta(u)D)$ and $\dot c$ the pinned representative of $c$. Pullback through $n$ rotations gives the ordered product
\[
 \Pi_0(u)=1,\qquad
 \Pi_n(u)=s_{\mathrm{as}}(u)^{-1}\cdots s_{\mathrm{as}}(q^{n-1}u)^{-1}
 \quad(n\ge1).
\]
The rightmost factor acts first, in the pullback from $q^nu$ to $q^{n-1}u$. Each factor lies in $N_G(H)$, the normalizer of the Cartan torus $H$. For $x\in N_G(H)$ and a vector $v$ of weight $\mu$, the identity $t(xv)=\mu(x^{-1}tx)\,xv$ for $t\in H$ says that $xv$ has weight $w\mu$, where $w=xH\in W=N_G(H)/H$. Since $h_{\mathrm{as}}(q^ju)\in H$, each factor of $\Pi_n(u)$ represents $c^{-1}$ in $W$. Thus $\Pi_n(u)$ carries the highest-weight line to the line of weight $c^{-n}\omega_i$. To recover the rotated solution from this line, we must also identify the scalar multiplying $\rho_i(\dot c^{-n})v_{\omega_i}$.

\begin{proposition}\label{prop:appendix-shifted-solution-columns}
For each fundamental representation $V_i$ and $0\le n\le d$, the rotated distinguished solution satisfies
\begin{equation}\label{eq:appendix-shifted-solution-columns}
 \Phi_n^{(i)}(z;u)=\mathcal Y_{\mathrm{as},V_i}(z;u)
                    \rho_i(\Pi_n(u))v_{\omega_i}.
\end{equation}
\end{proposition}
\begin{proof}
At $n=0$, the identity is the definition of $\Phi_0^{(i)}$. For $1\le n\le d$, we compare both sides on a sector where the normalizations at $u,qu,\ldots,q^nu$ all apply. The sectorial frame is normalized on $a-\pi<\arg z/h<b$, where $-\pi/h<a<b<0$. A rotation increases $\arg z/h$ by $2\pi/h$, so the common interval is
\[
 a-\pi<\arg z/h<b-2\pi n/h .
\]
This interval is nonempty since $2\pi n/h\le\pi$ and $a<b$. Each pair of consecutive rotations then lies in the overlap of Hypothesis~\ref{hyp:common-orbit}. Only finitely many shifts occur, so the comparisons are uniform on closed subsectors and for $u$ in compact subsets of the spectral domain.

Applying the phase comparison at $u,qu,\ldots,q^{n-1}u$ gives
\[
 S(e^{2\pi in}z;q^nu)=\gamma_h^n S(z;u)+\delta_n(u)+o(1),
 \qquad
 \delta_n(u)=\sum_{j=0}^{n-1}\gamma_h^{\,n-1-j}\delta(q^ju).
\]
The accumulated change of phase satisfies $\delta_{n+1}(u)=\gamma_h\delta_n(u)+\delta(q^nu)$ with $\delta_0=0$. The same recurrence determines the torus factor of $\Pi_n(u)$. At $n=0$ we have $\Pi_0=1$. If $\Pi_n(u)=\dot c^{-n}e^{-\delta_n(u)D}$, multiplication on the right by $s_{\mathrm{as}}(q^nu)^{-1}$ and the relation $\dot cD\dot c^{-1}=\gamma_hD$ give
\[
 \Pi_{n+1}(u)
 =\dot c^{-n}e^{-\delta_n(u)D}\dot c^{-1}e^{-\delta(q^nu)D}
 =\dot c^{-(n+1)}e^{-(\gamma_h\delta_n(u)+\delta(q^nu))D}.
\]
By the recurrence for $\delta_n$, induction gives
\[
 \Pi_n(u)=\dot c^{-n}e^{-\delta_n(u)D},\qquad
 \rho_i(\Pi_n(u))v_{\omega_i}
   =e^{-\delta_n(u)\nu_i}\rho_i(\dot c^{-n})v_{\omega_i},
 \qquad \nu_i=\omega_i(D).
\]
Here $\nu_i$ is the eigenvalue of the infinitesimal action of $D$ on the highest-weight vector.

We now express both leading vectors in the gauge $a_\Lambda^{-1}T_\infty(z;u)$ at the original point. In the gauge $a_\Lambda^{-1}T_\infty(e^{2\pi in}z;q^nu)$ at the rotated point, multiplication by $e^{\nu_iS(e^{2\pi in}z;q^nu)}$ makes the distinguished solution tend to $v_{\omega_i}$ by \eqref{eq:general-sectorial-columns}. The change from the rotated gauge to the original one has limit
\[
 a_\Lambda^{-1}T_\infty(z;u)\,T_\infty(e^{2\pi in}z;q^nu)^{-1}a_\Lambda
       \longrightarrow \dot c^{-n}.
\]
Indeed, the $n$ successive changes of gauge in \eqref{eq:analytic-gauge-rotation} each tend to $\dot c^{-1}$, and the intermediate gauges cancel in their product. Applying this change of gauge to the normalized column and using the phase comparison gives
\begin{equation}\label{eq:appendix-shifted-solution-asymptotics}
 \begin{aligned}
 &\rho_i(a_\Lambda^{-1}T_\infty(z;u))\Phi_n^{(i)}(z;u)\\
 &\qquad=e^{-\gamma_h^n\nu_i S(z;u)}
    \bigl(e^{-\delta_n(u)\nu_i}\rho_i(\dot c^{-n})v_{\omega_i}+o(1)\bigr).
 \end{aligned}
\end{equation}
For the sectorial solution on the right of \eqref{eq:appendix-shifted-solution-columns}, the coefficient vector $\rho_i(\Pi_n(u))v_{\omega_i}$ has weight $c^{-n}\omega_i$. The eigenvalue of $D$ on that vector is $(c^{-n}\omega_i)(D)=\omega_i(c^nD)=\gamma_h^n\nu_i$. By the normalization at the unshifted point,
\[
 \rho_i(a_\Lambda^{-1}T_\infty(z;u))\mathcal Y_{\mathrm{as},V_i}(z;u)
      =F_{V_i}(z;u)e^{-S(z;u)\rho_i(D)},\qquad F_{V_i}(z;u)\longrightarrow I_{V_i}.
\]
Applying this frame to $\rho_i(\Pi_n(u))v_{\omega_i}$ gives precisely the exponential and the leading vector in \eqref{eq:appendix-shifted-solution-asymptotics}.

To show that this leading term determines the solution, we use a ray on which the exponential for the weight $c^{-n}\omega_i$ decays faster than those for all other weights of $V_i$. For $\theta\in(a,b)$, the real Cartan element $H_\theta=\operatorname{Re}(e^{i\theta}D)$ lies in the open dominant chamber. Indeed, \eqref{eq:cyclic-Cartan-PF-vector} gives
\[
 \alpha_j(H_\theta)=
 \begin{cases}
 |\alpha_j(D)|\sin(\theta+\pi/h),& j\in I_0,\\
 -|\alpha_j(D)|\sin\theta,& j\in I_1.
 \end{cases}
\]
Here $I_0$ and $I_1$ are the two color classes of the Dynkin diagram. Both expressions are positive because $-\pi/h<\theta<0$. On the ray $\arg z/h=\theta-2\pi n/h$, which lies in the common sector since $2\pi n/h\le\pi$, the corresponding Cartan element is $H_{\theta-2\pi n/h}=c^{-n}H_\theta$. This follows from $c^{-n}D=\gamma_h^{-n}D$ and the real action of $c$ on the span of the coroots. For any weight $\mu\ne c^{-n}\omega_i$ of $V_i$, its translate $c^n\mu$ is again a weight of $V_i$ and is distinct from the highest weight. Hence $\omega_i-c^n\mu$ is a nonzero sum of simple roots with nonnegative integral coefficients, and
\[
 (c^{-n}\omega_i-\mu)(H_{\theta-2\pi n/h})
      =(\omega_i-c^n\mu)(H_\theta)>0 .
\]
For fixed $u$, the difference of the two solutions is $\mathcal Y_{\mathrm{as},V_i}(z;u)\mathbf b$ with $\mathbf b$ independent of $z$. Its normalized leading vector vanishes by \eqref{eq:appendix-shifted-solution-asymptotics}. Multiplying the difference in the normalizing gauge by $e^{\gamma_h^n\nu_iS(z;u)}F_{V_i}(z;u)^{-1}$ and taking its component in the space of weight $\mu$ therefore gives
\[
 e^{((c^{-n}\omega_i)(D)-\mu(D))S(z;u)}\mathbf b_\mu\longrightarrow0 .
\]
On the chosen ray $S(z;u)=re^{i(\theta-2\pi n/h)}+o(r)$, where $r=h|z|^{1/h}\to\infty$. For $\mu\ne c^{-n}\omega_i$, the strict inequality above makes the modulus of the exponential tend to infinity, forcing $\mathbf b_\mu=0$. For $\mu=c^{-n}\omega_i$ the exponential is one, so this component also vanishes. The solutions agree on the sectorial overlap. By Hypothesis~\ref{hyp:ccm}, both sides continue as single-valued solutions on the connected differential domain $\mathscr X$, which is invariant under the rotations. Their equality extends throughout $\mathscr X$ and holds as a meromorphic identity in $u$.
\end{proof}

We can now compute the Frobenius coefficients of $\Phi_n^{(i)}$ at $u$. At $q^nu$ the distinguished solution has coefficient vector $\rho_i(C(q^nu))v_{\omega_i}$, since $\mathcal Y_{\mathrm{as},V_i}=\mathcal Y_{\mathrm F,V_i}\rho_i(C)$. Under pullback to $u$, the Frobenius frame becomes $\mathcal Y_{\mathrm F,V_i}(z;u)\rho_i(Z^{-n})$, where $Z=e^{2\pi i\ell}$ and $\ell$ is the Frobenius exponent, by \eqref{eq:Frobenius-Symanzik-frame}. Alternatively, Proposition~\ref{prop:appendix-shifted-solution-columns} gives the column of $\Phi_n^{(i)}$ in the sectorial frame at $u$. Multiplication by $\rho_i(C(u))$ expresses this column in the same Frobenius frame. The two expressions give the shifted $(G,q)$-Wronskian equations
\begin{equation}\label{eq:appendix-shifted-Wronskian}
 \rho_i\bigl(Z^{-n}C(q^nu)\bigr)v_{\omega_i}
   =\rho_i\bigl(C(u)\Pi_n(u)\bigr)v_{\omega_i},
 \qquad i\in I,\quad 0\le n\le d .
\end{equation}

\begin{corollary}\label{cor:appendix-differential-reconstruction}
The central connection $C$ is the universal $(G,q)$-Wronskian in its orbit $[C]_{\mathrm{Wr}}$. In the original sectorial frame, the $q$-connection from $u$ to $qu$ has Steinberg form
\begin{equation}\label{eq:appendix-differential-reconstruction}
 \mathscr M(u)=C(qu)^{-1}ZC(u)\in U_c(K)s_{\mathrm{as}}(u),
 \qquad U_c=N_+\cap\dot cN_-\dot c^{-1},
\end{equation}
where $K$ is the field of meromorphic functions on the spectral domain and $N_\pm$ are the unipotent radicals of the two fixed opposite Borel subgroups.
\end{corollary}
\begin{proof}
The equation for $n=1$ makes $C$ a generalized $(G,q)$-Wronskian. By Definitions~\ref{def:Wronskian-right-orbit} and~\ref{def:universal-Wronskian-notation}, the universal representative is the solution of all the shifted equations with the prescribed fundamental highest-weight columns $\rho_i(C(u))v_{\omega_i}$. Equation~\eqref{eq:appendix-shifted-Wronskian} gives such a solution for every $1\le n\le d$.

Uniqueness follows from the common stabilizer of the vectors $v_{\omega_i}$, which is $N_+$ by Lemma~\ref{lem:unipotent-weight-scalar}. Indeed, preserving their lines places an element in $B_+=HN_+$. Fixing the vectors requires every fundamental weight to take value one on its torus factor. Since $G$ is simply connected, the fundamental weights form a basis of the character lattice of $H$, so that factor is the identity. A second solution with the same columns therefore has the form $C(u)b(u)$ with $b(u)\in N_+(K)$. At $n=d$, the two left sides of \eqref{eq:appendix-shifted-Wronskian} agree because $b(q^du)$ fixes the highest-weight vectors. Comparing the right sides shows that $b(u)$ also fixes every $\rho_i(\Pi_d(u))v_{\omega_i}$. Their common stabilizer is $\Pi_d(u)N_+\Pi_d(u)^{-1}$, since conjugating by $\Pi_d(u)^{-1}$ brings the stabilizer back to that of the highest-weight vectors. The element $\Pi_d(u)$ represents $c^{-d}=w_0$ in $N_G(H)/H$. Conjugation by a representative of $w_0$ takes the positive root subgroups generating $N_+$ to the negative root subgroups generating $N_-$. Thus $b(u)$ belongs to $N_+(K)\cap N_-(K)=\{1\}$, which proves uniqueness.

The same stabilizer argument yields the Steinberg form, as in Corollary~\ref{cor:universal-Wronskian-Steinberg}. The frame identity $C(qu)\mathscr M(u)=ZC(u)$ gives the first equality in \eqref{eq:appendix-differential-reconstruction}. Multiplying the equation at $n=1$ by $\rho_i(C(qu)^{-1}Z)$ shows that $\mathscr M(u)s_{\mathrm{as}}(u)^{-1}$ fixes every $v_{\omega_i}$, so it belongs to $N_+(K)$. The equation at $n=d-1$ based at $qu$ and the equation at $n=d$ based at $u$ both involve the Frobenius coefficient at $q^du$. Comparing them gives a second family of vectors fixed by $\mathscr M(u)s_{\mathrm{as}}(u)^{-1}$. Since $\Pi_d(u)=s_{\mathrm{as}}(u)^{-1}\Pi_{d-1}(qu)$, we obtain
\[
 \begin{aligned}
 \rho_i(Z^{-1}C(qu))\,\rho_i(\Pi_{d-1}(qu))v_{\omega_i}
    &=\rho_i(Z^{-d}C(q^du))v_{\omega_i}\\
    &=\rho_i(C(u)s_{\mathrm{as}}(u)^{-1})\,\rho_i(\Pi_{d-1}(qu))v_{\omega_i}.
 \end{aligned}
\]
Multiplication by $\rho_i(C(qu)^{-1}Z)$ now shows that $\mathscr M(u)s_{\mathrm{as}}(u)^{-1}$ also fixes the vectors $\rho_i(\Pi_{d-1}(qu))v_{\omega_i}$.

The common stabilizer of this second family of vectors is $\Pi_{d-1}(qu)N_+\Pi_{d-1}(qu)^{-1}$. For $1\le j\le d-1$, the $j$th factor of $\Pi_{d-1}(qu)$ is $s_{\mathrm{as}}(q^ju)^{-1}$ and represents $c^{-1}$ in the quotient by $H$. Therefore
\[
 \Pi_{d-1}(qu)H=c^{-(d-1)}=c(c^d)^{-1}=cw_0
 \qquad\text{in }W=N_G(H)/H .
\]
Here $c^d=w_0$ and $w_0^{-1}=w_0$. In particular, $\rho_i(\Pi_{d-1}(qu))v_{\omega_i}$ has weight $cw_0\omega_i$ for every $u$. With $\dot w_0=\dot c^{\,d}$ as representative of $w_0$, the equality of cosets means that $\Pi_{d-1}(qu)=\dot c\dot w_0\,t(u)$ for some $t(u)\in H$. The torus normalizes $N_+$, so the two elements conjugate $N_+$ to the same subgroup. Hence
\[
 \Pi_{d-1}(qu)N_+\Pi_{d-1}(qu)^{-1}
    =\dot c\bigl(\dot w_0N_+\dot w_0^{-1}\bigr)\dot c^{-1}
    =\dot cN_-\dot c^{-1}.
\]
Intersecting this subgroup with $N_+$ gives $\mathscr M(u)s_{\mathrm{as}}(u)^{-1}\in U_c(K)$, as claimed. In rank one $d=1$, so $\Pi_{d-1}(qu)=1$ and $cw_0=1$. Here $U_c=N_+$, and the first shifted equation already gives the stated form.
\end{proof}

\begin{remark}\label{rem:forward-opposite-normal-forms}
The equation at $n=d$ gives an explicit change from the sectorial frame to a frame in which the inverse connection has the normalization used in the proof of \cite[Proposition~4.20]{KZ2021}. For a universal $(G,q)$-Wronskian $\mathscr W$, write $A_{\mathscr W}(u)=\mathscr W(qu)^{-1}Z\mathscr W(u)$ for its forward connection. A second meromorphic frame in the reference bundle is specified by
\[
 g(q^du)=Z^d\mathscr W(u)\Pi_d(u).
\]
This defines $g$ on the spectral domain because $u\mapsto q^du$ is a bijection of that domain. The equation for $\mathscr W$ at $n=d$ says that $g(q^du)$ and $\mathscr W(q^du)$ agree on every fundamental highest-weight vector. Their quotient $\mathscr W(q^du)^{-1}g(q^du)$ therefore lies in $N_+(K)$, so the two frames define the same $B_+$-reduction. The reference connection has matrix $Z$. Its matrix from $q^du$ to $q^{d+1}u$ in the new frame is
\begin{equation}\label{eq:BW-shifted-Steinberg-frame}
 g(q^{d+1}u)^{-1}Zg(q^du)=\Pi_d(qu)^{-1}A_{\mathscr W}(u)\Pi_d(u).
\end{equation}
The two products on the right can be written as
\[
 \Pi_d(u)=s_{\mathrm{as}}(u)^{-1}\Pi_{d-1}(qu),
 \qquad
 \Pi_d(qu)=\Pi_{d-1}(qu)s_{\mathrm{as}}(q^du)^{-1}.
\]
After multiplication on the left by $s_{\mathrm{as}}(q^du)^{-1}$, the right-hand side of \eqref{eq:BW-shifted-Steinberg-frame} becomes
\[
 \Pi_{d-1}(qu)^{-1}\bigl(A_{\mathscr W}(u)s_{\mathrm{as}}(u)^{-1}\bigr)\Pi_{d-1}(qu).
\]
The middle factor belongs to $U_c$ by Corollary~\ref{cor:universal-Wronskian-Steinberg}. Conjugation by $\Pi_{d-1}(qu)^{-1}$ takes $\dot cN_-\dot c^{-1}$ back to $N_+$ by the stabilizer calculation above. It takes $N_+$ to $\dot c^{-1}N_-\dot c$ because the inverse product represents $(cw_0)^{-1}=c^{-1}w_0$. Here $c$ commutes with $w_0=c^d$. Consequently it sends their intersection $U_c$ to $U_{c^{-1}}=N_+\cap\dot c^{-1}N_-\dot c$. The forward connection in \eqref{eq:BW-shifted-Steinberg-frame} therefore lies in $s_{\mathrm{as}}(q^du)U_{c^{-1}}(K)$. Its inverse belongs to $U_{c^{-1}}(K)s_{\mathrm{as}}(q^du)^{-1}$, as in the opposite normalization. When $h_{\mathrm{as}}\equiv1$, both products $\Pi_d(u)$ and $\Pi_d(qu)$ equal $\dot c^{-d}$, so the comparison in \eqref{eq:BW-shifted-Steinberg-frame} is conjugation by $\dot c^{\,d}$.
\end{remark}

\end{document}